\documentclass{amsart}
\usepackage{amssymb}
\usepackage{mathrsfs}
\usepackage{stmaryrd} 
\usepackage{bbm}
\usepackage{oldgerm}
\usepackage[francais]{babel}
\usepackage[T1]{fontenc}
\usepackage[latin1]{inputenc}
\usepackage[all]{xy}
\usepackage{hyperref}

\makeindex

\newtheorem{teo}[subsection]{Théorème}
\newtheorem{prop}[subsection]{Proposition}
\newtheorem{cor}[subsection]{Corollaire}
\newtheorem{lem}[subsection]{Lemme}

\theoremstyle{definition}

\newtheorem{defi}[subsection]{Définition}
\newtheorem{rema}[subsection]{Remarque}
\newtheorem{remas}[subsection]{Remarques}
\newtheorem{exemple}[subsection]{Exemple}

\numberwithin{equation}{subsection}

\newcommand{\gtimes}{\stackrel{\leftarrow}{\times}}

\newcommand{\mN}{{\mathbb N}}

\newcommand{\mZ}{{\mathbb Z}}

\newcommand{\mU}{{\mathbb U}}
\newcommand{\mV}{{\mathbb V}}

\newcommand{\bB}{{\bf B}}
\newcommand{\bD}{{\bf D}}

\newcommand{\Et}{{\bf \acute{E}t}}
\newcommand{\Sch}{{\bf Sch}}
\newcommand{\Ens}{{\bf Ens}}
\newcommand{\Pt}{{\bf Pt}}

\newcommand{\Top}{{\bf Top}}
\newcommand{\bHom}{{\bf Hom}}
\newcommand{\bMod}{{\bf Mod}}

\newcommand{\et}{{\rm \acute{e}t}}
\newcommand{\fet}{{\rm f\acute{e}t}}

\newcommand{\coh}{{\rm coh}}
\newcommand{\cart}{{\rm cart}}
\newcommand{\scoh}{{\rm scoh}}

\newcommand{\rf}{{\rm f}}

\newcommand{\Fl}{{\rm Fl}}
\newcommand{\ob}{{\rm Ob}}

\newcommand{\pr}{{\rm pr}}

\newcommand{\coker}{{\rm coker}}

\newcommand{\id}{{\rm id}}

\newcommand{\rb}{{\rm b}}

\newcommand{\Hom}{{\rm Hom}}

\newcommand{\rE}{{\rm E}}
\newcommand{\rH}{{\rm H}}

\newcommand{\rR}{{\rm R}}

\newcommand{\rp}{{\rm p}}

\newcommand{\oK}{{\overline{K}}}

\newcommand{\ox}{{\overline{x}}}
\newcommand{\oy}{{\overline{y}}}
\newcommand{\oz}{{\overline{z}}}

\newcommand{\uE}{{\underline{E}}}
\newcommand{\uX}{{\underline{X}}}
\newcommand{\uY}{{\underline{Y}}}
\newcommand{\uf}{{\underline{f}}}

\newcommand{\ubeta}{{\underline{\beta}}}
\newcommand{\urho}{{\underline{\rho}}}

\newcommand{\bvA}{{\breve{A}}}

\newcommand{\hD}{{\widehat{D}}}
\newcommand{\hE}{{\widehat{E}}}
\newcommand{\hI}{{\widehat{I}}}

\newcommand{\hPhi}{\widehat{\Phi}}

\newcommand{\cC}{{\mathscr C}}
\newcommand{\cD}{{\mathscr D}}
\newcommand{\cE}{{\mathscr E}}
\newcommand{\cF}{{\mathscr F}}
\newcommand{\cG}{{\mathscr G}}

\newcommand{\cP}{{\mathscr P}}
\newcommand{\co}{{\mathscr O}}
\newcommand{\cR}{{\mathscr R}}

\newcommand{\cT}{{\mathscr T}}
\newcommand{\cH}{{\mathscr H}}

\newcommand{\cQ}{{\mathscr Q}}

\newcommand{\cHom}{{\mathscr Hom}}

\newcommand{\fD}{{\mathfrak D}}
\newcommand{\fE}{{\mathfrak E}}
\newcommand{\fF}{{\mathfrak F}}
\newcommand{\fM}{{\mathfrak M}}

\newcommand{\fR}{{\mathfrak R}}

\newcommand{\hY}{{\widehat{Y}}}

\newcommand{\hcC}{{\widehat{\cC}}}

\newcommand{\tC}{{\widetilde{C}}}
\newcommand{\tD}{{\widetilde{D}}}
\newcommand{\tE}{{\widetilde{E}}}

\newcommand{\tuE}{{\widetilde{\uE}}}
\newcommand{\tI}{{\widetilde{I}}}

\newcommand{\tS}{{\widetilde{S}}}

\newcommand{\tX}{{\widetilde{X}}}
\newcommand{\tY}{{\widetilde{Y}}}
\newcommand{\tZ}{{\widetilde{Z}}}

\newcommand{\tx}{{\widetilde{x}}}
\newcommand{\ty}{{\widetilde{y}}}

\newcommand{\tfE}{{\widetilde{\fE}}}

\newcommand{\ttau}{{\widetilde{\tau}}}

\newcommand{\tcC}{{\widetilde{\cC}}}

\begin{document}

\title{Topos co-évanescents et généralisations}
\author{Ahmed Abbes et Michel Gros}
\address{A.A. Laboratoire Alexander Grothendieck, ERL 9216 du CNRS, Institut des Hautes \'Etudes Scientifiques, 35 route de Chartres, 91440 Bures-sur-Yvette, France}
\address{M.G. CNRS UMR 6625, IRMAR, Université de Rennes 1,
Campus de Beaulieu, 35042 Rennes cedex, France}
\email{abbes@ihes.fr}
\email{michel.gros@univ-rennes1.fr}

\maketitle

\setcounter{tocdepth}{1}
\tableofcontents

\section{Introduction}

\subsection{}\label{higgs2-introd1}
Le but de cet article est de consolider les fondements topologiques nécessaires à l'approche de Faltings  
en théorie de Hodge $p$-adique \cite{faltings1,faltings2,faltings3}. 
Sa genèse a été motivée par notre travail sur la correspondance de Simpson $p$-adique \cite{ag0,ag1,ag3}. 
L'approche de Faltings se compose schématiquement de deux étapes. La première, de nature locale,
est une généralisation des techniques galoisiennes de Tate, Sen et Fontaine à certains schémas 
affines au-dessus de corps locaux $p$-adiques. Elle utilise d'une façon essentielle sa théorie des extensions 
presque-étales \cite{tsuji2}. La seconde étape, de nature plus globale, relie la cohomologie étale $p$-adique 
de certains schémas sur des corps locaux $p$-adiques à la cohomologie galoisienne étudiée dans la première étape. 
Pour ce faire, Faltings utilise d'une part la notion de schémas $K(\pi,1)$ \cite{faltings1},
et d'autre part un nouveau topos \cite{faltings2}.
Ce dernier n'a été explicitement défini qu'assez tardivement par rapport au reste de la théorie
(\cite{faltings2} page 214), et ne nous semble pas avoir reçu l'attention qu'il mérite.

\subsection{}\label{higgs2-introd2}
Nous avons tiré profit   
d'une lettre de Deligne à Illusie \cite{deligne1} (antérieure à \cite{faltings2}) et de notes d'Illusie \cite{illusie3}. 
Dans cette lettre, Deligne suggère que le topos dont a besoin Faltings devrait être en réalité un 
topos {\em co-évanescent}, autrement dit, un cas particulier du produit orienté 
de topos qu'il avait introduit pour développer le formalisme des cycles évanescents en dimensions supérieures \cite{illusie2}. 
Le premier auteur du présent article a remarqué en 2008 que la topologie co-évanescente de Deligne diffère 
en général de celle considérée par Faltings dans (\cite{faltings2} page 214); toutefois, ce dernier utilise une caractérisation
des faisceaux qui vaut pour les topos co-évanescents (cf. \ref{higgs2-introd6}).\footnote{Ce problème ne semble
pas avoir été observé antérieurement.} Il s'avère que cette caractérisation 
n'est pas en général vérifiée par le topos considéré initialement par Faltings (cf. \ref{higgs2-tf7} pour un contre exemple). 
En fait, l'un des principaux faisceaux utilisés par Faltings, à savoir le faisceau structural, 
n'est pas en général un faisceau pour la topologie qu'il a définie dans (\cite{faltings2} page 214); 
mais c'est un faisceau pour la topologie co-évanescente (cf. \cite{ag3} 8.16 et 8.18). 
Nous nous proposons dans cet article de corriger  
la définition du topos introduit par Faltings et de le développer suivant la suggestion de Deligne. 
C'est ce nouveau topos que nous appellerons {\em topos de Faltings}~; 
celui introduit  dans (\cite{faltings2} page 214) semble avoir très peu d'intérêt.

\subsection{}\label{higgs2-introd3}
Les résultats de cet article étant de nature assez technique, nous allons maintenant en présenter un résumé détaillé.
Commençons par rappeler la définition des produits orientés
de topos due à Deligne. Soient $f\colon \tX \rightarrow \tS$ et $g\colon \tY\rightarrow \tS$ deux morphismes de $\mU$-topos,
où $\mU$ est un univers fixé. Le {\em produit orienté} $\tX\gtimes_{\tS}\tY$ est un $\mU$-topos muni de deux morphismes 
\begin{equation}\label{higgs2-introd3a}
\rp_1\colon \tX\gtimes_{\tS}\tY\rightarrow \tX\ \ \ {\rm et}\ \ \ 
\rp_2\colon \tX\gtimes_{\tS}\tY\rightarrow \tY
\end{equation}
et d'un $2$-morphisme 
\begin{equation}
\tau\colon g\rp_2\rightarrow f\rp_1,
\end{equation}
tels que le quadruplet $(\tX\gtimes_{\tS}\tY,\rp_1,\rp_2,\tau)$ soit universel dans la $2$-catégorie des $\mU$-topos. 
On peut lui construire explicitement un site sous-jacent $C$ à partir 
de la donnée de $\mU$-sites $X$, $Y$ et $S$ sous-jacents à $\tX$, $\tY$ et $\tS$, respectivement, dans lesquels 
les limites projectives finies sont représentables, et de deux foncteurs 
continus et exacts à gauche  $f^+\colon S \rightarrow X$ et $g^+\colon S\rightarrow Y$ définissant $f$ et $g$ 
(cf. \ref{higgs2-topfl1} et \ref{higgs2-topfl6}).
Suivant Illusie, nous appelons $\tX\gtimes_{\tS}\tS$ le {\em topos évanescent} de $f$, 
et $\tS\gtimes_{\tS}\tY$ le {\em topos co-évanescent} de $g$. Le premier topos a été utilisé par Deligne 
pour développer le formalisme des cycles évanescents de $f$, 
et a été étudié par Gabber, Illusie \cite{illusie2}, Laumon \cite{laumon} et Orgogozo \cite{orgogozo}. 
Le second topos est le prototype du topos de Faltings. On peut lui construire explicitement un autre site sous-jacent 
$D$, plus simple que $C$, que nous appelons site {\em co-évanescent} de $g^+$ 
 (cf. \ref{higgs2-co-ev1} et \ref{higgs2-co-ev101}).

\subsection{}\label{higgs2-introd4}
\`A proprement parler, le topos de Faltings n'est en fait pas un topos co-évanescent, 
mais un cas particulier d'une notion plus générale que nous développons dans cet article, et que nous 
baptisons {\em topos co-évanescent généralisé}. Soient $I$ un $\mU$-site,  $\tI$ le topos des faisceaux 
de $\mU$-ensembles sur $I$,
\begin{equation}\label{higgs2-introd4a}
\pi\colon E\rightarrow I
\end{equation}
une catégorie fibrée, clivée et normalisée au-dessus de la catégorie sous-jacente à $I$. 
On suppose les conditions suivantes satisfaites~:
\begin{itemize}
\item[(i)] Les produits fibrés sont représentables dans $I$.  
\item[(ii)]  Pour tout $i\in \ob(I)$, la catégorie fibre $E_i$ de $E$ au-dessus de $i$ 
est munie d'une topologie faisant de celle-ci un $\mU$-site, 
et les limites projectives finies sont représentables dans $E_i$. 
On note $\tE_i$ le topos des faisceaux de $\mU$-ensembles sur $E_i$.
\item[(iii)] Pour tout morphisme $f\colon i\rightarrow j$ de $I$, le foncteur image inverse $f^+\colon E_j\rightarrow E_i$
est continu et exact à gauche. Il définit donc un morphisme de topos que l'on note aussi (abusivement) 
$f\colon \tE_i\rightarrow \tE_j$. 
\end{itemize}
Le foncteur $\pi$ est en fait un $\mU$-site fibré (cf. \ref{higgs2-tcevg1}). On appelle topologie {\em co-évanescente} de $E$ 
la topologie engendrée par les familles de recouvrements $(V_n\rightarrow V)_{n\in \Sigma}$
des deux types suivants~:
\begin{itemize}
\item[(v)] Il existe $i\in \ob(I)$ tel que $(V_n\rightarrow V)_{n\in \Sigma}$ soit 
une famille couvrante de $E_i$.
\item[(c)] Il existe une famille couvrante de morphismes $(f_n\colon i_n\rightarrow i)_{n\in \Sigma}$ de $I$ 
telle que $V_n$ soit isomorphe à $f_n^+(V)$ pour tout $n\in \Sigma$.
\end{itemize}
Les recouvrements du type (v) sont dits {\em verticaux}, et ceux du type (c) sont dits {\em cartésiens}. 
Le site ainsi défini est appelé site {\em co-évanescent} associé au site fibré $\pi$; c'est un $\mU$-site. 
On appelle topos {\em co-évanescent} associé au site fibré $\pi$, et l'on note $\tE$, 
le topos des faisceaux de $\mU$-ensembles sur $E$ (cf. \ref{higgs2-tcevg3}). 
Lorsque $I$ est munie de la topologie chaotique, c'est-à-dire de la topologie la moins fine
parmi toutes les topologies de $I$, 
on retrouve la topologie totale sur $E$ relative au site fibré $\pi$ (cf.~\ref{higgs2-tcevg41}). 

\subsection{}\label{higgs2-introd6}
On montre  \eqref{higgs2-tcevg5} que la donnée d'un faisceau 
$F$ sur $E$ est équivalente à la donnée pour tout objet $i$ de $I$ d'un faisceau $F_i$ sur 
$E_i$ et pour tout morphisme $f\colon i\rightarrow j$ de $I$ d'un morphisme 
$F_j\rightarrow f_*(F_i)$, ces morphismes étant soumis à des relations de compatibilité, 
tels que pour toute famille couvrante $(f_n\colon i_n\rightarrow i)_{n\in \Sigma}$ de $I$, 
si pour tout $(m,n)\in \Sigma^2$, on pose $i_{mn}=i_m\times_ii_n$ et on note $f_{mn}\colon i_{mn}\rightarrow i$
le morphisme canonique, la suite de morphismes de faisceaux sur $E_i$
\begin{equation}\label{higgs2-introd6a}
F_i\rightarrow \prod_{n\in \Sigma}(f_{n})_*(F_{i_n})\rightrightarrows \prod_{(m,n)\in \Sigma^2} (f_{mn})_*(F_{i_{mn}})
\end{equation}
soit exacte. On identifiera dans la suite $F$ au foncteur $\{i\mapsto F_i\}$ qui lui est associé.

On étudie la fonctorialité des sites et topos co-évanescents 
par rapport au foncteur fibrant $\pi$ (\ref{higgs2-tcevg85} et \ref{higgs2-tcevg30}) et par changement de base \eqref{higgs2-tcevg9}. 
On établit aussi quelques propriétés de cohérence. On montre \eqref{higgs2-tcevg155},
entre autres, que s'il existe une sous-catégorie pleine, $\mU$-petite et topologiquement génératrice $I'$ de $I$, 
formée d'objets quasi-compacts, stable par produits fibrés, telle que pour tout $i\in \ob(I')$, le topos $\tE_i$ 
soit cohérent et que pour tout morphisme 
$f\colon i\rightarrow j$ de $I'$,
le morphisme de topos $f\colon \tE_i\rightarrow \tE_j$ soit cohérent, alors le topos $\tE$ est localement cohérent.
Si de plus, la catégorie $I$ admet un objet final qui appartient à $I'$, alors le topos $\tE$ est cohérent.

\subsection{}\label{higgs2-introd5}
Le lien avec le topos co-évanescent de Deligne se fait de la façon suivante. 
Soient $X$ et $Y$ deux $\mU$-sites dans lesquels les 
limites projectives finies sont représentables, $f^+\colon X\rightarrow Y$ un foncteur continu et exact à gauche. 
On note $\tX$ et $\tY$ les topos des faisceaux de $\mU$-ensembles
sur $X$ et $Y$, respectivement,  et $f\colon \tY\rightarrow \tX$ le morphisme de topos défini par $f^+$. 
Considérons la catégorie $\Fl(Y)$ des morphismes de $Y$, et le ``foncteur but'' 
\begin{equation}\label{higgs2-introd5a}
\Fl(Y)\rightarrow Y,
\end{equation}
qui fait de $\Fl(Y)$ une catégorie fibrée, clivée et normalisée au-dessus de $Y$; la catégorie fibre au-dessus de 
tout $V\in \ob(Y)$ est canoniquement équivalente à la catégorie $Y_{/V}$.  
Munissant chaque fibre $Y_{/V}$ de la topologie induite par celle de $Y$,
$\Fl(Y)/Y$ devient un $\mU$-site fibré, vérifiant les conditions de \eqref{higgs2-introd4}. Soit 
\begin{equation}\label{higgs2-introd5b}
\pi\colon E\rightarrow X
\end{equation}
le site fibré déduit de $\Fl(Y)/Y$ par changement de base par le foncteur $f^+$. 
Le site co-évanescent $E$ associé à $\pi$ \eqref{higgs2-introd4} est canoniquement équivalent au site co-évanescent $D$ 
associé au foncteur $f^+$ \eqref{higgs2-introd3}; d'où la terminologie. Le topos co-évanescent $\tE$ 
associé à $\pi$ est donc canoniquement équivalent au topos co-évanescent $\tX\gtimes_{\tX}\tY$ associé à $f$ (cf. \ref{higgs2-tcevg40}).

\subsection{}\label{higgs2-introd7}
Conservons les hypothèses de \eqref{higgs2-introd4}, et supposons, de plus, que les limites projectives finies 
soient représentables dans $I$; compte tenu de \ref{higgs2-introd4}(i), il revient au même de demander que $I$ 
admette un objet final $\iota$, que nous supposons fixé dans la suite. On peut alors définir pour $\tE$ des analogues 
des projections canoniques du produit orienté \eqref{higgs2-introd3a}. 
D'une part, le foncteur d'injection canonique $\alpha_{\iota!}\colon E_\iota\rightarrow E$
est continu et exact à gauche (cf. \ref{higgs2-tcevg18}). Il définit donc un morphisme de topos 
\begin{equation}
\beta\colon \tE\rightarrow \tE_\iota,
\end{equation}
analogue de la seconde projection $\rp_2$ \eqref{higgs2-introd3a}.  D'autre part, on se donne un objet final $e$ de $E_\iota$,
qui existe d'après \ref{higgs2-introd4}(ii) et 
que nous supposons fixé dans la suite. Il existe alors essentiellement 
une unique section cartésienne de $\pi$ \eqref{higgs2-introd4a}
\begin{equation}
\sigma^+\colon I\rightarrow E
\end{equation}
telle que $\sigma^+(\iota)=e$. Pour tout $i\in \ob(I)$, $\sigma^+(i)$ est un objet final de $E_i$.  
On vérifie facilement que $\sigma^+$ est continu et exact à gauche (cf. \ref{higgs2-tcevg18}). Il définit donc un morphisme de topos 
\begin{equation}
\sigma\colon \tE\rightarrow \tI,
\end{equation}
analogue de la première projection $\rp_1$ \eqref{higgs2-introd3a}. 

\subsection{}\label{higgs2-introd8}
Conservons les hypothèses de \eqref{higgs2-introd7}, et soient, de plus, 
$V$ un objet de $E$, $c=\pi(V)$, 
\begin{equation}
\varpi\colon E_{/V}\rightarrow I_{/c}
\end{equation}
le foncteur induit par $\pi$. Pour tout morphisme $f\colon i\rightarrow c$ de $I$, la catégorie fibre de $\varpi$
au-dessus de l'objet $(i,f)$ de $I_{/c}$ est canoniquement équivalente à la catégorie $(E_i)_{/f^+(V)}$. 
Munissant $I_{/c}$ de la topologie induite par celle de $I$,
et chaque fibre $(E_i)_{/f^+(V)}$ de la topologie induite par celle de $E_i$, 
$\varpi$ devient un site fibré vérifiant les conditions de \eqref{higgs2-introd4}.
On montre \eqref{higgs2-tcevg71} que la topologie co-évanescente de $E_{/V}$ est induite
par celle de $E$ au moyen du foncteur canonique $E_{/V}\rightarrow E$. 
En particulier, le topos des faisceaux de $\mU$-ensembles
sur le site co-évanescent $E_{/V}$ est canoniquement équivalent à $\tE_{/\varepsilon(V)}$, où 
$\varepsilon \colon E\rightarrow \tE$ est le foncteur canonique.

\subsection{}\label{higgs2-introd9}
Une des principales caractéristiques du topos co-évanescent de Deligne est l'existence d'un morphisme 
des {\em cycles co-proches} (cf. \ref{higgs2-co-ev14}). 
Celle-ci s'étend aussi au topos co-évanescent généralisé introduit dans \eqref{higgs2-introd4}. 
Conservons les hypothèses de \eqref{higgs2-introd7} et considérons, de plus, un $\mU$-site $X$ et un foncteur continu
et exact à gauche $\Psi^+\colon E\rightarrow X$. 
On désigne par $\tX$ le topos des faisceaux de $\mU$-ensembles sur $X$ et par 
\begin{equation}\label{higgs2-introd9a}
\Psi\colon \tX\rightarrow \tE
\end{equation}
le morphisme de topos associé à $\Psi^+$. On pose $u^+=\Psi^+\circ \sigma^+\colon I\rightarrow X$ et 
\begin{equation}\label{higgs2-introd9b}
u=\sigma\Psi\colon \tX\rightarrow \tI.
\end{equation} 
Pour tout $i\in \ob(I)$, le foncteur $\Psi^+$ induit un foncteur $\Psi_i^+\colon E_i\rightarrow X_{/u^+(i)}$. 
Munissant $X_{/u^+(i)}$ de la topologie induite par celle de $X$, 
$\Psi^+_i$ est exact à gauche et continu (cf. \ref{higgs2-fccp1}).
Il définit donc un morphisme de topos 
\begin{equation}\label{higgs2-introd9c}
\Psi_i\colon \tX_{/u^*(i)}\rightarrow \tE_i.
\end{equation}
Le morphisme $\Psi_\iota$ n'est autre que le composé $\beta\Psi\colon \tX\rightarrow \tE_\iota$. 
De la relation $\Psi^*\beta^*=\Psi_\iota^*$, on déduit par adjonction un morphisme 
\begin{equation}\label{higgs2-introd9d}
\beta^*\rightarrow \Psi_*\Psi_\iota^*.
\end{equation}
Généralisant une propriété importante du topos co-évanescent de Deligne, 
on montre \eqref{higgs2-fccp4} que si pour tout $i\in \ob(I)$, le morphisme d'adjonction $\id\rightarrow \Psi_{i*}\Psi_i^*$
est un isomorphisme, alors les morphismes d'adjonction $\id\rightarrow \beta_*\beta^*$ et 
$\beta^*\rightarrow \Psi_*\Psi_\iota^*$ sont des isomorphismes. 

\subsection{}\label{higgs2-introd10}
Conservons les hypothèses de \eqref{higgs2-introd9}, 
supposons de plus que les limites projectives finies soient représentables dans $X$.
On note $\varpi\colon D\rightarrow I$ le site fibré associé au foncteur $u^+\colon I\rightarrow X$ défini dans \eqref{higgs2-introd5},
et on munit $D$ de la topologie co-évanescente relative à $\varpi$. On obtient ainsi le site co-évanescent 
associé au foncteur $u^+$ \eqref{higgs2-introd3}, dont le topos des faisceaux de $\mU$-ensembles est $\tI\gtimes_{\tI}\tX$. 
D'après la propriété universelle des produits orientés, 
les morphismes $u\colon \tX\rightarrow \tI$ et $\id_{\tX}$ et le $2$-morphisme $\id_u$ 
définissent un morphisme de topos, dit morphisme des {\em cycles co-proches} (cf.~\ref{higgs2-co-ev14})
\begin{equation}
\Psi_D\colon \tX\rightarrow \tI\gtimes_{\tI}\tX.
\end{equation}
D'autre part, les foncteurs $\Psi_i^+$ pour tout $i\in \ob(I)$ définissent
un $I$-foncteur cartésien $\rho^+\colon E\rightarrow D$. 
Celui-ci est continu et exact à gauche (cf. \ref{higgs2-fccp5}). Il définit donc un morphisme de topos 
\begin{equation}
\rho\colon \tI\gtimes_{\tI}\tX\rightarrow \tE.
\end{equation}
On vérifie aussitôt que les carrés du diagramme
\begin{equation}
\xymatrix{
{\tI}\ar@{=}[d]&{\tI\gtimes_{\tI}\tX}\ar[l]_-(0.4){\rp_1}\ar[d]^{\rho}\ar[r]^-(0.5){\rp_2}&{\tX}\ar[d]^-(0.5){\Psi_\iota}\\
{\tI}&{\tE}\ar[l]_-(0.5){\sigma}\ar[r]^-(0.5){\beta}&{\tE_\iota}}
\end{equation}
et le diagramme 
\begin{equation}
\xymatrix{
{\tX}\ar[r]^-(0.5){\Psi_D}\ar[rd]_{\Psi}&{\tI\gtimes_{\tI}\tX}\ar[d]^{\rho}\\
&{\tE}}
\end{equation}
sont commutatifs à isomorphismes canoniques près.

\subsection{}
Conservons les hypothèses de \eqref{higgs2-introd7} et soit $R$ un anneau de $\tE$; 
il revient au même de se donner, pour tout $i\in \ob(I)$, un anneau $R_i$ de $\tE_i$, 
et pour tout morphisme $f\colon i\rightarrow j$ de $I$, un homomorphisme d'anneaux
$R_j\rightarrow f_*(R_i)$, ces homomorphismes étant soumis à des relations de compatibilité
et de recollement \eqref{higgs2-introd6a}. 
Nous développons dans § \ref{higgs2-tcea} quelques sorites sur la catégorie des $R$-modules de $\tE$, 
en particulier, sur le produit tensoriel \eqref{higgs2-tcea9} et le faisceau des morphismes \eqref{higgs2-tcea8}.
Nous étudions aussi des invariants cohomologiques élémentaires du topos annelé $(\tE,R)$. 
Pour ce faire, nous reprenons dans § 7
le formalisme des topos totaux annelés (\cite{sga4} VI 8.6) qui nous servent d'intermédiaire. 
Supposons pour simplifier que la catégorie $I$ soit $\mU$-petite (cf. \ref{higgs2-tcea4} pour le cas général).
On désigne par $\Top(E)$ le topos total associé au site fibré $\pi$, c'est-à-dire le topos des faisceaux de 
$\mU$-ensembles sur le site total $E$ (cf. \ref{higgs2-tta1}). On a alors un plongement canonique de topos \eqref{higgs2-tcevg22a}
\begin{equation}
\delta\colon \tE\rightarrow \Top(E)
\end{equation}
tel que le foncteur $\delta_*$ soit le foncteur d'inclusion canonique de $\tE$ dans $\Top(E)$. On le
considère comme un morphisme de topos annelés (respectivement, par $R$ et $\delta_*(R)$),
et on calcule dans \ref{higgs2-tcea5} les foncteurs dérivés droits du foncteur $\delta_*$. 
La suite spectrale de Cartan-Leray associée à $\delta$ fournit alors une suite spectrale qui calcule 
les images directes supérieures d'un morphisme de topos co-évanescents généralisés annelés \eqref{higgs2-tcea6}. 
On notera que la suite spectrale de Cartan-Leray associée à $\delta$
qui calcule la cohomologie d'un $R$-module de $\tE$ est la même que celle associée à $\beta$  \eqref{higgs2-tcea22}.

\subsection{}\label{higgs2-introd11}
Le reste de l'article est consacré à l'étude d'un cas particulier du topos co-évanescent généralisé, 
à savoir le topos de Faltings. En préambule, on développe dans § \ref{higgs2-Kp} quelques sorites 
sur le topos fini étale que nous n'avons pu trouver dans la littérature.  
Pour tout schéma $X$, on désigne par $\Et_{/X}$ son site étale, c'est à dire la catégorie des 
schémas étales au-dessus de $X$ (éléments de $\mU$), munie de la topologie étale, et par $X_\et$
le topos des faisceaux de $\mU$-ensembles sur $\Et_{/X}$.
On appelle site {\em fini étale} de $X$ et l'on note $\Et_{\rf/X}$ la sous-catégorie pleine de $\Et_{/X}$
formée des revêtements étales de $X$ (c'est-à-dire, des schémas étales finis sur $X$), 
munie de la topologie induite par celle de $\Et_{/X}$. 
On appelle topos  {\em fini étale} de $X$ et l'on note $X_\fet$ le topos des faisceaux de $\mU$-ensembles sur $\Et_{\rf/X}$.
Le foncteur d'injection canonique $\Et_{\rf/X}\rightarrow \Et_{/X}$ induit un morphisme de topos
\begin{equation}\label{higgs2-introd11a}
\rho_X\colon X_\et\rightarrow X_\fet.
\end{equation}
On montre \eqref{higgs2-Kpp4} que si $X$ est un schéma cohérent, ayant un nombre fini de composantes connexes, 
alors le morphisme d'adjonction $\id\rightarrow \rho_{X*}\rho_X^*$ est un isomorphisme; 
en particulier, le foncteur $\rho^*_X\colon X_\fet\rightarrow X_\et$ est pleinement fidèle. 
Son image essentielle est formée des limites inductives filtrantes de faisceaux localement constants et constructibles 
\eqref{higgs2-Kpp5}; on établit aussi une variante pour les faisceaux abéliens de torsion. 

On désigne par $\Sch$ la catégorie des schémas éléments de $\mU$,
par $\cR$ la catégorie des revêtements étales 
({\em i.e.}, la sous-catégorie pleine de la catégorie des morphismes de $\Sch$, formée des revêtements étales) et par 
\begin{equation}\label{higgs2-introd11b}
\cR\rightarrow \Sch
\end{equation}
le ``foncteur but'', qui fait de $\cR$ une catégorie fibrée clivée et normalisée au-dessus de $\Sch$~; 
la catégorie fibre au-dessus d'un schéma $X$ est canoniquement équivalente à la catégorie $\Et_{\rf/X}$. 
On considère $\cR/\Sch$ comme un $\mU$-site fibré en munissant chaque fibre de la topologie étale.

\subsection{}\label{higgs2-introd12}
Soit $f\colon Y\rightarrow X$ un morphisme de schémas. On appelle {\em site fibré de Faltings} associé à $f$
le $\mU$-site fibré 
\begin{equation}\label{higgs2-introd12a}
\pi\colon E\rightarrow \Et_{/X}
\end{equation}
déduit de $\cR/\Sch$ \eqref{higgs2-introd11b} 
par changement de base par le foncteur 
\begin{equation}\label{higgs2-introd12b}
\Et_{/X}\rightarrow \Sch, \ \ \ U\mapsto U\times_XY.
\end{equation} 
On peut décrire explicitement la catégorie $E$ de la façon suivante (cf. \ref{higgs2-tf1}). Les objets de $E$ 
sont les morphismes de schémas $V\rightarrow U$ au-dessus de $f\colon Y\rightarrow X$ tels que le morphisme
$U\rightarrow X$ soit étale et que le morphisme $V\rightarrow U_Y=U\times_XY$ soit étale fini. 
Soient $(V'\rightarrow U')$, $(V\rightarrow U)$ deux objets de $E$. Un morphisme 
de $(V'\rightarrow U')$ dans $(V\rightarrow U)$ est la donnée d'un $X$-morphisme $U'\rightarrow U$ et 
d'un $Y$-morphisme $V'\rightarrow V$ tels que le diagramme
\begin{equation}\label{higgs2-introd12c}
\xymatrix{
V'\ar[r]\ar[d]&U'\ar[d]\\
V\ar[r]&U}
\end{equation}
soit commutatif. Le foncteur $\pi$ est alors défini pour tout $(V\rightarrow U)\in \ob(E)$, par
$\pi(V\rightarrow U)=U$. Celui-ci vérifie clairement les conditions de \eqref{higgs2-introd4} et \eqref{higgs2-introd7}. 
On peut donc lui appliquer les constructions développées plus haut. On munit $E$ de la topologie 
co-évanescente relative à $\pi$. Le site co-évanescent ainsi défini est appelé {\em site de Faltings} associé à $f$;
c'est un $\mU$-site. On appelle {\em topos de Faltings} associé à $f$ et l'on note $\tE$
le topos des faisceaux de $\mU$-ensembles 
sur $E$. En fait, Faltings se limite au cas suivant~: 
soient $K$ un corps de valuation discrète complet de caractéristique $0$, de corps résiduel parfait de 
caractéristique $p>0$, $\co_K$ l'anneau de valuation de $K$, $\oK$ une clôture algébrique de $K$,
$X$ un $\co_K$-schéma séparé de type fini, $X^\circ$ un ouvert de $X$. Faltings considère exclusivement 
le cas où $f$ est le morphisme canonique $X^\circ_\oK\rightarrow X$. 

De nombreux résultats sur les site et topos de Faltings dans § \ref{higgs2-tf} 
sont des applications directes de ceux établis dans §§ \ref{higgs2-tcevg} et \ref{higgs2-fccp}
pour les sites et topos co-évanescents généralisés.
On y trouve en particulier une étude de la fonctorialité par rapport à $f$ \eqref{higgs2-tf19}, 
et de la localisation par rapport à un objet de $E$ \eqref{higgs2-tf17}. On résume dans la suite de cette introduction 
d'autres propriétés plus spécifiques. 

\subsection{}\label{higgs2-introd13}
Conservons les hypothèses de \eqref{higgs2-introd12}, et supposons de plus $X$ quasi-séparé. 
On désigne par $\Et_{\coh/X}$ (resp.
$\Et_{\scoh/X}$) la sous-catégorie pleine de $\Et_{/X}$ formée des schémas étales 
de présentation finie sur $X$ (resp. étales, séparés et de présentation finie sur $X$), 
munie de la topologie induite par celle de $\Et_{/X}$; ce sont des sites $\mU$-petits. 
Notons par $\star$ l'un des deux symboles ``$\coh$'' ou ``$\scoh$''. 
On rappelle que le foncteur de restriction de $X_\et$ dans le topos des faisceaux de $\mU$-ensembles sur 
$\Et_{\star/X}$ est une équivalence de catégories. On désigne par 
\begin{equation}\label{higgs2-introd13a}
\pi_\star\colon E_\star\rightarrow \Et_{\star/X}
\end{equation}
le site fibré déduit de $\pi$ par changement de base par le foncteur d'injection canonique 
\begin{equation}\label{higgs2-introd13b}
\Et_{\star/X}\rightarrow \Et_{/X}, 
\end{equation} 
et par $\Phi\colon E_\star\rightarrow E$ la projection canonique. On montre \eqref{higgs2-tcevg10} 
que si l'on munit $E_\star$ de la topologie co-évanescente définie par $\pi_\star$
et si l'on note $\tE_\star$ le topos des faisceaux de $\mU$-ensembles sur $E_\star$, 
le foncteur $\Phi$ induit par restriction une équivalence de catégories 
$\tE\stackrel{\sim}{\rightarrow} \tE_\star$. De plus, la topologie co-évanescente de $E_\star$ 
est induite par celle de $E$ au moyen du foncteur $\Phi$ \eqref{higgs2-tcevg11}. 

La sous-catégorie $E_\coh$ permet d'appliquer les résultats de cohérence 
établis précédemment au topos de Faltings. 
On montre \eqref{higgs2-tf4} par exemple que si $X$ et $Y$ sont cohérents, alors le topos $\tE$ est cohérent~; 
en particulier, il a suffisamment de points. Nous verrons plus loin l'intérêt d'introduire la sous-catégorie $E_\scoh$ \eqref{higgs2-introd16}.

\subsection{}\label{higgs2-introd14}
Conservons les hypothèses de \eqref{higgs2-introd12}, et munissons $\Et_{/X}$ de l'objet final $X$ et 
$E$ de l'objet final $(Y\rightarrow X)$. Les foncteurs $\alpha_{X!}$ et $\sigma^+$ introduits dans \ref{higgs2-introd7} 
sont explicitement définis par 
\begin{eqnarray}
\alpha_{X!}\colon \Et_{\rf/Y}\rightarrow E,&& V\mapsto (V\rightarrow X),\label{higgs2-introd14a}\\
\sigma^+\colon \Et_{/X}\rightarrow E,&& U\mapsto (U_Y\rightarrow U).\label{higgs2-introd14b}
\end{eqnarray}
Ceux-ci sont exacts à gauche et continus. Ils définissent donc deux morphismes de topos 
\begin{eqnarray}
\beta\colon \tE\rightarrow Y_\fet,\label{higgs2-introd14c}\\
\sigma\colon \tE\rightarrow X_\et.\label{higgs2-introd14d}
\end{eqnarray} 
D'autre part, le foncteur 
\begin{equation}\label{higgs2-introd14e}
\Psi^+\colon E\rightarrow \Et_{/Y},\ \ \ (V\rightarrow U)\mapsto V
\end{equation}
est continu et exact à gauche \eqref{higgs2-tf13}~; il définit donc un morphisme de topos 
\begin{equation}\label{higgs2-introd14f}
\Psi\colon Y_\et\rightarrow \tE.
\end{equation}
On a $f_\et=\sigma\Psi$. Pour tout objet $U$ de $\Et_{/X}$, on peut identifier le morphisme 
$\Psi_U$ défini dans \eqref{higgs2-introd9c}
au morphisme canonique $\rho_{U_Y}\colon (U_Y)_\et\rightarrow (U_Y)_\fet$ \eqref{higgs2-introd11a}; 
en particulier, on a $\beta\Psi=\rho_Y$. De l'isomorphisme $\Psi^*\beta^*=\rho^*_Y$, 
on déduit par adjonction un morphisme
\begin{equation}\label{higgs2-introd14g}
\beta^*\rightarrow\Psi_* \rho^*_Y. 
\end{equation} 

Supposons de plus $X$ quasi-séparé et $Y$ cohérent et étale-localement connexe ({\em i.e.}, pour tout 
morphisme étale $X'\rightarrow X$, toute composante connexe de $X'$ est un ensemble ouvert dans $X'$).
On montre \eqref{higgs2-tf15} que les morphismes d'adjonction $\id\rightarrow \beta_*\beta^*$ et 
$\beta^*\rightarrow \Psi_*\rho^*_Y$ sont des isomorphismes. 

\subsection{}\label{higgs2-introd15}
Conservons les hypothèses de \eqref{higgs2-introd12}. On note $\varpi\colon D\rightarrow \Et_{/X}$ 
le site fibré associé au foncteur image inverse 
$f^+\colon \Et_{/X}\rightarrow \Et_{/Y}$ défini dans \ref{higgs2-introd5}, et  
l'on munit $D$ de la topologie co-évanescente associée à $\varpi$. 
On obtient ainsi le site co-évanescent associé au foncteur $f^+$, dont le topos des faisceaux de $\mU$-ensembles 
est $X_\et\gtimes_{X_\et}Y_\et$. 
Tout objet de $E$ est naturellement un objet de $D$. On a donc un foncteur pleinement fidèle et exact à gauche 
$\rho^+\colon E\rightarrow D$, 
qui n'est autre que le foncteur portant le même nom défini dans le cadre plus général de \eqref{higgs2-introd10}. 
Celui-ci étant continu et exact à gauche, il définit un morphisme de topos 
\begin{equation}
\rho\colon X_\et\gtimes_{X_\et}Y_\et\rightarrow \tE.
\end{equation}
On renvoie à \ref{higgs2-tf21} pour plus de détails. 
On montre \eqref{higgs2-tf251} que si $X$ et $Y$ sont cohérents, 
alors la famille des points de $\tE$ images par $\rho$ des points de $X_\et\gtimes_{X_\et}Y_\et$ 
est conservative. On notera que la donnée d'un point de $X_\et\gtimes_{X_\et}Y_\et$ 
est équivalente à la donnée d'une paire de points géométriques $\ox$ de $X$ et $\oy$ de $Y$
et d'une flèche de spécialisation de $f(\oy)$ vers $\ox$, c'est-à-dire, d'un $X$-morphisme 
$\oy\rightarrow X_{(\ox)}$, où $X_{(\ox)}$ désigne le localisé strict de $X$ en $\ox$ (cf. \ref{higgs2-tf24}). 

\subsection{}\label{higgs2-introd16}
Conservons les hypothèses de \eqref{higgs2-introd12}, supposons de plus $X$ strictement local, de point fermé $x$. 
Pour tout $X$-schéma étale, séparé et de présentation finie $U$, 
on désigne par $U^\rf$ la somme disjointe des localisés stricts de $U$ en les points de $U_x$;
c'est un sous-schéma ouvert et fermé de $U$, qui est fini sur $X$ (cf. \ref{higgs2-tf26}). Considérons le site fibré 
\begin{equation}\label{higgs2-introd16a}
\pi_\scoh\colon E_\scoh\rightarrow \Et_{\scoh/X}
\end{equation}
défini dans \eqref{higgs2-introd13a}, et munissons $E_\scoh$ de la topologie co-évanescente associée à $\pi_\scoh$.  
Pour tout objet $(V\rightarrow U)$ de $E_{\scoh}$, $V\times_UU^\rf=V\times_{U_Y}U^\rf_Y$ 
est un revêtement étale de $Y$. On obtient ainsi un foncteur 
\begin{equation}\label{higgs2-introd16b}
\theta^+\colon E_{\scoh}\rightarrow \Et_{\rf/Y}, \ \ \ (V\rightarrow U)\mapsto V\times_UU^\rf. 
\end{equation}
Celui-ci est continu et exact à gauche \eqref{higgs2-tf261}. Il définit donc un morphisme de topos 
\begin{equation}\label{higgs2-introd16c}
\theta\colon Y_\fet\rightarrow \tE, 
\end{equation}
que l'on étudie suivant l'approche introduite dans \eqref{higgs2-introd9} (cf. \ref{higgs2-tf262}). 
On a un isomorphisme canonique 
$\beta\theta \stackrel{\sim}{\rightarrow} \id_{Y_\fet}$, qui induit un morphisme de changement de base 
\begin{equation}\label{higgs2-introd16e}
\beta_*\rightarrow \theta^*.
\end{equation}
On montre \eqref{higgs2-tf265} que c'est un isomorphisme~; en particulier, le foncteur $\beta_*$ est exact . 
On en déduit \eqref{higgs2-tf266} que pour tout faisceau $F$ de $\tE$, l'application canonique 
\begin{equation}
\Gamma(\tE,F)\rightarrow \Gamma(Y_\fet, \theta^*F)
\end{equation}
est bijective, et que pour tout faisceau abélien $F$ de $\tE$, l'application canonique 
\begin{equation}
\rH^i(\tE,F)\rightarrow \rH^i(Y_\fet, \theta^*F)
\end{equation}
est bijective pour tout $i\geq 0$.

\subsection{}\label{higgs2-introd17}
Conservons les hypothèses de \eqref{higgs2-introd12} ; supposons, de plus, que $f$ soit cohérent. 
Soient $\ox$ un point géométrique de $X$, $\uX$ le localisé strict de $X$ en $\ox$, 
$\uY=Y\times_X\uX$, $\uf\colon \uY\rightarrow \uX$ la projection canonique. 
On désigne par $\uE$ le site de Faltings associé à $\uf$,
par $\tuE$ le topos des faisceaux de $\mU$-ensembles sur $\uE$ et par
\begin{equation}\label{higgs2-introd17a}
\theta\colon \uY_\fet\rightarrow \tuE
\end{equation}
le morphisme de topos défini dans \eqref{higgs2-introd16c}. Le foncteur  
\begin{equation}\label{higgs2-introd17b}
\Phi^+\colon E\rightarrow \uE, \ \ \ (V\rightarrow U)\mapsto (V\times_Y\uY\rightarrow U\times_X\uX)
\end{equation}
étant continu et exact à gauche (cf. \ref{higgs2-tf19}), il définit un morphisme de topos 
\begin{equation}\label{higgs2-introd17c}
\Phi\colon \tuE\rightarrow \tE.
\end{equation}
On montre \eqref{higgs2-tf268} que pour tout faisceau $F$ de $\tE$, on a un isomorphisme canonique fonctoriel 
\begin{equation}\label{higgs2-introd17d}
\sigma_*(F)_\ox\stackrel{\sim}{\rightarrow}\Gamma(\uY_\fet,\theta^*(\Phi^*F));
\end{equation}
et pour tout faisceau abélien $F$ de $\tE$ et tout entier $i\geq 0$, on a un isomorphisme canonique fonctoriel 
\begin{equation}\label{higgs2-introd17e}
\rR^i\sigma_*(F)_\ox\stackrel{\sim}{\rightarrow}\rH^i(\uY_\fet,\theta^*(\Phi^*F)).
\end{equation}
La preuve de ce résultat repose sur le calcul d'une limite projective de topos de Faltings \eqref{higgs2-lptf3} 
et ses conséquences cohomologiques \eqref{higgs2-lptf5}. 

On donne dans \ref{higgs2-tf28} une variante globale de l'isomorphisme \eqref{higgs2-introd17e}. 

\subsection{}
Signalons enfin que le topos co-évanescent généralisé est susceptible d'autres applications, dont des variantes rigide
et hensélienne du topos de Faltings. 
Nous tenons à remercier L. Illusie de nous avoir communiqué la lettre de Deligne \cite{deligne1},
ses notes \cite{illusie3} et son article \cite{illusie2} qui ont été la source principale d'inspiration de ce travail. 
Après que nous lui ayons communiqué une première 
version de ce travail, O. Gabber nous a transmis une copie d'un courrier électronique qu'il avait adressé à L. Illusie 
en 2006  dans lequel il définit la topologie co-évanescente pour un site fibré au-dessus d'un site et   esquisse 
des énoncés qui recoupent certains de nos résultats.  
Nous lui sommes très reconnaissants pour cet échange et les perspectives de développements ainsi offertes. 
Une grande partie de cet article a été écrite lors d'un séjour
du premier auteur (A.A.) à l'Université de Tokyo pendant l'automne 2010 et l'hiver 2011. Il
souhaite remercier cette Université pour son hospitalité.

\section{Notations et conventions}

{\it Tous les anneaux considérés dans cet article possèdent un élément unité~; les homomorphismes
d'anneaux sont toujours supposés transformer l'élément unité en l'élément unité.}

\subsection{}\label{higgs2-not1}
Dans tout cet article, on fixe un univers $\mU$ possédant un élément de cardinal infini. 
On appelle catégorie des $\mU$-ensembles et l'on note $\Ens$, 
la catégorie des ensembles qui se trouvent dans $\mU$. 
C'est un $\mU$-topos ponctuel que l'on note encore $\Pt$ (\cite{sga4} IV 2.2). 
Sauf mention explicite du contraire, il sera sous-entendu que les schémas 
envisagés dans cet article sont éléments de l'univers $\mU$.
On désigne par $\Sch$ la catégorie des schémas éléments de $\mU$.

\subsection{}\label{higgs2-not3}
Pour une catégorie $\cC$, nous notons $\ob(\cC)$ l'ensemble de ses objets,
$\cC^\circ$ la catégorie opposée, et pour $X,Y\in \ob(\cC)$, 
$\Hom_\cC(X,Y)$ (ou $\Hom(X,Y)$ lorsqu'il n'y a aucune ambiguïté) 
l'ensemble des morphismes de $X$ dans $Y$. 

Si $\cC$ et $\cC'$ sont deux catégories, nous désignons par $\Hom(\cC,\cC')$ 
l'ensemble des foncteurs de $\cC$ dans $\cC'$, et  
par $\bHom(\cC,\cC')$ la catégorie des foncteurs de $\cC$ dans $\cC'$. 

Soient $I$ une catégorie, $\cC$ et $\cC'$ deux catégories sur $I$ (\cite{sga1} VI 2). 
Nous notons $\Hom_{I}(\cC,\cC')$ l'ensemble des $I$-foncteurs de $\cC$ dans $\cC'$ 
et $\Hom_{\cart/I}(\cC,\cC')$ l'ensemble des foncteurs cartésiens (\cite{sga1} VI 5.2).
Nous désignons par $\bHom_{I}(\cC,\cC')$ la catégorie des $I$-foncteurs de $\cC$ dans $\cC'$ et 
par $\bHom_{\cart/I}(\cC,\cC')$ la sous-catégorie pleine formée des foncteurs cartésiens.

\subsection{}\label{higgs2-not6}
\index{100000207@$\cC_{/F}$}
Soit $\cC$ une catégorie. On désigne par $\hcC$ la catégorie des préfaisceaux 
de $\mU$-ensembles sur $\cC$, c'est-à-dire la catégorie des foncteurs
contravariants sur $\cC$ à valeurs dans $\Ens$ (\cite{sga4} I 1.2). 
Si $\cC$ est munie d'une topologie (\cite{sga4} II 1.1), on désigne par $\tcC$ le topos 
des faisceaux de $\mU$-ensembles sur $\cC$ (\cite{sga4} II 2.1). 

Pour $F$ un objet de $\hcC$, on note $\cC_{/F}$ la catégorie 
suivante (\cite{sga4} I 3.4.0). Les objets de $\cC_{/F}$
sont les couples formés d'un objet $X$ de $\cC$ 
et d'un morphisme $u$ de $X$ dans $F$. Si $(X,u)$ et $(Y,v)$ sont deux objets, 
un morphisme de $(X,u)$ vers $(Y,v)$ est un morphisme $g\colon X\rightarrow Y$
tel que $u=v\circ g$.

\subsection{}\label{higgs2-not4}
Si $(\cE,R)$ est un topos annelé, on note $\bMod(R)$ ou $\bMod(R,\cE)$ la catégorie des 
$R$-modules\footnote{Pour fixer les idées, nous prendrons les modules à gauche.}  de $\cE$,
$\bD(R)$ sa catégorie dérivée,  $\bD^-(R)$, $\bD^+(R)$ et $\bD^\rb(R)$ les sous-catégories
pleines de $\bD(R)$ formées des complexes à cohomologie bornée
supérieurement, inférieurement et des deux côtés, respectivement.

\subsection{}\label{higgs2-not5}\index{Coherent@Cohérent (quasi-compact et quasi-séparé)}
Suivant les conventions de (\cite{sga4} VI), nous utilisons l'adjectif {\em cohérent} 
comme synonyme de quasi-compact et quasi-séparé.

\section{Produits orientés de topos}

La notion de produit orienté de topos, rappelée ci-dessous, est due à Deligne. Elle a été étudiée par Gabber, 
Illusie, Laumon et Orgogozo \cite{illusie2,laumon,orgogozo}. 

\subsection{}\label{higgs2-topfl1}
Dans cette section, $X$, $Y$ et $S$ désignent des $\mU$-sites (\cite{sga4} II 3.0.2) 
dans lesquels les limites projectives finies sont représentables, 
\begin{equation}
f^+\colon S\rightarrow X\ \ \ {\rm et}\ \ \ g^+\colon S\rightarrow Y
\end{equation}
deux foncteurs continus et exacts à gauche. On désigne par 
$\tX$, $\tY$ et $\tS$ les topos des faisceaux de $\mU$-ensembles sur $X$, $Y$ et $S$, respectivement, par 
\begin{equation}
f\colon \tX\rightarrow \tS \ \ \ {\rm et}\ \ \ g\colon \tY\rightarrow \tS
\end{equation}
les morphismes de topos définis par $f^+$ et $g^+$ (\cite{sga4} IV 4.9.2), respectivement, 
et par $\varepsilon_X\colon X\rightarrow \tX$, 
$\varepsilon_Y\colon Y\rightarrow \tY$ et $\varepsilon_S\colon S\rightarrow \tS$  les foncteurs canoniques. 
Soient $e_X$, $e_Y$ et $e_S$ des objets finaux de $X$, $Y$ et $S$, respectivement, qui
existent par hypothèse. Comme les foncteurs canoniques sont exacts à gauche, 
$\varepsilon_X(e_X)$, $\varepsilon_Y(e_Y)$ et $\varepsilon_S(e_S)$ sont des objets finaux de 
$\tX$, $\tY$ et $\tS$, respectivement. 

On désigne par $C$ la catégorie des triplets 
\[
(W, U\rightarrow f^+(W), V\rightarrow g^+(W)),
\]
où $W$ est un objet de $S$, 
$U\rightarrow f^+(W)$ est un morphisme de $X$ et $V\rightarrow g^+(W)$ est un morphisme 
de $Y$; un tel objet sera noté $(U\rightarrow W\leftarrow V)$. 
Soient $(U\rightarrow W\leftarrow V)$ et $(U'\rightarrow W'\leftarrow V')$
deux objets de $C$. Un morphisme de $(U'\rightarrow W'\leftarrow V')$ dans $(U\rightarrow W\leftarrow V)$
est la donnée de trois morphismes $U\rightarrow U'$, $V\rightarrow V'$ et $W\rightarrow W'$ de
$X$, $Y$ et $S$, respectivement, tels que les diagrammes 
\begin{equation}
\xymatrix{
U'\ar[r]\ar[d]&{f^+(W')}\ar[d]\\
{U}\ar[r]&{f^+(W)}}\ \ \ \ \ \ \ \ 
\xymatrix{
V'\ar[d]\ar[r]&{g^+(W')}\ar[d]\\
V\ar[r]&{g^+(W)}}
\end{equation}
soient commutatifs. 

Il résulte aussitôt de la définition et du fait que les foncteurs $f^+$ 
et $g^+$ sont exacts à gauche que les limites projectives finies dans $C$ sont représentables.

On munit $C$ de la topologie engendrée par les recouvrements 
\[
\{(U_i\rightarrow W_i\leftarrow V_i)\rightarrow (U\rightarrow W\leftarrow V)\}_{i\in I}
\] 
des trois types suivants~:
\begin{itemize}
\item[(a)] $V_i=V$, $W_i=W$ pour tout $i\in I$, et $(U_i\rightarrow U)_{i\in I}$ est une famille couvrante.
\item[(b)] $U_i=U$, $W_i=W$ pour tout $i\in I$, et $(V_i\rightarrow V)_{i\in I}$ est une famille couvrante.
\item[(c)] $I=\{'\}$, $U'=U$ et le morphisme $V'\rightarrow V\times_{g^+(W)}g^+(W')$ est un isomorphisme
(il n'y a aucune condition sur le morphisme $W'\rightarrow W$). 
\end{itemize}
On notera que chacune de ces familles est stable par changement de base. 
On désigne par $\tC$ le topos des faisceaux de $\mU$-ensembles sur $C$. Si $F$ est un préfaisceau de $C$,
on note $F^a$ le faisceau associé.

\begin{lem}\label{higgs2-topfl2}
Pour qu'un préfaisceau $F$ sur $C$ soit un faisceau, il faut et il suffit que les conditions suivantes soient remplies~:
\begin{itemize}
\item[{\rm (i)}] Pour tout famille couvrante $(Z_i\rightarrow Z)_{i\in I}$ de $C$ du type {\rm (a)} ou {\rm (b)}, 
la suite 
\begin{equation}
F(Z)\rightarrow \prod_{i\in I}F(Z_i)\rightrightarrows \prod_{(i,j)\in I\times J}F(Z_i\times_ZZ_j)
\end{equation}
est exacte. 
\item[{\rm (ii)}] Pour tout recouvrement $(U\rightarrow W'\leftarrow V')\rightarrow (U\rightarrow W\leftarrow V)$ 
du type {\rm (c)}, l'application 
\begin{equation}\label{higgs2-topfl2a}
F(U\rightarrow W\leftarrow V) \rightarrow F(U\rightarrow W'\leftarrow V')
\end{equation}
est bijective.
\end{itemize}
\end{lem}

En effet, quitte à élargir l'univers $\mU$, on peut supposer les catégories $X$ et $Y$ petites (\cite{sga4} II 2.7(2)). 
Pour tout recouvrement $(U\rightarrow W'\leftarrow V')\rightarrow (U\rightarrow W\leftarrow V)$ 
du type (c), le morphisme diagonal 
\[
(U\rightarrow W'\leftarrow V')\rightarrow (U\rightarrow W'\times_WW'\leftarrow V'\times_VV')
\]
est un recouvrement du type (c) qui égalise les deux projections canoniques
\[
(U\rightarrow W'\times_WW'\leftarrow V'\times_VV')\rightrightarrows 
(U\rightarrow W'\leftarrow V').
\]
La proposition résulte donc de (\cite{sga4} II 2.3, I 3.5 et I 2.12).

\begin{rema}\label{higgs2-topfl3}
Il résulte de \ref{higgs2-topfl2} que le morphisme de faisceaux associé à un recouvrement du type (c)
est un isomorphisme; en particulier, la topologie de $C$ n'est pas toujours moins 
fine que la topologie canonique.
\end{rema} 

\subsection{}\label{higgs2-topfl4}\index{Projections canoniques d'un produit oriente de topos@Projections canoniques d'un produit orienté de topos}
Les foncteurs 
\begin{eqnarray}
\rp_1^+\colon X&\rightarrow& C,\ \ \ U\mapsto (U\rightarrow e_Z\leftarrow e_Y),\label{higgs2-topfl4a}\\
\rp_2^+\colon Y&\rightarrow& C,\ \ \ V\mapsto (e_X\rightarrow e_Z\leftarrow V),\label{higgs2-topfl4b}
\end{eqnarray}
sont exacts à gauche et continus (\cite{sga4} III 1.6). Ils définissent donc deux morphismes de topos (\cite{sga4} IV 4.9.2)
\begin{eqnarray}
\rp_1\colon \tC&\rightarrow& \tX,\\
\rp_2\colon \tC&\rightarrow& \tY.
\end{eqnarray}
Par ailleurs, on a un $2$-morphisme 
\begin{equation}\label{higgs2-topfl4e}
\tau\colon g \rp_2\rightarrow f \rp_1,
\end{equation}
donné par le morphisme de foncteurs $(g \rp_{2})_*\rightarrow (f \rp_{1})_*$ suivant~: 
pour tout faisceau $F$ sur $C$ et tout $W\in \ob(S)$, 
\begin{equation}\label{higgs2-topfl4f}
g_*(\rp_{2*}(F))(W)\rightarrow f_*(\rp_{1*}(F))(W)
\end{equation}
est l'application composée 
\[
F(e_X\rightarrow e_Z\leftarrow g^+(W))\rightarrow F(f^+(W)\rightarrow W\leftarrow g^+(W))
\rightarrow F(f^+(W)\rightarrow e_Z\leftarrow e_Y)
\]
où la première flèche est le morphisme canonique et la seconde flèche est l'inverse de l'isomorphisme \eqref{higgs2-topfl2a}.

\begin{rema}\label{higgs2-topfl41}
Pour tout $W\in \ob(S)$, le morphisme $\tau\colon (f\rp_1)^*(W)\rightarrow (g\rp_2)^*(W)$ \eqref{higgs2-topfl4e} est le morphisme
\begin{equation}\label{higgs2-topfl41a}
(\rp_1^+(f^+W))^a\rightarrow (\rp_2^+(g^+W))^a
\end{equation}
composé de 
\[
(f^+W\rightarrow e_Z\leftarrow e_Y)^a\rightarrow (f^+W\rightarrow W\leftarrow g^+W)^a \rightarrow 
(e_X\rightarrow e_Z\leftarrow g^+W)^a,
\]
où la première flèche est l'inverse de l'isomorphisme canonique \eqref{higgs2-topfl3} et la seconde flèche est le morphisme
canonique. 
\end{rema}

\begin{lem}\label{higgs2-topfl5}
Pour tout objet $Z=(U\rightarrow W\leftarrow V)$ de $C$, on a un diagramme cartésien  
\begin{equation}\label{higgs2-topfl5a}
\xymatrix{
Z^a\ar[r]\ar[d]&{\rp_2^*(V)}\ar[d]\\
{\rp_1^*(U)}\ar[r]^-(0.5)i&{(g\rp_2)^*(W)}}
\end{equation}
où les flèches non libellées sont les morphismes canoniques et $i$ est le composé du morphisme  
$\rp_1^*(U)\rightarrow (f\rp_1)^*(W)$ et du morphisme 
$\tau\colon (f\rp_1)^*(W)\rightarrow (g\rp_2)^*(W)$ \eqref{higgs2-topfl4e}.
\end{lem}

On a un diagramme canonique commutatif à carrés cartésiens de $C$
\begin{equation}
\xymatrix{
Z\ar[d]\ar[r]&{(f^+W\rightarrow W\leftarrow V)}\ar[r]\ar[d]&{\rp_2^+V}\ar[d]\\
{(U\rightarrow W\leftarrow g^+W)}\ar[r]\ar[d]&{(f^+W\rightarrow W\leftarrow g^+W)}\ar[d]_-(0.5){u}\ar[r]^-(0.5){v}&
{\rp_2^+(g^+W)}\\
{\rp_1^+U}\ar[r]&{\rp_1^+(f^+W)}&}
\end{equation}
D'autre part, $u^a$ est un isomorphisme \eqref{higgs2-topfl3} et $\tau\colon (f\rp_1)^*W\rightarrow (g\rp_2)^*W$ est égal à 
$v^a\circ (u^a)^{-1}$ \eqref{higgs2-topfl41}. Comme le foncteur canonique $C\rightarrow \tC$ est exact à gauche,
on en déduit que le diagramme \eqref{higgs2-topfl5a} est cartésien.

\begin{teo}\label{higgs2-topfl6}
Soient $T$ un $\mU$-topos, $a\colon T\rightarrow \tX$, $b\colon T\rightarrow \tY$ 
deux morphismes de topos, $t\colon gb\rightarrow fa$ un $2$-morphisme. Alors il existe un triplet 
\[
(h\colon T\rightarrow \tC, \alpha\colon \rp_1h\stackrel{\sim}{\rightarrow}a,  
\beta\colon \rp_2h\stackrel{\sim}{\rightarrow}b),
\] 
unique à isomorphisme unique près, formé d'un morphisme de topos $h$ et de deux isomorphismes 
de morphismes de topos $\alpha$ et $\beta$, tel que le diagramme 
\begin{equation}\label{higgs2-topfl6a}
\xymatrix{
{g\rp_2h}\ar[r]^{\tau*h}\ar[d]_{g*\beta}&{f\rp_1h}\ar[d]^{f*\alpha}\\
{gb}\ar[r]^t&{fa}}
\end{equation}
soit commutatif.
\end{teo}

L'unicité de $(h,\alpha,\beta)$ est claire. En effet, d'après \ref{higgs2-topfl5},  
la ``restriction'' $h^+\colon C\rightarrow T$ du foncteur $h^*$ à $C$
est nécessairement donnée, pour tout objet $Z=(U\rightarrow W\leftarrow V)$ de $C$, par 
\begin{equation}\label{higgs2-topfl6b}
h^+(Z)=a^*(U)\times_{(gb)^*(W)}b^*(V),
\end{equation}
où le morphisme $a^*(U)\rightarrow (gb)^*(W)$ est le composé 
\[
a^*(U)\rightarrow (fa)^*(W)\stackrel{t}{\rightarrow} (gb)^*(W). 
\]
Montrons que le foncteur $h^+$ ainsi défini est un morphisme de sites 
et que le morphisme de topos associé répond à la question.  
Le foncteur $h^+$ est clairement exact à gauche et transforme 
les familles couvrantes de $C$ du type (a) ou (b)  en familles couvrantes de $T$. 
D'autre part, si 
\[
(U\rightarrow W'\leftarrow V')\rightarrow (U\rightarrow W\leftarrow V)
\] 
est un recouvrement de $C$ du type (c), le carré 
\begin{equation}\label{higgs2-topfl6c}
\xymatrix{
{b^*V'}\ar[r]\ar[d]&{b^*V}\ar[d]\\
{b^*(g^+W')}\ar[r]&{b^*(g^+W)}}
\end{equation}
est cartésien, et par suite le morphisme 
\begin{equation}\label{higgs2-topfl6d}
h^+((U\rightarrow W'\leftarrow V'))\rightarrow 
h^+((U\rightarrow W\leftarrow V))
\end{equation}
est un isomorphisme. Donc le foncteur $h^+$ est continu en vertu de \ref{higgs2-topfl2}.
On en déduit que $h^+$ est un morphisme de sites  (\cite{sga4} IV 4.9.4); 
il définit donc un morphisme de topos $h\colon T\rightarrow \tC$. 

On a des isomorphismes canoniques $\alpha\colon a^*\stackrel{\sim}{\rightarrow}h^*\rp_1^*$
et $\beta\colon b^*\stackrel{\sim}{\rightarrow}h^*\rp_2^*$ dont les ``restrictions'' à $X$ et $Y$, respectivement,
sont les isomorphismes tautologiques. Pour vérifier que le diagramme \eqref{higgs2-topfl6a} est commutatif, 
il suffit de montrer que sa ``restriction'' à $S$ l'est. Pour tout $W\in \ob(S)$, considérons le diagramme  
\begin{equation}\label{higgs2-topfl6e}
\xymatrix{
{h^+((f^+W\rightarrow W\leftarrow g^+W))}\ar[r]^-(0.5)u\ar@{=}[d]&{(fa)^*(W)}\ar[r]^{t}\ar[d]^{\alpha(f^*W)}&
{(gb)^*(W)}\ar[d]^{\beta(g^*W)}\\
{h^+((f^+W\rightarrow W\leftarrow g^+W))}\ar[r]^-(0.5)v&{h^*((f\rp_1)^*W)}\ar[r]^{h^*(\tau)}&{h^*((g\rp_2)^*W)}}
\end{equation}
où $u$ est la projection déduite de la formule \eqref{higgs2-topfl6b} et $v$ est le morphisme canonique. Par définition, 
on a $v=\alpha(f^*W) \circ u$. D'autre part, $u'=t\circ u$ est la projection déduite de la formule \eqref{higgs2-topfl6b}, 
et $v'=h^*(\tau)\circ v$ est le morphisme canonique. Par définition, on a $v'=\beta(g^*W) \circ u'$.
Comme $u$ et $v$ sont des isomorphismes, on en déduit que le carré de droite dans \eqref{higgs2-topfl6e} 
est commutatif. Donc le diagramme \eqref{higgs2-topfl6a} est commutatif.

\subsection{}\label{higgs2-topfl10}
Soient $X'$, $Y'$, $S'$ trois $\mU$-sites dans lesquels les limites projectives finies sont représentables, 
$f'^+\colon S'\rightarrow X'$, $g'^+\colon S'\rightarrow Y'$ deux foncteurs continus et exacts à gauche. On note
$\tX'$, $\tY'$ et $\tS'$ les topos de faisceaux de $\mU$-ensembles sur $X'$, $Y'$ et $S'$, respectivement, et
\begin{equation}
f'\colon \tX'\rightarrow \tS' \ \ \ {\rm et}\ \ \ g'\colon \tY'\rightarrow \tS'
\end{equation}
les morphismes de topos définis par $f'^+$ et $g'^+$. 
On désigne par $C'$ le site associé aux foncteurs $(f'^+,g'^+)$ défini dans \eqref{higgs2-topfl1},
par $\tC'$ le topos des faisceaux de $\mU$-ensembles sur $C'$, par 
$\rp'_1\colon \tC'\rightarrow \tX'$ et $\rp'_2\colon \tC'\rightarrow \tY'$ les projections canoniques 
et par $\tau'\colon g'\rp'_2\rightarrow f'\rp'_1$ le $2$-morphisme canonique \eqref{higgs2-topfl4}. 
Considérons un diagramme de morphismes de topos 
\begin{equation}\label{higgs2-topfl10a}
\xymatrix{
{\tX'}\ar[d]_u\ar[r]^-(0.5){f'}&{\tS'}\ar[d]^w&{\tY'}\ar[d]^v\ar[l]_-(0.5){g'}\\
{\tX}\ar[r]^-(0.5){f}&{\tS}&{\tY}\ar[l]_-(0.5){g}}
\end{equation}
et deux $2$-morphismes 
\begin{equation}\label{higgs2-topfl10b}
a\colon wf'\rightarrow fu\ \ \ {\rm et}\ \ \
b\colon gv\rightarrow wg'.
\end{equation}
D'après \ref{higgs2-topfl6}, les morphismes de topos $u\rp'_1\colon \tC'\rightarrow \tX$ et $v\rp'_2\colon \tC'\rightarrow \tY$ 
et le $2$-morphisme composé $t$ 
\begin{equation}\label{higgs2-topfl10c}
\xymatrix{
{gv\rp'_2}\ar[r]^-(0.4){b*\rp'_2}&{wg'\rp'_2}\ar[r]^-(0.5){h*\tau'}&{wf'\rp'_1}\ar[r]^-(0.4){a*\rp'_1}&{fu\rp'_1}},
\end{equation}
définissent un morphisme de topos 
\begin{equation}\label{higgs2-topfl10d}
h\colon \tC'\rightarrow \tC
\end{equation}
et des $2$-isomorphismes $\alpha\colon \rp_1 h\stackrel{\sim}{\rightarrow} u\rp'_1$ et 
$\beta\colon \rp_2 h\stackrel{\sim}{\rightarrow} v\rp'_2$ rendant commutatif le diagramme 
\begin{equation}\label{higgs2-topfl10e}
\xymatrix{
{g\rp_2h}\ar[r]^{\tau*h}\ar[d]_{g*\beta}&{f\rp_1h}\ar[d]^{f*\alpha}\\
{gv\rp'_2}\ar[r]^t&{fu\rp'_1}}
\end{equation}

\begin{cor}\label{higgs2-topfl11}
Sous les hypothèses de \eqref{higgs2-topfl10}, si $u$, $v$ et $w$ sont des équivalences de topos et si $a$ et $b$
sont des $2$-isomorphismes, alors $h$ est une équivalence de topos. 
\end{cor}

Cela résulte de \ref{higgs2-topfl6}. 

\vspace{2mm}

Il résulte de \ref{higgs2-topfl11} que le topos $\tC$ ne dépend que du couple de morphismes de topos
$(f,g)$, à équivalence canonique près. 
Ceci justifie la terminologie et la notation suivantes~: 

\begin{defi}\label{higgs2-topfl7}\index{Produit oriente de topos@Produit orienté de topos}\index{100000310@$\tX\gtimes_\tS\tY$ (produit orienté de topos)}
Le topos $\tC$ est appelé {\em produit orienté}  de $\tX$ et $\tY$ au-dessus de $\tS$, 
et noté $\tX\gtimes_\tS\tY$.
\end{defi}

Sous  les hypothèses de \eqref{higgs2-topfl10}, on note le morphisme $h$ \eqref{higgs2-topfl10d} par 
\begin{equation}\label{higgs2-topfl7a}
u\gtimes_wv\colon \tX'\gtimes_{\tS'}\tY'\rightarrow \tX\gtimes_\tS\tY.
\end{equation}

\begin{cor}\label{higgs2-topfl8}
La donnée d'un point de $\tX\gtimes_\tS\tY$ est équivalente à 
la donnée d'une paire de points $x\colon \Pt\rightarrow \tX$ et $y\colon \Pt\rightarrow \tY$ 
et d'un $2$-morphisme $u\colon gy\rightarrow fx$. 
\end{cor}

Cela résulte de \ref{higgs2-topfl6}. 

\begin{defi}[\cite{illusie2} 4.1]\label{higgs2-topfl14}\index{Topos!evanescent@évanescent}\index{Topos!co-evanescent@co-évanescent}
Le topos $\tX\gtimes_\tS\tS$ est appelé topos {\em évanescent} de $f$, 
et le topos $\tS\gtimes_\tS\tY$ est appelé topos {\em co-évanescent} de $g$. 
\end{defi}

Nous donnerons dans \eqref{higgs2-co-ev101} une description plus simple de $\tS\gtimes_\tS\tY$ due à Deligne (\cite{illusie2} 4.6).

\subsection{}\label{higgs2-topfl9}
Considérons $\tX$, $\tY$ et $\tS$ comme des $\mU$-sites munis des topologies canoniques~; 
les $\mU$-topos associés s'identifient alors à $\tX$, $\tY$ et $\tS$, respectivement (\cite{sga4} IV 1.2). 
Par ailleurs, les foncteurs $f^*\colon \tS\rightarrow \tX$ et $g^*\colon \tS\rightarrow \tY$ sont clairement 
continus et exacts à gauche. On peut donc considérer le site $C^\dagger$ associé à $(f^*,g^*)$ défini 
dans \eqref{higgs2-topfl1}.
Notons provisoirement $\tC^\dagger$ le topos des faisceaux de $\mU$-ensembles sur $C^\dagger$,
$\pi_1^+\colon \tX\rightarrow C^\dagger$ et $\pi_2^+\colon \tY\rightarrow C^\dagger$ 
les foncteurs définis dans \eqref{higgs2-topfl4a} et \eqref{higgs2-topfl4b},  
\begin{eqnarray}
\pi_1\colon \tC^\dagger&\rightarrow& \tX,\\
\pi_2\colon \tC^\dagger&\rightarrow& \tY,
\end{eqnarray}
les morphismes de topos associés, et $\nu\colon g \pi_2\rightarrow f \pi_1$
le $2$-morphisme défini dans \eqref{higgs2-topfl4e}.

Les foncteurs canoniques $\varepsilon_X$, $\varepsilon_Y$ et $\varepsilon_S$
induisent un foncteur exact à gauche
\begin{equation}\label{higgs2-topfl9a}
\varphi^+ \colon C\rightarrow C^\dagger.
\end{equation}
Celui-ci transforme les recouvrements de $C$ du type (a) (resp. (b), resp. (c)) 
en des recouvrements de $C^\dagger$ du même type. Il résulte alors de \ref{higgs2-topfl2} 
que pour tout faisceau $F$ sur $C^\dagger$, $F\circ \varphi^+$ est un faisceau sur $C$. Par suite, $\varphi^+$
est continu. Il définit donc un morphisme de topos
\begin{equation}\label{higgs2-topfl9b}
\varphi \colon \tC^\dagger\rightarrow \tC.
\end{equation}
On a des isomorphismes canoniques 
\[
\varphi^+\circ \rp_1^+\stackrel{\sim}{\rightarrow}\pi_1^{+}\circ \varepsilon_X\ \ \ {\rm et} \ \ \ 
\varphi^+\circ \rp_2^+\stackrel{\sim}{\rightarrow}\pi_2^{+}\circ \varepsilon_Y.
\]
On en déduit des isomorphismes 
\begin{equation}
\pi_1\stackrel{\sim}{\rightarrow} \rp_1\varphi \ \ \ {\rm et} \ \ \
\pi_2\stackrel{\sim}{\rightarrow} \rp_2\varphi.
\end{equation}
De plus, il résulte aussitôt de la définition \eqref{higgs2-topfl4e} que $\tau*\varphi$ s'identifie à $\nu$.  
Par suite, $\varphi$ est une équivalence de topos en vertu de \ref{higgs2-topfl6}. En fait, $\varphi$ est le morphisme
de topos \eqref{higgs2-topfl10d} défini dans \eqref{higgs2-topfl10} en prenant pour $u$, $v$ et $w$ 
les morphismes identité de $\tX$, $\tY$ et $\tS$, respectivement. Nous identifions dans la suite $\tC^\dagger$
au topos $\tX\gtimes_{\tS}\tY$ au moyen de l'équivalence $\varphi$, le morphisme $\pi_1$ 
(resp. $\pi_2$) à $\rp_1$ (resp. $\rp_2$) et le $2$-morphisme $\nu$ à $\tau$.

\subsection{}\label{higgs2-topfl12}
Soit $(F\rightarrow H\leftarrow G)$ un objet de $C^\dagger$ \eqref{higgs2-topfl9}. 
Rappelons (\cite{sga4} IV 5.1) que la catégorie $\tX_{/F}$ est un $\mU$-topos, dit topos induit 
sur $F$ par $\tX$, et qu'on a un morphisme canonique $j_F\colon \tX_{/F}\rightarrow \tX$,
dit morphisme de localisation de $\tX$ en $F$. De même, on a des morphismes de localisation
$j_G\colon \tY_{/G}\rightarrow \tY$ et $j_H\colon \tS_{/H}\rightarrow \tS$. Notons $f'$ le morphisme 
composé 
\begin{equation}\label{higgs2-topfl12b}
\xymatrix{
{\tX_{/F}}\ar[r]&{\tX_{/f^*(H)}}\ar[r]^{f_{/H}}&{\tS_{/H}}}
\end{equation}
où la première flèche est le morphisme de localisation de $\tX_{/f^*(H)}$
en $F\rightarrow f^*(H)$ (\cite{sga4} IV 5.5), et la seconde flèche est le morphisme déduit de $f$ (\cite{sga4} IV 5.10).
On définit de même le morphisme $g'\colon \tY_{/G}\rightarrow \tS_{/H}$. Les carrés du diagramme 
\begin{equation}\label{higgs2-topfl12c}
\xymatrix{
{\tX_{/F}}\ar[d]_{j_F}\ar[r]^{f'}&{\tS_{/H}}\ar[d]^{j_H}&{\tY_{/G}}\ar[l]_{g'}\ar[d]^{j_G}\\
{\tX}\ar[r]^f&{\tS}&{\tY}\ar[l]_g}
\end{equation}
sont commutatifs à isomorphismes canoniques près. On désigne par $\tX_{/F}\gtimes_{\tS_{/H}}\tY_{/G}$
le produit orienté de $\tX_{/F}$ et $\tY_{/G}$ au-dessus de $\tS_{/H}$.
Le diagramme \eqref{higgs2-topfl12c} induit alors un morphisme canonique \eqref{higgs2-topfl7a}
\begin{equation}\label{higgs2-topfl12ab}
j_{F}\gtimes_{j_{H}}j_{G}\colon \tX_{/F}\gtimes_{\tS_{/H}}\tY_{/G}\rightarrow \tX\gtimes_{\tS}\tY.
\end{equation}

On désigne par $F\gtimes_HG$
le faisceau de $\tX\gtimes_{\tS}\tY$ associé à $(F\rightarrow H\leftarrow G)$ (cf. \ref{higgs2-topfl9}). 
D'après \ref{higgs2-topfl5}, on a un isomorphisme canonique
\begin{equation}\label{higgs2-topfl12a}
F\gtimes_HG\stackrel{\sim}{\rightarrow}\rp_1^*(F)\times_{(g\rp_2)^*(H)} \rp_2^*(G),
\end{equation}
où le morphisme $\rp_1^*(F)\rightarrow (g\rp_2)^*(H)$ est le composé du morphisme $\rp_1^*(F)\rightarrow (f\rp_1)^*(H)$ 
et du morphisme $\tau\colon (f\rp_1)^*(H)\rightarrow (g\rp_2)^*(H)$ \eqref{higgs2-topfl4e}. On note 
\begin{equation}\label{higgs2-topfl12d}
j_{F\gtimes_HG}\colon (\tX\gtimes_\tS\tY)_{/(F\gtimes_HG)} \rightarrow \tX\gtimes_\tS\tY
\end{equation}
le morphisme de localisation de $\tX\gtimes_\tS\tY$ en $F\gtimes_HG$. 

Par des constructions analogues à \eqref{higgs2-topfl12b}, les morphismes $\rp_1$  et $\rp_2$ induisent des morphismes 
\begin{eqnarray}
q_1\colon (\tX\gtimes_\tS\tY)_{/(F\gtimes_HG)} \rightarrow \tX_{/F},\\
q_2\colon (\tX\gtimes_\tS\tY)_{/(F\gtimes_HG)} \rightarrow \tY_{/G},
\end{eqnarray}
qui s'insèrent dans un diagramme à carrés commutatifs à isomorphismes canoniques près
\begin{equation}\label{higgs2-topfl12e}
\xymatrix{
{\tX_{/F}}\ar[d]_{j_F}&{(\tX\gtimes_\tS\tY)_{/(F\gtimes_HG)}}\ar[l]_-(0.4){q_1}
\ar[r]^-(0.4){q_2}\ar[d]_{j_{F\gtimes_HG}}&{\tY_{/G}}\ar[d]^{j_G}\\
{\tX}&{\tX\gtimes_\tS\tY}\ar[r]^-(0.5){\rp_2}\ar[l]_-(0.5){\rp_1}&{\tY}}
\end{equation}
Le $2$-morphisme $\tau\colon g\rp_2\rightarrow f\rp_1$ \eqref{higgs2-topfl4e} et le diagramme commutatif 
\begin{equation}\label{higgs2-topfl12ee}
\xymatrix{
{F\gtimes_HG}\ar[rr]\ar[d]&&{\rp_2^*G}\ar[d]\\
{\rp_1^*F}\ar[r]&{\rp_1^*(f^*H)}\ar[r]^\tau&{\rp_2^*(g^*H)}}
\end{equation}
induisent, pour tout $L\in \ob(\tS_{/H})$, un morphisme fonctoriel
\begin{equation}
\rp_1^*(f^*L)\times_{\rp_1^*(f^*H)}(F\gtimes_HG)\rightarrow \rp_2^*(g^*L)\times_{\rp_2^*(g^*H)}(F\gtimes_HG).
\end{equation}
Comme $\rp_1^*$ et $\rp_2^*$ sont exacts à gauche, on obtient un $2$-morphisme 
\begin{equation}\label{higgs2-topfl12f}
\tau'\colon g' q_2\rightarrow f'q_1
\end{equation}
tel que $j_H*\tau'=\tau*j_{F\gtimes_HG}$. 
D'après \ref{higgs2-topfl6}, $q_1$, $q_2$ et $\tau'$ définissent un morphisme 
\begin{equation}\label{higgs2-topfl12g}
m\colon (\tX\gtimes_\tS\tY)_{/(F\gtimes_HG)}\rightarrow \tX_{/F}\gtimes_{\tS_{/H}}\tY_{/G}.
\end{equation}

\begin{prop}\label{higgs2-topfl13}
Le morphisme $m$ est une équivalence de topos, et on a un isomorphisme canonique 
\begin{equation}\label{higgs2-topfl13a}
(j_{F}\gtimes_{j_H}j_G)\circ m\stackrel{\sim}{\rightarrow} j_{F\gtimes_HG}.
\end{equation}
\end{prop}

On notera d'abord que l'isomorphisme \eqref{higgs2-topfl13a} résulte de \ref{higgs2-topfl6}, compte tenu 
de \eqref{higgs2-topfl12e} et de la relation $j_H*\tau'=\tau*j_{F\gtimes_HG}$. 
Notons 
\begin{equation}
j_{(F\rightarrow H\leftarrow G)}\colon C^\dagger_{/(F\rightarrow H\leftarrow G)}\rightarrow C^\dagger
\end{equation} 
le foncteur canonique, $\cT$ la topologie de $C^\dagger$ \eqref{higgs2-topfl9} et $\cT_1$ la topologie
de $ C^\dagger_{/(F\rightarrow H\leftarrow G)}$ induite par $\cT$ via $j_{(F\rightarrow H\leftarrow G)}$.
Une famille $(L_i\rightarrow L)_{i\in I}$ 
de morphismes de $C^\dagger_{/(F\rightarrow H\leftarrow G)}$ est couvrante 
pour $\cT_1$ si et seulement si son image par $j_{(F\rightarrow H\leftarrow G)}$  est couvrante 
dans $C^\dagger$ pour $\cT$  (\cite{sga4} III 5.2(1)).
En vertu de  (\cite{sga4} III 5.4), $(\tX\gtimes_{\tS}\tY)_{/(F\gtimes_HG)}$ est canoniquement équivalent au topos des 
faisceaux de $\mU$-ensembles sur le site $(C^\dagger_{/(F\rightarrow H\leftarrow G)},\cT_1)$. 

Pour définir un site sous-jacent au topos $\tX_{/F}\gtimes_{\tS_{/H}}\tY_{/G}$, 
nous considérons  $\tX_{/F}$, $\tY_{/G}$ et $\tS_{/H}$ comme des $\mU$-sites munis des topologies canoniques 
(cf. \ref{higgs2-topfl9}). 
Pour tous $F'\in \ob(\tX_{/F})$ et $H'\in \ob(\tS_{/H})$, se donner un morphisme $F'\rightarrow f'^*(H')$ de $\tX_{/F}$
revient à se donner un morphisme $F'\rightarrow f^*(H')$ de $\tX$ au-dessus de $F\rightarrow f^*(H)$.  
Par suite, le site associé au couple de foncteurs $f'^*\colon \tS_{/H}\rightarrow \tX_{/F}$
et $g'^*\colon \tS_{/H}\rightarrow \tY_{/G}$ défini dans \eqref{higgs2-topfl1},
s'identifie canoniquement à la catégorie $C^\dagger_{/(F\rightarrow H\leftarrow G)}$,
munie d'une topologie $\cT_2$, qui est a priori moins fine que la topologie $\cT_1$. 
Le foncteur identique de $C^\dagger_{/(F\rightarrow H\leftarrow G)}$ définit alors un morphisme de sites 
\begin{equation}\label{higgs2-topfl13b}
\id \colon (C^\dagger_{/(F\rightarrow H\leftarrow G)}, \cT_1)\rightarrow (C^\dagger_{/(F\rightarrow H\leftarrow G)},\cT_2).
\end{equation}
Montrons que $m$ est le morphisme de topos associé à \eqref{higgs2-topfl13b}.

Il résulte de la preuve de \ref{higgs2-topfl6}, en particulier de \eqref{higgs2-topfl6b}, que la restriction 
\begin{equation}
m^+\colon C^\dagger_{/(F\rightarrow H\leftarrow G)}\rightarrow (\tX\gtimes_\tZ\tY)_{/(F\gtimes_HG)}
\end{equation}
du foncteur $m^*$ est donnée, pour tout objet $(F'\rightarrow H'\leftarrow G')$ de 
$C^\dagger_{/(F\rightarrow H\leftarrow G)}$, par 
\begin{equation}
m^+((F'\rightarrow H'\leftarrow G'))=q_1^*(F')\times_{(g'q_2)^*(H')}q_2^*(G'),
\end{equation}
où le morphisme $q_1^*(F')\rightarrow (g'q_2)^*(H')$ est le composé du morphisme 
$q_1^*(F')\rightarrow (f'q_1)^*(H')$ et du morphisme $\tau'\colon (f'q'_1)^*(H')\rightarrow (g'q_2)^*(H')$ \eqref{higgs2-topfl12f}.
On a des isomorphismes canoniques
\begin{eqnarray}
q_1^*(F')&\simeq&\rp_1^*(F')\times_{\rp_1^*(F)}(F\gtimes_HG),\\
q_2^*(G')&\simeq&\rp_2^*(G')\times_{\rp_2^*(G)}(F\gtimes_HG),\\
(g'q_2)^*(H')&\simeq&(g\rp_2)^*(H')\times_{(g\rp_2)^*(H)}(F\gtimes_HG).
\end{eqnarray}
De plus, compte tenu de \eqref{higgs2-topfl12ee}, le morphisme $q_1^*(F')\rightarrow (g'q_2)^*(H')$ provient du morphisme composé 
\[
\rp_1^*(F')\rightarrow (f\rp_1)^*(H')\stackrel{\tau}{\rightarrow}(g\rp_2)^*(H').
\] 
On en déduit un isomorphisme (\cite{sga4} I 2.5.0)
\begin{equation}
m^+((F'\rightarrow H'\leftarrow G'))\simeq \rp_1^*(F')\times_{(g\rp_2)^*(H')}\rp_2^*(G').
\end{equation}
Par suite, en vertu de \ref{higgs2-topfl5} et (\cite{sga4} III 5.5), le diagramme 
\begin{equation}\label{higgs2-topfl13c}
\xymatrix{
{\tX_{/F}\gtimes_{\tS_{/H}}\tY_{/G}}\ar[r]^-(0.5){m^*}&{(\tX\gtimes_\tS\tY)_{/(F\gtimes_HG)}}\\
{C^\dagger_{/(F\rightarrow H\leftarrow G)}}\ar[u]\ar@{=}[r]&{C^\dagger_{/(F\rightarrow H\leftarrow G)}}\ar[u]}
\end{equation}
où les flèches verticales sont les foncteurs canoniques, est commutatif. Donc $m$ est le morphisme de topos associé à \eqref{higgs2-topfl13b}. 

Pour montrer que $m$ est une équivalence de topos, il suffit de montrer que $\cT_1=\cT_2$, ou
que $\cT_1$ est moins fine que $\cT_2$, ou encore que le foncteur canonique 
\begin{equation}
j_{(F\rightarrow H\leftarrow G)}\colon (C^\dagger_{/(F\rightarrow H\leftarrow G)},\cT_2)\rightarrow (C^\dagger,\cT)
\end{equation}
est cocontinu (\cite{sga4} III 2.1).
Le foncteur $j_{(F\rightarrow H\leftarrow G)}$ est un adjoint à gauche du foncteur 
\begin{equation}
j_{(F\rightarrow H\leftarrow G)}^+\colon 
\begin{array}[t]{clcr}
(C^\dagger,\cT) &\rightarrow &(C^\dagger_{/(F\rightarrow H\leftarrow G)},\cT_2),\\
L&\mapsto& L\times (F\rightarrow H\leftarrow G).
\end{array}
\end{equation}
On rappelle que les limites projectives finies sont représentables dans $C^\dagger$. 
On montre comme plus haut que le diagramme 
\begin{equation}\label{higgs2-topfl13d}
\xymatrix{
{\tX\gtimes_\tS\tY}\ar[rr]^-(0.5){(j_F\gtimes_{j_H}j_G)^*}&&{\tX_{/F}\gtimes_{\tS_{/H}}\tY_{/G}}\\
{C^\dagger}\ar[u]\ar[rr]^-(0.5){j_{(F\rightarrow H\leftarrow G)}^+}&&{C^\dagger_{/(F\rightarrow H\leftarrow G)}}\ar[u]}
\end{equation}
où les flèches verticales sont les foncteurs canoniques, est commutatif. Par suite, $j_{(F\rightarrow H\leftarrow G)}^+$
est continu en vertu de (\cite{sga4} III 1.6), et donc $j_{(F\rightarrow H\leftarrow G)}$ est cocontinu d'après 
(\cite{sga4} III 2.5).

\section{Topos co-évanescents} \label{higgs2-co-ev}

\subsection{}\label{higgs2-co-ev1}
Dans cette section, $X$ et $Y$ désignent deux $\mU$-sites dans lesquels les 
limites projectives finies sont représentables, et $f^+\colon X\rightarrow Y$ un foncteur continu et exact à gauche. 
On désigne par $\tX$ et $\tY$ les topos des faisceaux de $\mU$-ensembles sur $X$ et $Y$, respectivement, 
par  $f\colon \tY\rightarrow \tX$ le morphisme de topos associé à $f^+$ et par 
$\varepsilon_X\colon X\rightarrow \tX$ et $\varepsilon_Y\colon Y\rightarrow \tY$ les foncteurs canoniques. 
Soient $e_X$ et $e_Y$ des objets finaux de $X$ et $Y$, respectivement, qui existent par hypothèse.
Comme les foncteurs canoniques sont exacts à gauche, 
$\varepsilon_X(e_X)$ et $\varepsilon_Y(e_Y)$ sont des objets finaux de $\tX$ et $\tY$, respectivement.

On désigne par $D$ la catégorie des paires $(U, V\rightarrow f^+(U))$, où $U$ est un objet de $X$ et  
$V\rightarrow f^+(U)$ est un morphisme de $Y$; un tel objet sera noté $(V\rightarrow U)$. 
Soient $(V\rightarrow U)$, $(V'\rightarrow U')$ deux objets de $D$. 
Un morphisme de $(V'\rightarrow U')$ dans $(V\rightarrow U)$
est la donnée de deux morphismes $V'\rightarrow V$ de $Y$ et $U'\rightarrow U$ de $X$, tels que le diagramme
\[
\xymatrix{
V'\ar[r]\ar[d]&{f^+(U')}\ar[d]\\
{V}\ar[r]&{f^+(U)}}
\] 
soit commutatif. Il résulte aussitôt de la définition et du fait que le foncteur $f^+$ 
est exact à gauche, que les limites projectives finies dans $D$ sont représentables.

On appelle topologie {\em co-évanescente} de $D$ la topologie engendrée par les recouvrements 
\[
\{(V_i\rightarrow U_i)\rightarrow (V\rightarrow U)\}_{i\in I}
\] 
des deux types suivants~:
\begin{itemize}
\item[$(\alpha)$] $U_i=U$ pour tout $i\in I$, et $(V_i\rightarrow V)_{i\in I}$ est une famille couvrante.
\item[$(\beta)$] $(U_i\rightarrow U)_{i\in I}$ est une famille couvrante, 
et pour tout $i\in I$, le morphisme canonique $V_i\rightarrow V\times_{f^+(U)}f^+(U_i)$ est un isomorphisme. 
\end{itemize}
On notera que chacune de ces familles est stable par changement de base. 
Le site ainsi défini est appelé site {\em co-évanescent} associé au foncteur $f^+$; c'est un $\mU$-site. 
On désigne par $\hD$ (resp. $\tD$) la catégorie des préfaisceaux (resp. le topos des faisceaux) de $\mU$-ensembles
sur $D$. On dit que $\tD$ est le topos {\em co-évanescent} associé au foncteur $f^+$. 
Nous montrerons dans \eqref{higgs2-co-ev101} que cette terminologie n'induit aucune confusion avec celle introduite 
dans \eqref{higgs2-topfl14}. Si $F$ est un préfaisceau sur $D$, on note $F^a$ le faisceau associé.

\begin{rema}\label{higgs2-co-ev2}
La topologie de $D$ est engendrée par les recouvrements 
\[
\{(V_{ij}\rightarrow U_i)\rightarrow (V\rightarrow U)\}_{(i,j)}
\] 
vérifiant les conditions suivantes~:
\begin{itemize}
\item[(i)] La famille $(U_i\rightarrow U)_i$ est couvrante.
\item[(ii)] Pour tout $i$, la famille $(V_{ij}\rightarrow V\times_{f^+(U)}f^+(U_i))_j$
est couvrante.
\end{itemize}
On prendra garde que la famille de ces recouvrements n'est pas en général stable par composition, et
donc ne forme pas une prétopologie.
\end{rema}

\subsection{}\label{higgs2-co-ev3}
On désigne par $\hY$ la catégorie des préfaisceaux de $\mU$-ensembles sur $Y$,
et par $\cQ$ la catégorie scindée des préfaisceaux de $\mU$-ensembles sur $Y$ (\cite{giraud2} I 2.6.1), 
c'est-à-dire, la catégorie fibrée sur $Y$ obtenue en associant à tout $V\in \ob(Y)$ 
la catégorie $(Y_{/V})^\wedge=\hY_{/V}$, et à tout morphisme $h\colon V'\rightarrow V$ de $Y$ le foncteur 
$h^*\colon \hY_{/V}\rightarrow \hY_{/V'}$ défini par composition avec le foncteur $Y_{/h}\colon Y_{/V'}\rightarrow Y_{/V}$. 
On notera que $h^*$ est aussi le changement de base dans $\hY$ par $h$. 
Comme les produits fibrés sont représentables dans $Y$, $h^*$ admet un adjoint à droite, à savoir 
le foncteur ``restriction de Weil'' $h_*\colon \hY_{/V'}\rightarrow \hY_{/V}$,
défini, pour tous $F\in \ob(\hY_{/V'})$ et $W\in \ob(Y_{/V})$, par  
\begin{equation}\label{higgs2-co-ev3c}
h_*(F)(W)=F(W\times_VV').
\end{equation}
On désigne par $\cQ^\vee$ la catégorie fibrée clivée et normalisée au-dessus de $Y^\circ$, 
obtenue en associant à tout $V\in \ob(Y)$ la catégorie $\hY_{/V}$, 
et à tout morphisme $h\colon V'\rightarrow V$ de $Y$ le foncteur $h_*$ (\cite{egr1} 1.1.2), et par 
\begin{eqnarray}
\cP&\rightarrow& X\label{higgs2-co-ev3d}\\
\cP^\vee&\rightarrow& X^\circ\label{higgs2-co-ev3e}
\end{eqnarray}
les catégories fibrées déduites de $\cQ$ et $\cQ^\vee$, respectivement, par changement de base 
par le foncteur $f^+\colon X\rightarrow Y$. 
D'après (\cite{sga1} VI 12; cf. aussi \cite{egr1} 1.1.2), on a une équivalence de catégories
\begin{eqnarray}\label{higgs2-co-ev3f}
\hD&\stackrel{\sim}{\rightarrow} &\bHom_{X^\circ}(X^\circ,\cP^\vee)\\
F&\mapsto&\{U\mapsto F_U\},\nonumber
\end{eqnarray}
définie, pour tout $(V\rightarrow U)\in \ob(D)$, par la relation  
\begin{equation}\label{higgs2-co-ev3g}
F_U(V)=F(V\rightarrow U).
\end{equation}
On identifiera dans la suite $F$ à la section $\{U\mapsto F_U\}$ qui lui est associée par cette équivalence.

\begin{prop}\label{higgs2-co-ev4}
Pour qu'un préfaisceau $F=\{U\mapsto F_U\}$ sur $D$ soit un faisceau, il faut et il suffit que les conditions suivantes 
soient remplies~:
\begin{itemize}
\item[{\rm (i)}] Pour tout $U\in \ob(X)$, $F_U$ est un faisceau sur $Y_{/f^+(U)}$. 
\item[{\rm (ii)}] Pour toute famille couvrante $(U_i\rightarrow U)_{i\in I}$ de $X$, si pour $(i,j)\in I^2$, 
on pose $U_{ij}=U_i\times_UU_j$ et on note $h_i\colon f^+(U_i)\rightarrow f^+(U)$ et 
$h_{ij}\colon f^+(U_{ij})\rightarrow f^+(U)$ les morphismes structuraux, alors la suite de morphismes de faisceaux sur 
$Y_{/f^+(U)}$
\begin{equation}\label{higgs2-co-ev4a}
F_U\rightarrow \prod_{i\in I}h_{i*}(F_{U_i})\rightrightarrows \prod_{(i,j)\in I^2}
h_{ij*}(F_{U_{ij}})
\end{equation}
est exacte.
\end{itemize}
\end{prop}

En effet, quitte à élargir l'univers $\mU$, on peut supposer la catégorie $X$ petite (\cite{sga4} II 2.7(2)). 
La proposition résulte alors de (\cite{sga4} II 2.3, I 3.5, I 2.12 et II 4.1(3)). 

\begin{rema}\label{higgs2-co-ev5}
La condition \ref{higgs2-co-ev4}(ii) revient à dire que, pour tout $(V\rightarrow U)\in \ob(D)$, 
si l'on pose $V_i=V\times_{f^+(U)}f^+(U_i)$ et $V_{ij}=V\times_{f^+(U)}f^+(U_{ij})$, la suite d'applications d'ensembles
\begin{equation}\label{higgs2-co-ev5a}
F_U(V)\rightarrow \prod_{i\in I}F_{U_i}(V_i)\rightrightarrows \prod_{(i,j)\in I^2}F_{U_{ij}}(V_{ij})
\end{equation}
est exacte. 
\end{rema}

\begin{remas}\label{higgs2-co-ev8}
(i)\ Pour tout objet $(V\rightarrow U)$ de $D$, le diagramme 
\begin{equation}\label{higgs2-co-ev8a}
\xymatrix{
{(V\rightarrow U)^a}\ar[r]\ar[d]&{(V\rightarrow e_X)^a}\ar[d]\\
{(f^+(U)\rightarrow U)^a}\ar[r]&{(f^+(U)\rightarrow e_X)^a}}
\end{equation}
est cartésien dans $\tD$. En effet, le foncteur canonique $D\rightarrow \tD$ est exact à gauche. 

(ii)\ Soient $W$ un objet de $X$, $F=\{U\mapsto F_U\}$ le préfaisceau sur $D$ défini 
par $(f^+(W)\rightarrow W)$. Pour tout $U\in \ob(X)$, $F_U$ est le préfaisceau constant
de valeur $\Hom_X(U,W)$ sur $Y_{/f^+(U)}$. En particulier, la topologie de $D$ n'est pas 
en général moins fine que la topologie canonique.

(iii)\ Soient $V$ un objet de $Y$, $F=\{U\mapsto F_U\}$ le préfaisceau sur $D$ défini 
par $(V\rightarrow e_X)$. Pour tout $U\in \ob(X)$, $F_U$ est le préfaisceau $V\times f^+(U)$ sur $Y_{/f^+(U)}$.
Si les topologies de $X$ et $Y$ sont moins fines que les topologies canoniques, 
$F$ est un faisceau sur $D$ en vertu de \ref{higgs2-co-ev4}.
\end{remas}

\subsection{}\label{higgs2-co-ev9}
Notons $C$ le site associé au couple de foncteurs $(\id_X,f^+)$ défini dans \eqref{higgs2-topfl1},
et considérons les foncteurs
\begin{eqnarray}
\iota^+\colon D&\rightarrow &C,\ \ \ (V\rightarrow U)\mapsto (U\rightarrow U\leftarrow V),\label{higgs2-co-ev9a}\\
\jmath^+\colon C&\rightarrow &D,\ \ \  (U\rightarrow W\leftarrow V)\mapsto (V\times_{f^+(W)}f^+(U)\rightarrow U).
\label{higgs2-co-ev9b}
\end{eqnarray}
Il est clair que $\iota^+$ est un adjoint à gauche de $\jmath^+$, 
que le morphisme d'adjonction $\id\rightarrow \jmath^+\circ \iota^+$ 
est un isomorphisme ({\em i.e.}, $\iota^+$ est pleinement fidèle), et que $\iota^+$ et $\jmath^+$ sont exacts à gauche.

\begin{prop}\label{higgs2-co-ev10}
{\rm (i)}\ Les foncteurs $\iota^+$ et $\jmath^+$ sont continus.

{\rm (ii)}\ La topologie de $D$ est induite par celle de $C$ au moyen du foncteur $\iota^+$.
\end{prop}
 
(i) Le foncteur $\iota^+$ transforme les familles couvrantes
de $D$ du type $(\alpha)$ en familles couvrantes de $C$ du type (b), 
et les familles couvrantes de $D$ du type $(\beta)$ en familles couvrantes de $C$~:
\begin{equation}\label{higgs2-co-ev10a}
\xymatrix{
{V_i}\ar[r]\ar[d]\ar@{}[rd]|{\Box}&{U_i}\ar[d]\\
V\ar[r]&U}
\ \ \ \mapsto \ \ \ 
\xymatrix{
{U_i}\ar[r]\ar[d]&{U_i}\ar[d]&{V_i}\ar[l]\ar[d]\ar@{}[ld]|{\Box}\\
{U_i}\ar[r]\ar[d]&U\ar[d]&V\ar[l]\ar[d]\\
U\ar[r]&U&V\ar[l]}
\end{equation}
Soient $G$ un préfaisceau sur $C$, $F=\{U\mapsto F_U\}=G\circ \iota^+$. 
Pour tout $(V\rightarrow U)\in \ob(D)$, on a 
\begin{equation}\label{higgs2-co-ev10b}
F_U(V)=G(U\rightarrow U\leftarrow V).
\end{equation}
Par suite, si $G$ est un faisceau sur $C$, 
$F$ est un faisceau sur $D$ en vertu de \ref{higgs2-topfl2}, \ref{higgs2-co-ev4} et \eqref{higgs2-co-ev10a}; donc $\iota^+$ est continu. 

Le foncteur $\jmath^+$ transforme les recouvrements de $C$ de type $(a)$ (resp. $(b)$)
en recouvrements de $D$ de type $(\beta)$ (resp. $(\alpha)$), et les recouvrements de $C$ de type $(c)$
en isomorphismes. Par suite, pour tout faisceau $F$ sur $D$, $F\circ \jmath^+$
est un faisceau sur $C$ en vertu de \ref{higgs2-topfl2}; donc $\jmath^+$ est continu. 

(ii) On sait que la topologie de $D$ est induite par la topologie canonique de $\tD$ (\cite{sga4} III 3.5),
autrement dit, la topologie de $D$ est la plus fine telle que tout $F\in \ob(\tD)$ soit un faisceau.
D'après (i), on peut considérer les foncteurs 
\begin{eqnarray}
\iota_s\colon \tC&\rightarrow& \tD,\ \ \  G\mapsto G\circ \iota^+,\\
\jmath_s\colon \tD&\rightarrow& \tC,\ \ \ F\mapsto F\circ \jmath^+.
\end{eqnarray}
L'isomorphisme d'adjonction $\id\stackrel{\sim}{\rightarrow} \jmath^+\circ \iota^+$ induit un isomorphisme 
$\iota_s\circ \jmath_s\stackrel{\sim}{\rightarrow} \id$. Le foncteur $\iota_s$ est donc essentiellement surjectif. 
Par suite, la topologie de $D$ est la plus fine telle que, pour tout $G\in \ob(\tC)$, 
$\iota_s(G)$ soit un faisceau sur $D$; d'où la proposition.

\subsection{}\label{higgs2-co-ev100}
Les foncteurs $\iota^+$ et $\jmath^+$ étant continus et exacts à gauche \eqref{higgs2-co-ev10}, ils définissent 
des morphismes de topos (\cite{sga4} IV 4.9.2)
\begin{eqnarray}
\iota\colon \tC&\rightarrow& \tD,\label{higgs2-co-ev100a}\\
\jmath\colon \tD&\rightarrow& \tC.\label{higgs2-co-ev100b}
\end{eqnarray}
Les morphismes d'adjonction $\id \rightarrow \jmath^+\circ \iota^+$ et $\iota^+ \circ \jmath^+\rightarrow \id$ 
induisent des morphismes $\iota_*\circ \jmath_*\rightarrow \id$ et $\id\rightarrow \jmath_*\circ \iota_*$
qui font de $\iota_*$ un adjoint à droite de $\jmath_*$.

\begin{prop}\label{higgs2-co-ev101}
Les morphismes d'adjonction $\iota_*\circ \jmath_*\rightarrow \id$ et $\id\rightarrow \jmath_*\circ \iota_*$ 
sont des isomorphismes. En particulier, $\iota$ \eqref{higgs2-co-ev100a} et $\jmath$ \eqref{higgs2-co-ev100b} 
sont des équivalences de topos quasi-inverses l'une de l'autre.
\end{prop}

En effet, comme le morphisme d'adjonction $\id \rightarrow \jmath^+\circ \iota^+$ est un isomorphisme, 
$\iota_*\circ \jmath_*\rightarrow \id$ est un isomorphisme.  D'autre part, le morphisme 
d'adjonction $\iota^+ \circ \jmath^+\rightarrow \id$ est défini, pour tout objet $(U\rightarrow W\leftarrow V)$ de $C$, 
par le morphisme canonique
\begin{equation}\label{higgs2-co-ev101a}
(U\rightarrow U\leftarrow V\times_{f^+(W)}f^+(U))\rightarrow (U\rightarrow W\leftarrow V),
\end{equation}
qui est un recouvrement de type (c). Comme le morphisme de faisceaux associé à \eqref{higgs2-co-ev101a} 
est un isomorphisme dans $\tC$ \eqref{higgs2-topfl3}, $\id\rightarrow \jmath_*\circ \iota_*$ est un isomorphisme.

\subsection{}\label{higgs2-co-ev11}
Les foncteurs 
\begin{eqnarray}
\rp_1^+\colon X&\rightarrow& D,\ \ \ U\mapsto (f^+(U)\rightarrow U),\label{higgs2-co-ev11a}\\
\rp_2^+\colon Y&\rightarrow& D,\ \ \ V\mapsto (V\rightarrow e_X),\label{higgs2-co-ev11b}
\end{eqnarray}
sont exacts à gauche et continus (\cite{sga4} III 1.6). Ils définissent donc deux morphismes de topos (\cite{sga4} IV 4.9.2)
\begin{eqnarray}
\rp_1\colon \tD&\rightarrow& \tX,\\
\rp_2\colon \tD&\rightarrow& \tY.
\end{eqnarray}
Pour tout $U\in \ob(X)$, le morphisme $(U\rightarrow U\leftarrow f^+(U))^a\rightarrow 
(U\rightarrow e_X\leftarrow e_Y)^a$ est un isomorphisme de $\tC$ \eqref{higgs2-topfl3}. 
Les morphismes $\rp_1\circ \iota$  et $\rp_2\circ \iota$ s'identifient donc aux morphismes 
$\rp_1\colon \tC\rightarrow \tX$ et $\rp_2\colon \tC\rightarrow \tY$ définis dans \eqref{higgs2-topfl4}, d'où la terminologie. 
Le $2$-morphisme \eqref{higgs2-topfl4e}
\begin{equation}\label{higgs2-co-ev11c}
\tau\colon f\rp_2\rightarrow \rp_1
\end{equation}
est alors défini par le morphisme de foncteurs $(f\rp_2)_*\rightarrow \rp_{1*}$ suivant~: pour tout faisceau $F$ sur $D$ et 
tout $U\in \ob(X)$, 
\begin{equation}\label{higgs2-co-ev11cc}
f_*(\rp_{2*}(F))(U)\rightarrow \rp_{1*}(F)(U)
\end{equation}
est l'application canonique
\[
F(f^+(U)\rightarrow e_X)\rightarrow F(f^+(U)\rightarrow U).
\]

Le $2$-morphisme $\tau$ induit un morphisme de changement de base 
\begin{equation}\label{higgs2-co-ev15a}
f_*\rightarrow \rp_{1*}\rp_2^*
\end{equation}
composé de 
\[
\xymatrix{
{f_*}\ar[r]&{f_*\rp_{2*}\rp_2^*}\ar[r]^-(0.4){\tau*\rp_2^*}&{\rp_{1*}\rp_2^*}},
\]
où le premier morphisme est déduit du morphisme d'adjonction. Pour tout anneau $\Lambda$, 
le morphisme \eqref{higgs2-co-ev15a} induit un morphisme de foncteurs de $\bD^+(\tY,\Lambda)$
dans $\bD^+(\tX,\Lambda)$   
\begin{equation}\label{higgs2-co-ev15b}
\rR f_*\rightarrow \rR\rp_{1*}\rp_2^*.
\end{equation}

\begin{prop}\label{higgs2-co-ev13}
{\rm (i)}\ Pour tout faisceau $F$ sur $X$, $\rp_1^*(F)$ est le faisceau associé au préfaisceau 
$\{U\mapsto F(U) \}$ sur $D$ \eqref{higgs2-co-ev3f}, où pour tout $U\in \ob(X)$, $F(U)$ est le préfaisceau
constant sur $Y_{/f^+(U)}$ de valeur $F(U)$. 

{\rm (ii)}\ Pour tout faisceau $F$ sur $Y$, $\rp_2^*(F)$ est le faisceau $\{U\mapsto F\times f^*(U)\}$.

{\rm (iii)}\ Pour tout faisceau $F$ sur $X$, le morphisme $\tau\colon \rp_1^* (F)\rightarrow (f\rp_2)^*(F)$ \eqref{higgs2-co-ev11c}
est induit par le morphisme de préfaisceaux sur $D$ défini, pour tout $U\in \ob(X)$, par le morphisme de 
préfaisceaux sur $f^+(U)$
\begin{equation}\label{higgs2-co-ev13a}
F(U)\rightarrow f^*(F\times U)
\end{equation}
donné, pour tout $(V\rightarrow U)\in \ob(D)$, par le composé 
\begin{equation}\label{higgs2-co-ev13b}
F(U)\rightarrow (f^*F)(f^+U)\rightarrow (f^*F)(V).
\end{equation}

{\rm (iv)}\ Le morphisme d'adjonction $\id\rightarrow \rp_{2*}\rp_2^*$ est un isomorphisme. 
\end{prop}

En effet, quitte à élargir $\mU$, on peut supposer que les catégories $X$ et $Y$ sont $\mU$-petites 
(\cite{sga4} II 3.6 et III 1.5). 

(i) D'après (\cite{sga4} I 5.1 et III 1.3), le faisceau $\rp_1^*(F)$ est le faisceau sur $D$ 
associé au préfaisceau $G$ défini pour $(V\rightarrow U)\in \ob(D)$ par 
\begin{equation}\label{higgs2-co-ev13c}
G(V\rightarrow U)=\underset{\underset{(P,u)\in I^\circ_{(V\rightarrow U)}}{\longrightarrow}}{\lim}\ F(P),
\end{equation}
où $I_{(V\rightarrow U)}$ est la catégorie des couples $(P,u)$ formés d'un objet $P$ de $X$ et 
d'un morphisme $u\colon (V\rightarrow U)\rightarrow (f^+(P)\rightarrow P)$ de $D$. 
Cette catégorie admet comme objet initial le couple formé de $U$ et du morphisme canonique
$(V\rightarrow U)\rightarrow (f^+(U)\rightarrow U)$. On a donc $G(V\rightarrow U)=F(U)$. 

(ii) Le faisceau $\rp_2^*(F)$ est le faisceau sur $D$ 
associé au préfaisceau $H$ défini pour $(V\rightarrow U)\in \ob(D)$ par 
\begin{equation}\label{higgs2-co-ev13d}
H(V\rightarrow U)=\underset{\underset{(Q,v)\in J^\circ_{(V\rightarrow U)}}{\longrightarrow}}{\lim}\ F(Q),
\end{equation}
où $J_{(V\rightarrow U)}$ est la catégorie des couples $(Q,v)$ formés d'un objet $Q$ de $Y$ et 
d'un morphisme $v\colon (V\rightarrow U)\rightarrow (Q\rightarrow e_X)$ de $D$. 
Cette catégorie admet comme objet initial le couple formé de $V$ et du morphisme canonique
$(V\rightarrow U)\rightarrow (V\rightarrow e_X)$. On a donc $H(V\rightarrow U)=F(V)$. Par suite, 
$H=\{U\mapsto F\times f^*(U)\}$, qui est en fait un faisceau  sur $D$
en vertu de \ref{higgs2-co-ev4}.

(iii) Notons $G$ le préfaisceau sur $D$ associé à $F$ défini dans \eqref{higgs2-co-ev13c} et 
$H$ le faisceau sur $D$ associé à $f^*(F)$ défini dans \eqref{higgs2-co-ev13d}. 
Pour tout objet $(V\rightarrow U)$ de $D$, on a un  foncteur 
\begin{equation}\label{higgs2-co-ev13f}
I_{(V\rightarrow U)}\rightarrow J_{(V\rightarrow U)}, 
\end{equation} 
défini par $(P,u)\mapsto (f^+(P),v)$, où $v$ est le morphisme composé 
\[
(V\rightarrow U)\stackrel{u}{\rightarrow}(f^+(P)\rightarrow  P) \rightarrow (f^+(P)\rightarrow e_X).
\] 
Alors l'application composée  
\[
\underset{\underset{(P,u)\in I^\circ_{(V\rightarrow U)}}{\longrightarrow}}{\lim}\ F(P)
\rightarrow \underset{\underset{(P,u)\in I^\circ_{(V\rightarrow U)}}{\longrightarrow}}{\lim}\ (f^*F)(f^+P)
\rightarrow \underset{\underset{(Q,v)\in J^\circ_{(V\rightarrow U)}}{\longrightarrow}}{\lim}\ (f^*F)(Q),
\]
où la première flèche est l'application canonique et la seconde flèche est induite par le foncteur \eqref{higgs2-co-ev13f}, 
est égale à l'application 
\begin{equation}\label{higgs2-co-ev13e}
G(V\rightarrow U)\rightarrow H(V\rightarrow U)
\end{equation} 
définie dans \eqref{higgs2-co-ev13b}. Par suite, le morphisme de faisceaux 
$\rp_1^*(F)\rightarrow (f\rp_2)^*(F)$ associé à \eqref{higgs2-co-ev13e} est l'adjoint du morphisme \eqref{higgs2-co-ev11cc}, d'où 
la proposition. 

(iv) En effet, pour tout faisceau $F$ sur $Y$ et tout $V\in \ob(Y)$, 
le morphisme d'adjonction $F(V)\rightarrow (\rp_{2}^*F)(V\rightarrow e_X)$ s'identifie au morphisme
identique de $F(V)$ en vertu de (ii).

\subsection{}\label{higgs2-co-ev14}
Le foncteur 
\begin{equation}
\Psi^+\colon D\rightarrow Y, \ \ \ (V\rightarrow U)\mapsto V
\end{equation}
est clairement exact à gauche. Pour tout faisceau $F$ sur $Y$, 
on a 
\begin{equation}
F\circ \Psi^+=\{U\mapsto F|f^+(U)\},
\end{equation}
où pour tout morphisme $g\colon U'\rightarrow U$ de $X$, si l'on pose $h= f^+(g)$, le morphisme de transition 
\begin{equation}
F|f^+(U)\rightarrow h_*(F|f^+(U'))
\end{equation}
est l'adjoint de l'isomorphisme canonique $h^*(F|f^+(U))\stackrel{\sim}{\rightarrow}F|f^+(U')$.
Il résulte de \ref{higgs2-co-ev4} que $F\circ \Psi^+$ est un faisceau sur $D$. Par suite, $\Psi^+$
est continu. Il définit donc un morphisme de topos
\begin{equation}\label{higgs2-co-ev14a}
\Psi\colon \tY\rightarrow \tD
\end{equation}
tel que $\rp_1\Psi=f$, $\rp_2\Psi=\id_\tY$ et $\tau*\Psi=\id_f$, où $\tau$ est le $2$-morphisme \eqref{higgs2-co-ev11c}.
Par suite, $\jmath \Psi $  est le morphisme 
défini par les morphismes $f\colon \tY\rightarrow \tX$ et $\id_\tY$ et le $2$-morphisme $\id_f$,
compte tenu de la propriété universelle des produits orientés \eqref{higgs2-topfl6}~: 
\begin{equation}\label{higgs2-co-ev14b}
\xymatrix{
&\tY\ar[dl]_f\ar[rd]^{\id_\tY}\ar[d]_{\jmath \Psi}&\\
{\tX}\ar[dr]&{\tX\gtimes_\tX\tY}\ar[r]^-(0.4){\rp_2}\ar[l]_-(0.4){\rp_1}&{\tY}\ar[ld]^f\\
&{\tX}&}
\end{equation}
Le morphisme $\Psi$ (ou $\jmath \Psi$) est appelé morphisme des {\em cycles co-proches}.
De la relation $\rp_{2*}\Psi_*=\id_{\tY}$, on obtient par adjonction un morphisme 
\begin{equation}\label{higgs2-co-ev14c}
\rp_2^*\rightarrow \Psi_*.
\end{equation}

\begin{prop}\label{higgs2-co-ev12}
Le morphisme $\rp_2^*\rightarrow \Psi_*$ \eqref{higgs2-co-ev14c} est un isomorphisme, en particulier, 
le foncteur $\Psi_*$ est exact. 
\end{prop}

En effet, pour tout faisceau $F$ sur $Y$ et tout $(V\rightarrow U)\in \ob(D)$, on a un diagramme commutatif 
\begin{equation}\label{higgs2-co-ev12c}
\xymatrix{
{\rp_2^*(F)(V\rightarrow e_X)}\ar[r]\ar[d]&{\Psi_*(F)(V\rightarrow e_X)}\ar[d]\\
{\rp_2^*(F)(V\rightarrow U)}\ar[r]&{\Psi_*(F)(V\rightarrow U)}}
\end{equation}
où les flèches horizontales sont les applications \eqref{higgs2-co-ev14c} et les flèches verticales sont les applications
canoniques. Ces derniers sont des isomorphismes en vertu de \ref{higgs2-co-ev13}(ii). 
D'autre part, la flèche horizontale supérieure est induite par le morphisme 
\begin{equation}\label{higgs2-co-ev12d}
\rp_{2*}\rp_2^*\rightarrow \rp_{2*}\Psi_*
\end{equation}
déduit de \eqref{higgs2-co-ev14c}. Le composé de \eqref{higgs2-co-ev12d} avec le morphisme d'adjonction $\id\rightarrow \rp_{2*}\rp_2^*$
est l'isomorphisme canonique $\id\stackrel{\sim}{\rightarrow}\rp_{2*}\Psi_*$. Il résulte alors de \ref{higgs2-co-ev13}(iv)
que \eqref{higgs2-co-ev12d} est un isomorphisme. Par suite, la flèche horizontale inférieure de \eqref{higgs2-co-ev12c} est un 
isomorphisme, d'où la proposition.

\begin{prop}\label{higgs2-co-ev16}
{\rm (i)}\ Pour tout faisceau d'ensembles $F$ sur $Y$, le morphisme \eqref{higgs2-co-ev15a}
\begin{equation}\label{higgs2-co-ev16a}
f_*(F)\rightarrow \rp_{1*}\rp_2^*(F)
\end{equation}
est un isomorphisme.

{\rm (ii)}\ Soit $\Lambda$ un anneau.
Pour tout complexe $F$ de $\bD^+(\tY,\Lambda)$, le morphisme \eqref{higgs2-co-ev15b}
\begin{equation}\label{higgs2-co-ev16b}
\rR f_*(F)\rightarrow \rR \rp_{1*}\rp_2^*(F)
\end{equation}
est un isomorphisme.
\end{prop}

(i) Considérons le diagramme commutatif 
\begin{equation}\label{higgs2-co-ev16c}
\xymatrix{
{f_*}\ar[r]^-(0.5)a\ar[rd]_-(0.5)d&{f_*\rp_{2*}\rp_2^*}\ar[r]^-(0.4){\tau*\rp_2^*}\ar[d]^b&{\rp_{1*}\rp_2^*}\ar[d]^-(0.5)c\\
&{f_*\rp_{2*}\Psi_*}\ar[r]^-(0.4){\tau*\Psi_*}&{\rp_{1*}\Psi_*}}
\end{equation}
où $a$ est induit par le morphisme d'adjonction et $b$ et $c$ sont induits par \eqref{higgs2-co-ev14c}. 
Comme $d=b\circ a$ est l'isomorphisme déduit de la relation $\rp_2\Psi=\id_{\tY}$, 
$(\tau*\Psi_*)\circ d$ s'identifie à l'isomorphisme déduit de la relation $\rp_1\Psi=f$ 
\eqref{higgs2-co-ev14}. D'autre part, $c$ est un isomorphisme en vertu de \ref{higgs2-co-ev12}, d'où la proposition.

(ii) Soit $\ttau\colon \rR f_* \rR \rp_{2*}\rightarrow \rR \rp_{1*}$ le morphisme 
induit par $\tau$ \eqref{higgs2-co-ev11c}. Comme $\Psi_*$ est exact en vertu de \ref{higgs2-co-ev12}, 
le diagramme \eqref{higgs2-co-ev16c} induit un diagramme commutatif 
\begin{equation}\label{higgs2-co-ev16d}
\xymatrix{
{\rR f_*}\ar[r]^-(0.5)\alpha\ar[rd]_-(0.5)\delta&{\rR f_* \rR \rp_{2*}\rp_2^*}\ar[r]^-(0.4){\ttau*\rp_2^*}\ar[d]^\beta&
{\rR \rp_{1*}\rp_2^*}\ar[d]^-(0.5)\gamma\\
&{\rR f_* \rR \rp_{2*}\Psi_*}\ar[r]^-(0.4){\ttau*\Psi_*}&{\rR \rp_{1*}\Psi_*}}
\end{equation}
D'autre part, $\Psi_*$ étant exact, $\delta$ est l'isomorphisme déduit de la relation $\rp_2\Psi=\id_{\tY}$,
et par suite $(\ttau*\Psi_*)\circ \delta$ s'identifie à l'isomorphisme déduit de la relation $\rp_1\Psi=f$ \eqref{higgs2-co-ev14}.
Comme $\gamma$ est un isomorphisme en vertu de \ref{higgs2-co-ev12}, la proposition s'ensuit. 

\begin{rema}\label{higgs2-co-ev18}
La proposition \ref{higgs2-co-ev16} et sa preuve sont dues à Gabber (\cite{illusie2} 4.9). C'est un cas particulier 
d'un théorème de changement de base pour les topos orientés (\cite{illusie2} 2.4), qui requiert par contre 
des hypothèses de cohérence plus restrictives.
\end{rema}

\subsection{}\label{higgs2-co-ev20}
Soit $(B\rightarrow A)$ un objet de $D$. On note $j_{A} \colon X_{/A}\rightarrow X$
et $j_B\colon Y_{/B}\rightarrow Y$ les foncteurs canoniques, et on munit $X_{/A}$
et $Y_{/B}$ des topologies induites par celles de $X$ et $Y$ via les foncteurs $j_A$ et $j_B$,
respectivement. Notons $(X_{/A})^\sim$ le topos des faisceaux de $\mU$-ensembles sur $X_{/A}$.
D'après (\cite{sga4} III 5.2), le foncteur $j_A$ est continu et cocontinu. 
Il induit donc une suite de trois foncteurs adjoints~:  
\begin{equation}\label{higgs2-co-ev20e}
j_{A!}\colon (X_{/A})^\sim \rightarrow \tX, \ \ \ j_A^*\colon  \tX \rightarrow (X_{/A})^\sim, 
\ \ \ j_{A*}\colon (X_{/A})^\sim \rightarrow \tX
\end{equation}
dans le sens que pour deux foncteurs consécutifs de la suite, celui
de droite est adjoint à droite de l'autre. Le foncteur $j_{A!}$ se factorise à travers une équivalence de catégories 
$(X_{/A})^\sim\stackrel{\sim}{\rightarrow} \tX_{/A^a}$, où $A^a=\varepsilon_X(A)$ (\cite{sga4} III 5.4), et le
couple $(j_A^*,j_{A*})$ définit un morphisme de topos $\tX_{/A^a}\rightarrow \tX$, dit morphisme 
de localisation de $\tX$ en $A^a$; et de même pour $j_B$.  

Les limites projectives finies sont représentables dans $X_{/A}$ et $Y_{/B}$. D'autre part, le foncteur 
\begin{equation}\label{higgs2-co-ev20d}
f'^+\colon  X_{/A}\rightarrow Y_{/B}, \ \ \ U\mapsto f^+(j_A(U))\times_{f^+(A)}B
\end{equation}
est exact à gauche et continu en vertu de (\cite{sga4} III 1.6 et 3.3). 
Le morphisme de topos 
\begin{equation}\label{higgs2-co-ev20c}
f'\colon \tY_{/B^a}\rightarrow \tX_{/A^a}
\end{equation}
associé à $f'^+$ s'identifie au morphisme composé 
\[
\xymatrix{
{\tY_{/B^a}}\ar[r]&{\tY_{/f^*(A)}}\ar[r]^-(0.5){f_{/A^a}}&{\tX_{/A^a}}},
\]
où la première flèche est le morphisme de localisation associé à $B^a\rightarrow f^*(A^a)$ (\cite{sga4} IV 5.5)
et la seconde flèche est le morphisme déduit de $f$ (\cite{sga4} IV 5.10).

On désigne par $D'$ (resp. $\tD'$) le site (resp. topos) co-évanescent associé au foncteur $f'^+$.
Les foncteurs $j_A$ et $j_B$ induisent un foncteur 
\begin{equation}\label{higgs2-co-ev20f}
j_{(B\rightarrow A)}\colon D'\rightarrow D,
\end{equation}
qui se factorise à travers une équivalence de catégories 
\begin{equation}\label{higgs2-co-ev20g}
n\colon D'\stackrel{\sim}{\rightarrow}D_{/(B\rightarrow A)}.
\end{equation}
 
\begin{prop}\label{higgs2-co-ev21}
Sous les hypothèses de \eqref{higgs2-co-ev20},
la topologie co-évanescente de $D'$ est induite par la topologie co-évanescente de $D$ au moyen du foncteur 
$j_{(B\rightarrow A)}$ \eqref{higgs2-co-ev20f}; en particulier, $n$ \eqref{higgs2-co-ev20g} induit une équivalence de topos 
\begin{equation}\label{higgs2-co-ev21a}
m\colon \tD_{/(B\rightarrow A)^a}\stackrel{\sim}{\rightarrow} \tD'.
\end{equation}
\end{prop}

Identifions les topos $\tD$ et $\tX\gtimes_{\tX}\tY$
et les topos $\tD'$ et $\tX_{/A^a}\gtimes_{\tX_{/A^a}}\tY_{/B^a}$ par les équivalences \eqref{higgs2-co-ev100a}. 
D'après \ref{higgs2-topfl13}, on a une équivalence de topos 
\begin{equation}\label{higgs2-co-ev21b}
m\colon (\tX\gtimes_\tX\tY)_{/(B\rightarrow A)^a}\stackrel{\sim}{\rightarrow} \tX_{/A^a}\gtimes_{\tX_{/A^a}}\tY_{/B^a}.
\end{equation}

A priori, la topologie co-évanescente de $D'$ est moins fine que la topologie 
induite par la topologie co-évanescente de $D$ au moyen du foncteur $j_{(B\rightarrow A)}$. 
Mais il résulte de la preuve de \ref{higgs2-topfl13}, en particulier de \eqref{higgs2-topfl13c}, que le diagramme 
\begin{equation}\label{higgs2-co-ev21c}
\xymatrix{
{\tX_{/A^a}\gtimes_{\tX_{/A^a}}\tY_{/B^a}}\ar[r]^-(0.5){m^*}&{(\tX\gtimes_{\tX}\tY)_{/(B\rightarrow A)^a}}\\
{D'}\ar[r]^-(0.5)n\ar[u]&{D_{/(B\rightarrow A)}}\ar[u]}
\end{equation}
où les flèches verticales sont les foncteurs canoniques, est commutatif. 
On en déduit, par (\cite{sga4} III 3.5), que 
la topologie co-évanescente de $D'$  est induite  
par la topologie co-évanescente de $D$ au moyen du foncteur $j_{(B\rightarrow A)}$. 
On notera que l'équivalence \eqref{higgs2-co-ev21a} induite par $n$ s'identifie à l'équivalence \eqref{higgs2-co-ev21b} 
en vertu de \eqref{higgs2-co-ev21c}, d'où la notation.

\begin{rema}
Nous pouvons donner une preuve directe de \ref{higgs2-co-ev21} qui ne passe pas par \ref{higgs2-topfl13}, mais qui utilise
les mêmes arguments. La preuve devient particulièrement simple lorsque $B=f^+(A)$. Nous traiterons 
directement ce cas dans un cadre plus général \eqref{higgs2-tcevg71}.
\end{rema}

\subsection{}\label{higgs2-co-ev17}
Rappelons \eqref{higgs2-topfl8} que la donnée d'un point de $\tX\gtimes_\tX\tY$ est équivalente à 
la donnée d'une paire de points $x\colon \Pt\rightarrow \tX$ et $y\colon \Pt\rightarrow \tY$ 
et d'un $2$-morphisme $u\colon f(y)\rightarrow x$. Un tel point sera noté $(y\rightarrow x)$,  
ou encore $(u\colon y\rightarrow x)$. 
Pour tous $F\in \ob(\tX)$ et $G\in \ob(\tY)$, on a des isomorphismes canoniques fonctoriels
\begin{eqnarray}
(\rp_1^*F)_{(y\rightarrow x)} &\stackrel{\sim}{\rightarrow}& F_x, \label{higgs2-co-ev17a}\\
(\rp_2^*G)_{(y\rightarrow x)} &\stackrel{\sim}{\rightarrow}& G_y. \label{higgs2-co-ev17b}
\end{eqnarray}
D'après \ref{higgs2-topfl6}, l'application 
\begin{equation}
(\rp_1^*F)_{(y\rightarrow x)} \rightarrow (\rp_2^*(f^*F))_{(y\rightarrow x)} 
\end{equation}
induite par $\tau$ \eqref{higgs2-co-ev11c}, s'identifie canoniquement au morphisme 
de spécialisation $F_x\rightarrow F_{f(y)}$ défini par $u$. 
Identifions les topos $\tD$ et $\tX\gtimes_{\tX}\tY$ par l'équivalence \eqref{higgs2-co-ev100b}. 
Les isomorphismes \eqref{higgs2-co-ev14c} et \eqref{higgs2-co-ev17b} induisent 
un isomorphisme canonique fonctoriel
\begin{equation}\label{higgs2-co-ev17c}
(\Psi_*G)_{(y\rightarrow x)} \stackrel{\sim}{\rightarrow} G_y.
\end{equation}
Il résulte \ref{higgs2-topfl6}  et \eqref{higgs2-co-ev14b} qu'on a un isomorphisme canonique fonctoriel  
de points de $\tX\gtimes_\tX\tY$
\begin{equation}\label{higgs2-co-ev17e}
\Psi(y) \stackrel{\sim}{\rightarrow} (y\rightarrow f(y)).
\end{equation}
D'après \ref{higgs2-co-ev8}(i), pour tout $(V\rightarrow U)\in \ob(D)$, on a un isomorphisme canonique fonctoriel
\begin{equation} \label{higgs2-co-ev17d}
(V\rightarrow U)^a_{(y\rightarrow x)} \stackrel{\sim}{\rightarrow} U^a_x\times_{U^a_{f(y)}}V^a_y,
\end{equation}
où l'exposant $^a$ désigne les faisceaux associés,  
l'application $V^a_y\rightarrow U^a_{f(y)}$ est induite par le morphisme structural 
$V\rightarrow f^+(U)$ et l'application $U^a_x\rightarrow U^a_{f(y)}$
est le morphisme de spécialisation défini par~$u$.

\section{Topos co-évanescents généralisés}\label{higgs2-tcevg}

\subsection{}\label{higgs2-tcevg1}\index{100000500@$\pi\colon E\rightarrow I$}
\index{100000503@$\alpha_{i}$}\index{100000505@$\cF$, $\cF^\vee$, $\cP^\vee$}
Dans cette section, $I$ désigne un $\mU$-site, $\tI$ le topos des faisceaux de $\mU$-ensembles sur $I$ et 
\begin{equation}\label{higgs2-tcevg1a}
\pi\colon E\rightarrow I
\end{equation}
une catégorie fibrée, clivée et normalisée au-dessus de la catégorie sous-jacente à $I$ (\cite{sga1} VI 7.1). 
On suppose les conditions suivantes satisfaites~:
\begin{itemize}
\item[(i)] Les produits fibrés sont représentables dans $I$. 
\item[(ii)]  Pour tout $i\in \ob(I)$, la catégorie fibre $E_i$ de $E$ au-dessus de $i$ 
est munie d'une topologie faisant de celle-ci un $\mU$-site, 
et les limites projectives finies sont représentables dans $E_i$. 
On note $\tE_i$ le topos des faisceaux de $\mU$-ensembles sur $E_i$. 
\item[(iii)] Pour tout morphisme $f\colon i\rightarrow j$ de $I$, le foncteur image inverse $f^+\colon E_j\rightarrow E_i$
est continu et exact à gauche. Il définit donc un morphisme de topos que l'on note aussi (abusivement) 
$f\colon \tE_i\rightarrow \tE_j$ (\cite{sga4} IV 4.9.2). 
\end{itemize}
La condition (i) sera renforcée à partir de \ref{higgs2-tcevg18}. Pour tout $i\in \ob(I)$, on note 
\begin{equation}\label{higgs2-tcevg1ab}
\alpha_{i!}\colon E_i\rightarrow E
\end{equation}
le foncteur d'inclusion canonique. 

Le foncteur $\pi$ est en fait un $\mU$-site fibré  (\cite{sga4} VI 7.2.1 et 7.2.4). On désigne par 
\begin{equation}\label{higgs2-tcevg1b}
\cF\rightarrow I
\end{equation}
le $\mU$-topos fibré associé à $\pi$ (\cite{sga4} VI 7.2.6). La catégorie fibre de $\cF$ au-dessus de tout $i\in \ob(I)$
est canoniquement équivalente au topos $\tE_i$, et le foncteur image inverse par tout  
morphisme $f\colon i\rightarrow j$ de $I$ s'identifie au foncteur image inverse $f^*\colon \tE_j\rightarrow \tE_i$ 
par le morphisme de topos $f\colon \tE_i\rightarrow \tE_j$. On désigne par
\begin{equation}\label{higgs2-tcevg1c}
\cF^\vee\rightarrow I^\circ
\end{equation}
la catégorie fibrée obtenue en associant à tout $i\in \ob(I)$ la catégorie $\tE_i$, et à tout morphisme 
$f\colon i\rightarrow j$ de $I$ le foncteur $f_*\colon \tE_i\rightarrow \tE_j$ image directe par le morphisme 
de topos $f\colon \tE_i\rightarrow \tE_j$. On désigne par
\begin{equation}\label{higgs2-tcevg1d}
\cP^\vee\rightarrow I^\circ
\end{equation}
la catégorie fibrée obtenue en associant à tout $i\in \ob(I)$ la catégorie $\hE_i$ des préfaisceaux de $\mU$-ensembles
sur $E_i$, et à tout morphisme $f\colon i\rightarrow j$ de $I$ le foncteur $f_*\colon \hE_i\rightarrow \hE_j$ obtenu en composant 
avec le foncteur image inverse $f^+\colon E_j\rightarrow E_i$.  Cette convention de notation ne suit 
pas celle de (\cite{sga4} I 5.0); elle est faite de sorte que le $I^\circ$-foncteur canonique 
$\cF^\vee\rightarrow \cP^\vee$ devient compatible aux foncteurs images inverses.

\subsection{}\label{higgs2-tcevg2}
\index{100000510@$\hE$}\index{100000511@$\{i\mapsto F_i\}$}
On note que $E$ est une $\mU$-catégorie. 
On désigne par $\hE$ la catégorie des préfaisceaux de $\mU$-ensembles sur $E$. On observera que 
la catégorie $\hE$ n'étant pas naturellement fibrée au-dessus de $I$, la notation $\hE_i$ des fibres de $\cP^\vee$
au-dessus de $I^\circ$ n'induit aucune confusion.  

D'après (\cite{sga1} VI 12; cf. aussi \cite{egr1} 1.1.2) et avec les notations de \eqref{higgs2-not3}, 
on a une équivalence de catégories
\begin{eqnarray}\label{higgs2-tcevg2a}
\hE&\stackrel{\sim}{\rightarrow} & \bHom_{I^\circ}(I^\circ,\cP^\vee)\\
F&\mapsto& \{i\mapsto F\circ \alpha_{i!}\},\nonumber
\end{eqnarray}
où $\alpha_{i!}$ est le foncteur \eqref{higgs2-tcevg1ab}.  
On identifiera dans la suite $F$ à la section $\{i\mapsto F\circ \alpha_{i!}\}$ qui lui est associée par cette équivalence.

\subsection{}\label{higgs2-tcevg3}\index{Topologie co-évanescente@Topologie co-evanescente}
\index{Topos!co-evanescent generalise@co-évanescent généralisé}
\index{Recouvrements!verticaux}\index{Recouvrements!cartesiens@cartésiens}
\index{100000512@$\tE$}\index{100000513@$\varepsilon$}
On appelle topologie {\em co-évanescente} de $E$ 
la topologie engendrée par les familles de recouvrements $(V_n\rightarrow V)_{n\in \Sigma}$
des deux types suivants~:
\begin{itemize}
\item[(v)] Il existe $i\in \ob(I)$ tel que $(V_n\rightarrow V)_{n\in \Sigma}$ soit 
une famille couvrante de $E_i$.
\item[(c)] Il existe une famille couvrante de morphismes $(f_n\colon i_n\rightarrow i)_{n\in \Sigma}$ de $I$ 
telle que $V_n$ soit isomorphe à $f_n^+(V)$ pour tout $n\in \Sigma$.
\end{itemize}
Les recouvrements du type (v) sont dits {\em verticaux}, et ceux du type (c) sont dits {\em cartésiens}. 
Le site ainsi défini est appelé site {\em co-évanescent} associé au site fibré $\pi$ \eqref{higgs2-tcevg1a}; c'est un $\mU$-site. 
On appelle topos {\em co-évanescent} associé au site fibré $\pi$, et l'on note $\tE$, 
le topos des faisceaux de $\mU$-ensembles sur $E$. On désigne par 
\begin{equation}
\varepsilon \colon E\rightarrow \tE
\end{equation}
le foncteur canonique. 

\begin{exemple}\label{higgs2-tcevg41}
Supposons que $I$ soit munie de la topologie grossière ou chaotique, c'est-à-dire de la topologie 
la moins fine parmi toutes les topologies de $I$ (\cite{sga4} II 1.1.4). On notera que sous cette hypothèse,
il revient au même de demander que $I$ soit un $\mU$-site ou que $I$ soit équivalente à une $\mU$-petite catégorie.  
La topologie totale sur $E$ relative au site fibré $\pi$ (\cite{sga4} VI 7.4.1)  
est engendrée par les recouvrements verticaux, en vertu de (\cite{sga4} VI 7.4.2(1)).
Elle est donc égale à la topologie co-évanescente sur $E$ . 
\end{exemple}

\begin{exemple}\label{higgs2-tcevg40}
Soient $X$ et $Y$ deux $\mU$-sites dans lesquels les 
limites projectives finies sont représentables, $f^+\colon X\rightarrow Y$ un foncteur continu et exact à gauche. 
On associe à $f^+$ un $\mU$-site fibré 
\begin{equation}\label{higgs2-tcevg40b}
\pi\colon E\rightarrow X
\end{equation} 
vérifiant les conditions de \eqref{higgs2-tcevg1} de la façon suivante. Considérons la catégorie $\Fl(Y)$ des morphismes de $Y$,
et le ``foncteur but'' 
\begin{equation}\label{higgs2-tcevg40a}
\Fl(Y)\rightarrow Y,
\end{equation}
qui fait de $\Fl(Y)$ une catégorie fibrée, clivée et normalisée au-dessus de $Y$~: 
la catégorie fibre au-dessus de 
tout objet $V$ de $Y$ est canoniquement équivalente à la catégorie $Y_{/V}$, et pour tout morphisme
$h\colon V'\rightarrow V$ de $Y$, le foncteur image inverse $h^*\colon Y_{/V}\rightarrow Y_{/V'}$ 
n'est autre que le foncteur de changement de base par $h$. 
Munissant chaque fibre $Y_{/V}$ de la topologie induite par celle de $Y$, $\Fl(Y)/Y$ devient un $\mU$-site fibré, 
vérifiant les conditions de \eqref{higgs2-tcevg1}. 
On prend alors pour $\pi$ le site fibré déduit de $\Fl(Y)/Y$ par changement de base par le foncteur $f^+$. 
Le site co-évanescent $E$ associé au site fibré $\pi$ \eqref{higgs2-tcevg3} 
est canoniquement équivalent au site co-évanescent $D$ associé au foncteur $f^+$ \eqref{higgs2-co-ev1}
en vertu de (\cite{sga4} III 5.2(1)), d'où la terminologie. 
\end{exemple}

\begin{lem}\label{higgs2-tcevg4}
{\rm (i)}\ Les produits fibrés sont représentables dans $E$.

{\rm (ii)}\ Les foncteurs $\pi$ et $\alpha_{i!}$ \eqref{higgs2-tcevg1ab}, pour tout $i\in \ob(I)$,
commutent aux produits fibrés.

{\rm (iii)}\ La famille des recouvrements verticaux (resp. cartésiens) de $E$ est stable par changement de base. 
\end{lem}

(i) Considérons un diagramme commutatif de $E$
\begin{equation}\label{higgs2-tcevg4a}
\xymatrix{
X\ar[d]\ar[r]&V\ar[d]\\
U\ar[r]&W}
\end{equation}
au-dessus d'un diagramme commutatif de $I$
\begin{equation}\label{higgs2-tcevg4b}
\xymatrix{
x\ar[d]\ar[r]&v\ar[d]\\
u\ar[r]& w}
\end{equation}
Alors $u\times_wv$ est représentable dans $I$.  
Notons $U'$, $V'$ et $W'$ les images inverses de $U$, $V$ et $W$ au-dessus de $u\times_wv$, 
par les morphismes canoniques de $u\times_wv$ dans $u$, $v$ et $w$, respectivement. 
Alors $U'\times_{W'}V'$ est représentable dans $E_{u\times_wv}$. 
Comme les foncteurs image inverse de $\pi$ sont exacts à gauche (\ref{higgs2-tcevg1}(iii)), le diagramme \eqref{higgs2-tcevg4a} 
détermine uniquement un morphisme $X\rightarrow U'\times_{W'}V'$ au-dessus du morphisme canonique 
$x\rightarrow u\times_wv$. Par suite, $U'\times_{W'}V'$ représente le produit fibré $U\times_WV$ dans $E$. 

(ii) \& (iii) Ceux-ci résultent aussitôt de la preuve de (i). 

\begin{rema}\label{higgs2-tcevg42}
On vérifie aisément que la famille des recouvrements verticaux (resp. horizontaux) forme une prétopologie \eqref{higgs2-tcevg4}.
Ce n'est pas le cas de leur union, ce qui est à l'origine de beaucoup de difficultés.  
\end{rema}

\begin{lem}\label{higgs2-tcevg43}
Soit $(V_m\rightarrow V)_{m\in M}$ un recouvrement vertical {\em fini} de $E$ et pour tout 
$m\in M$, soit $(V_{m,n}\rightarrow V_m)_{n\in N_m}$ un recouvrement cartésien de $E$. 
Alors, il existe un recouvrement cartésien $(W_\ell\rightarrow V)_{\ell\in L}$
tel que pour tous $m\in M$ et $\ell\in L$, il existe $n_\ell\in N_m$ et un $V_m$-morphisme 
$V_m\times_VW_\ell\rightarrow V_{m,n_\ell}$; en particulier, le recouvrement 
$(V_m\times_VW_\ell\rightarrow V)_{m\in M,\ell\in L}$ raffine le recouvrement $(V_{m,n}\rightarrow V)_{m\in M,n\in N_m}$. 
\end{lem}

En effet, pour tout $m\in M$, $(\pi(V_{m,n})\rightarrow \pi(V))_{n\in N_m}$ est un recouvrement de $I$. 
Comme $M$ est fini et que $I$ est stable par produits fibrés, 
il existe un recouvrement $(f_\ell\colon i_\ell\rightarrow \pi(V))_{\ell\in L}$ 
de $I$ tel que pour tous $m\in M$ et $\ell\in L$, il existe $n_\ell\in N_m$ et un $\pi(V)$-morphisme
$g_{m,\ell}\colon i_\ell\rightarrow \pi(V_{m,n_\ell})$ de $I$. Pour tout $\ell\in L$, 
posons $W_\ell=f_\ell^+(V)$. Pour tout $m\in M$, le morphisme $g_{m,\ell}$ induit 
alors un $V_m$-morphisme $V_m\times_VW_\ell\rightarrow V_{m,n_\ell}$; d'où la proposition.

\begin{prop}[\cite{tsuji3} 6.1.3]\label{higgs2-tcevg44}
Supposons que pour tout $i\in \ob(I)$, tout objet de $E_i$ soit quasi-compact. Alors, pour qu'un crible $R$ d'un objet $V$ de $E$ 
soit couvrant, il faut et il suffit qu'il existe un recouvrement cartésien 
$(V_n\rightarrow V)_{n\in N}$ et pour tout $n\in N$, un recouvrement vertical fini
$(V_{n,m}\rightarrow V_n)_{n\in M_n}$ tels que pour tous $n\in N$ et $m\in M_n$, on ait un $V$-morphisme
$V_{n,m}\rightarrow R$.
\end{prop}

Pour tout objet $V$ de $E$, notons $J(V)$ l'ensemble des cribles $R$ de $V$ dans $E$ vérifiant la propriété requise.
Par définition, tout crible de $J(V)$ est couvrant pour la topologie co-évanescente, 
et tout crible engendré par un recouvrement cartésien (resp. vertical) de $V$ appartient à $J(V)$. Il suffit donc de 
montrer que les $J(V)$, pour $V\in \ob(E)$, définissent une topologie (\cite{sga4} II 1.1). 
Il est clair que $V$ appartient à $J(V)$ (axiome (T3) de {\em loc. cit.}).
La stabilité par changement de base (axiome (T1) de {\em loc. cit.}) résulte de \ref{higgs2-tcevg4}(iii). 
Il reste à établir le caractère local (axiome (T2) de {\em loc. cit.}). Soient $V\in \ob(E)$, 
$R$ et $R'$ deux cribles de $V$ tels que $R\in J(V)$ et que pour tout $W\in \ob(E)$ et tout morphisme $W\rightarrow R$, 
le crible $R'\times_VW$ appartienne à $J(W)$. Montrons que $R'$ appartient à $J(V)$. Par hypothèse, 
il existe un recouvrement cartésien $(V_n\rightarrow V)_{n\in N}$ et pour tout $n\in N$, un recouvrement vertical fini
$(V_{n,m}\rightarrow V_n)_{n\in M_n}$ tels que pour tous $n\in N$ et $m\in M_n$, 
le morphisme composé $V_{n,m}\rightarrow V$ appartienne à $R$.
De plus, pour tous $n\in N$ et $m\in M_n$, il existe un recouvrement cartésien 
$(V_{n,m}^\alpha\rightarrow V_{n,m})_{\alpha\in A_{n,m}}$ et pour tout $\alpha\in A_{n,m}$, 
un recouvrement vertical fini $(V_{n,m}^{\alpha,\beta}\rightarrow V^\alpha_{n,m})_{\beta\in B^\alpha_{n,m}}$ tels que 
pour tous $\alpha\in A_{n,m}$ et $\beta\in B^\alpha_{n,m}$, 
le morphisme composé $V_{m,n}^{\alpha,\beta}\rightarrow V$ appartienne à $R'$.
D'après \ref{higgs2-tcevg43}, pour tout $n\in N$, il existe un recouvrement cartésien $(W_{n,\ell}\rightarrow V_n)_{\ell\in L_n}$
tel que pour tous $m\in M_n$ et $\ell\in L_n$, il existe $\alpha_{n,m,\ell}\in A_{n,m}$ et un $V_{n,m}$-morphisme
$p_{n,m,\ell}\colon V_{n,m}\times_{V_n}W_{n,\ell}\rightarrow V_{n,m}^{\alpha_{n,m,\ell}}$. 
Pour tous $n\in N$, $m\in M_n$ et $\ell\in L_n$, on note 
$(W_{n,\ell}^{m,\beta}\rightarrow V_{n,m}\times_{V_n}W_{n,\ell})_{\beta\in B_{n,m}^{\alpha_{n,m,\ell}}}$ l'image inverse par 
$p_{n,m,\ell}$ du recouvrement vertical 
$(V_{n,m}^{\alpha_\ell,\beta}\rightarrow V^{\alpha_{n,m,\ell}}_{n,m})_{\beta\in B^{\alpha_{n,m,\ell}}_{n,m}}$. 
Il est clair que $(W_{n,\ell}\rightarrow V)_{n\in N,\ell\in L_n}$ est un recouvrement cartésien, et que pour tous $n\in N$ et 
$\ell\in L_n$, $(W_{n,\ell}^{m,\beta}\rightarrow W_{n,\ell})_{m\in M_n, \beta\in B^{\alpha_{n,m,\ell}}_{n,m}}$ 
est un recouvrement vertical fini.
Pour tous $n\in N$, $\ell\in L_n$, $m\in M_n$ et $\beta\in B^{\alpha_{n,m,\ell}}_{n,m}$, le morphisme composé
$W_{n,\ell}^{\beta,m}\rightarrow V$ appartient à $R'$; d'où la propriété requise.

\begin{prop}\label{higgs2-tcevg5}
Pour qu'un préfaisceau $F=\{i\mapsto F_i\}$ sur $E$ soit un faisceau, 
il faut et il suffit que les conditions suivantes soient remplies~:
\begin{itemize}
\item[{\rm (i)}] Pour tout $i\in \ob(I)$, $F_i$ est un faisceau sur $E_i$. 
\item[{\rm (ii)}] Pour toute famille couvrante $(f_n\colon i_n\rightarrow i)_{n\in \Sigma}$ de $I$, si 
pour tout $(m,n)\in \Sigma^2$, on pose $i_{mn}=i_m\times_ii_n$
et on note $f_{mn}\colon i_{mn} \rightarrow i$ le morphisme canonique, alors la suite de morphismes de faisceaux sur $E_i$
\begin{equation}\label{higgs2-tcevg5a}
F_i\rightarrow \prod_{n\in \Sigma}(f_{n})_*(F_{i_n})\rightrightarrows \prod_{(m,n)\in \Sigma^2} (f_{mn})_*(F_{i_{mn}})
\end{equation}
est exacte.
\end{itemize}
\end{prop}

En effet, quitte à élargir l'univers $\mU$, on peut supposer la catégorie $I$ petite (\cite{sga4} II 2.7(2)). 
La proposition résulte alors de  \ref{higgs2-tcevg4} et (\cite{sga4} II 2.3, I 3.5, I 2.12 et II 4.1(3)).

\begin{cor}\label{higgs2-tcevg6}
Le foncteur \eqref{higgs2-tcevg2a} induit une équivalence de catégories entre $\tE$ et la sous-catégorie 
pleine de $\bHom_{I^\circ}(I^\circ,\cF^\vee)$ formée des sections $\{i\mapsto F_i\}$
telles que pour toute famille couvrante $(f_n\colon i_n\rightarrow i)_{n\in \Sigma}$ de $I$, si
pour tout $(m,n)\in \Sigma^2$,  on pose $i_{mn}=i_m\times_ii_n$
et on note $f_{mn}\colon i_{mn} \rightarrow i$ le morphisme canonique, la suite de faisceaux sur $E_i$
\begin{equation}\label{higgs2-tcevg6a}
F_i\rightarrow \prod_{n\in \Sigma}(f_{n})_*(F_{i_n})\rightrightarrows \prod_{(m,n)\in \Sigma^2} (f_{mn})_*(F_{i_{mn}})
\end{equation}
soit exacte.
\end{cor} 

\begin{cor}\label{higgs2-tcevg66}
Pour tout $i\in \ob(I)$, le foncteur $\alpha_{i!}\colon E_i\rightarrow E$ \eqref{higgs2-tcevg1ab} est continu. 
\end{cor}

\begin{rema}\label{higgs2-tcevg7}
La condition \ref{higgs2-tcevg5}(ii) est équivalente à la condition suivante~: 
\begin{itemize}
\item[(ii')] Pour tout $i\in \ob(I)$ et pour tout crible couvrant $\fR$ de $i$ dans $I$, le morphisme canonique 
\begin{equation}\label{higgs2-tcevg7a}
F_i\rightarrow \underset{\underset{(i',u)\in \fR^\circ}{\longleftarrow}}{\lim}\ u_*(F_{i'}),
\end{equation}
où $u\colon i'\rightarrow i$ désigne le morphisme structural, est un isomorphisme. 
\end{itemize}
En effet, quitte à élargir l'univers $\mU$, on peut supposer la catégorie $I$ petite. 
L'assertion résulte alors de (\cite{sga4} I 2.12) appliqué au foncteur $\fR^\circ \rightarrow \tE_i$, $(i',u)\mapsto u_*(F_{i'})$.
\end{rema}

\begin{rema}\label{higgs2-tcevg77}
Supposons que tout objet de $I$ soit quasi-compact. Alors la proposition \ref{higgs2-tcevg5} reste vraie si l'on se limite dans (ii)
aux familles couvrantes {\em finies} $(f_n\colon i_n\rightarrow i)_{n\in \Sigma}$ de $I$.
En effet, tout recouvrement cartésien de $E$ admet une sous-famille couvrante finie.  
Par suite, la topologie de $E$ est engendrée par  les recouvrements cartésiens {\em finis} 
et les recouvrements verticaux, d'où l'assertion. 
\end{rema}

\begin{rema}\label{higgs2-tcevg777}
Pour que les foncteurs $\alpha_{i!}$ soient cocontinus pour tout $i\in \ob(I)$,  
il faut et il suffit que la topologie co-évanescente de $E$ soit égale à sa topologie totale relative à $\pi$ \eqref{higgs2-tcevg41}.
En effet, cette dernière est par définition la topologie la moins fine qui
rend continus les foncteurs $\alpha_{i!}$ pour tout $i\in \ob(I)$ (\cite{sga4} VI 7.4.1), 
et est aussi la topologie la plus fine qui rend cocontinus les foncteurs $\alpha_{i!}$ pour tout $i\in \ob(I)$
(\cite{sga4} VI 7.4.3(2)). L'assertion résulte donc de \ref{higgs2-tcevg66}. 
On notera que les topologies co-évanescente et totale de $E$
ne sont pas en général égales, compte tenu de \ref{higgs2-tcevg5}.   
\end{rema}

\subsection{}\label{higgs2-tcevg22}\index{Topos!total}
\index{100000520@$\Top(E)$}\index{100000521@$\delta$}
Supposons que $I$ soit équivalente à une $\mU$-petite catégorie. 
Il résulte de \ref{higgs2-tcevg5} et (\cite{sga4} VI 7.4.7)
que le foncteur identique $\id_E\colon E\rightarrow E$ est continu lorsque l'on munit la source de la topologie totale
relative au site fibré $\pi$  \eqref{higgs2-tcevg41} et le but de la topologie co-évanescente.   
Notons $\Top(E)$ le {\em topos total} associé au site fibré $\pi$, c'est-à-dire le topos des faisceaux de $\mU$-ensembles sur 
le site total $E$. On a donc un morphisme canonique de topos (\cite{sga4} IV 4.9.2)
\begin{equation}\label{higgs2-tcevg22a}
\delta\colon \tE\rightarrow \Top(E)
\end{equation}
tel que le foncteur $\delta_*$ soit le foncteur d'inclusion canonique de $\tE$ dans $\Top(E)$; 
c'est un plongement (\cite{sga4} IV 9.1.1). On notera que le diagramme 
\begin{equation}\label{higgs2-tcevg22b}
\xymatrix{
{\tE}\ar[d]_{\delta_*}\ar[rd]&\\
{\Top(E)}\ar[r]^-(0.4)\sim&{\bHom_{I^\circ}(I^\circ,\cF^\vee)}}
\end{equation}
où la flèche horizontale est l'équivalence canonique de catégories (\cite{sga4} VI 7.4.7) et la flèche oblique est induite par 
le foncteur \eqref{higgs2-tcevg2a}, est commutatif à isomorphisme canonique près. 
D'autre part, pour tout objet $F$ de $\Top(E)$, $\delta^*(F)$ est canoniquement isomorphe 
au faisceau associé au préfaisceau $F$ sur le site co-évanescent $E$. En effet, quitte à élargir $\mU$, on peut supposer 
que la catégorie $E$ est $\mU$-petite (\cite{sga4} II 3.6 et III 1.5), auquel cas l'assertion résulte de (\cite{sga4} I 5.1 et III 1.3).

\begin{lem}\label{higgs2-tcevg8}
Soit $F=\{i\mapsto F_i\}$ un préfaisceau sur $E$. Pour chaque $i\in \ob(I)$,
notons $F_i^a$ le faisceau de $\tE_i$ associé à $F_i$. Alors $\{i\mapsto F_i^a\}$ est un 
préfaisceau sur $E$ et on a un morphisme canonique 
$\{i\mapsto F_i\}\rightarrow \{i\mapsto F_i^a\}$ de $\hE$, 
induisant un isomorphisme entre les faisceaux associés. 
\end{lem}

Pour tout morphisme $f\colon i'\rightarrow i$ de $I$, notons 
$\gamma_f\colon F_i\rightarrow f_*(F_{i'})$ le morphisme de transition de $F$ associé à $f$ 
\eqref{higgs2-tcevg2a}. Le morphisme canonique $F_{i'}\rightarrow F_{i'}^a$ induit un morphisme 
$f_*(F_{i'}) \rightarrow f_*(F_{i'}^a)$. Comme $f_*(F_{i'}^a)$ est un faisceau 
de $\tE_i$, $\gamma_f$ induit un morphisme 
\begin{equation}\label{higgs2-tcevg8a}
\gamma_f^a\colon F_i^a\rightarrow f_*(F_{i'}^a).
\end{equation}
Les morphismes $\gamma_f$ satisfont des relations de cocycles du type (\cite{egr1} (1.1.2.2)), 
déduites de la composition des morphismes dans $I$. 
Celles-ci induisent des relations de cocycles analogues pour les morphismes $\gamma_f^a$. 
Par suite, $\{i\mapsto F_i^a\}$ est une section de la catégorie fibrée $\cP^\vee$, 
et est donc un préfaisceau sur $E$ \eqref{higgs2-tcevg2a} (cf. \cite{sga1} VI 12 ou \cite{egr1} 1.1.2). 
De plus, on a un morphisme canonique 
\begin{equation}\label{higgs2-tcevg8b}
\{i\mapsto F_i\}\rightarrow \{i\mapsto F_i^a\}. 
\end{equation}
Soit $G=\{i\mapsto G_i\}$ un faisceau sur $E$. Pour tout $i\in \ob (I)$, 
$G_i$ est un faisceau de $\tE_i$ en vertu de \ref{higgs2-tcevg5}. On en déduit aussitôt que l'application 
\begin{equation}\label{higgs2-tcevg8c}
\Hom_{\hE}(\{i\mapsto F_i^a\},\{i\mapsto G_i\})\rightarrow 
\Hom_{\hE}(\{i\mapsto F_i\},\{i\mapsto G_i\})
\end{equation}
induite par \eqref{higgs2-tcevg8b} est un isomorphisme. Par suite, le morphisme \eqref{higgs2-tcevg8b} induit un isomorphisme
entre les faisceaux associés.

\subsection{}\label{higgs2-tcevg85}
Soient $\pi'\colon E'\rightarrow I$ un $\mU$-site fibré, clivé et normalisé vérifiant les conditions de \eqref{higgs2-tcevg1}, 
\begin{equation}\label{higgs2-tcevg85a}
\Phi\colon E'\rightarrow E
\end{equation}
un $I$-foncteur cartésien \eqref{higgs2-not3}. On munit $E'$ de la topologie co-évanescente 
définie par $\pi'$, et on note $\tE'$ le topos des faisceaux de $\mU$-ensembles sur $E'$.
On associe à $\pi'$ des objets analogues à ceux associés à $\pi$,
et on les note par les mêmes lettres affectées d'un prime $'$. 
Pour tout $i\in \ob(I)$, on note $\Phi_i\colon E'_i\rightarrow E_i$
le foncteur induit par $\Phi$ sur les catégories fibres en $i$, 
et $\hPhi_{i}^*\colon \hE_i\rightarrow \hE'_i$ le foncteur obtenu en composant avec $\Phi_i$. 
Pour tout morphisme $f\colon j\rightarrow i$ de $I$, on a un isomorphisme de foncteurs 
\begin{equation}\label{higgs2-tcevg85b}
\Phi_j \circ f_{E'}^+ \stackrel{\sim}{\rightarrow} f_E^+\circ \Phi_i,
\end{equation}
où $f_E^+$ et $f_{E'}^+$ sont les foncteurs image inverse de $E$ et $E'$, respectivement. 
Il induit un isomorphisme de foncteurs 
\begin{equation}\label{higgs2-tcevg85dd}
\hPhi_i^* \circ f_{E*} \stackrel{\sim}{\rightarrow}  f_{E'*} \circ \hPhi_j^*, 
\end{equation}
où $f_{E*}$ et $f_{E'*}$ sont les foncteurs image inverse de $\cP^\vee$ et $\cP'^\vee$, respectivement \eqref{higgs2-tcevg1d}. 
Les isomorphismes \eqref{higgs2-tcevg85b} vérifient une relation de cocycle du type (\cite{egr1} (1.1.2.2)), 
qui induit une relation analogue pour les isomorphismes \eqref{higgs2-tcevg85dd}. 
D'après (\cite{sga1} VI 12; cf. aussi \cite{egr1} 1.1.2), 
les foncteurs $\hPhi_i^*$ définissent donc un $I^\circ$-foncteur cartésien
\begin{equation}\label{higgs2-tcevg85de}
\cP^\vee\rightarrow \cP'^\vee.
\end{equation}
On vérifie aussitôt que le diagramme de foncteurs 
\begin{equation}\label{higgs2-tcevg85e}
\xymatrix{
{\hE}\ar[r]^-(0.5)\sim\ar[d]_{\hPhi^*}&{\bHom_{I^\circ}(I^\circ,\cP^\vee)}\ar[d]\\
{\hE'}\ar[r]^-(0.5)\sim&{\bHom_{I^\circ}(I^\circ,\cP'^\vee)}}
\end{equation}
où $\hPhi^*$ est le foncteur défini par la composition avec $\Phi$, 
la flèche verticale de droite est définie par la composition avec \eqref{higgs2-tcevg85de}, 
et les flèches horizontales sont les équivalences de catégories \eqref{higgs2-tcevg2a}, est commutatif
à isomorphisme canonique près. Par suite, pour tout préfaisceau $F=\{i\mapsto F_i\}$ sur $E$, on a 
\begin{equation}\label{higgs2-tcevg85f}
\hPhi^*(F)=\{i\mapsto \hPhi_i^*(F_{i})\}.
\end{equation}

Supposons que pour tout $i\in \ob(I)$, le foncteur $\Phi_i$ soit continu. Alors le foncteur $\hPhi_i^*$ induit 
un foncteur $\Phi_{i,s}\colon \tE_i\rightarrow \tE'_i$, qui commute aux limites projectives (\cite{sga4} III 1.2). 
Il résulte de \ref{higgs2-tcevg5}  et \eqref{higgs2-tcevg85e} que pour tout faisceau $F$ sur $E$,
$\hPhi^*(F)$ est un faisceau sur $E'$, et par suite, que 
$\Phi$ est continu pour les topologies co-évanescentes de $E$ et $E'$. Il induit donc un foncteur 
(\cite{sga4} III 1.1.1)
\begin{equation}\label{higgs2-tcevg85g}
\Phi_s\colon \tE\rightarrow \tE'.
\end{equation}
Si de plus, pour tout $i\in \ob(I)$, le foncteur $\Phi_{i,s}\colon \tE_i\rightarrow \tE'_i$ est une équivalence 
de catégories, alors $\Phi_s$ est une équivalence de catégories en vertu de \ref{higgs2-tcevg6}. 

\begin{rema}\label{higgs2-tcevg86}
Considérons le topos fibré $\cF/I$ \eqref{higgs2-tcevg1b} comme un $\mU$-site fibré, en munissant chaque fibre de 
la topologie canonique. Celui-ci vérifie clairement les conditions \eqref{higgs2-tcevg1}. Notons 
\begin{equation}\label{higgs2-tcevg86a}
\varepsilon_I\colon E\rightarrow \cF
\end{equation}
le $I$-foncteur cartésien canonique, qui induit sur les fibres les foncteurs canoniques
$\varepsilon_i\colon E_i\rightarrow \tE_i$ (\cite{sga4} VI (7.2.6.7)). Il résulte de \ref{higgs2-tcevg85} 
que $\varepsilon_I$ induit une équivalence entre les topos co-évanescents associés à $E/I$ et $\cF/I$.  
\end{rema}

\subsection{}\label{higgs2-tcevg9}
Soient $I'$ un $\mU$-site dans lequel les produits fibrés sont représentables, 
$\varphi\colon I'\rightarrow I$ un foncteur continu, qui commute aux produits fibrés. On désigne par 
\begin{equation}\label{higgs2-tcevg9a}
\pi'\colon E' \rightarrow I'
\end{equation}
le changement de base de $\pi$ \eqref{higgs2-tcevg1a} par $\varphi$, et par 
\begin{equation}\label{higgs2-tcevg9d}
\Phi\colon E' \rightarrow E
\end{equation}
la projection canonique (\cite{sga1} VI § 3). Alors $E'/I'$ est un site fibré vérifiant les conditions de \eqref{higgs2-tcevg1}. 
Le $\mU$-topos fibré $\cF'\rightarrow I'$ associé à $\pi'$ est canoniquement
$I'$-équivalent au topos fibré déduit de $\cF/I$ par changement de base par $\varphi$. 
On note
\begin{eqnarray}
\cF'^\vee&\rightarrow &I'^\circ,\label{higgs2-tcevg9b}\\
\cP'^\vee&\rightarrow &I'^\circ,\label{higgs2-tcevg9c}
\end{eqnarray}
les catégories fibrées associées à $\pi'$, définies dans \eqref{higgs2-tcevg1c} et \eqref{higgs2-tcevg1d}, respectivement. 
Elles sont canoniquement $I'^\circ$-équivalentes  aux catégories fibrées déduites de 
$\cF^\vee/I^\circ$ et $\cP^\vee/I^\circ$ par changement de base par $\varphi^\circ$.

On note $\hE'$ la catégorie des préfaisceaux de $\mU$-ensembles sur $E'$.
On vérifie aussitôt que le diagramme de foncteurs 
\begin{equation}\label{higgs2-tcevg9e}
\xymatrix{
{\hE}\ar[r]^-(0.5)\sim\ar[d]_{\hPhi^*}&{\bHom_{I^\circ}(I^\circ,\cP^\vee)}\ar[d]\\
{\hE'}\ar[r]^-(0.5)\sim&{\bHom_{I'^\circ}(I'^\circ,\cP'^\vee)}}
\end{equation}
où les flèches horizontales sont les équivalences de catégories \eqref{higgs2-tcevg2a}, 
$\hPhi^*$ est le foncteur défini par la composition avec $\Phi$ 
et la flèche verticale de droite est le foncteur canonique (\cite{sga1} VI § 3), est commutatif
à isomorphisme canonique près. Par suite, pour tout préfaisceau $F=\{i\mapsto F_i\}$ sur $E$, on a 
\begin{equation}\label{higgs2-tcevg9f}
\hPhi^*(F)=\{i'\mapsto F_{\varphi(i')}\}, 
\end{equation}
où pour tout $i'\in \ob(I')$, on a identifié les catégories fibres $E'_{i'}$ et $E_{\varphi(i')}$. 

On munit $E'$ de la topologie co-évanescente associée au site fibré $\pi'$, 
et on note $\tE'$ le topos des faisceaux de $\mU$-ensembles sur $E'$.  
Il résulte aussitôt de \ref{higgs2-tcevg5} et \eqref{higgs2-tcevg9e} que pour tout faisceau $F$ sur $E$,
$\hPhi^*(F)$ est un faisceau sur $E'$, et par suite, que le foncteur $\Phi$ est continu. Il induit donc un foncteur 
(\cite{sga4} III 1.1.1)
\begin{equation}\label{higgs2-tcevg9g}
\Phi_s\colon \tE\rightarrow \tE'.
\end{equation}

\begin{prop}\label{higgs2-tcevg10}
Conservons les hypothèses de \eqref{higgs2-tcevg9}, supposons de plus les conditions suivantes satisfaites~:
\begin{itemize}
\item[{\rm (i)}] La catégorie $I'$ est $\mU$-petite et le foncteur $\varphi$ est pleinement fidèle. 
\item[{\rm (ii)}] La topologie de $I'$ est induite par celle de $I$ au moyen du foncteur $\varphi$.
\item[{\rm (iii)}] Tout objet de $I$ peut être recouvert par des objets provenants de $I'$.
\end{itemize}
Alors le foncteur $\Phi_s\colon \tE\rightarrow \tE'$ \eqref{higgs2-tcevg9g} est une équivalence de catégories. 
\end{prop}

Pour tout objet $i$ de $I$, on désigne par $I'_{/i}$ la catégorie des couples $(i',u)$, où 
$i'\in \ob(I')$ et $u\colon \varphi(i')\rightarrow i$ est un morphisme de $I$. Soient $(i',u)$, $(i'',v)$
deux objets de $I'_{/i}$. Un morphisme de $(i',u)$ dans $(i'',v)$ est un morphisme $w\colon i'\rightarrow i''$
de $I'$ tel que $u=v\circ \varphi(w)$. On notera que le foncteur 
$I'_{/i}\rightarrow I_{/i}$,  $(i',u)\mapsto (\varphi(i'),u)$ est pleinement fidèle.

Montrons d'abord que pour tout objet $F=\{i\mapsto F_i\}$ de $\tE$ et tout objet $i$ de $I$, 
le morphisme de faisceaux sur $E_i$ 
\begin{equation}\label{higgs2-tcevg10aa}
F_i\rightarrow \underset{\underset{(i',u)\in I'^\circ_{/i}}{\longleftarrow}}{\lim}\ u_*(F_{\varphi(i')})
\end{equation}
est un isomorphisme. En effet, il résulte de \ref{higgs2-tcevg5} et des hypothèses que la suite 
\[
F_i\rightarrow \prod_{(i',u)\in \ob(I'_{/i})}u_*(F_{\varphi(i')})\rightrightarrows \prod_{((i',u),(i'',v))\in \ob(I'_{/i})^2}\ 
\prod_{(j,w)\in \ob(I'_{/\varphi(i')\times_i\varphi(i'')})}w_*(F_{\varphi(j)})
\]
est exacte. D'autre part, le morphisme canonique 
\begin{eqnarray*}
\lefteqn{\underset{\underset{(i',u)\in I'^\circ_{/i}}{\longleftarrow}}{\lim}\ u_*(F_{\varphi(i')})
\rightarrow}\\
&& \ker\left(\prod_{(i',u)\in \ob(I'_{/i})}u_*(F_{\varphi(i')})\rightrightarrows \prod_{((i',u),(i'',v))\in \ob(I'_{/i})^2}\ 
\prod_{(j,w)\in \ob(I'_{/\varphi(i')\times_i\varphi(i'')})}w_*(F_{\varphi(j)})\right)
\end{eqnarray*}
est un monomorphisme d'après (\cite{sga4} II 4.1(3)), d'où l'assertion. 

L'isomorphisme \eqref{higgs2-tcevg10aa} montre que le foncteur $\Phi_s$ est pleinement fidèle. 
Montrons que $\Phi_s$ est essentiellement surjectif. Soit $F'=\{i'\mapsto F'_{i'}\}$ un faisceau de $\tE'$. 
Pour tout $i\in \ob(I)$, la catégorie $I'_{/i}$ étant $\mU$-petite, 
on définit le faisceau $F_i$ de $\tE_i$ par la formule (\cite{sga4} II 4.1)
\begin{equation}\label{higgs2-tcevg10a}
F_i= \underset{\underset{(i',u)\in I'^\circ_{/i}}{\longleftarrow}}{\lim}\ u_*(F'_{i'}).
\end{equation} 
Pour tout morphisme $f\colon j\rightarrow i$ de $I$, le foncteur $f_*\colon \tE_j\rightarrow \tE_i$ 
commute aux limites projectives. 
On en déduit que $\{i\mapsto F_i\}$ est une section de la catégorie fibrée $\cP^\vee/I^\circ$ \eqref{higgs2-tcevg1d}
et est donc un préfaisceau sur $E$ \eqref{higgs2-tcevg2a}. On a clairement un isomorphisme canonique 
$\Phi^*(F)\stackrel{\sim}{\rightarrow}F'$ \eqref{higgs2-tcevg9f}.  
Il suffit donc de montrer que $F$ est un faisceau sur $E$, ou encore d'après \ref{higgs2-tcevg5} et \ref{higgs2-tcevg7}, 
que pour tout $i\in \ob(I)$ et tout crible couvrant $\fR$ de $i$ dans $I$, le morphisme canonique 
\begin{equation}\label{higgs2-tcevg10b}
F_i\rightarrow \underset{\underset{(j,v)\in \fR^\circ}{\longleftarrow}}{\lim}\ v_*(F_j),
\end{equation}
où $v\colon j\rightarrow i$ désigne le morphisme structural, est un isomorphisme. 
Pour tout $j\in \ob(I)$, notons $J_{/j}$ la sous-catégorie pleine de $I'_{/j}$ 
formée des couples $(j',v)$ tels que $\varphi(j')$ soit un objet de $\fR$.
On a $J_{/j}=I'_{/j}$ si $j\in \ob(\fR)$. 
Comme les foncteurs $u_*$ commutent aux limites projectives, 
il suffit de montrer que le morphisme canonique 
\begin{equation}\label{higgs2-tcevg10c}
F_i\rightarrow \underset{\underset{(i',u)\in J_{/i}^\circ}{\longleftarrow}}{\lim}\ u_*(F'_{i'})
\end{equation}
est un isomorphisme. 
Pour tout $(i',u)\in \ob(I'_{/i})$, $J_{/\varphi(i')}$ est un crible couvrant de $i'$ dans $I'$ (\cite{sga4} III 1.6). 
Par suite, le morphisme canonique 
\begin{equation}\label{higgs2-tcevg10d}
F'_{i'}
\rightarrow \underset{\underset{(j',v)\in J_{/\varphi(i')}^\circ}{\longleftarrow}}{\lim}\ v_*(F'_{j'})
\end{equation}
est un isomorphisme \eqref{higgs2-tcevg7a}. Appliquant le foncteur $u_*$ et 
prenant la limite projective suivant la catégorie $I'^\circ_{/i}$, 
on en déduit que le morphisme canonique 
\begin{equation}\label{higgs2-tcevg10e}
\underset{\underset{(i',u)\in I'^\circ_{/i}}{\longleftarrow}}{\lim}\ u_*(F'_{i'})
\rightarrow \underset{\underset{(j',v)\in J_{/i}^\circ}{\longleftarrow}}{\lim}\ v_*(F'_{j'})
\end{equation}
est un isomorphisme, et par suite que \eqref{higgs2-tcevg10c} est un isomorphisme. 

\begin{prop}\label{higgs2-tcevg11}
Les hypothèses étant celles de \eqref{higgs2-tcevg10}, supposons de plus que, pour tout $i'\in \ob(I')$, 
la catégorie $E'_{i'}$ soit $\mU$-petite.  Alors la topologie co-évanescente de $E'$ est induite par celle de $E$
au moyen du foncteur $\Phi$ \eqref{higgs2-tcevg9d}. 
\end{prop}

Le foncteur $\Phi_s$ \eqref{higgs2-tcevg9g} admet un adjoint à gauche $\Phi^s$
qui prolonge $\Phi$, {\em i.e.}, qui s'insère dans un diagramme commutatif à isomorphisme près
\begin{equation}\label{higgs2-tcevg11a}
\xymatrix{
{E}\ar[r]^{\varepsilon}&{\tE}\\
{E'}\ar[u]^{\Phi}\ar[r]^{\varepsilon'}&{\tE'}\ar[u]_{\Phi^s}}
\end{equation}
où les flèches horizontales sont les foncteurs canoniques (\cite{sga4} III 1.2). Comme  
$\Phi_s$ est une équivalence de catégories en vertu de \ref{higgs2-tcevg10}, $\Phi^s$ est un quasi-inverse de $\Phi_s$. 
On en déduit un isomorphisme $\varepsilon'\stackrel{\sim}{\rightarrow} \Phi_s\circ \varepsilon \circ \Phi$. 
Par suite, la topologie co-évanescente de $E'$ est induite par la topologie canonique de $\tE$ au moyen du foncteur 
$\varepsilon\circ \Phi$ (\cite{sga4} III 3.5). 

D'autre part, il résulte des hypothèses que la catégorie $E'$ est $\mU$-petite et que le foncteur $\Phi$ est pleinement fidèle. 
Donc en vertu du lemme de comparaison (\cite{sga4} III 4.1), si l'on munit $E'$ de la topologie induite par celle de $E$, 
le foncteur de restriction de $\tE$ dans la catégorie des faisceaux de $\mU$-ensembles
sur $E'$ est une équivalence de catégories.  Le raisonnement précédent montre alors que la topologie de $E'$
induite par celle de $E$ est aussi induite par la topologie canonique de $\tE$ au moyen du foncteur 
$\varepsilon\circ \Phi$,  d'où la proposition. 

\begin{rema}\label{higgs2-tcevg12}
Les hypothèses étant celles de \eqref{higgs2-tcevg10}, supposons de plus que, pour tout $i'\in \ob(I')$,  
la catégorie $E'_{i'}$ soit $\mU$-petite. 
On peut alors déduire \ref{higgs2-tcevg10} de \ref{higgs2-tcevg11} et (\cite{sga4} III 4.1), 
bien que nous ayons procédé inversement.
\end{rema}

\begin{lem}\label{higgs2-tcevg125}
\`A toute $\mU$-petite catégorie $J$ et tout foncteur 
\begin{equation}\label{higgs2-tcevg125a}
\phi\colon J\rightarrow \tE,\ \ \ j\mapsto F_j=\{i\mapsto F_{j,i}\},
\end{equation}
sont canoniquement associés les données suivantes~:
\begin{itemize}
\item[{\rm (i)}] un faisceau $\{i\mapsto \underset{\underset{j\in J}{\longleftarrow}}{\lim}\ F_{j,i}\}$ sur $E$, 
et un isomorphisme canonique 
\begin{equation}\label{higgs2-tcevg125c}
\{i\mapsto \underset{\underset{j\in J}{\longleftarrow}}{\lim}\ F_{j,i}\}
\stackrel{\sim}{\rightarrow} \underset{\underset{J}{\longleftarrow}}{\lim}\ \phi; 
\end{equation} 
\item[{\rm (ii)}]  un préfaisceau 
$\{i\mapsto \underset{\underset{j\in J}{\longrightarrow}}{\lim}\ F_{j,i}\}$ sur $E$, où les limites sont prises dans $\tE_i$,
et un isomorphisme canonique
\begin{equation}\label{higgs2-tcevg125d}
\underset{\underset{J}{\longrightarrow}}{\lim}\ \phi\stackrel{\sim}{\rightarrow}
\{i\mapsto \underset{\underset{j\in J}{\longrightarrow}}{\lim}\ F_{j,i}\}^a,
\end{equation}
où le terme de droite est le faisceau sur $E$ associé au préfaisceau 
$\{i\mapsto \underset{\underset{j\in J}{\longrightarrow}}{\lim}\ F_{j,i}\}$.
\end{itemize}
En particulier, pour tout $i\in I$, le foncteur canonique 
\begin{equation}\label{higgs2-tcevg125b}
\tE\rightarrow \tE_i, \ \ \ \{i\mapsto F_i\}\mapsto F_i
\end{equation}
commute aux $\mU$-limites projectives.
\end{lem}

Il résulte aussitôt de \eqref{higgs2-tcevg2a} et de (\cite{sga4} I 3.1) que le foncteur 
\begin{equation}\label{higgs2-tcevg125e}
\hE\rightarrow \hE_i, \ \ \ \{i\mapsto G_i\}\mapsto G_i
\end{equation}
commute aux $\mU$-limites projectives et inductives. De plus, 
à toute $\mU$-petite catégorie $J$ et tout foncteur 
\begin{equation}\label{higgs2-tcevg125f}
\psi\colon J\rightarrow \hE,\ \ \ j\mapsto G_j=\{i\mapsto G_{j,i}\},
\end{equation}
sont canoniquement associés deux préfaisceaux  $\{i\mapsto \underset{\underset{j\in J}{\longrightarrow}}{\lim}\ G_{j,i}\}$ et 
$\{i\mapsto \underset{\underset{j\in J}{\longleftarrow}}{\lim}\ G_{j,i}\}$ sur $E$, et deux morphismes canoniques 
\begin{eqnarray}
\{i\mapsto \underset{\underset{j\in J}{\longleftarrow}}{\lim}\ G_{j,i}\}
\rightarrow \underset{\underset{J}{\longleftarrow}}{\lim}\ \psi, \label{higgs2-tcevg125g} \\
\underset{\underset{J}{\longrightarrow}}{\lim}\ \psi\rightarrow  
\{i\mapsto \underset{\underset{j\in J}{\longrightarrow}}{\lim}\ G_{j,i}\}.\label{higgs2-tcevg125h}
\end{eqnarray}
Ces derniers sont donc des isomorphismes.

Supposons que $\psi$ soit induit par un foncteur $\phi$ comme dans \eqref{higgs2-tcevg125a}. 
D'après (\cite{sga4} II 4.1(3)), l'isomorphisme \eqref{higgs2-tcevg125g} induit l'isomorphisme \eqref{higgs2-tcevg125c};
en particulier, $\{i\mapsto \underset{\underset{j\in J}{\longleftarrow}}{\lim}\ G_{j,i}\}$ est un faisceau sur $E$. 
D'autre part, en vertu de \ref{higgs2-tcevg8} et (\cite{sga4} II 4.1), 
$\{i\mapsto \underset{\underset{j\in J}{\longrightarrow}}{\lim}\ F_{j,i}\}$ est un préfaisceau sur $E$, et 
l'isomorphisme \eqref{higgs2-tcevg125h} induit l'isomorphisme \eqref{higgs2-tcevg125d}. 

\begin{prop}\label{higgs2-tcevg13}
Supposons les conditions suivantes satisfaites~:
\begin{itemize}
\item[{\rm (i)}] Tout objet de $I$ est quasi-compact. 
\item[{\rm (ii)}] Pour tout objet $i$ de $I$, le topos $\tE_i$ est algébrique {\rm (\cite{sga4} VI 2.3)}.
\item[{\rm (iii)}] Pour tout morphisme $f\colon i\rightarrow j$ de $I$, 
le morphisme de topos $f\colon \tE_i\rightarrow \tE_j$ est cohérent {\rm (\cite{sga4} VI 3.1)}. 
\end{itemize}
Alors pour toute $\mU$-petite catégorie filtrante $J$ et tout foncteur 
\begin{equation}\label{higgs2-tcevg13a}
\phi\colon J\rightarrow \tE,\ \ \ j\mapsto F_j=\{i\mapsto F_{j,i}\},
\end{equation}
$\{i\mapsto \underset{\underset{j\in J}{\longrightarrow}}{\lim}\ F_{j,i}\}$ est un faisceau sur $E$,
et on a un isomorphisme canonique
\begin{equation}\label{higgs2-tcevg13b}
\underset{\underset{J}{\longrightarrow}}{\lim}\ \phi\stackrel{\sim}{\rightarrow}
\{i\mapsto \underset{\underset{j\in J}{\longrightarrow}}{\lim}\ F_{j,i}\}.
\end{equation}
En particulier, pour tout $i\in I$, le foncteur canonique 
\begin{equation}\label{higgs2-tcevg13aa}
\tE\rightarrow \tE_i, \ \ \ \{i\mapsto F_i\}\mapsto F_i
\end{equation}
commute aux $\mU$-limites inductives filtrantes.
\end{prop}

Compte tenu de \ref{higgs2-tcevg125}(ii), il suffit de montrer que
$\{i\mapsto \underset{\underset{j\in J}{\longrightarrow}}{\lim}\ F_{j,i}\}$ est un faisceau sur $E$. 
Soit $(f_n\colon i_n\rightarrow i)_{n\in \Sigma}$ une famille couvrante {\em finie} de $I$. 
Pour tout $(m,n)\in \Sigma^2$,  on pose $i_{mn}=i_m\times_ii_n$
et on note $f_{mn}\colon i_{mn} \rightarrow i$ le morphisme canonique. 
Pour tout $j\in J$, la suite 
\begin{equation}\label{higgs2-tcevg13d}
F_{j,i}\rightarrow \prod_{n\in \Sigma} (f_n)_*(F_{j,i_n})
\rightrightarrows \prod_{(m,n)\in \Sigma^2} (f_{mn})_*(F_{j,i_{mn}})
\end{equation}
est exacte. Comme $\Sigma$ est fini, on en déduit par passage à la limite inductive sur $J$  que la suite  
\begin{equation}\label{higgs2-tcevg13e}
\underset{\underset{j\in J}{\longrightarrow}}{\lim}\ F_{j,i}\rightarrow 
\prod_{n\in \Sigma} \underset{\underset{j\in J}{\longrightarrow}}{\lim}\ (f_n)_*(F_{j,i_n})
\rightrightarrows \prod_{(m,n)\in \Sigma^2} \underset{\underset{j\in J}{\longrightarrow}}{\lim}\ 
(f_{mn})_*(F_{j,i_{mn}})
\end{equation}
est exacte (\cite{sga4} II 4.3(4)).  Comme les foncteurs $(f_n)_*$ et $(f_{mn})_*$
commutent aux limites inductives filtrantes de faisceaux d'ensembles en vertu de (\cite{sga4} VI 5.1 et VII 5.14), 
on en déduit que la suite 
\begin{equation}\label{higgs2-tcevg13c}
\underset{\underset{j\in J}{\longrightarrow}}{\lim}\ F_{j,i}\rightarrow \prod_{n\in \Sigma}
(f_n)_*(\underset{\underset{j\in J}{\longrightarrow}}{\lim}\ F_{j,i_n})
\rightrightarrows \prod_{(m,n)\in \Sigma^2} 
(f_{mn})_*(\underset{\underset{j\in J}{\longrightarrow}}{\lim}\ F_{j,i_{mn}})
\end{equation}
est exacte. L'assertion recherchée s'ensuit compte tenu de \ref{higgs2-tcevg77}.

\begin{cor}\label{higgs2-tcevg14}
Les hypothèses étant celles de \eqref{higgs2-tcevg13}, soient, de plus, $V$ un objet de $E$, $v$ son image dans $I$. 
Alors pour que $V$ soit quasi-compact dans $E$, il faut et il suffit qu'il le soit dans $E_v$. 
\end{cor}

Il résulte aussitôt de la définition de la topologie co-évanescente \eqref{higgs2-tcevg3} que si $V$ 
est quasi-compact dans $E$, il l'est aussi dans $E_v$. 
Montrons l'implication inverse. Notons $\varepsilon_v\colon E_v\rightarrow \tE_v$ le foncteur canonique. 
D'après \ref{higgs2-tcevg13}, pour toute $\mU$-petite catégorie filtrante $J$ et tout foncteur 
\begin{equation}
J\rightarrow \tE,\ \ \ j\mapsto F_j=\{i\mapsto F_{j,i}\},
\end{equation}
on peut identifier les applications canoniques 
\begin{eqnarray}
\underset{\underset{j\in J}{\longrightarrow}}{\lim}\ \Hom_{\tE}(\varepsilon(V),F_j)&\rightarrow& 
\Hom_{\tE}(\varepsilon(V),\underset{\underset{j\in J}{\longrightarrow}}{\lim}\ F_j),\label{higgs2-tcevg14b}\\
\underset{\underset{j\in J}{\longrightarrow}}{\lim}\ \Hom_{\tE_v}(\varepsilon_v(V),F_{j,v})&\rightarrow& 
\Hom_{\tE_v}(\varepsilon_v(V),\underset{\underset{j\in J}{\longrightarrow}}{\lim}\ F_{j,v}).\label{higgs2-tcevg14c}
\end{eqnarray}
Si $V$ est quasi-compact dans $E_v$, l'application \eqref{higgs2-tcevg14c} est injective d'après (\cite{sga4} VI 1.2 et 1.23(i)). 
Il en est donc de même de l'application \eqref{higgs2-tcevg14b}. Par suite, $V$ est quasi-compact dans $E$
en vertu de (\cite{sga4} VI 1.23(i)).

\begin{prop}\label{higgs2-tcevg15}
Supposons que tout objet de $I$ soit quasi-compact et que pour tout $i\in \ob(I)$, 
tout objet de $E_i$ soit quasi-compact. Alors~:
\begin{itemize}
\item[{\rm (i)}] Pour tout objet $i$ de $I$, le topos $\tE_i$ est cohérent.
\item[{\rm (ii)}] Pour tout morphisme $f\colon i\rightarrow j$ de $I$, 
le morphisme de topos $f\colon \tE_i\rightarrow \tE_j$ est cohérent. 
\item[{\rm (iii)}] Pour tout objet $V$ de $E$, $\varepsilon(V)$ est un objet cohérent de $\tE$~; 
en particulier, le topos $\tE$ est localement cohérent. 
\item[{\rm (iv)}] Si, de plus, la catégorie $E$ admet un objet final $e$, alors le topos $\tE$ est cohérent. 
\end{itemize}
\end{prop}

(i) Notons $\varepsilon_i\colon E_i\rightarrow \tE_i$ le foncteur canonique.
Pour tout $U\in \ob(E_i)$, $\varepsilon_i(U)$ est un objet cohérent de $\tE_i$ en vertu de 
\ref{higgs2-tcevg1}(ii) et (\cite{sga4} VI 2.1). D'autre part, $E_i$ admet un objet final $e_i$ (\ref{higgs2-tcevg1}(ii)).  
Comme $\varepsilon _i$ est exact à gauche, $\varepsilon_i(e_i)$ est l'objet final de $\tE_i$.
Pour tout $U\in \ob(E_i)$, le morphisme diagonal 
$\delta\colon \varepsilon_i(U)\rightarrow \varepsilon_i(U)\times_{\varepsilon_i(e_i)}\varepsilon_i(U)$ 
est l'image par $\varepsilon_i$ du morphisme diagonal $U\rightarrow U\times_{e_i}U$.
Donc $\delta$ est quasi-compact puisque sa source et son but sont cohérents, autrement dit, $\varepsilon_i(U)$
est quasi-séparé sur $\varepsilon_i(e_i)$. Par suite, le topos $\tE_i$ est cohérent (\cite{sga4} VI 2.3). 

(ii) Cela résulte de (\cite{sga4} VI 3.3). 

(iii) Tout objet de $E$ est quasi-compact d'après \ref{higgs2-tcevg44} (ou d'après (i), (ii) et \ref{higgs2-tcevg14}). 
Comme les produits fibrés sont représentables dans $E$ \eqref{higgs2-tcevg4}, la proposition résulte de (\cite{sga4} VI 2.1). 

(iv) Comme le foncteur canonique $\varepsilon \colon E\rightarrow \tE$ est exact à gauche, $\varepsilon(e)$
est l'objet final de $\tE$. Pour tout $V\in \ob(E)$, le morphisme diagonal 
$\delta\colon \varepsilon(V)\rightarrow \varepsilon(V)\times_{\varepsilon(e)}\varepsilon(V)$ est l'image par $\varepsilon$
du morphisme diagonal $V\rightarrow V\times_eV$ \eqref{higgs2-tcevg4}. Donc $\delta$ est quasi-compact 
puisque sa source et son but sont cohérents d'après (iii), autrement dit, $\varepsilon(V)$
est quasi-séparé sur $\varepsilon(e)$. La proposition s'ensuit compte tenu de~(iii).

\begin{prop}\label{higgs2-tcevg155}
Supposons les conditions suivantes satisfaites~:
\begin{itemize}
\item[{\rm (i)}] Il existe une sous-catégorie pleine, $\mU$-petite et topologiquement génératrice $I'$ de $I$, 
formée d'objets quasi-compacts, et stable par produits fibrés. 
\item[{\rm (ii)}] Pour tout objet $i$ de $I'$, le topos $\tE_i$ est cohérent.
\item[{\rm (iii)}] Pour tout morphisme $f\colon i\rightarrow j$ de $I'$, 
le morphisme de topos $f\colon \tE_i\rightarrow \tE_j$ est cohérent. 
\end{itemize}
Alors le topos $\tE$ est localement cohérent. 
Si, de plus, la catégorie $I$ admet un objet final qui appartient à $I'$, 
alors le topos $\tE$ est cohérent.
\end{prop}

D'après \ref{higgs2-tcevg10}, on peut se borner au cas où $I=I'$. 
En vertu de \ref{higgs2-tcevg86}, il suffit de montrer la proposition analogue pour le site fibré $\cF/I$. 
Notons $\cC$ la sous-catégorie pleine de $\cF$ formée des objets qui sont cohérents dans leurs fibres. 
D'après la condition (iii), le foncteur structural $\cC\rightarrow I$ fait de $\cC/I$
une catégorie fibrée, clivée et normalisée. La fibre $\cC_i$ de $\cC$ au-dessus d'un objet $i$ de $I$
s'identifie à la sous-catégorie pleine de $\tE_i$ formée des objets cohérents de $\tE_i$. 
Celle-ci est génératrice dans $\tE_i$, et est stable par produits fibrés d'après (\cite{sga4} VI  2.2). 
Comme l'objet final de $\tE_i$ est un objet de $\cC_i$ par hypothèse, les limites projectives finies sont 
représentables dans $\cC_i$ (\cite{sga4} I 2.3.1). 
Munissant chaque fibre $\cC_i$ de la topologie induite par la topologie canonique de $\tE_i$, 
$\cC/I$ devient un site fibré vérifiant les conditions de \eqref{higgs2-tcevg1}. 
Les topos co-évanescents associés à $\cC/I$ et $\cF/I$ sont équivalents d'après \ref{higgs2-tcevg85}. 
D'autre part, si $I$ admet un objet final $\iota$, alors l'objet final $e$ de $\cC_\iota$, 
qui existe d'après \ref{higgs2-tcevg1}(ii), est un objet final de $\cC$ (cela résulte facilement de \ref{higgs2-tcevg1}(iii); 
cf. \ref{higgs2-tcevg18} ci-dessous). 
On se réduit ainsi à l'énoncé de la proposition \ref{higgs2-tcevg15}.

\begin{cor}\label{higgs2-tcevg16}
Sous les hypothèses de \eqref{higgs2-tcevg155}, le topos $\tE$ a suffisamment de points.
\end{cor}

Cela résulte de (\cite{sga4} VI § 9).

\begin{cor}\label{higgs2-tcevg17}
Soient $\tX$, $\tY$ deux $\mU$-topos cohérents, $f\colon \tY\rightarrow \tX$ un morphisme cohérent. 
Alors le produit orienté $\tX\gtimes_{\tX}\tY$ est cohérent. En particulier, il a suffisamment de points. 
\end{cor}

D'après \ref{higgs2-co-ev102} ci-dessous, il existe une sous-catégorie pleine, $\mU$-petite et génératrice $X$ (resp. $Y$)  
de $\tX$ (resp. $\tY$), formée d'objets cohérents de $\tX$ (resp. $\tY$) et stable par limites projectives finies, 
telles que pour tout $U\in \ob(X)$, $f^*(U)$ appartienne à $Y$.   
Le produit orienté $\tX\gtimes_{\tX}\tY$ est équivalent au topos co-évanescent associé au foncteur 
$f^*\colon X\rightarrow Y$ \eqref{higgs2-co-ev101}. Par ailleurs, pour tout objet $V$
de $Y$, le topos $\tY_{/V}$ est cohérent (\cite{sga4} VI 2.4.2), et pour tout morphisme $V'\rightarrow V$ de $Y$,
le morphisme de localisation $\tY_{/V'}\rightarrow \tY_{/V}$ est cohérent (\cite{sga4} VI 3.3). 
La proposition résulte donc de \ref{higgs2-tcevg40} et \ref{higgs2-tcevg155}.

\begin{lem}\label{higgs2-co-ev102}
Pour tout morphisme de $\mU$-topos $f\colon \tY\rightarrow \tX$,  
il existe une sous-catégorie pleine, $\mU$-petite et génératrice $X$ (resp. $Y$)  
de $\tX$ (resp. $\tY$), stable par limites projectives finies telles que pour tout $U\in \ob(X)$, 
$f^*(U)$ appartienne à $Y$. Si, de plus, $\tX$ et $\tY$ sont cohérents
et si $f$ est cohérent, on peut supposer $X$ (resp. $Y$) formée d'objets cohérents. 
\end{lem}

En effet, il existe une sous-catégorie pleine, $\mU$-petite et génératrice $X$ (resp. $Y$)  
de $\tX$ (resp. $\tY$) telle que pour tout $U\in \ob(X)$, $f^*(U)$ appartienne à $Y$. 
D'après le premier argument de (\cite{sga4} IV 1.2.3), quitte à augmenter $X$, puis $Y$, 
on peut les supposer stables par limites projectives finies.  

Supposons ensuite $\tX$ et $\tY$ cohérents et $f$ cohérent.  
D'après (\cite{sga4} VI 2.1), il existe une sous-catégorie pleine, $\mU$-petite et génératrice $X$ (resp. $Y$)  
de $\tX$ (resp. $\tY$), formée d'objets cohérents, telle que pour tout $U\in \ob(X)$, $f^*(U)$ appartienne à $Y$.
Par ailleurs, la sous-catégorie pleine de $\tX$ (resp. $\tY$) formée des objets cohérents de $\tX$ (resp. $\tY$) 
est stable par limites projectives finies (\cite{sga4} VI 2.4.4).  L'argument de (\cite{sga4} IV 1.2.3) montre alors 
que quitte à augmenter $X$ puis $Y$, on peut les supposer stables par limites projectives finies.

\subsection{}\label{higgs2-tcevg18}
\index{100000525@$\sigma$, $\beta$}
Dans la suite de cette section, en plus des hypothèses faites dans \ref{higgs2-tcevg1}, nous supposons 
la condition suivante satisfaite~:
\begin{itemize}
\item[{\rm (i')}] Les limites projectives finies sont représentables dans $I$. 
\end{itemize}
Compte tenu de \ref{higgs2-tcevg1}(i), il revient au même de demander que $I$ admette un objet final (\cite{sga4} I 2.3.1).
On se donne un objet final  $\iota$ de $I$ et un objet final $e$ de $E_\iota$ 
(qui existe d'après \ref{higgs2-tcevg1}(ii)), que nous supposons fixés dans la suite. Pour tout $i\in \ob(I)$, on note 
\begin{equation}\label{higgs2-tcevg18h}
f_i\colon i\rightarrow \iota
\end{equation} 
le morphisme canonique. D'après \ref{higgs2-tcevg1}(iii), $f^+_i(e)$ est un objet final de $E_i$. 
Par suite, $e$ est un objet final de $E$, 
et les limites projectives finies sont représentables dans $E$ d'après \ref{higgs2-tcevg4}(i) et (\cite{sga4} I 2.3.1).

Le foncteur d'injection canonique $\alpha_{\iota!}\colon E_\iota\rightarrow E$ 
commute aux produits fibrés \eqref{higgs2-tcevg4},
et transforme objet final en objet final~; il est donc exact à gauche.
Par ailleurs, il est continu \eqref{higgs2-tcevg66}. Il définit alors un morphisme de topos (\cite{sga4} IV 4.9.2)
\begin{equation}\label{higgs2-tcevg18b}
\beta\colon \tE\rightarrow \tE_\iota.
\end{equation}
Pour tout faisceau $F=\{i\mapsto F_i\}$ sur $E$, on a un isomorphisme canonique
\begin{equation}\label{higgs2-tcevg18c}
\beta_*(F)\stackrel{\sim}{\rightarrow}F_\iota.
\end{equation}

Il existe essentiellement une unique section cartésienne de $\pi$ \eqref{higgs2-tcevg1a}
\begin{equation}\label{higgs2-tcevg18d}
\sigma^+\colon I\rightarrow E
\end{equation}
telle que $\sigma^+(\iota)=e$. 
Elle est définie, compte tenu de (\cite{egr1} 1.1.2), pour tout objet $i$ de $I$, par  $\sigma^+(i)=f^+_i(e)$ 
et pour tout morphisme $f\colon j\rightarrow i$ de $I$ par l'isomorphisme canonique 
\begin{equation}\label{higgs2-tcevg18e}
\sigma^+(j)\stackrel{\sim}{\rightarrow}f^+(\sigma^+(i)).
\end{equation}
Pour tout $i\in \ob(I)$, $\sigma^+(i)$ est un objet final de $E_i$, et 
on a un morphisme canonique de foncteurs 
\begin{equation}\label{higgs2-tcevg18g}
\id_E\rightarrow \sigma^+\circ \pi.
\end{equation}
Il résulte aussitôt de la preuve de \ref{higgs2-tcevg4} que $\sigma^+$ commute aux produits fibrés, et est donc exact à gauche. 
D'autre part, $\sigma^+$ est continu en vertu de (\cite{sga4} III 1.6). Il définit donc un 
morphisme de topos 
\begin{equation}\label{higgs2-tcevg18f}
\sigma\colon \tE\rightarrow \tI.
\end{equation}

\begin{rema}\label{higgs2-tcevg188}
Les morphismes $\beta$ \eqref{higgs2-tcevg18b} et $\sigma$ \eqref{higgs2-tcevg18f} 
s'identifient canoniquement à ceux définis à partir du site fibré $\cF/I$ \eqref{higgs2-tcevg86}.
\end{rema} 

\begin{lem}\label{higgs2-tcevg19}
{\rm (i)}\ Pour tout faisceau $F$ de $\tE_\iota$, $\beta^*(F)$ est le faisceau associé au préfaisceau
$\{i\mapsto f_i^*F\}$ sur $E$, où  pour tout morphisme $f\colon j\rightarrow i$ de $I$, le morphisme de transition 
\begin{equation}\label{higgs2-tcevg19a}
f_i^*F\rightarrow f_*(f_j^*F)
\end{equation}
est l'adjoint de l'isomorphisme canonique $f^*(f_i^*F)\stackrel{\sim}{\rightarrow} f_j^*F$.
De plus, le morphisme d'adjonction $F\rightarrow \beta_*(\beta^*F)$
se factorise à travers le morphisme de préfaisceaux sur $E_\iota$ 
\begin{equation}\label{higgs2-tcevg19b}
F\rightarrow \{i\mapsto f_i^*F\}\circ \alpha_{\iota!}
\end{equation}
défini par l'identité de $F$. 

{\rm (ii)}\ Pour tout faisceau $F$ de $\tI$, $\sigma^*(F)$ est le faisceau associé au préfaisceau $\{i\mapsto F(i)\}$ 
sur $E$, où pour tout $i\in \ob(I)$, on a noté $F(i)$ le préfaisceau constant sur $E_i$ de valeur $F(i)$.  
\end{lem}

En effet, quitte à élargir $\mU$, on peut supposer que les catégories $I$ et $E_\iota$ sont $\mU$-petites 
(\cite{sga4} II 3.6 et III 1.5). 

(i) D'après (\cite{sga4} I 5.1 et III 1.3), $\beta^*(F)$ est le faisceau  
associé au préfaisceau $G$ sur $E$ défini, pour tout $V\in \ob(E)$, par 
\begin{equation}\label{higgs2-tcevg19c}
G(V)=\underset{\underset{(U,u)\in A^\circ_{V}}{\longrightarrow}}{\lim}\ F(U),
\end{equation}
où $A_{V}$ est la catégorie des couples $(U,u)$ formés d'un objet $U$ de $E_\iota$ 
et d'un morphisme $u\colon V\rightarrow U$ de $E$. Si $i$ désigne l'image de $V$ dans $I$, 
$A_V$ s'identifie à la catégorie des couples $(U,u)$ formés d'un objet $U$ de $E_\iota$ 
et d'un morphisme $u\colon V\rightarrow f_i^+(U)$ de $E_i$. 
Il résulte alors de \ref{higgs2-tcevg8} que $\beta^*(F)$ est le faisceau associé au préfaisceau 
$\{i\mapsto f_i^*F\}$ sur $E$ défini par les morphismes de transition \eqref{higgs2-tcevg19a}.
La dernière assertion résulte aussitôt de ce qui précède. 

(ii) En effet, $\sigma^*(F)$ est le faisceau associé au préfaisceau $H$ sur $E$
défini, pour tout $V\in \ob(E)$, par 
\begin{equation}\label{higgs2-tcevg19d}
H(V)=\underset{\underset{(i,u)\in B^\circ_{V}}{\longrightarrow}}{\lim}\ F(i),
\end{equation}
où $B_{V}$ est la catégorie des couples $(i,u)$ formés d'un objet $i$ de $I$ et 
d'un morphisme $u\colon V\rightarrow \sigma^+(i)$ de $E$. 
Cette catégorie admet comme objet initial le couple formé de $\pi(V)$ et du morphisme canonique
$V\rightarrow \sigma^+(\pi(V))$ \eqref{higgs2-tcevg18g}. On a donc $H(V)=F(\pi(V))$. 

\begin{prop}\label{higgs2-tcevg199}
Supposons les conditions suivantes satisfaites~:
\begin{itemize}
\item[{\rm (i)}] Il existe une sous-catégorie pleine, $\mU$-petite et topologiquement génératrice $I'$ de $I$, 
formée d'objets quasi-compacts, et stable par limites projectives finies. 
\item[{\rm (ii)}] Pour tout objet $i$ de $I'$, le topos $\tE_i$ est cohérent.
\item[{\rm (iii)}] Pour tout morphisme $f\colon i\rightarrow j$ de $I'$, 
le morphisme de topos $f\colon \tE_i\rightarrow \tE_j$ est cohérent. 
\end{itemize}
Alors les morphismes $\beta$ \eqref{higgs2-tcevg18b} et $\sigma$ \eqref{higgs2-tcevg18f} sont cohérents. 
\end{prop}

Compte tenu de \ref{higgs2-tcevg188}, procédant comme dans la preuve de \ref{higgs2-tcevg155}, on peut se réduire au cas où 
tout objet de $I$ est quasi-compact et pour tout $i\in \ob(I)$, tout objet de $E_i$ est quasi-compact.
Par suite, tout objet de $E$ est quasi-compact en vertu de \ref{higgs2-tcevg15}.
La proposition résulte alors de (\cite{sga4} VI 3.3). 

\subsection{}\label{higgs2-tcevg30}
Soient $\pi'\colon E'\rightarrow I$ un $\mU$-site fibré, clivé et normalisé vérifiant les conditions de \eqref{higgs2-tcevg1}, 
\begin{equation}\label{higgs2-tcevg30a}
\Phi^+\colon E'\rightarrow E
\end{equation}
un $I$-foncteur cartésien. On munit $E'$ de la topologie co-évanescente 
définie par $\pi'$, et on note $\tE'$ le topos des faisceaux de $\mU$-ensembles sur $E'$.
On associe à $\pi'$ des objets analogues à ceux associés à $\pi$,
et on les note par les mêmes lettres affectées d'un prime $'$.  
Supposons que pour tout $i\in \ob(I)$, le foncteur $\Phi_i^+\colon E'_i\rightarrow E_i$, 
induit par $\Phi^+$ sur les catégories fibres en $i$, soit continu et exact à gauche. Celui-ci définit alors
un morphisme de topos 
\begin{equation}\label{higgs2-tcevg30b}
\Phi_i\colon \tE_i\rightarrow \tE'_i,
\end{equation}
caractérisé par $\Phi_{i*}(F)=F\circ \Phi_i$ et $\Phi_i^*$ prolonge $\Phi_i^+$  (\cite{sga4} IV 4.9.4). 
On sait \eqref{higgs2-tcevg85} que le foncteur $\Phi^+$ est continu pour les topologies co-évanescentes de $E$ et $E'$. 
Par ailleurs, il commute aux produits fibrés (cf. la preuve de \ref{higgs2-tcevg4}(i)), et
$\Phi^+(e)=\Phi_\iota^+(e)$ est un objet final de $E'_\iota$ et donc de $E'$ \eqref{higgs2-tcevg18}. 
Par suite, $\Phi^+$ est exact à gauche. Il définit donc un morphisme de topos 
\begin{equation}\label{higgs2-tcevg30d}
\Phi\colon \tE\rightarrow \tE'
\end{equation}
caractérisé par $\Phi_*(F)=F\circ \Phi^+$ et $\Phi^*$ prolonge $\Phi^+$. D'après \eqref{higgs2-tcevg85f},
pour tout faisceau $F=\{i\mapsto F_i\}$ sur $E$, on a 
\begin{equation}\label{higgs2-tcevg30e}
\Phi_*(F)=\{i\mapsto \Phi_{i*}(F_i)\}.
\end{equation}

\subsection{}\label{higgs2-tcevg70}
Soient $V$ un objet de $E$, $c=\pi(V)$. On désigne par 
\begin{equation}\label{higgs2-tcevg70a}
\varpi\colon E_{/V}\rightarrow I_{/c}
\end{equation}
le foncteur induit par $\pi$, et par 
\begin{equation}\label{higgs2-tcevg70b}
\gamma_V\colon E_{/V}\rightarrow E
\end{equation}
le foncteur canonique. Pour tout morphisme $f\colon i\rightarrow c$ de $I$, la catégorie fibre de $\varpi$
au-dessus de l'objet $(i,f)$ de $I_{/c}$ est canoniquement équivalente à la catégorie $(E_i)_{/f^+(V)}$. 
Munissant $I_{/c}$ de la topologie induite par celle de $I$,
et chaque fibre $(E_i)_{/f^+(V)}$ de la topologie induite par celle de $E_i$, 
$\varpi$ devient un site fibré vérifiant les conditions de \eqref{higgs2-tcevg1}. On munit alors $E_{/V}$
de la topologie co-évanescente relative à $\varpi$.

\begin{prop}\label{higgs2-tcevg71}
Sous les hypothèses de \eqref{higgs2-tcevg70}, la topologie co-évanescente de $E_{/V}$ est induite par celle de $E$
au moyen du foncteur $\gamma_V$. En particulier, le topos des faisceaux de $\mU$-ensembles sur $E_{/V}$ 
est canoniquement équivalent à $\tE_{/\varepsilon(V)}$. Le foncteur restriction à $\tE_{/\varepsilon(V)}$
est isomorphe au foncteur
\begin{equation}\label{higgs2-tcevg71a}
\tE\rightarrow \tE_{/\varepsilon(V)}, \ \ \ \{i\mapsto F_i\}\mapsto \{(i,f)\mapsto F_i\times f^*(V)\},
\end{equation}
où $(i,f)$ désigne un objet de $I_{/c}$, autrement dit, $f\colon i\rightarrow c$ est un morphisme de $I$. 
\end{prop}

Soit $F=\{i\mapsto F_i\}$ un faisceau sur $E$. 
D'après (\cite{sga4} III 5.4), on a un isomorphisme canonique de préfaisceaux sur $E_{/V}$ 
\begin{equation}\label{higgs2-tcevg71b}
F\circ \gamma_V \stackrel{\sim}{\rightarrow} \{(i,f)\mapsto F_i\times f^*(V)\}.
\end{equation}
Pour tout morphisme $g\colon (i,f)\rightarrow (j,h)$ de $I_{/c}$, le diagramme de morphismes de topos 
\begin{equation}\label{higgs2-tcevg71c}
\xymatrix{
{(\tE_i)_{/f^*(V)}}\ar[r]^{g_\varpi}\ar[d]&{(\tE_j)_{/h^*(V)}}\ar[d]\\
{\tE_i}\ar[r]^{g_\pi}&{\tE_j}}
\end{equation}
où $g_\varpi$ et $g_\pi$ désignent les morphismes de topos associés à $g$ par $\varpi$ et $\pi$, respectivement,
et les flèches verticales sont les morphismes de localisation, est commutatif à isomorphisme canonique près 
(\cite{sga4} IV 5.10). Le morphisme de transition
\begin{equation}\label{higgs2-tcevg71d}
F_j\times h^*(V)\rightarrow g_{\varpi *}(F_i\times f^*(V))
\end{equation}
du préfaisceau $F\circ \gamma_V$ est le composé 
\begin{equation}\label{higgs2-tcevg71e}
F_j\times h^*(V)\rightarrow g_{\pi*}(F_i)\times h^*(V)\stackrel{\sim}{\rightarrow} g_{\varpi *}(F_i\times f^*(V)),
\end{equation}
où la première flèche est induite par le morphisme de transition de $F$, et la seconde flèche
est le morphisme de changement de base relativement à \eqref{higgs2-tcevg71c} (\cite{egr1} (1.2.2.2)),
qui est en fait un isomorphisme. 

Soient $(i,f)$ un objet de $I_{/c}$, $(g_n\colon i_n\rightarrow i)_{n\in \Sigma}$ une famille couvrante de $I$.
Pour tout $(m,n)\in \Sigma^2$, on pose $i_{mn}=i_m\times_ii_n$
et on note $g_{mn}\colon i_{mn} \rightarrow i$ le morphisme canonique. Posons 
$f_n=f\circ g_n$ et $f_{mn}=f\circ g_{mn}$. Comme le foncteur restriction $\tE_i\rightarrow (\tE_i)_{/f^*(V)}$ 
admet un adjoint à gauche, il commute aux limites projectives. 
Par suite, compte tenu de \eqref{higgs2-tcevg71e}, la suite exacte de morphismes de faisceaux sur $E_i$
\begin{equation}\label{higgs2-tcevg71f}
F_i\rightarrow \prod_{n\in \Sigma}(g_{n})_{\pi *}(F_{i_n})\rightrightarrows \prod_{(m,n)\in \Sigma^2} (g_{mn})_{\pi *}(F_{i_{mn}})
\end{equation}
induit une suite exacte de morphismes de faisceaux sur $(\tE_i)_{/f^*(V)}$
\begin{equation}\label{higgs2-tcevg71g}
F_i\times f^*(V)\rightarrow \prod_{n\in \Sigma}(g_{n})_{\varpi *}(F_{i_n}\times f_n^*(V))\rightrightarrows 
\prod_{(m,n)\in \Sigma^2} (g_{mn})_{\varpi*}(F_{i_{mn}}\times f_{mn}^*(V)).
\end{equation}
On en déduit que $F\circ \gamma_V$ est un faisceau pour la topologie co-évanescente de $E_{/V}$. Donc $\gamma_V$
est continu, autrement dit, la topologie co-évanescente de $E_{/V}$ est moins fine que sa topologie induite 
par celle de $E$. Pour montrer la première proposition, il suffit donc de montrer que $\gamma_V$ est cocontinu 
(\cite{sga4} III 2.1). 

Le foncteur $\gamma_V$ est un adjoint à gauche du foncteur 
\begin{equation}\label{higgs2-tcevg71h}
\Phi^+\colon E\rightarrow E_{/V}, \ \ \ W\mapsto W\times V.
\end{equation}
Soit $G=\{(i,f)\mapsto G_{(i,f)}\}$ un faisceau sur $E_{/V}$ (où $(i,f)\in \ob(I_{/c})$). Pour tout $i\in \ob(I)$, 
on note $p_i\colon i\times c\rightarrow c$ et $q_i\colon i\times c\rightarrow i$ les projections canoniques, et 
\begin{equation}\label{higgs2-tcevg71i}
q'_i\colon (\tE_{i\times c})_{/p_i^*(V)}\rightarrow \tE_i
\end{equation}
le morphisme composé
\[
(\tE_{i\times c})_{/p_i^*(V)}\longrightarrow \tE_{i\times c} \stackrel{q_i}{\longrightarrow} \tE_i,
\]
où la première flèche est le morphisme de localisation. On a alors un isomorphisme canonique 
\begin{equation}\label{higgs2-tcevg71j}
G\circ \Phi^+\simeq \{i\mapsto q'_{i*}(G_{i\times c})\}.
\end{equation}
Pour tout morphisme $f\colon j\rightarrow i$ de $I$, le morphisme de transition
\begin{equation}\label{higgs2-tcevg71k}
q'_{i*}(G_{i\times c})\rightarrow f_{\pi *}(q'_{j*}(G_{j\times c}))
\end{equation}
est le composé 
\[
q'_{i*}(G_{i\times c})\rightarrow q'_{i*}((f\times\id_c)_{\varpi *}(G_{j\times c}))\stackrel{\sim}{\rightarrow}f_{\pi*}
(q'_{j*}(G_{j\times c})),
\]
où la première flèche provient du morphisme de transition de $F$ et la seconde flèche est l'isomorphisme 
canonique \eqref{higgs2-tcevg71c}. Soit $(i_n\rightarrow i)_{n\in \Sigma}$ une famille couvrante de $I$. 
Comme le foncteur $q'_{i*}$
commute aux limites projectives, la relation de recollement \eqref{higgs2-tcevg5a} 
pour $G$ relativement au recouvrement $(i_n\times c\rightarrow i\times c)_{n\in \Sigma}$
de $I_{/c}$  implique la relation analogue pour $G\circ \Phi^+$ relativement au recouvrement
$(i_n\rightarrow i)_{n\in \Sigma}$. On en déduit que 
$G\circ \Phi^+$ est un faisceau sur $E$ d'après \ref{higgs2-tcevg5}; donc $\Phi^+$ est continu.
Par suite, $\gamma_V$ est cocontinu en vertu de (\cite{sga4} III 2.5), d'où la première proposition.
La deuxième proposition s'ensuit en vertu de (\cite{sga4} III 5.4). La dernière proposition
résulte de ce qui précède, en particulier de \eqref{higgs2-tcevg71b}.

\subsection{}\label{higgs2-tcevg20}
Soient $c$ un objet de $I$, $e_c=\sigma^+(c)=f_c^+(e)$ \eqref{higgs2-tcevg18h}. 
Le site fibré 
\begin{equation}\label{higgs2-tcevg20a}
\varpi\colon E_{/e_c}\rightarrow I_{/c}
\end{equation} 
est déduit du site fibré $\pi$ par changement de base par le foncteur canonique 
$I_{/c}\rightarrow I$. Munissant $I_{/c}$ de la topologie induite par celle de $I$,
le site fibré $\varpi$ vérifie alors les conditions de \eqref{higgs2-tcevg1}. En vertu de \ref{higgs2-tcevg71}, 
la topologie co-évanescente de $E_{/e_c}$ est induite par celle de $E$ au moyen du foncteur
canonique $\gamma_{c}\colon E_{/e_c}\rightarrow E$. 
En particulier, le topos des faisceaux de $\mU$-ensembles sur $E_{/e_c}$ 
est canoniquement équivalent à $\tE_{/\sigma^*(c)}$.
L'objet final $\id_c$ de $I_{/c}$ définit alors un morphisme de topos \eqref{higgs2-tcevg18b}
\begin{equation}\label{higgs2-tcevg23a}
\beta_c\colon \tE_{/\sigma^*(c)}\rightarrow \tE_c.
\end{equation}
Notons encore
\begin{equation}\label{higgs2-tcevg23b}
\gamma_c\colon \tE_{/\sigma^*(c)}\rightarrow \tE
\end{equation}
le morphisme de localisation de $\tE$ en $\sigma^*(c)$.
D'après \ref{higgs2-tcevg71} et \eqref{higgs2-tcevg18c}, pour tout faisceau $F=\{i\mapsto F_i\}$ sur $E$, 
on a un isomorphisme canonique 
\begin{equation}\label{higgs2-tcevg23c}
\beta_{c*}(\gamma_c^*(\{i\mapsto F_i\}))\stackrel{\sim}{\rightarrow}F_c.
\end{equation}

\section{Morphismes à valeurs dans un topos co-évanescent généralisé}\label{higgs2-fccp}

\subsection{}\label{higgs2-fccp1}
Les notations et conventions du § \ref{higgs2-tcevg}, en particulier, celles introduites dans \ref{higgs2-tcevg18},
sont en vigueur dans cette section. Plus précisément, $I$ désigne un $\mU$-site,
et $\pi\colon E\rightarrow I$ un $\mU$-site fibré, clivé et normalisé,
tels que les conditions suivantes soient satisfaites~:
\begin{itemize}
\item[(i)] Les limites projectives finies sont représentables dans $I$. 
\item[(ii)]  Pour tout $i\in \ob(I)$, les limites projectives finies sont représentables dans $E_i$. 
\item[(iii)] Pour tout morphisme $f\colon i\rightarrow j$ de $I$, le foncteur image inverse 
$f^+\colon E_j\rightarrow E_i$ est continu et exact à gauche. 
\end{itemize}
On munit $E$ de la topologie co-évanescente définie par $\pi$, et on note $\tE$ le topos des faisceaux
de $\mU$-ensembles sur $E$. On se donne, de plus, un $\mU$-site $X$ et un foncteur continu et exact à gauche
\begin{equation}\label{higgs2-fccp1a}
\Psi^+\colon E\rightarrow X.
\end{equation}
On désigne par $\tX$ le topos des faisceaux de $\mU$-ensembles sur $X$ et par 
\begin{equation}\label{higgs2-fccp1b}
\Psi\colon \tX\rightarrow \tE
\end{equation}
le morphisme de topos associé à $\Psi^+$  (\cite{sga4} IV 4.9.2). On pose 
\begin{equation}\label{higgs2-fccp1c}
u^+=\Psi^+\circ \sigma^+\colon I \rightarrow X,\\
\end{equation}
où $\sigma^+$ est le foncteur défini dans \eqref{higgs2-tcevg18d}. 
Le foncteur $u^+$ étant continu et exact à gauche, on note
\begin{equation}\label{higgs2-fccp1d}
u=\sigma\Psi\colon \tX \rightarrow \tI
\end{equation}
le morphisme de topos associé. 

Soit $i$ un objet de $I$. Comme $\sigma^+(i)$ est un objet final de $E_i$, $\Psi^+$
induit un foncteur 
\begin{equation}\label{higgs2-fccp1aa}
\Psi_i^+\colon E_i\rightarrow X_{/u^+(i)}. 
\end{equation}
Celui-ci commute aux produits fibrés \eqref{higgs2-tcevg4} et transforme objet final en objet final~; il est donc 
exact à gauche. D'autre part, lorsque l'on munit $X_{/u^+(i)}$ de la topologie induite par celle de $X$,
$\Psi_i^+$ transforme famille couvrante en famille couvrante d'après (\cite{sga4} III 1.6 et 5.2);
il est donc continu. On désigne par
\begin{equation}\label{higgs2-fccp1ac}
\Psi_i\colon \tX_{/u^*(i)}\rightarrow \tE_i
\end{equation}
le morphisme de topos défini par $\Psi_i^+$, et par
\begin{equation}\label{higgs2-fccp1ab}
\jmath_i\colon \tX_{/u^*(i)}\rightarrow \tX
\end{equation}
le morphisme de localisation de $\tX$ en $u^*(i)$ (\cite{sga4} IV 5.1). Pour tout morphisme $f\colon i'\rightarrow i$ de $I$, 
on note 
\begin{equation}\label{higgs2-fccp1ba}
\jmath_f\colon \tX_{/u^*(i')}\rightarrow \tX_{/u^*(i)}
\end{equation}
le morphisme de localisation associé à $u^+(f)\colon u^+(i')\rightarrow u^+(i)$ (\cite{sga4} IV  5.5). 

Rappelons que nous avons fixé un objet final $\iota$ de $I$ \eqref{higgs2-tcevg18}. Comme $u^+(\iota)$ est un objet final de $X$, 
$\Psi_\iota$ n'est autre que le morphisme composé
\begin{equation}\label{higgs2-fccp1ad}
\Psi_\iota=\beta\Psi\colon \tX\rightarrow \tE_\iota.
\end{equation}
Pour tout $i\in \ob(I)$, on note $f_i\colon i\rightarrow \iota$ est le morphisme canonique \eqref{higgs2-tcevg18h}.
On identifie $\jmath_{i}$ et $\jmath_{f_i}$.

Des relations $\Psi^*\sigma^* =u^*$ et $\Psi^*\beta^*=\Psi_\iota^*$, on obtient par adjonction des morphismes 
\begin{eqnarray}
\sigma^*&\rightarrow& \Psi_* u^*, \label{higgs2-fccp1e}\\
\beta^*&\rightarrow&\Psi_* \Psi_\iota^*. \label{higgs2-fccp1f}
\end{eqnarray}

\begin{lem}\label{higgs2-fccp2}
\begin{itemize}
\item[{\rm (i)}] Pour tout morphisme $f\colon i'\rightarrow i$ de $I$, le diagramme de morphismes de topos
\begin{equation}\label{higgs2-fccp2a}
\xymatrix{
{\tX_{/u^*(i')}}\ar[r]^{\Psi_{i'}}\ar[d]_{\jmath_f}&{\tE_{i'}}\ar[d]^{f}\\
{\tX_{/u^*(i)}}\ar[r]^{\Psi_i}&{\tE_i}}
\end{equation}
est commutatif à un isomorphisme canonique près
\begin{equation}\label{higgs2-fccp2b}
\Psi_i \jmath_f\stackrel{\sim}{\rightarrow} f \Psi_{i'}.
\end{equation}

\item[{\rm (ii)}] Pour tout objet $F$ de $\tX$, on a un isomorphisme canonique de $\tE$
\begin{equation}\label{higgs2-fccp2c}
\Psi_*(F)\stackrel{\sim}{\rightarrow} \{i\mapsto \Psi_{i*}(\jmath_i^*F)\},
\end{equation}
où pour tout morphisme $f\colon i'\rightarrow i$ de $I$, le morphisme de transition 
\begin{equation}\label{higgs2-fccp2d}
\Psi_{i*}(j_i^*F)\rightarrow f_*(\Psi_{i'*}(\jmath_{i'}^*F))
\end{equation}
est le composé 
\begin{equation}\label{higgs2-fccp2e}
\Psi_{i*}(\jmath_i^*F)\rightarrow \Psi_{i*}(\jmath_{f*}(\jmath^*_f(\jmath_i^*F)))\stackrel{\sim}{\rightarrow} 
\Psi_{i*}(\jmath_{f*}(\jmath_{i'}^*F))
\stackrel{\sim}{\rightarrow} f_*(\Psi_{i'*}(\jmath_{i'}^*F)),
\end{equation}
dans lequel le premier morphisme provient du morphisme d'adjonction $\id\rightarrow \jmath_{f*}\jmath^*_f$, 
le second morphisme est l'isomorphisme canonique et le dernier morphisme provient de \eqref{higgs2-fccp2b}.

\item[{\rm (iii)}]Pour tout faisceau $F$ de $\tE_\iota$, le morphisme d'adjonction 
$\beta^*(F)\rightarrow \Psi_*(\Psi^*_\iota F)$ \eqref{higgs2-fccp1f} est induit par le morphisme de préfaisceaux sur $E$
\begin{equation}\label{higgs2-fccp2f}
\{i\mapsto f_i^*F\}\rightarrow \{i\mapsto \Psi_{i*}(\jmath_i^*(\Psi_\iota^*F))\}
\end{equation}
défini, pour tout $i\in \ob(I)$, par le morphisme de $\tE_i$
\begin{equation}\label{higgs2-fccp2g}
f_i^*F \rightarrow \Psi_{i*}(\jmath_i^*(\Psi_\iota^*F))
\end{equation}
adjoint de l'isomorphisme \eqref{higgs2-fccp2b}. 
\end{itemize}
\end{lem}

(i) Pour tout objet $V$ de $E_i$, on a un isomorphisme canonique dans $E$
\begin{equation}\label{higgs2-fccp2h}
f^+(V)\stackrel{\sim}{\rightarrow}V\times_{\sigma^+(i)}\sigma^+(i').
\end{equation}
Comme $\Psi^+$ est exact à gauche, on en déduit que le diagramme de morphismes de foncteurs 
\begin{equation}\label{higgs2-fccp2i}
\xymatrix{
{E_i}\ar[r]^-(0.5){\Psi_i^+}\ar[d]_{f^+}&{X_{/u^+(i)}}\ar[d]^{\jmath^+_f}\\
{E_{i'}}\ar[r]^-(0.5){\Psi_{i'}^+}&{X_{/u^+(i')}}}
\end{equation}
où $\jmath_f^+$ est le foncteur de changement de base par $u^+(f)$, est commutatif à isomorphisme canonique près.  
La proposition s'ensuit compte tenu de l'interprétation du foncteur $\jmath_f^*$ comme un foncteur 
de changement de base par le morphisme $u^*(f)$ (\cite{sga4} III 5.4). 

(ii) Cela résulte aussitôt des définitions.

(iii) Posons $\beta^*(F)=\{i\mapsto G_i\}$. D'après (ii), le morphisme 
$\beta^*(F)\rightarrow \Psi_*(\Psi^*_\iota F)$ \eqref{higgs2-fccp1f} est défini, pour tout $i\in \ob(I)$, 
par un morphisme de $\tE_i$ 
\begin{equation}\label{higgs2-fccp2j}
s_i\colon G_i\rightarrow  \Psi_{i*}(\jmath_i^*(\Psi_\iota^*F)).
\end{equation}
En vertu de \ref{higgs2-tcevg19}(i), on a un morphisme canonique de préfaisceaux sur $E$
\begin{equation}\label{higgs2-fccp2k}
\{i\mapsto f_i^*F\}\rightarrow \{i\mapsto G_i\},
\end{equation}
défini, pour tout $i\in \ob(I)$, 
par un morphisme $t_i\colon f_i^*F\rightarrow  G_i$ de $\tE_i$. Le diagramme
\begin{equation}\label{higgs2-fccp2l}
\xymatrix{
{f_i^*F}\ar[rr]^{f_i^*(t_\iota)}\ar[rrd]_{t_i}&&{f_i^*(G_\iota)}\ar[rr]^-(0.5){f_i^*(s_\iota)}\ar[d]&&
{f_i^*(\Psi_{\iota*}(\Psi_\iota^*F))}\ar[d]\\
&&{G_i}\ar[rr]^-(0.5){s_i}&&{\Psi_{i*}(\jmath_i^*(\Psi_\iota^*F))}}
\end{equation}
où les flèches verticales sont les adjoints des morphismes de transition, est commutatif.  
D'une part, $s_\iota$ s'identifie au morphisme 
\begin{equation}\label{higgs2-fccp2m}
\beta_*(\beta^*F)\rightarrow \beta_*(\Psi_*(\Psi_\iota^*F))
\end{equation}
induit par \eqref{higgs2-fccp1f}. D'autre part, $t_\iota$ s'identifie au morphisme d'adjonction $F\rightarrow \beta_*(\beta^*F)$ 
en vertu de \ref{higgs2-tcevg19}(i). Par suite, compte tenu de la définition de \eqref{higgs2-fccp1f},  
$s_\iota\circ t_\iota$ est le morphisme d'adjonction $F\rightarrow \Psi_{\iota*}(\Psi_\iota^*F)$. 
Il résulte alors de (ii) et \eqref{higgs2-fccp2l} que le morphisme 
$s_i\circ t_i$ est le composé 
\begin{equation}
f_i^*F\rightarrow f_i^*(\Psi_{\iota*}(\Psi_\iota^*F))\rightarrow \Psi_{i*}(\jmath_{i}^*(\Psi_\iota^*F)),
\end{equation}
où la première flèche provient du morphisme d'adjonction $\id \rightarrow \Psi_{\iota*}\Psi_\iota^*$
et la seconde flèche est le morphisme de changement de base relativement au diagramme \eqref{higgs2-fccp2a} (pour $f=f_i$). 
D'après (\cite{sga4} XVII 2.1.3), l'adjoint de $s_i\circ t_i$ est le composé 
\begin{equation}
\Psi_{i}^*(f_i^*F)\rightarrow \Psi_{i}^*(f_i^*(\Psi_{\iota*}(\Psi_\iota^*F)))\stackrel{\sim}{\longrightarrow} 
\jmath_i^*(\Psi_\iota^*(\Psi_{\iota*}(\Psi_\iota^*F)))
\rightarrow \jmath_{i}^*(\Psi_\iota^*F),
\end{equation}
où la première flèche provient du morphisme d'adjonction $\id \rightarrow \Psi_{\iota*}\Psi_\iota^*$, 
la seconde flèche est l'isomorphisme \eqref{higgs2-fccp2b} et la troisième flèche provient du 
morphisme d'adjonction $\Psi_{\iota}^*\Psi_{\iota*}\rightarrow \id$. Ce morphisme composé est égal au composé 
\begin{equation}
\Psi_{i}^*(f_i^*F)\stackrel{\sim}{\longrightarrow} 
\jmath_i^*(\Psi_\iota^*F)\stackrel{u}{\longrightarrow} \jmath_i^*(\Psi_\iota^*(\Psi_{\iota*}(\Psi_\iota^*F)))
\stackrel{v}{\longrightarrow} \jmath_i^*(\Psi_\iota^*F),
\end{equation}
où la première flèche est l'isomorphisme \eqref{higgs2-fccp2b},
$u$ provient du morphisme d'adjonction $\id \rightarrow \Psi_{\iota*}\Psi_\iota^*$
et $v$ provient du morphisme d'adjonction $\Psi_{\iota}^*\Psi_{\iota*}\rightarrow \id$. 
Comme $v\circ u$ est l'identité, la proposition s'ensuit.

\begin{prop}\label{higgs2-fccp4}
Supposons que, pour tout $i\in \ob(I)$, 
le morphisme d'adjonction $\id\rightarrow \Psi_{i*}\Psi_i^*$ soit un isomorphisme. Alors~:
\begin{itemize}
\item[{\rm (i)}] Pour tout faisceau $F$ de $\tE_\iota$, $\beta^*(F)$ est le faisceau sur $E$ défini 
par $\{i\mapsto f_i^*F\}$.
\item[{\rm (ii)}] Le morphisme d'adjonction $\id\rightarrow \beta_*\beta^*$ est un isomorphisme. 
\item[{\rm (iii)}] Le morphisme d'adjonction $\beta^*\rightarrow \Psi_*\Psi^*_\iota$ \eqref{higgs2-fccp1f} est un isomorphisme.
\end{itemize}
\end{prop}

(i) En effet, pour tout $i\in \ob(I)$, le morphisme 
\begin{equation}\label{higgs2-fccp4a}
f_i^*F \rightarrow \Psi_{i*}(\jmath_i^*(\Psi_\iota^*F))
\end{equation}
défini dans \eqref{higgs2-fccp2g}, est le composé 
\[
f_i^*F \rightarrow \Psi_{i*} (\Psi_{i}^*(f_i^*F))\stackrel{\sim}{\rightarrow}
\Psi_{i*}(\jmath_i^*(\Psi_{\iota}^*F)),
\]
où le premier morphisme provient du morphisme d'adjonction $\id\rightarrow \Psi_{i*} \Psi_{i}^*$,
et le second de \eqref{higgs2-fccp2b}. C'est donc un isomorphisme.
Par suite, le morphisme \eqref{higgs2-fccp2f} est un isomorphisme de  $\hE$
\begin{equation}\label{higgs2-fccp4b}
\{i\mapsto f_i^*F\}\stackrel{\sim}{\rightarrow} \{i\mapsto \Psi_{i*}(\jmath_i^*(\Psi_{\iota}^*F))\}.
\end{equation}
Comme le but de ce morphisme est un faisceau sur $E$ \eqref{higgs2-fccp2c}, il en est de même de sa source. 
La proposition résulte alors de \ref{higgs2-tcevg19}(i). 

(ii) Cela résulte aussitôt de (i) et \ref{higgs2-tcevg19}(i).

(iii) Cela résulte de \ref{higgs2-fccp2}(iii) et de la preuve de (i).

\subsection{}\label{higgs2-fccp5}
Dans la suite de cette section, en plus des hypothèses générales faites dans \eqref{higgs2-fccp1}, 
on suppose que les limites projectives finies sont représentables dans $X$.
On désigne par 
\begin{equation}\label{higgs2-fccp5a}
\varpi\colon D\rightarrow I
\end{equation}
le $\mU$-site fibré associé au foncteur $u^+\colon I\rightarrow X$ \eqref{higgs2-fccp1c} 
défini dans \ref{higgs2-tcevg40}~: la fibre de $D$ au-dessus de tout objet $i$ de $I$
est la catégorie $X_{/u^+(i)}$, et pour tout morphisme $f\colon i'\rightarrow i$ de $I$, le 
foncteur image inverse $\jmath_f^+\colon X_{/u^+(i)}\rightarrow X_{/u^+(i')}$ est le foncteur 
de changement de base par $u^+(f)$. On munit $D$ de la topologie co-évanescente 
associée à $\varpi$. On obtient ainsi le site co-évanescent associé au foncteur $u^+$ défini dans \eqref{higgs2-co-ev1}, 
dont le topos des faisceaux de $\mU$-ensembles est $\tI\gtimes_{\tI}\tX$ \eqref{higgs2-co-ev101}. 

D'après (\cite{sga1} VI 12; cf. aussi \cite{egr1} 1.1.2), les foncteurs $\Psi_i^+$ \eqref{higgs2-fccp1aa} pour tout $i\in \ob(I)$
et les isomorphismes \eqref{higgs2-fccp2b} définissent un $I$-foncteur cartésien 
\begin{equation}\label{higgs2-fccp5b}
\rho^+\colon E\rightarrow D. 
\end{equation}
Pour tout $V\in \ob(E)$, le morphisme canonique $V\mapsto \sigma^+(\pi(V))$ \eqref{higgs2-tcevg18g}
induit un morphisme $\Psi^+(V)\rightarrow u^+(\pi(V))$, et on a  
\begin{equation}\label{higgs2-fccp5c}
\rho^+(V)= (\Psi^+(V)\rightarrow \pi(V)). 
\end{equation}
Le foncteur $\rho^+$ transforme objet final en objet final et commute aux produits fibrés~; 
il est donc exact à gauche. Il est d'autre part continu en vertu de \ref{higgs2-tcevg85}. 
Il définit donc un morphisme de topos (\cite{sga4} IV 4.9.2) 
\begin{equation}\label{higgs2-fccp5d}
\rho\colon \tI\gtimes_{\tI}\tX\rightarrow \tE.
\end{equation}
Pour tout faisceau $F=\{i\mapsto F_i\}$ sur $D$, on a 
\begin{equation}\label{higgs2-fccp5e}
\rho_*(F)=\{i\mapsto \Psi_{i*}(F_i)\}. 
\end{equation}

Il résulte aussitôt des définitions que les carrés du diagramme 
\begin{equation}\label{higgs2-tcevg24f}
\xymatrix{
{\tI}\ar@{=}[d]&{\tI\gtimes_{\tI}\tX}\ar[l]_-(0.4){\rp_1}\ar[d]^{\rho}\ar[r]^-(0.5){\rp_2}&{\tX}\ar[d]^-(0.5){\Psi_\iota}\\
{\tI}&{\tE}\ar[l]_-(0.5){\sigma}\ar[r]^-(0.5){\beta}&{\tE_\iota}}
\end{equation}
et le diagramme 
\begin{equation}\label{higgs2-fccp5g}
\xymatrix{
{\tX}\ar[r]^-(0.5){\Psi_D}\ar[rd]_{\Psi}&{\tI\gtimes_{\tI}\tX}\ar[d]^{\rho}\\
&{\tE}}
\end{equation}
où $\Psi_D$ est le morphisme \eqref{higgs2-co-ev14a}, 
sont commutatifs à isomorphismes canoniques près. 

\begin{prop}\label{higgs2-fccp6}
Pour tout faisceau $F=\{i\mapsto F_i\}$ de $\tE$, 
$\rho^*(F)$ est le faisceau associé au préfaisceau sur $D$ défini par $\{i\mapsto \Psi_{i}^*(F_i)\}$,
où pour tout morphisme $f\colon i'\rightarrow i$ de $I$, le morphisme de transition 
\begin{equation}\label{higgs2-fccp6a}
\Psi_{i}^*(F_i)\rightarrow \jmath_{f*}(\Psi_{i'}^*F_{i'})
\end{equation}
est l'adjoint du morphisme composé 
\[
\jmath_f^*(\Psi_{i}^*F_i)\stackrel{\sim}{\rightarrow} \Psi_{i'}^*(f^*F_i),
\rightarrow\Psi_{i'}^*F_{i'},
\]
où la première flèche est l'isomorphisme \eqref{higgs2-fccp2b}, et 
la seconde flèche est induite par le morphisme de transition de $F$. 
\end{prop}

En effet, quitte à élargir $\mU$, on peut supposer que la catégorie $E$ est $\mU$-petite (\cite{sga4} II 3.6 et et III 1.5). 
D'après (\cite{sga4} I 5.1 et III 1.3), le faisceau $\rho^*(F)$ est le faisceau sur $D$ 
associé au préfaisceau $G$ défini, pour tout $U\in \ob(D)$, par 
\begin{equation}\label{higgs2-fccp6b}
G(U)=\underset{\underset{(V,u)\in A^\circ_{U}}{\longrightarrow}}{\lim}\ 
F(V),
\end{equation}
où $A_{U}$ est la catégorie des couples $(V,u)$ formés d'un objet $V$
de $E$ et d'un morphisme $u\colon U \rightarrow \rho^+(V)$ de $D$.
Posons $i=\varpi(U)$, considérons $U$ comme un objet de $X_{/u^+(i)}$,
et notons $B_{U}$ la catégorie des couples $(W,v)$ formés d'un objet $W$ de 
$E_i$ et d'un $u^+(i)$-morphisme $v\colon U\rightarrow \Psi_i^+(W)$ de $X$. 
Les catégories $A_{U}$ et $B_{U}$ sont clairement cofiltrantes. 
Tout objet $(W,v)$ de $B_{U}$ peut être naturellement considéré comme un objet de $A_{U}$. 
On définit ainsi un foncteur pleinement fidèle
\[
\varphi\colon B_{U}\rightarrow A_{U}.
\]
Pour tout objet $(V,u)$ de $A_{U}$, $u$ induit un morphisme 
$f\colon i\rightarrow \pi(V)$ de $I$ et un $u^+(i)$-morphisme 
$v\colon U\rightarrow \Psi^+_i(f^+V)$ de $X$, se sorte que
$(f^+(V), v)$ est un objet de $B_{U}$. De plus, on a un morphisme canonique  
$\varphi((f^+(V), v))\rightarrow (V,u)$ de $A_{U}$. 
Il résulte alors de  (\cite{sga4} I 8.1.3(c)) que le foncteur  $\varphi^\circ$ est cofinal. Par suite, 
$\varphi$ induit un isomorphisme 
\begin{equation}\label{higgs2-fccp6d}
G(U)\simeq \underset{\underset{(W,v)\in B^\circ_{U}}{\longrightarrow}}{\lim}\ F_i(W).
\end{equation}
Donc en vertu de \ref{higgs2-tcevg8} , $\rho^*(F)$ est le faisceau sur $D$ associé au préfaisceau 
$\{i\mapsto \Psi_i^*(F_i)\}$ défini par les morphismes de transition \eqref{higgs2-fccp6a}.

\section{Topos total annelé}\label{higgs2-tta}

\subsection{}\label{higgs2-tta1}
Dans cette section, $I$ désigne une catégorie équivalente à une $\mU$-petite catégorie, et 
\begin{equation}\label{higgs2-tta1a}
\pi\colon E\rightarrow I
\end{equation}
un $\mU$-site fibré clivé et normalisé  (\cite{sga4} VI 7.2.1). Pour tout $i\in \ob(I)$, on note $E_i$ la catégorie fibre de 
$E$ au-dessus de $i$, que l'on considère toujours comme munie de la topologie donnée par $\pi$, 
$\tE_i$ le topos des faisceaux de $\mU$-ensembles sur $E_i$, et
\begin{equation}\label{higgs2-tta1b}
\alpha_{i!}\colon E_i\rightarrow E
\end{equation}
le foncteur d'inclusion canonique. On désigne par 
\begin{equation}\label{higgs2-tta1c}
\cF\rightarrow I
\end{equation}
le $\mU$-topos fibré associé à $\pi$ (\cite{sga4} VI 7.2.6) et par
\begin{equation}\label{higgs2-tta1cc}
\cF^\vee\rightarrow I^\circ
\end{equation}
la catégorie fibrée obtenue en associant à tout $i\in \ob(I)$ la catégorie $\tE_i$, et à tout morphisme 
$f\colon i\rightarrow j$ de $I$ le foncteur $f_*\colon \tE_i\rightarrow \tE_j$ image directe par le morphisme 
de topos $f\colon \tE_i\rightarrow \tE_j$. 

On munit $E$ de la topologie totale relative à $\pi$ 
(\cite{sga4} VI 7.4.1); c'est un $\mU$-site (\cite{sga4} VI 7.4.3(3)). 
On note $\Top(E)$ le topos des faisceaux de $\mU$-ensembles sur $E$, dit topos total associé à $\pi$.  
Pour tout $i\in \ob(I)$, le foncteur $\alpha_{i!}$ étant cocontinu (\cite{sga4} 7.4.2), il définit un morphisme de topos 
(\cite{sga4} 4.7)
\begin{equation}\label{higgs2-tta1d}
\alpha_i\colon \tE_i\rightarrow \Top(E). 
\end{equation}
Comme, de plus, $\alpha_{i!}$ est continu, le foncteur $\alpha_i^*$ admet un adjoint à gauche noté encore
\begin{equation}\label{higgs2-tta1e}
\alpha_{i!}\colon \tE_i\rightarrow \Top(E),
\end{equation}
qui prolonge le foncteur $\alpha_{i!}\colon E_i\rightarrow E$. 
On rappelle (\cite{sga4} VI 7.4.5) que pour tout morphisme $f\colon i\rightarrow j$
de $I$, le diagramme 
\begin{equation}\label{higgs2-tta1f}
\xymatrix{
{\tE_i}\ar[d]_f\ar[r]^-(0.4){\alpha_i}&{\Top(E)}\\
{\tE_j}\ar[ru]_{\alpha_j}&}
\end{equation}
n'est pas commutatif en général, mais il existe un $2$-morphisme de topos 
\begin{equation}\label{higgs2-tta1g}
\alpha_i\rightarrow \alpha_j f,
\end{equation} 
autrement dit un morphisme de foncteurs $f^*\circ \alpha_j^*\rightarrow \alpha_i^*$, 
ou encore par adjonction un morphisme de foncteurs 
\begin{equation}\label{higgs2-tta1h}
\alpha_j^*\rightarrow f_*\circ \alpha_i^*.
\end{equation} 
Ces derniers vérifient une relation de cocycle du type (\cite{egr1} (1.1.2.2)). Ils induisent donc un foncteur
\begin{eqnarray}
\Top(E)&\rightarrow& \bHom_{I^\circ}(I^\circ,\cF^\vee)\label{higgs2-tta1i}\\
F&\mapsto& \{i\mapsto \alpha_{i}^*(F)\},\nonumber
\end{eqnarray}
qui est en fait une équivalence de catégories (\cite{sga4} VI 7.4.7).
On identifiera dans la suite $F$ à la section $\{i\mapsto F_i\}$ qui lui est associée. 

\begin{defi}\label{higgs2-tta10}
On dit qu'un objet $F$ de $\Top(E)$ est {\em cartésien} si la section $\{i\mapsto F_i\}$ de 
$\cF^\vee\rightarrow I^\circ$ qui lui correspond par l'équivalence de catégories \eqref{higgs2-tta1i}
est cartésienne, autrement dit, si pour tout morphisme $f\colon i\rightarrow j$ de $I$, 
le morphisme de transition $F_j\rightarrow f_*(F_i)$ est un isomorphisme. 
\end{defi}

Si la catégorie $I$ est cofiltrante et $\mU$-petite, alors la limite projective du topos fibré
$\cF\rightarrow I$ existe (\cite{sga4} VI 8.2.3), et la catégorie sous-jacente est canoniquement équivalente à 
la sous-catégorie des objets cartésiens de $\Top(E)$ (\cite{sga4} VI 8.2.9). 

\subsection{}\label{higgs2-tta11}
Soient $F=\{i\mapsto F_i\}$ et $G=\{i\mapsto G_i\}$ deux objets de $\Top(E)$ tels que $G$ soit cartésien. 
On a alors un isomorphisme canonique 
\begin{equation}\label{higgs2-tta11a}
\Hom_{\Top(E)}(F,G)\stackrel{\sim}{\rightarrow}\underset{\underset{i\in I}{\longleftarrow}}{\lim}\ \Hom_{E_i}(F_i,G_i),
\end{equation}
où pour tout morphisme $f\colon i\rightarrow j$ de $I$, le morphisme de transition 
\[
d_f\colon \Hom_{E_i}(F_i,G_i)\rightarrow \Hom_{E_j}(F_j,G_j)
\]
du système projectif qui apparaît dans \eqref{higgs2-tta11a} associe à tout morphisme 
$u\colon M_i\rightarrow N_i$ le morphisme $d_f(u)$ composé de 
\[
\xymatrix{
{M_j}\ar[r]&{f_*(M_i)}\ar[r]^-(0.5){f_*(u)}&{f_*(N_i)}\ar[r]^-(0.5)\sim&{N_j}}
\]
où la première flèche est le morphisme de transition de $M$ et la dernière flèche
est l'inverse de l'isomorphisme de transition $N_j\stackrel{\sim}{\rightarrow}f_*(N_i)$ de $N$.

\subsection{}\label{higgs2-tta2}
En plus des données fixées plus haut, on se donne un faisceau d'anneaux $R$ de $\Top(E)$,
autrement dit, une structure annelée sur le topos fibré $\cF$ selon la terminologie de (\cite{sga4} VI 8.6.1). 
Il revient au même de se donner, pour tout $i\in \ob(I)$, un anneau $R_i$ de $\tE_i$, 
et pour tout morphisme $f\colon i\rightarrow j$ de $I$, un homomorphisme d'anneaux
$R_j\rightarrow f_*(R_i)$, ces homomorphismes étant soumis à des relations de compatibilité \eqref{higgs2-tta1i}. 
Les morphismes de topos $f\colon \tE_i\rightarrow \tE_j$
sont donc des morphismes de topos annelés (respectivement, par $R_i$ et $R_j$). 
Nous utilisons pour les modules la notation $f^{-1}$ pour désigner l'image
inverse au sens des faisceaux abéliens et nous réservons la notation 
$f^*$ pour l'image inverse au sens des modules.

Par ailleurs, pour tout $i\in \ob(I)$, le morphisme de topos $\alpha_i\colon \tE_i\rightarrow \Top(E)$ \eqref{higgs2-tta1d} 
est un morphisme de topos annelés (respectivement, par $R_i$ et $R$). On notera que 
comme $R_i=\alpha_i^*(R)$, il n'y a pas de différence pour les modules entre l'image
inverse au sens des faisceaux abéliens et l'image inverse au sens des modules.

La donnée d'une structure de $R$-module à gauche (resp. à droite) sur un faisceau $M=\{i\mapsto M_i\}$ de $\Top(E)$ 
est équivalente à la donnée, pour tout $i\in \ob(I)$, d'une structure de $R_i$-module à gauche (resp. à droite) sur $M_i$
telle que pour tout morphisme $f\colon i\rightarrow j$ de $I$, le morphisme de transition 
$M_j\rightarrow f_*(M_i)$ soit $R_j$-linéaire. 

\begin{lem}\label{higgs2-tta8}
Soient $M=\{i\mapsto M_i\}$ un $R$-module à droite de $\Top(E)$, $N=\{i\mapsto N_i\}$ un $R$-module à gauche 
de $\Top(E)$. Alors on a un isomorphisme canonique bifonctoriel de faisceaux abéliens sur $E$
\begin{equation}
M\otimes_RN\stackrel{\sim}{\rightarrow} \{i\mapsto  M_i\otimes_{R_i}N_i\}.
\end{equation}
\end{lem}
En effet, pour tout $i\in \ob(I)$, on a un isomorphisme canonique (\cite{sga4} IV 13.4)
\begin{equation}
\alpha_i^*(M\otimes_RN)\stackrel{\sim}{\rightarrow} \alpha_i^*(M)\otimes_{\alpha_i^*(R)}\alpha_i^*(N);
\end{equation}
d'où la proposition. 

\subsection{}\label{higgs2-tta3}
Soient $\pi'\colon E'\rightarrow I$ un $\mU$-site fibré, clivé et normalisé,
\begin{equation}\label{higgs2-tta3a}
\Phi^+\colon E'\rightarrow E
\end{equation}
un $I$-foncteur cartésien. On munit $E'$ de la topologie totale 
relative à $\pi'$, et on note $\Top(E')$ le topos des faisceaux de $\mU$-ensembles sur $E'$. 
On associe à $\pi'$ des objets analogues à ceux associés à $\pi$ \eqref{higgs2-tta1},
et on les note par les mêmes lettres affectées d'un prime $'$.
Supposons que pour tout $i\in \ob(I)$, le foncteur $\Phi_i^+\colon E'_i\rightarrow E_i$, 
induit par $\Phi^+$ sur les catégories fibres en $i$, soit un morphisme du site $E_i$ dans le site $E'_i$, et notons 
\begin{equation}\label{higgs2-tta3b}
\Phi_i\colon \tE_i\rightarrow \tE'_i
\end{equation}
le morphisme de topos associé. Les morphismes $\Phi_i$ définissent un morphisme cartésien 
de topos fibrés (\cite{sga4} VI 7.1.5 et 7.1.7)
\begin{equation}\label{higgs2-tta3c}
\cF\rightarrow \cF'.
\end{equation}
En vertu de (\cite{sga4} 7.4.10), $\Phi^+$ est un morphisme du site total $E$ dans le site total $E'$.   
Il définit donc un morphisme de topos 
\begin{equation}\label{higgs2-tta3d}
\Psi\colon \Top(E)\rightarrow \Top(E')
\end{equation}
tel que pour tout faisceau $F=\{i\mapsto F_i\}$ sur $E$, on ait
\begin{equation}\label{higgs2-tta3e}
\Psi_*(F)=\{i\mapsto \Phi_{i*}(F_i)\}.
\end{equation}

Soient $R'=\{i\mapsto R'_i\}$ un anneau de $\Top(E')$, $h\colon R'\rightarrow \Psi_*(R)$ un homomorphisme d'anneaux,
de sorte que $\Psi\colon \Top(E)\rightarrow \Top(E')$ est un morphisme de topos annelés (respectivement, par $R$ et $R'$). 
La donnée de $h$ est équivalente à la donnée, pour tout $i\in \ob(I)$, d'un homomorphisme d'anneaux
$h_i\colon R'_i\rightarrow  \Phi_{i*}(R_i)$ vérifiant des relations de compatibilité.
En particulier, pour tout $i\in \ob(I)$, $\Phi_i\colon \tE_i\rightarrow \tE'_i$
est un morphisme de topos annelés (respectivement, par $R_i$ et $R'_i$).
On désigne par $\rR^n\Psi_*$ et $\rR^n\Phi_{i*}$ ($n\in \mN$) les foncteurs dérivés droits des foncteurs (cf. \ref{higgs2-not4})
\begin{eqnarray*}
\Psi_*\colon \bMod(R,\Top(E))&\rightarrow& \bMod(R',\Top(E')),\\
\Phi_{i*}\colon \bMod(R_i,\tE_i)&\rightarrow& \bMod(R'_i,\tE'_i).
\end{eqnarray*}

\begin{prop}\label{higgs2-tta4}
Les hypothèses étant celles de \eqref{higgs2-tta3}, soient, de plus, $M=\{i\mapsto M_i\}$ un $R$-module 
de $\Top(E)$, $n$ un entier $\geq 0$. On a alors un $R'$-isomorphisme canonique et fonctoriel 
\begin{equation}\label{higgs2-tta4a}
\rR^n\Psi_*(M)\stackrel{\sim}{\rightarrow}\{i\mapsto \rR^n\Phi_{i*}(M_i)\}.
\end{equation}
\end{prop}

En effet, pour tout $c\in \ob(I)$, le foncteur 
\begin{equation}\label{higgs2-tta4ab}
\bMod(R,\tE)\rightarrow \bMod(R_c,\tE_c),\ \ \ M=\{i\mapsto M_i\}\mapsto M_c=\alpha_c^*(M)
\end{equation}
est additif et exact. 
D'autre part, pour tout $R$-module injectif $M=\{i\mapsto M_i\}$, $M_c$
est flasque en vertu de (\cite{sga4} VI 8.7.2). 
La proposition résulte donc de \eqref{higgs2-tta3e}. 

Explicitons les morphismes de transition du faisceau $\{i\mapsto \rR^n\Phi_{i*}(M_i)\}$. 
Le morphisme de topos fibrés annelés $(\cF,R)\rightarrow (\cF',R')$ induit 
pour tout morphisme $f\colon i\rightarrow j$ de $I$, un diagramme de morphismes de topos annelés,
commutatif à isomorphisme canonique près (\cite{egr1} 1.2.3),
\begin{equation}\label{higgs2-tta4b}
\xymatrix{
{(\tE_i,R_i)}\ar[r]^{\Phi_i}\ar[d]_{f_E}&{(\tE'_i,R'_i)}\ar[d]^{f_{E'}}\\
{(\tE_j,R_j)}\ar[r]^{\Phi_j}&{(\tE'_j,R'_j)}}
\end{equation}
où les flèches verticales, notées habituellement $f$, 
ont été affectées d'un indice $E$ ou $E'$ pour les distinguer. 
On laissera le soin au lecteur de vérifier que le morphisme de transition associé à $f$
est le composé 
\begin{eqnarray}
\lefteqn{\rR^n\Phi_{j*}(M_j)\rightarrow \rR^n\Phi_{j*}(f_{E*}M_i)\rightarrow} \label{higgs2-tta4c}\\
&& \rR^n(\Phi_jf_E)_*(M_i)\simeq
\rR^n(f_{E'}\Phi_i)_*(M_i)\rightarrow f_{E'*}(\rR^n\Phi_{i*}(M_i)),\nonumber
\end{eqnarray}
où la première flèche provient du morphisme de transition de $M$, la deuxième et la dernière flèches
sont induites par la suite spectrale de Cartan-Leray (\cite{sga4} V 5.4) et le troisième isomorphisme 
provient de \eqref{higgs2-tta4b}.

\begin{prop}[\cite{sga4} VI 7.4.15]\label{higgs2-tta9}
Pour tout $R$-module $M=\{i\mapsto M_i\}$ de $\tE$, on a une suite spectrale canonique et fonctorielle
\begin{equation}\label{higgs2-tta9a}
\rE_2^{p,q}=\rR^p \underset{\underset{i\in I^\circ}{\longleftarrow}}{\lim}\  \rH^q(\tE_i,M_i)
\Rightarrow \rH^{p+q}(\Top(E),M).
\end{equation}
\end{prop}

Considérons le topos fibré constant $\varpi\colon \Ens\times I\rightarrow I$ de fibre le topos ponctuel $\Ens$ \eqref{higgs2-not1}. 
D'après \eqref{higgs2-tta1i}, le topos des faisceaux de $\mU$-ensembles sur le site total $\Ens\times I$ associé à $\varpi$
est équivalent à la catégorie $\hI$ des préfaisceaux de $\mU$-ensembles sur $I$. On a un morphisme 
cartésien de topos fibrés 
\begin{equation}\label{higgs2-tta9c}
\Phi\colon \cF\rightarrow \Ens\times I
\end{equation}
dont la fibre au-dessus de $i\in \ob(I)$ est le morphisme canonique de topos $\Phi_i\colon \tE_i\rightarrow \Ens$
(\cite{sga4} IV 2.2). D'après (\cite{sga4} VI 7.4.10), $\Phi$ définit un morphisme de topos 
\begin{equation}\label{higgs2-tta9d}
\Psi\colon \Top(E)\rightarrow \hI
\end{equation}
tel que pour tout faisceau $F=\{i\mapsto F_i\}$ sur $E$, on ait
\begin{equation}\label{higgs2-tta9e}
\Psi_*(F)=\{i\mapsto \Gamma(\tE_i,F_i)\}.
\end{equation}
Considérons $\Psi$ comme un morphisme de topos annelés (respectivement, par $R$ et $\Psi_*(R)$). 
D'après \ref{higgs2-tta4}, pour tout $R$-module $M=\{i\mapsto M_i\}$ 
de $\Top(E)$ et tout entier $q\geq 0$, on a un isomorphisme canonique et fonctoriel 
\begin{equation}\label{higgs2-tta9f}
\rR^q\Psi_*(M)\stackrel{\sim}{\rightarrow}\{i\mapsto \rH^q(\tE_i,M_i)\}.
\end{equation}
Par ailleurs, pour tout groupe abélien $N$ de $\hI$  et tout entier $p\geq 0$, 
on a un isomorphisme canonique fonctoriel (\cite{sga4} V 2.3.1)
\begin{equation}\label{higgs2-tta9g}
\rH^p(\hI,N)=\rR^p \underset{\underset{i\in I^\circ}{\longleftarrow}}{\lim}\  N(i).
\end{equation}
On prend alors pour \eqref{higgs2-tta9a} la suite spectrale de Cartan-Leray relative à $\Psi$ (\cite{sga4} V 5.3).

\begin{cor}\label{higgs2-tta14}
Supposons que $I$ admette un objet final $\iota$. Alors pour tout $R$-module $M=\{i\mapsto M_i\}$ de $\tE$ 
et tout entier $n\geq 0$, on a un isomorphisme canonique et fonctoriel
\begin{equation}\label{higgs2-tta14a}
\rH^n(\Top(E),M)\stackrel{\sim}{\rightarrow} \rH^n(\tE_\iota,M_\iota).
\end{equation}
\end{cor}

En effet, avec les notations de la preuve de \ref{higgs2-tta9}, le foncteur 
$\hI\rightarrow \Ens$, $N\mapsto \Gamma(\hI,N)=N(\iota)$ est exact. 
Par suite, $\rR^p \underset{\underset{i\in I^\circ}{\longleftarrow}}{\lim} =0$ pour tout $p\geq 1$. 
La proposition résulte alors de \ref{higgs2-tta9}.

\begin{cor}\label{higgs2-tta12}
Supposons que $I$ soit la catégorie filtrante définie par 
l'ensemble ordonné des entiers naturels $\mN$.  Alors pour tout $R$-module $M=\{i\mapsto M_i\}$ de $\tE$
et tout entier $n\geq 0$, on a une suite exacte canonique et fonctorielle 
\begin{equation}
0\rightarrow \rR^1\underset{\underset{i\in \mN^\circ}{\longleftarrow}}{\lim}\ \rH^{n-1}(\tE_i,M_i)\rightarrow 
\rH^n(\Top(E),M)\rightarrow \underset{\underset{i\in \mN^\circ}{\longleftarrow}}{\lim}\ \rH^{n}(\tE_i,M_i)\rightarrow 0,
\end{equation}
où l'on a posé $\rH^{-1}(\tE_n,M_n)=0$ pour tout $n\in \mN$. 
\end{cor}
Cela résulte de \ref{higgs2-tta9} et du fait que $\rR^p\underset{\underset{n\in \mN^\circ}{\longleftarrow}}{\lim} =0$ pour tout $p\geq 2$
(\cite{jannsen} 1.4 et \cite{roos} 2.1).

\begin{defi}\label{higgs2-tta5}\index{Module co-cartesien d'un topos@Module co-cartésien d'un topos!1@total annelé}
On dit qu'un $R$-module $M=\{i\mapsto M_i\}$ de $\Top(E)$ est {\em co-cartésien} si pour tout morphisme 
$f\colon i\rightarrow j$ de $I$, le morphisme de transition $f^*(M_j)\rightarrow M_i$ de $M$
est un isomorphisme, où $f^*$ désigne l'image inverse au sens des modules \eqref{higgs2-tta2}.
\end{defi}

\begin{exemple}\label{higgs2-tta6}
Soient $X$ un $\mU$-site, $\tX$ le topos des faisceaux de $\mU$-ensembles sur $X$, 
$A$ un anneau commutatif de $\tX$, $J$ un idéal de $A$. 
On note encore $\mN$ la catégorie filtrante définie par l'ensemble ordonné des entiers naturels $\mN$. 
Le topos total associé au $\mU$-site fibré constant $X\times \mN\rightarrow \mN$ de fibre $X$,
est canoniquement équivalent à la catégorie $\bHom(\mN^\circ,\tX)$ \eqref{higgs2-tta1i}, que l'on note aussi $\tX^{\mN^\circ}$.
On le munit de l'anneau $\bvA=\{n\mapsto A/J^{n+1}\}$. 
Pour qu'un $\bvA$-module $M=\{n\mapsto M_n\}$ de $\tX^{\mN^\circ}$ soit co-cartésien, 
il faut et il suffit que pour tout entier $n\geq 0$, le morphisme $M_{n+1}\otimes_A(A/J^{n+1})\rightarrow M_n$
déduit du morphisme de transition $M_{n+1}\rightarrow M_n$,
soit un isomorphisme, autrement dit, que le système projectif $(M_n)_{n\in \mN}$ 
soit $J$-adique dans le sens de (\cite{sga5} V 3.1.1). 
\end{exemple}

\subsection{}\label{higgs2-tta7}
Soient $M=\{i\mapsto M_i\}$ et $N=\{i\mapsto N_i\}$ deux $R$-modules de $\Top(E)$ 
tels que $M$ soit co-cartésien. On a alors un isomorphisme canonique 
\begin{equation}\label{higgs2-tta7a}
\Hom_{R}(M,N)\stackrel{\sim}{\rightarrow}\underset{\underset{i\in I^\circ}{\longleftarrow}}{\lim}\ \Hom_{R_i}(M_i,N_i),
\end{equation}
où pour tout morphisme $f\colon i\rightarrow j$ de $I$, le morphisme de transition 
\[
d_f\colon \Hom_{R_j}(M_j,N_j)\rightarrow \Hom_{R_i}(M_i,N_i)
\]
du système projectif qui apparaît dans \eqref{higgs2-tta7a} associe à tout $R_j$-morphisme 
$u\colon M_j\rightarrow N_j$ le morphisme $d_f(u)$ composé de 
\[
\xymatrix{
{M_i}\ar[r]^(0.4)\sim&{f^*(M_j)}\ar[r]^-(0.5){f^*(u)}&{f^*(N_j)}\ar[r]&{N_i}}
\]
où la première flèche est l'inverse de l'isomorphisme de transition $f^*(M_j)\stackrel{\sim}{\rightarrow}M_i$ 
de $M$ et la dernière flèche est le morphisme de transition de $N$.
En particulier, si $I$ admet un objet final $\iota$, on a un isomorphisme canonique 
\begin{equation}\label{higgs2-tta7b}
\Hom_{R}(M,N)\stackrel{\sim}{\rightarrow}\Hom_{R_\iota}(M_\iota,N_\iota).
\end{equation}

\section{Topos co-évanescent annelé}\label{higgs2-tcea}

\subsection{}\label{higgs2-tcea1}
Dans cette section, $I$ désigne un $\mU$-site, et $\pi\colon E\rightarrow I$ un $\mU$-site fibré, clivé et normalisé
au-dessus de la catégorie sous-jacente à $I$
tels que les conditions suivantes soient satisfaites~:
\begin{itemize}
\item[(i)] Les limites projectives finies sont représentables dans $I$. 
\item[(ii)]  Pour tout $i\in \ob(I)$, les limites projectives finies sont représentables dans $E_i$. 
\item[(iii)] Pour tout morphisme $f\colon i\rightarrow j$ de $I$, le foncteur image inverse 
$f^+\colon E_j\rightarrow E_i$ est continu et exact à gauche. 
\end{itemize}
On munit $E$ de la topologie co-évanescente définie par $\pi$, et on note $\tE$ le topos des faisceaux
de $\mU$-ensembles sur $E$. On fixe un objet final  $\iota$ de $I$ et un objet final $e$ de $E_\iota$,  
et on reprend les notations introduites dans § \ref{higgs2-tcevg}, en particulier, celles introduites dans \eqref{higgs2-tcevg18}.

On se donne aussi un anneau $R$ de $\tE$. 
D'après \ref{higgs2-tcevg5}, il revient au même de se donner, pour tout $i\in \ob(I)$, un anneau $R_i$ de $\tE_i$, 
et pour tout morphisme $f\colon i\rightarrow j$ de $I$, un homomorphisme d'anneaux
$R_j\rightarrow f_*(R_i)$, ces homomorphismes étant soumis à des relations de compatibilité \eqref{higgs2-tcevg2a}
et de recollement \eqref{higgs2-tcevg5a}. Les morphismes de topos $f\colon \tE_i\rightarrow \tE_j$
sont donc des morphismes de topos annelés (respectivement, par $R_i$ et $R_j$). 
Nous utilisons pour les modules la notation $f^{-1}$ pour désigner l'image
inverse au sens des faisceaux abéliens et nous réservons la notation 
$f^*$ pour l'image inverse au sens des modules. 

La donnée d'une structure de $R$-module sur un faisceau $M=\{i\mapsto M_i\}$ de $\tE$ 
est équivalente à la donnée pour tout $i\in \ob(I)$, d'une structure de $R_i$-module sur $M_i$
telle que pour tout morphisme $f\colon i\rightarrow j$ de $I$, le morphisme de transition 
$M_j\rightarrow f_*(M_i)$ soit $R_j$-linéaire. 
Pour tout $c\in \ob(I)$, le foncteur 
\begin{equation}\label{higgs2-tcea1a}
\bMod(R,\tE)\rightarrow \bMod(R_c,\tE_c),\ \ \ \{i\mapsto M_i\}\mapsto M_c
\end{equation}
est clairement additif.

\begin{lem}\label{higgs2-tcea100}
Soient $u\colon \{i\mapsto M_i\}\rightarrow \{i\mapsto N_i\}$ un morphisme de $R$-modules de $\tE$, 
$F$ son noyau, $G$ son conoyau. Pour tout $i\in \ob(I)$, notons $u_i\colon M_i\rightarrow N_i$ 
le $R_i$-morphisme induit par $u$, et $F_i$ (resp. $G_i$) son noyau (resp. conoyau). 
Alors $\{i\mapsto F_i\}$ est un $R$-module de $\tE$,  $\{i\mapsto G_i\}$ est un $R$-module de $\hE$,
et on a des $R$-isomorphismes canoniques
\begin{eqnarray}
F&\stackrel{\sim}{\rightarrow}&\{i\mapsto F_i\},\label{higgs2-tcea100a}\\
G&\stackrel{\sim}{\rightarrow}&\{i\mapsto G_i\}^a,\label{higgs2-tcea100b}
\end{eqnarray}
où $\{i\mapsto G_i\}^a$ est le faisceau sur $E$ associé au préfaisceau $\{i\mapsto G_i\}$. 
En particulier, le foncteur \eqref{higgs2-tcea1a} est exact à gauche. 
\end{lem}

En effet, notons $\bMod(R,\hE)$ la catégorie des $R$-modules de $\hE$, 
$j_R\colon \bMod(R,\tE)\rightarrow \bMod(R,\hE)$ le foncteur d'inclusion canonique 
et $a_R\colon \bMod(R,\hE)\rightarrow \bMod(R,\tE)$ le foncteur ``faisceau associé'' (\cite{sga4} II 6.4). 
Alors $a_R$ est exact à gauche et commute aux limites inductives, et $j_R$ commute aux limites projectives 
(\cite{sga4} II 4.1). On en déduit des $R$-isomorphismes
canoniques 
\begin{eqnarray}
j_R(F)&\stackrel{\sim}{\rightarrow}& \ker(j_R(u)),\\
G&\stackrel{\sim}{\rightarrow}& a_R(\coker(j_R(u)));
\end{eqnarray}
et de même pour les morphismes $u_i$ pour tout $i\in \ob (I)$. 
Compte tenu de \eqref{higgs2-tcevg2a} et (\cite{sga4} I 3.1), on a $\ker(j_R(u))=\{i\mapsto F_i\}$. 
D'autre part, il résulte de \ref{higgs2-tcevg8} que $\{i\mapsto G_i\}$ est un $R$-module de $\hE$ et 
qu'on a un $R$-isomorphisme canonique 
\begin{equation}
a_R(\{i\mapsto G_i\}) \stackrel{\sim}{\rightarrow} a_R(\coker(j_R(u))),
\end{equation}
d'où la proposition. 

\subsection{}\label{higgs2-tcea3}
Pour tout $c\in \ob(I)$, on désigne par 
\begin{equation}\label{higgs2-tcea3a}
\gamma_c\colon \tE_{/\sigma^*(c)}\rightarrow \tE
\end{equation}
le morphisme de localisation de $\tE$ en $\sigma^*(c)$, et par
\begin{equation}\label{higgs2-tcea3b}
\beta_c\colon \tE_{/\sigma^*(c)}\rightarrow \tE_c
\end{equation}
le morphisme défini dans \eqref{higgs2-tcevg23a}. 
D'après \eqref{higgs2-tcevg23c}, le foncteur \eqref{higgs2-tcea1a} s'identifie au foncteur composé 
$\beta_{c*}\circ \gamma_c^*$. Par suite, $\beta_c$ est un morphisme de topos annelés
respectivement par $\gamma_c^*(R)$ et $R_c$. Nous utilisons pour les modules la notation $\beta_c^{-1}$ 
pour désigner l'image inverse au sens des faisceaux abéliens et nous réservons la notation 
$\beta_c^*$ pour l'image inverse au sens des modules. On désigne par $\rR^n\beta_{c*}$ 
($n\in \mN$) les foncteurs dérivés droits du foncteur 
\begin{equation}\label{higgs2-tcea3c}
\beta_{c*}\colon \bMod(\gamma_c^*(R),\tE_{/\sigma^*(c)})\rightarrow \bMod(R_c,\tE_c).
\end{equation}
Le $n$-ième foncteur dérivé droit du foncteur \eqref{higgs2-tcea1a} s'identifie
alors au foncteur $(\rR^n\beta_{c*}) \circ \gamma_c^*$ en vertu de (\cite{sga4} V 4.11).

On appelle que $\beta=\beta_\iota$, que l'on considère donc un morphisme de topos annelés
(respectivement par $R$ et $R_\iota$). Nous utilisons pour les modules la notation $\beta^{-1}$ 
pour désigner l'image inverse au sens des faisceaux abéliens et nous réservons la notation 
$\beta^*$ pour l'image inverse au sens des modules.

\subsection{}\label{higgs2-tcea4}
Soit $\mV$ un univers  tel que $\mU\subset \mV$ et $I\in \mV$. La catégorie $E$ munie 
de la topologie totale relative au site fibré $\pi$ est un $\mV$-site~; mais ce n'est pas en général un $\mU$-site. 
On désigne par $\tE_\mV$ le topos des faisceaux de $\mV$-ensembles sur le site co-évanescent $E$, par 
\begin{equation}\label{higgs2-tcea4aa}
\jmath\colon \tE\rightarrow \tE_\mV
\end{equation}
le foncteur d'inclusion canonique
et par $\Top_\mV(E)$ le topos des faisceaux de $\mV$-ensembles sur le site total $E$ \eqref{higgs2-tta1}.
Considérons le morphisme canonique \eqref{higgs2-tcevg22a}
\begin{equation}\label{higgs2-tcea4a}
\delta\colon \tE_\mV\rightarrow \Top_\mV(E)
\end{equation} 
comme un morphisme de topos annelés (respectivement, par $R$ et $\delta_*(R)$). 
Comme $\delta$ est un plongement, le morphisme d'adjonction $\delta^*\delta_*\rightarrow \id$ est un isomorphisme. 
Il n'y a donc pas de différence pour les modules entre l'image
inverse au sens des faisceaux abéliens et l'image inverse au sens des modules.
On note $\rR^n\delta_*$ $(n\in \mN)$ les foncteurs dérivés droits du foncteur 
\begin{equation}\label{higgs2-tcea4b}
\delta_*\colon \bMod(R,\tE_\mV)\rightarrow \bMod(\delta_*(R),\Top_\mV(E)).
\end{equation}  

On rappelle que le foncteur $\jmath$ est exact et pleinement fidèle 
sur les catégories des modules et transforme les modules injectifs en modules injectifs (\cite{sga4} V 1.9). 
Par suite, le $n$-ième foncteur dérivé droit de $\delta_*\circ \jmath$ est 
canoniquement isomorphe à $(\rR^n\delta_*)\circ \jmath$.

\begin{prop}\label{higgs2-tcea5}
Les hypothèses étant celles de \eqref{higgs2-tcea4}, soit, de plus, $M$ un $R$-module de $\tE$. 
\begin{itemize}
\item[{\rm (i)}] Pour tout entier $n\geq 0$, on a un isomorphisme canonique fonctoriel de $\delta_*(R)$-modules
\begin{equation}
\rR^n\delta_*(M)\stackrel{\sim}{\rightarrow}\{i\mapsto \rR^n\beta_{i*}(\gamma_i^*(M))\}.
\end{equation}
\item[{\rm (ii)}] Si $M$ est un $R$-module injectif, alors $\delta_*(M)$ est un $\delta_*(R)$-module injectif. 
\item[{\rm (iii)}] Pour tout entier $n>0$, $\delta^*(\rR^n\delta_*(M))=0$.
\end{itemize}
\end{prop}

(i) Pour tout $i\in \ob(I)$, on désigne par $\tE_{i,\mV}$ le topos des faisceaux de $\mV$-ensembles sur $E_i$,
par $\jmath_i\colon \tE_i\rightarrow\tE_{i,\mV}$ le foncteur d'inclusion canonique et par 
\begin{equation}
\alpha_i\colon \tE_{i,\mV} \rightarrow \Top_\mV(E)
\end{equation} 
le morphisme canonique \eqref{higgs2-tta1d}.
Ce dernier est un morphisme de topos annelés (respectivement, par $R_i$ et $\delta_*(R)$);
on rappelle \eqref{higgs2-tta2} qu'il n'y a pas de différence pour les modules entre l'image
inverse au sens des faisceaux abéliens et l'image inverse au sens des modules. 
On a donc un isomorphisme canonique fonctoriel
\begin{equation}
\alpha_i^*(\rR^n\delta_*(M))\stackrel{\sim}{\rightarrow} \rR^n(\alpha_i^*\circ \delta_*)(M).
\end{equation}
D'autre part, d'après \eqref{higgs2-tcevg22b} et \eqref{higgs2-tcevg23c}, on a un isomorphisme canonique de foncteurs
\begin{equation}
\alpha_i^*\circ \delta_*\circ \jmath \stackrel{\sim}{\rightarrow}\jmath_i\circ \beta_{i*}\circ \gamma_i^*.
\end{equation}
La proposition s'ensuit compte tenu de (\cite{sga4} V 1.9 et 4.11(1)). On laissera le soin au lecteur d'expliciter 
les morphismes de transition du faisceau $\{i\mapsto \rR^n\beta_{i*}(\gamma_i^*(M))\}$ de $\Top_\mV(E)$. 

(ii) Le module $\jmath(M)$ étant injectif, on peut se borner au cas où $\mU=\mV$. 
La proposition résulte alors de (\cite{sga4} V 0.2) puisque le foncteur $\delta_*$ admet un adjoint à gauche exact 
$\delta^*=\delta^{-1}$ \eqref{higgs2-tcea4}.

(iii) Comme le foncteur $\delta^*=\delta^{-1}$ est exact, on a un isomorphisme canonique fonctoriel
\begin{equation}
\delta^*(\rR^n\delta_*(M))\stackrel{\sim}{\rightarrow} \rR^n(\delta^*\circ \delta_*)(M).
\end{equation}
D'autre part, $\delta$ étant un plongement, 
le morphisme d'adjonction $\delta^*\delta_*\rightarrow \id$ est un isomorphisme~; d'où la proposition.  

\begin{rema}\label{higgs2-tcea22}
Pour tout $R$-module $M$ de $\tE$, la suite spectrale de Cartan-Leray associée à $\delta$
qui calcule la cohomologie de $M$ (\cite{sga4} V 5.3) se réduit à celle associée à $\beta$ 
\begin{equation}
\rE_2^{p,q}=\rH^p(\tE_\iota,\rR^q\beta_*(M))\Rightarrow \rH^{p+q}(\tE,M), 
\end{equation}
en vertu de \ref{higgs2-tta14} et \ref{higgs2-tcea5}. 
\end{rema}

\subsection{}\label{higgs2-tcea2}
Soient $\pi'\colon E'\rightarrow I$ un $\mU$-site fibré, clivé et normalisé vérifiant les conditions (ii) et (iii) de \ref{higgs2-tcea1},
\begin{equation}\label{higgs2-tcea2a}
\Phi^+\colon E'\rightarrow E
\end{equation}
un $I$-foncteur cartésien. On munit $E'$ de la topologie co-évanescente 
définie par $\pi'$, et on note $\tE'$ le topos des faisceaux de $\mU$-ensembles sur $E'$. 
On associe à $\pi'$ des objets analogues à ceux associés à $\pi$, 
et on les note par les mêmes lettres affectées d'un prime $'$. 
Supposons que pour tout $i\in \ob(I)$, le foncteur $\Phi_i^+\colon E'_i\rightarrow E_i$, 
induit par $\Phi^+$ sur les catégories fibres en $i$, soit continu et exact à gauche, et notons 
\begin{equation}\label{higgs2-tcea2b}
\Phi_i\colon \tE_i\rightarrow \tE'_i
\end{equation}
le morphisme de topos associé. En vertu de \ref{higgs2-tcevg30}, $\Phi^+$ est continu et exact à gauche.
Il définit donc un morphisme de topos 
\begin{equation}\label{higgs2-tcea2c}
\Phi\colon \tE\rightarrow \tE'
\end{equation}
tel que pour tout faisceau $F=\{i\mapsto F_i\}$ sur $E$, on ait 
\begin{equation}\label{higgs2-tcea2d}
\Phi_*(F)=\{i\mapsto \Phi_{i*}(F_i)\}.
\end{equation}

Soient $R'=\{i\mapsto R'_i\}$ un anneau de $\tE'$, $h\colon R'\rightarrow \Phi_*(R)$ un homomorphisme d'anneaux,
de sorte que $\Phi\colon \tE\rightarrow \tE'$ est un morphisme de topos annelés (respectivement, par $R$ et $R'$). 
La donnée de $h$ est équivalente à la donnée, pour tout $i\in \ob(I)$, d'un homomorphisme d'anneaux
$h_i\colon R'_i\rightarrow  \Phi_{i*}(R_i)$ vérifiant des relations de compatibilité \eqref{higgs2-tcevg85dd}.
En particulier, pour tout $i\in \ob(I)$, $\Phi_i\colon \tE_i\rightarrow \tE'_i$
est un morphisme de topos annelés (respectivement, par $R_i$ et $R'_i$).
On désigne par $\rR^n\Phi_*$ et $\rR^n\Phi_{i*}$ $(n\in \mN)$ les foncteurs dérivés droits des foncteurs 
\begin{eqnarray*}
\Phi_*\colon \bMod(R,\tE)&\rightarrow& \bMod(R',\tE'),\\
\Phi_{i*}\colon \bMod(R_i,\tE_i)&\rightarrow& \bMod(R'_i,\tE'_i).
\end{eqnarray*}

\begin{prop}\label{higgs2-tcea6}
Sous les hypothèses de \eqref{higgs2-tcea2}, pour tout $R$-module $M$, on a  
une suite spectrale canonique et fonctorielle 
\begin{equation}\label{higgs2-tcea6a}
\rE_2^{p,q}=\{i\mapsto \rR^p\Phi_{i*} (\rR^q\beta_{i*}(\gamma_i^*M))\}^a\Rightarrow \rR^{p+q}\Phi_*(M),
\end{equation}
où la source désigne le faisceau associé au préfaisceau $\{i\mapsto \rR^p\Phi_{i*} (\rR^q\beta_{i*}(\gamma_i^*M))\}$ sur $E'$.
\end{prop}

Quitte à élargir $\mU$, on peut se borner au cas où $I\in \mU$ (\cite{sga4} II 3.6 et V 1.9). 
On désigne par $\Top(E)$ et $\Top(E')$ les topos des faisceaux de $\mU$-ensembles sur les sites totaux $E$ 
et $E'$, respectivement  \eqref{higgs2-tta1}. Les morphismes $\Phi_i$ définissent un morphisme cartésien de topos fibrés 
$\cF\rightarrow \cF'$. 
Par ailleurs, $\Phi^+$ est un morphisme du site total $E$ dans le site total $E'$ (\cite{sga4} 7.4.10).   
Il définit donc un morphisme de topos 
\begin{equation}\label{higgs2-tcea6b}
\Psi\colon \Top(E)\rightarrow \Top(E')
\end{equation}
tel que pour tout faisceau $F=\{i\mapsto F_i\}$ sur $E$, on ait 
\begin{equation}\label{higgs2-tcea6c}
\Psi_*(F)=\{i\mapsto \Phi_{i*}(F_i)\}.
\end{equation}
En particulier, le diagramme de morphismes de topos 
\begin{equation}\label{higgs2-tcea6d}
\xymatrix{
{\tE}\ar[d]_{\Phi}\ar[r]^-(0.5){\delta} &{\Top(E)}\ar[d]^{\Psi}\\
{\tE'}\ar[r]^-(0.5){\delta'} &{\Top(E')}}
\end{equation}
est commutatif à isomorphisme canonique près.  
L'homomorphisme $h\colon R'\rightarrow \Phi_*(R)$ induit un homomorphisme 
$h'\colon \delta'_*(R')\rightarrow \Psi_*(\delta_*(R))$ qui fait de $\Psi\colon \Top(E)\rightarrow \Top(E')$ 
un morphisme de topos annelés (respectivement, par $\delta_*(R)$ et $\delta'_*(R')$). 

Comme $\delta'$ est un plongement, le morphisme d'adjonction $\delta'^*\delta'_*\Phi_*\rightarrow \Phi_*$ 
est un isomorphisme. Compte tenu de \eqref{higgs2-tcea6d}, on en déduit un isomorphisme 
$\delta'^*\Psi_*\delta_*\stackrel{\sim}{\rightarrow}\Phi_*$. Comme le foncteur $\delta'^*=\delta'^{-1}$ est exact  
\eqref{higgs2-tcea4}, la suite spectrale de Cartan-Leray (\cite{sga4} V 5.4)
induit alors une suite spectrale fonctorielle en $M$
\begin{equation}\label{higgs2-tcea6e}
\rE_2^{p,q}=\delta'^*(\rR^p\Psi_*(\rR^q\delta_*(M)))\Rightarrow \rR^{p+q}\Phi_*(M).
\end{equation}
En vertu de \ref{higgs2-tta4} et \ref{higgs2-tcea5}(i), on a un isomorphisme canonique de $\delta'_*(R')$-modules
\begin{equation}
\rR^p\Psi_*(\rR^q\delta_*(M)) \stackrel{\sim}{\rightarrow} \{i\mapsto \rR^p\Phi_{i*}(\rR^q\beta_{i*}(\gamma_i^*M))\}.
\end{equation}
D'autre part, pour tout objet $G$ de $\Top(E')$, $\delta'^*(G)$ est canoniquement isomorphe 
au faisceau sur le site co-évanescent $E'$ associé au préfaisceau $G$ (cf. \ref{higgs2-tcevg22}); d'où la proposition.

\begin{lem}\label{higgs2-tcea99}
Soient  $A=\{i\mapsto A_i\}$ un anneau de $\hE$ \eqref{higgs2-tcevg2}, 
$M=\{i\mapsto M_i\}$ un $A$-module à droite de $\hE$, $N=\{i\mapsto N_i\}$ un $A$-module à gauche 
de $\hE$, $A^a$, $M^a$ et $N^a$ les faisceaux sur $E$ associés à $A$, $M$ et $N$, respectivement.  
Pour tout $i\in \ob(I)$, on désigne par $M_i\otimes_{A_i}N_i$ le groupe abélien produit tensoriel de $M_i$
et $N_i$ dans $\hE_i$. 
Alors $\{i\mapsto M_i\otimes_{A_i}N_i\}$ est un préfaisceau sur $E$, et notant
$\{i\mapsto M_i\otimes_{A_i}N_i\}^a$ le faisceau associé sur $E$,
on a un isomorphisme canonique bifonctoriel de faisceaux abéliens sur $E$
\begin{equation}\label{higgs2-tcea99a}
M^a\otimes_{A^a}N^a\stackrel{\sim}{\rightarrow} \{i\mapsto  M_i\otimes_{A_i}N_i\}^a.
\end{equation}
\end{lem}
Cela résulte de (\cite{sga4} IV 12.10).

\begin{lem}\label{higgs2-tcea9}
Soient $M=\{i\mapsto M_i\}$ un $R$-module à droite de $\tE$, $N=\{i\mapsto N_i\}$ un $R$-module à gauche 
de $\tE$. Alors on a un isomorphisme canonique bifonctoriel de faisceaux abéliens sur $E$
\begin{equation}\label{higgs2-tcea9a}
M\otimes_{R}N\stackrel{\sim}{\rightarrow} \{i\mapsto  M_i\otimes_{R_i}N_i\}^a,
\end{equation}
où pour tout $i\in \ob(I)$, $M_i\otimes_{R_i}N_i$ est le groupe abélien produit tensoriel de $M_i$
et $N_i$ dans $\tE_i$, et $\{i\mapsto M_i\otimes_{R_i}N_i\}^a$ le faisceau sur $E$
associé  au préfaisceau $\{i\mapsto M_i\otimes_{R_i}N_i\}$. 
\end{lem}
Cela résulte de \ref{higgs2-tcevg8} et \ref{higgs2-tcea99}.

\begin{defi}\label{higgs2-tcea7}\index{Module co-cartesien d'un topos@Module co-cartésien d'un topos!2@co-évanescent annelé}
On dit qu'un $R$-module $M=\{i\mapsto M_i\}$ de $\tE$ est {\em co-cartésien} si 
le $\delta_*(R)$-module $\delta_*(M)$ est co-cartésien dans le sens de \eqref{higgs2-tta5}, autrement dit, 
si pour tout morphisme $f\colon i\rightarrow j$ de $I$, le morphisme de transition $f^*(M_j)\rightarrow M_i$ de $M$
est un isomorphisme, où $f^*$ désigne l'image inverse au sens des modules \eqref{higgs2-tcea1}.
\end{defi}

\begin{prop}\label{higgs2-tcea8}
Soient $M=\{i\mapsto M_i\}$ et $N=\{i\mapsto N_i\}$ deux $R$-modules de $\tE$ tels que $M$ soit co-cartésien. 
On a alors un $R$-isomorphisme canonique 
\begin{equation}\label{higgs2-tcea8a}
\cHom_R(M,N)\stackrel{\sim}{\rightarrow}\{i\mapsto \cHom_{R_i}(M_i,N_i)\}.
\end{equation}
\end{prop}

En effet, d'après \eqref{higgs2-tta7b}, on a un isomorphisme canonique 
\begin{equation}\label{higgs2-tcea8b}
\Hom_R(M,N)\stackrel{\sim}{\rightarrow}\Hom_{R_\iota}(M_\iota,N_\iota).
\end{equation}
Soient $V$ un objet de $E$, $c=\pi(V)$,
\begin{equation}\label{higgs2-tcea8c}
\varpi\colon E_{/V}\rightarrow I_{/c}
\end{equation}
le foncteur induit par $\pi$. Pour tout morphisme $f\colon i\rightarrow c$ de $I$, la catégorie fibre de $\varpi$
au-dessus de l'objet $(i,f)$ de $I_{/c}$ est canoniquement équivalente à la catégorie $(E_i)_{/f^+(V)}$. 
Munissant $I_{/c}$ de la topologie induite par celle de $I$,
et chaque fibre $(E_i)_{/f^+(V)}$ de la topologie induite par celle de $E_i$, 
$\varpi$ devient un site fibré vérifiant les conditions de \eqref{higgs2-tcevg1}. 
En vertu de \ref{higgs2-tcevg71}, la topologie co-évanescente de $E_{/V}$ relative à $\varpi$ est induite par celle de $E$
au moyen du foncteur canonique $\gamma_V\colon E_{/V}\rightarrow E$. En particulier, le topos des faisceaux de 
$\mU$-ensembles sur $E_{/V}$ est canoniquement équivalent à $\tE_{/\varepsilon(V)}$.
De plus, on a un isomorphisme canonique fonctoriel
\begin{equation}\label{higgs2-tcea8d}
M|\varepsilon(V)\stackrel{\sim}{\rightarrow} \{(i,f)\mapsto M_i|f^*(V)\};
\end{equation}
et de même pour $N$. Par suite, le $R|\varepsilon(V)$-module $M|\varepsilon(V)$ est co-cartésien. 
Notant $\varepsilon_c\colon E_c\rightarrow \tE_c$ le foncteur canonique,
les isomorphismes \eqref{higgs2-tcea8b} et \eqref{higgs2-tcea8d} induisent un isomorphisme 
\begin{equation}\label{higgs2-tcea8e}
\Hom_{R|\varepsilon(V)}(M|\varepsilon(V),N|\varepsilon(V))\stackrel{\sim}{\rightarrow} 
\Hom_{R_c|\varepsilon_c(V)}(M_c|\varepsilon_c(V),N_c|\varepsilon_c(V)), 
\end{equation}
qui est clairement fonctoriel en $V\in \ob(E_c)$; d'où la proposition (\cite{sga4} IV 12.1). 

\begin{rema}\label{higgs2-tcea10}
Sous les hypothèses de \ref{higgs2-tcea8}, il résulte aussitôt de \ref{higgs2-tta7} que pour tout morphisme $f\colon i\rightarrow j$ de $I$,
le morphisme de transition du $R$-module $\cHom_R(M,N)$ est le composé  
\[
f^*(\cHom_{R_j}(M_j,N_j))\rightarrow 
\cHom_{R_i}(f^*(M_j),f^*(N_j))\stackrel{\sim}{\rightarrow}\cHom_{R_i}(M_i,f^*(N_j))\rightarrow \cHom_{R_i}(M_i,N_i),
\]
où $f^*$ désigne l'image inverse au sens des modules, la première flèche est le morphisme canonique,
la deuxième flèche est induite par l'isomorphisme de transition $f^*(M_j)\stackrel{\sim}{\rightarrow} M_i$
de $M$, et la dernière flèche est induite par le morphisme de transition $f^*(N_j)\rightarrow N_i$ de $N$.  
\end{rema}

\begin{prop}\label{higgs2-tcea11}
Pour tout $R_\iota$-module localement projectif de type fini $M_\iota$ de $\tE_\iota$ (i.e.,
$M_\iota$ est localement un facteur direct d'un $R_\iota$-module libre de type fini), 
on a un isomorphisme canonique fonctoriel 
\begin{equation}\label{higgs2-tcea11a}
\beta^*(M_\iota)\stackrel{\sim}{\rightarrow} \{i\mapsto f_i^*(M_\iota)\},
\end{equation}
où $\beta^*$ désigne l'image inverse au sens des modules \eqref{higgs2-tcea3}, et pour tout $i\in \ob(I)$, 
$f_i\colon i\rightarrow \iota$ est le morphisme canonique, et $f_i^*$ 
désigne l'image inverse au sens des modules \eqref{higgs2-tcea1}. En particulier, $\beta^*(M_\iota)$ est un $R$-module 
co-cartésien. 
\end{prop}

En vertu de \ref{higgs2-tcevg19}(i), \ref{higgs2-tcea99} et  \ref{higgs2-tcevg8}, on a un isomorphisme canonique 
\begin{equation}
\beta^*(M_\iota)\stackrel{\sim}{\rightarrow} \{i\mapsto f_i^*(M_\iota)\}^a,
\end{equation}
où le terme de droite est le faisceau sur $E$ associé au préfaisceau $\{i\mapsto f_i^*(M_\iota)\}$.
Il reste à montrer que $M=\{i\mapsto f_i^*(M_\iota)\}$ est un faisceau sur $E$. 
Notons $\varepsilon_\iota\colon E_\iota\rightarrow \tE_\iota$ le foncteur canonique. 
Il existe un recouvrement vertical $(V_n\rightarrow e)_{n\in \Sigma}$ de $E$
({\em i.e.}, une famille couvrante de $E_\iota$) tel que pour tout $n\in \Sigma$, 
$M_\iota|\varepsilon_\iota(V_n)$ soit un facteur direct d'un $(R_\iota|\varepsilon_\iota(V_n))$-module libre de type fini. 
Il suffit de montrer que pour tout $n\in \Sigma$, la restriction de $M$ à la catégorie $E_{/V_n}$ 
est un faisceau pour la topologie induite par celle de $E$ (\cite{giraud2} II 3.4.4). Compte tenu de \ref{higgs2-tcevg71},
on peut se réduire au cas où $V_n=e$, auquel cas il existe un entier $d\geq 0$, un $R_\iota$-module $N_\iota$ 
de $\tE_\iota$ et un $R_\iota$-isomorphisme $R_\iota^d\stackrel{\sim}{\rightarrow} M_\iota\oplus N_\iota$. 
On a alors un isomorphisme de $\hE$
\begin{equation}
R^d\stackrel{\sim}{\rightarrow} M\oplus \{i\mapsto f_i^*(N_\iota)\},
\end{equation} 
ce qui implique que $M$ et $\{i\mapsto f_i^*(N_\iota)\}$ sont des faisceaux sur $E$.

\section{Site et topos finis étales d'un schéma}\label{higgs2-Kp}

\subsection{}\label{higgs2-not2}\index{100000911@$\Et_{/X}$, $\Et_{\coh/X}$, $\Et_{\scoh/X}$}\index{100000912@$X_\et$ (topos étale)}
Pour tout schéma $X$, on désigne par $\Et_{/X}$ le {\em site étale} de $X$, 
c'est-à-dire, la sous-catégorie pleine de $\Sch_{/X}$  \eqref{higgs2-not1} formée des schémas étales sur $X$,
munie de la topologie étale; c'est un $\mU$-site. 
On note $X_\et$ le {\em topos étale} de $X$, c'est-à-dire le topos des faisceaux de $\mU$-ensembles sur $\Et_{/X}$.
On désigne par $\Et_{\coh/X}$ (resp.
$\Et_{\scoh/X}$) la sous-catégorie pleine de $\Et_{/X}$ formée des schémas étales 
de présentation finie sur $X$ (resp. étales, séparés et de présentation finie sur $X$), 
munie de la topologie induite par celle de $\Et_{/X}$; ce sont des sites $\mU$-petits. 
Si $X$ est quasi-séparé, le foncteur de restriction de $X_\et$ dans le topos des faisceaux de $\mU$-ensembles sur 
$\Et_{\coh/X}$ (resp. $\Et_{\scoh/X}$) est une équivalence de catégories (\cite{sga4} VII 3.1 et 3.2).

\subsection{}\label{higgs2-Kp1}\index{100000921@$\Et_{\rf/X}$, $X_\fet$ (site et topos finis étales)}
\index{100000923@$\rho_X\colon X_\et\rightarrow X_\fet$}
Pour tout schéma $X$, on appelle {\em site fini étale} de $X$ et l'on note $\Et_{\rf/X}$ 
la sous-catégorie pleine de $\Et_{/X}$ formée des revêtements étales
de $X$ (c'est-à-dire, des schémas étales finis sur $X$), munie de la topologie induite par celle de $\Et_{/X}$; 
c'est un site $\mU$-petit. 
On appelle {\em topos fini étale} de $X$ et l'on note $X_\fet$ le topos des faisceaux de $\mU$-ensembles sur $\Et_{\rf/X}$. 
La topologie étale sur $\Et_{\rf/X}$ est clairement moins fine que la topologie canonique. 
L'injection canonique $\Et_{\rf/X}\rightarrow \Et_{/X}$ induit un morphisme de topos 
\begin{equation}\label{higgs2-Kp1a}
\rho_X\colon X_\et\rightarrow X_\fet.
\end{equation}

\subsection{}\label{higgs2-Kpp31}
Soit $f\colon Y\rightarrow X$ un morphisme de schémas. Le foncteur image inverse 
\begin{equation}\label{higgs2-Kpp31b}
f^\bullet\colon \Sch_{/X}\rightarrow \Sch_{/Y}, \ \ \ X'\mapsto X'\times_XY
\end{equation}
induit deux morphismes de topos que l'on note (pour les distinguer)  
\begin{eqnarray}
f_\et\colon Y_\et&\rightarrow& X_\et,\label{higgs2-Kpp31c}\\
f_\fet\colon Y_\fet&\rightarrow& X_\fet.\label{higgs2-Kpp31d}
\end{eqnarray}
On omettra les indices des notations $f_\et$ et $f_\fet$ lorsqu'il n'y a aucun risque d'ambiguïté. 
On vérifie aussitôt que le diagramme de morphismes de topos
\begin{equation}\label{higgs2-Kpp31a}
\xymatrix{
{Y_\et}\ar[r]^{f_\et}\ar[d]_{\rho_{Y}}&{X_\et}\ar[d]^{\rho_{X}}\\
{Y_\fet}\ar[r]^{f_\fet}&{X_\fet}}
\end{equation}
est commutatif à isomorphisme canonique près~: 
\begin{equation}
\rho_Xf_\et\stackrel{\sim}{\rightarrow}f_\fet\rho_Y.
\end{equation}

\subsection{}\label{higgs2-Kpp32}
Soient $f\colon Y\rightarrow X$ un revêtement étale, 
\begin{equation}\label{higgs2-Kpp32a}
\tau_f\colon \Et_{\rf/Y}\rightarrow \Et_{\rf/X}
\end{equation} 
le foncteur défini par composition à gauche avec $f$. 
Alors $\tau_f$ induit une équivalence de catégories 
$\Et_{\rf/Y}\stackrel{\sim}{\rightarrow} (\Et_{\rf/X})_{/(Y,f)}$.  
La topologie étale sur $\Et_{\rf/Y}$ est induite par celle sur $\Et_{\rf/X}$ au moyen du foncteur $\tau_f$. 
En vertu de (\cite{sga4} III 5.2), $\tau_f$ est continu et cocontinu. Il induit donc une suite de trois foncteurs adjoints~:
\begin{equation}\label{higgs2-Kpp32b}
\tau_{f!}\colon Y_\fet \rightarrow X_\fet,\ \ \ 
\tau_f^*\colon X_\fet\rightarrow Y_\fet,\ \ \ 
\tau_{f*}\colon  Y_\fet\rightarrow X_\fet,
\end{equation}
dans le sens que pour deux foncteurs consécutifs de la suite, celui de droite est
adjoint à droite de l'autre.  D'après (\cite{sga4} III 5.4), le foncteur $\tau_{f!}$ se factorise 
à travers une équivalence de catégories $(Y_\fet)\stackrel{\sim}{\rightarrow}(X_\fet)_{/Y}$. 
Le couple de foncteurs $(\tau_f^*,\tau_{f*})$ définit 
un morphisme de topos $Y_\fet\rightarrow X_\fet$, dit morphisme de localisation
de $X_\fet$ en $Y$. Par ailleurs, $\tau_{f}^*$ est adjoint à gauche du foncteur 
$f_{\fet *}$ image directe par le morphisme de topos $f_\fet$. 
Ce dernier s'identifie donc canoniquement au morphisme de localisation de $X_\fet$ en $Y$.
On renvoie à (\cite{sga4} VII 1.6) pour les énoncés analogues relatifs aux sites et topos étales.

\subsection{}\label{higgs2-Kp6}
On désigne par $\cR$ la catégorie des revêtements étales ({\em i.e.}, la sous-catégorie pleine de la catégorie 
des morphismes de $\Sch$, formée des revêtements étales) et par 
\begin{equation}\label{higgs2-Kp6a}
\cR\rightarrow \Sch
\end{equation}
le ``foncteur but'', qui fait de $\cR$ une catégorie fibrée clivée et normalisée au-dessus de $\Sch$~:  
la catégorie fibre au-dessus de tout schéma $X$ est canoniquement équivalente à
la catégorie $\Et_{\rf/X}$, et pour tout morphisme de schémas 
$f\colon Y\rightarrow X$, le foncteur image inverse $f^+\colon \Et_{\rf/X}\rightarrow \Et_{\rf/Y}$
n'est autre que le foncteur changement de base par $f$. 
Munissant chaque fibre de la topologie étale, $\cR/\Sch$ devient un $\mU$-site fibré (\cite{sga4} VI 7.2.1). 
On désigne par
\begin{equation}\label{higgs2-Kp6b}
\cG\rightarrow \Sch
\end{equation}
le $\mU$-topos fibré associé à $\cR/\Sch$ (\cite{sga4} VI 7.2.6)~: la catégorie fibre de $\cG$ 
au-dessus de tout schéma $X$ est le topos $X_\fet$, et pour tout morphisme de schémas 
$f\colon Y\rightarrow X$ le foncteur image inverse $f^*\colon X_\fet\rightarrow Y_\fet$
est le foncteur image inverse par le morphisme de topos $f_\fet$ \eqref{higgs2-Kpp31d}. On note 
\begin{equation}\label{higgs2-Kp6c}
\cG^\vee\rightarrow \Sch^\circ
\end{equation}
la catégorie fibrée obtenue en associant à tout schéma $X$ la catégorie $\cG_X=X_\fet$, et à 
tout morphisme de schémas $f\colon Y\rightarrow X$ le foncteur $f_{\fet*}\colon Y_\fet\rightarrow X_\fet$
image directe par le morphisme de topos $f_\fet$.

\begin{lem}\label{higgs2-Kpp20}
Soient $X$ un schéma, $\rho_X\colon X_\et\rightarrow X_\fet$ le morphisme canonique \eqref{higgs2-Kp1a}. 
Alors, la famille des foncteurs fibres de $X_\fet$ associés aux points $\rho_X(\ox)$, lorsque $\ox$ décrit la famille des points géométriques de $X$,  
est conservative {\rm (\cite{sga4} IV 6.4.0)}. 
\end{lem}

En effet, pour tout point géométrique $\ox$ de $X$, la donnée d'un voisinage du point $\rho_X(\ox)$ de $X_\fet$ dans le site $\Et_{\rf/X}$ est équivalente à la donnée 
d'un revêtement étale $\ox$-pointé de $X$ (\cite{sga4} IV 6.8.2 et VIII 3.9). Ces objets forment naturellement une catégorie cofiltrante, que l'on note $\cP_\ox$.
Pour tout objet $F$ de $X_\fet$, notant $F_{\ox}$ la fibre de $F$ en $\rho_X(\ox)$, on a un isomorphisme canonique fonctoriel 
\begin{equation}\label{higgs2-Kpp20a}
F_{\ox}\stackrel{\sim}{\rightarrow}\underset{\underset{(U,\xi)\in \cP_\ox^\circ}{\longrightarrow}}{\lim}\ F(U).
\end{equation}

Soit $u\colon F\rightarrow G$ un morphisme de $X_\fet$ tel que pour tout point géométrique $\ox$ de $X$, le morphisme correspondant 
$u_{\ox}\colon F_{\ox}\rightarrow G_{\ox}$ soit un monomorphisme. Montrons que $u$ est un monomorphisme.
Il s'agit de montrer que pour tous $X'\in \ob(\Et_{\rf/X})$ et $a,b \in F(X')$ tels que $u(a)=u(b)$, on a $a=b$. 
Compte tenu de \ref{higgs2-Kpp32} et \eqref{higgs2-Kpp31a}, on peut supposer $X'=X$. 
Pour tout point géométrique $\ox$ de $X$, on a $a_{\ox}=b_{\ox}$ puisque 
$u_{\ox}$ est un monomorphisme. Compte tenu de \eqref{higgs2-Kpp20a}, il existe un objet 
$(U_\ox,\zeta_{\ox})$ de $\cP_\ox$ tel que $a$ et $b$ aient même image dans $F(U_\ox)$. 
La famille des morphismes $(U_\ox\rightarrow X)_\ox$ étant clairement couvrante, on en déduit que $a=b$. Par suite, $u$ est un monomorphisme. 

Supposons, de plus, que pour tout point géométrique $\ox$ de $X$, le morphisme 
$u_{\ox}$ soit un épimorphisme et montrons qu'il en est de même de $u$. 
Il suffit de montrer que pour tous $X'\in \ob(\Et_{\rf/X})$ et $b \in G(X')$, 
il existe $a\in F(X')$ tel que $b=u(a)$. On peut encore supposer $X'=X$.
D'après \eqref{higgs2-Kpp20a}, pour tout point $\ox$ de $X$, 
il existe un objet $(U_\ox,\xi_\ox)$ de $\cP_\ox$ et une section $a_\ox\in F(U_\ox)$
dont l'image par $u$ dans $G(U_\ox)$ est la restriction de $b$. 
Comme $u$ est un monomorphisme, les sections $a_\ox$
coïncident sur $U_\ox\times_XU_{\ox'}$, pour tous les points géométriques $\ox$
et $\ox'$ de $X$. Ils proviennent donc d'une section $a\in F(X)$, et on a $u(a)=b$ puisque les restrictions aux $U_\ox$ coïncident.

\begin{defi}\label{higgs2-elc}\index{Schema@Schéma!(etale) localement connexe@(étale) localement connexe}
Soit $X$ un schéma. On dit que $X$ est {\em localement connexe} si l'espace topologique sous-jacent 
est localement connexe, et que $X$ est {\em étale-localement connexe} si pour tout 
morphisme étale $X'\rightarrow X$, toute composante connexe de $X'$ est un ensemble ouvert dans $X'$. 
\end{defi}

On peut faire les remarques suivantes~:

\subsubsection{}\label{higgs2-elc1}
Pour qu'un schéma $X$ soit étale-localement connexe, il faut et il suffit que tout $X$-schéma étale 
soit localement connexe (\cite{tg} I §11.6 prop.~11). 

\subsubsection{}\label{higgs2-elc0}
Pour que l'ensemble des composantes connexes d'un espace topologique $X$ soit localement fini, 
il faut et il suffit que toute composante connexe de $X$ soit ouverte dans $X$. 
C'est le cas si $X$ est localement connexe. 

\subsubsection{}\label{higgs2-elc2}
Soit $X$ un schéma dont l'ensemble des composantes irréductibles est localement fini
(par exemple un schéma dont l'espace topologique sous-jacent est localement noethérien). Alors~:
\begin{itemize}
\item[{\rm (i)}] Pour tout morphisme étale $f\colon X'\rightarrow X$, l'ensemble des composantes irréductibles de $X'$ 
est localement fini. En effet, la question étant locale sur $X$ et $X'$, on peut se borner au cas où $X$ et $X'$ sont affines, 
auquel cas, $f$ est quasi-compact et donc quasi-fini. L'assertion résulte alors de (\cite{ega4} 2.3.6(iii)). 

\item[{\rm (ii)}] Il résulte de (i) et (\cite{ega1n} 0.2.1.5(i)) que $X$ est étale-localement connexe.
\end{itemize}

\subsection{}\label{higgs2-Kp2}\index{Revetement universel normalise@Revêtement universel normalisé}
\index{Foncteur fibre d'un topos fini etale@Foncteur fibre d'un topos fini étale}\index{100000970@$\nu_\ox\colon X_\fet\rightarrow \bB_{\pi_1(X,\ox)}$}
Soient $X$ un schéma connexe, $\ox$ un point géométrique de $X$.
On désigne par 
\begin{equation}\label{higgs2-Kp2a}
\omega_\ox\colon \Et_{\rf/X}\rightarrow \Ens
\end{equation}
le foncteur fibre en $\ox$, qui à tout revêtement étale $Y$ de $X$ associe l'ensemble des points géométriques de 
$Y$ au-dessus de $\ox$, par $\pi_1(X,\ox)$ le groupe fondamental de $X$ en $\ox$ et 
par $\bB_{\pi_1(X,\ox)}$ le topos classifiant du groupe profini $\pi_1(X,\ox)$, 
c'est-à-dire la catégorie des $\mU$-ensembles discrets munis d'une action continue à gauche de $\pi_1(X,\ox)$ 
(\cite{sga4} IV 2.7). Alors $\omega_\ox$ induit un foncteur pleinement fidèle 
\begin{equation}\label{higgs2-Kp2aa}
\mu_\ox^+\colon \Et_{\rf/X}\rightarrow \bB_{\pi_1(X,\ox)}
\end{equation}
d'image essentielle la sous-catégorie pleine $\cC(\pi_1(X,\ox))$  
de $\bB_{\pi_1(X,\ox)}$ formée des ensembles finis (\cite{sga1} V  § 4 et § 7).
D'autre part, une famille 
$(Y_\lambda\rightarrow Y)_{\lambda\in \Lambda}$ de $\Et_{\rf/X}$ est couvrante pour la topologie étale, 
si et seulement si son image dans $\bB_{\pi_1(X,\ox)}$ est surjective, 
ou ce qui revient au même, couvrante pour la topologie canonique de $\bB_{\pi_1(X,\ox)}$. 
Par suite, la topologie étale sur $\Et_{\rf/X}$ est induite par la topologie canonique de 
$\bB_{\pi_1(X,\ox)}$ (\cite{sga4} III 3.3). 
Comme les objets de $\cC(\pi_1(X,\ox))$ forment une famille génératrice de $\bB_{\pi_1(X,\ox)}$, le foncteur 
\begin{equation}\label{higgs2-Kp2b}
\mu_\ox\colon \bB_{\pi_1(X,\ox)}\rightarrow X_\fet
\end{equation}
qui à tout objet $G$ de $\bB_{\pi_1(X,\ox)}$ (vu comme faisceau représentable) associe sa restriction à $\Et_{\rf/X}$,
est une équivalence de catégories en vertu de (\cite{sga4} IV 1.2.1).

Soit $(X_i)_{i\in I}$ un système projectif sur un ensemble ordonné filtrant $I$ dans $\Et_{\rf/X}$
qui pro-représente $\omega_\ox$, normalisé par le fait que les morphismes de transition $X_i\rightarrow X_j$
$(i\geq j)$ sont des épimorphismes et tout épimorphisme $X_i\rightarrow X'$ de $\Et_{\rf/X}$ 
est équivalent à un épimorphisme $X_i\rightarrow X_j$ $(j\leq i)$ convenable. 
Un tel pro-objet est essentiellement unique.
Il est appelé le {\em revêtement universel normalisé de $X$ en $\ox$} 
ou le {\em pro-objet fondamental normalisé de $\Et_{\rf/X}$ en $\ox$}. 
On notera que l'ensemble $I$ est $\mU$-petit.
Considérons le foncteur 
\begin{equation}\label{higgs2-Kp2c}
\nu_\ox\colon X_\fet\rightarrow \bB_{\pi_1(X,\ox)},
\ \ \ F\mapsto \underset{\underset{i\in I}{\longrightarrow}}{\lim}\ F(X_i).
\end{equation}
Par définition, la restriction de $\nu_\ox$ à $\Et_{\rf/X}$ est canoniquement isomorphe au foncteur $\mu_\ox^+$,
et on a un isomorphisme canonique de foncteurs
\begin{equation}\label{higgs2-Kp2f}
\nu_\ox\circ \mu_\ox\stackrel{\sim}{\rightarrow} \id.
\end{equation}  
Comme $\mu_\ox$ est une équivalence de catégories, $\nu_\ox$ 
est une équivalence de catégories quasi-inverse de $\mu_\ox$. 
On l'appelle {\em le foncteur fibre} de $X_\fet$ en $\ox$. 

Le foncteur $\nu_\ox$ induit une équivalence de catégories
entre la catégorie des faisceaux de groupes (resp. des faisceaux abéliens) de $X_\fet$
et la catégorie des $\pi_1(X,\ox)$-groupes discrets (resp. des $\pi_1(X,\ox)$-$\mZ$-modules discrets) 
(cf. \cite{ag1} 3.1). 
En particulier, pour tout faisceau abélien (resp. de groupes) $F$ de $X_\fet$ et tout entier $n\geq 0$
(resp. $n\in \{0,1\}$), 
on a un isomorphisme canonique fonctoriel
\begin{equation}\label{higgs2-Kp2e}
\rH^n(X_\fet,F)\simeq \rH^n(\pi_1(X,\ox),\nu_\ox(F)).
\end{equation}

\subsection{}\label{higgs2-Kpp17}
Soient $X$ un schéma  dont l'ensemble des composantes connexes est localement fini, 
$\ox$ un point géométrique de $X$, $\varphi_\ox$ le foncteur fibre de $X_\et$ associé à $\ox$, 
 $Y$ la composante connexe de $X$ contenant $\ox$, qui est alors une partie ouverte de $X$ \eqref{higgs2-elc0}.  
On munit $Y$ de la structure de schéma induite par celle de $X$
et on note $j\colon Y\rightarrow X$ l'injection canonique,
$\nu_\ox\colon  Y_\fet \rightarrow \bB_{\pi_1(Y,\ox)}$ le foncteur fibre de $Y_\fet$ en $\ox$ \eqref{higgs2-Kp2c}
et $\gamma\colon  \bB_{\pi_1(Y,\ox)} \rightarrow \Ens$
le foncteur d'oubli de l'action de $\pi_1(Y,\ox)$. Alors on a un isomorphisme canonique de foncteurs
\begin{equation}\label{higgs2-Kpp17c}
\varphi_\ox \circ \rho_X^*\stackrel{\sim}{\rightarrow} \gamma \circ  \nu_\ox \circ  j^*_\fet.
\end{equation}
En effet, pour tout $V\in \ob(\Et_{\rf/X})$, on a $\varphi_\ox(\rho_X^*V)=\varphi_\ox(V)$, et cet ensemble s'identifie 
à la fibre géométrique $V_\ox$ de $V$ au-dessus de $\ox$. 
Par suite, $V$ est un voisinage du point $\rho_X(\ox)$ de $X_\fet$ (\cite{sga4} IV 6.8.2)
si et seulement si $V_\ox$ est non vide, ou ce qui revient au même,
si et seulement si $Y$ est contenu dans l'image de $V$ dans $X$. 
Soit $(Y_i)_{i\in I}$ le revêtement universel normalisé de $Y$ en $\ox$. 
Pour tout $i\in I$, de par sa définition, $\varphi_\ox(Y_i)$ contient un élément distingué.  
On peut donc canoniquement considérer $Y_i$ comme un voisinage de $\rho_X(\ox)$ dans le site $\Et_{\rf/X}$.
De plus, il résulte de l'équivalence de catégories \eqref{higgs2-Kp2a} que  
le système $(Y_i)_{i\in I}$ est cofinal dans la catégorie (opposée de la catégorie) des voisinages 
de $\rho_X(\ox)$ dans le site $\Et_{\rf/X}$. 
L'isomorphisme \eqref{higgs2-Kpp17c} résulte alors de (\cite{sga4} IV 6.8.3), \ref{higgs2-Kpp32} et de la définition \eqref{higgs2-Kp2c}.

\subsection{}\label{higgs2-sli10}
Soient $X$ un schéma connexe, $\Pi(X)$ le groupoïde fondamental de $X$ (on rappelle qu'un groupoïde
est une catégorie dont l'ensemble des morphismes sont des isomorphismes). Les objets de $\Pi(X)$ sont 
les points géométriques de $X$. Pour tout point géométrique $\ox$ de $X$, on note 
$\omega_\ox\colon \Et_{\rf/X}\rightarrow \Ens$ le foncteur fibre correspondant \eqref{higgs2-Kp2a}.  
Si $\ox$ et $\ox'$ sont deux points géométriques de $X$, 
l'ensemble $\pi_1(X,\ox,\ox')$ des morphismes  de $\ox$ vers $\ox'$ dans $\Pi(X)$ est l'ensemble des morphismes 
(ou ce qui revient au même des isomorphismes) $\omega_{\ox}\rightarrow \omega_{\ox'}$ de foncteurs fibres associés. 
D'après (\cite{sga1} V 5.8), le foncteur 
\begin{equation}\label{higgs2-sli10a}
\mu^+\colon \Et_{\rf/X}\rightarrow \bHom(\Pi(X),\Ens), \ \ \ Y\mapsto (\ox\mapsto \omega_\ox(Y)),
\end{equation}
induit une équivalence entre la catégorie $\Et_{\rf/X}$ et la catégorie des foncteurs $\psi$ de $\Pi(X)$ dans 
$\Ens$ tels que pour tout point géométrique $\ox$ de $X$, $\psi(\ox)$ soit un ensemble fini muni d'une action continue 
de $\pi_1(X,\ox)$. 

Pour tout point géométrique $\ox$ de $X$, le foncteur  fibre 
$\omega_\ox$ est exact à gauche et transforme familles couvrantes en familles surjectives. 
Il se prolonge donc en un foncteur fibre $\phi_\ox\colon X_\fet\rightarrow \Ens$ (\cite{sga4} IV 6.3). 
Celui-ci se déduit du foncteur $\nu_\ox$ défini dans \eqref{higgs2-Kp2c} par oubli de l'action de $\pi_1(X,\ox)$. 
Compte tenu de \eqref{higgs2-Kpp17c}, il correspond donc au point $\rho_X(\ox)$ de $X_\fet$. 
Interprétant $\Pi(X)$ comme la catégorie opposée à la catégorie des pro-objets fondamentaux
normalisés de $\Et_{\rf/X}$ (\cite{sga1} V 5.7), pour tous points géométriques $\ox$ et $\ox'$ de $X$, 
tout morphisme de $\pi_1(X,\ox,\ox')$ induit un morphisme 
$\phi_\ox\rightarrow \phi_{\ox'}$ des foncteurs fibres associés de $X_\fet$. On en déduit un foncteur 
\begin{equation}\label{higgs2-sli10b}
\Pi(X)\rightarrow \Pt(X_\fet), \ \ \ \ox\mapsto \phi_\ox,
\end{equation}
où $\Pt(X_\fet)$ est la catégorie des points de $X_\fet$. Celui-ci est une équivalence de catégories 
en vertu de (\cite{sga4} IV 4.9.4 et 7.2.5).

\begin{lem}\label{higgs2-sli11}
Sous les hypothèses de \eqref{higgs2-sli10}, le foncteur 
\begin{equation}\label{higgs2-sli11a}
\nu\colon X_\fet\rightarrow \bHom(\Pi(X),\Ens), \ \ \ F\mapsto (\ox\mapsto \phi_\ox(F)),
\end{equation}
déduit de \eqref{higgs2-sli10b} induit une équivalence entre la catégorie $X_\fet$ et la catégorie $\Phi$ des foncteurs 
$\varphi$ de $\Pi(X)$ dans $\Ens$ tels que pour tout point géométrique $\ox$ de $X$, 
$\varphi(\ox)$ soit un $\mU$-ensemble discret muni d'une action continue à gauche de $\pi_1(X,\ox)$. 
\end{lem}

En effet, soient $\ox$ un point géométrique de $X$, $\bB_{\pi_1(X,\ox)}$ le topos classifiant du groupe profini $\pi_1(X,\ox)$. 
Le foncteur 
\begin{equation}\label{higgs2-sli11b}
\Phi\rightarrow \bB_{\pi_1(X,\ox)}, \ \ \ \varphi\mapsto \varphi(\ox),
\end{equation}
est clairement une équivalence de catégories. D'autre part, le composé de ce dernier et du foncteur $\nu$ est 
l'équivalence de catégories
\begin{equation}\label{higgs2-sli11c}
\nu_\ox\colon X_\fet\stackrel{\sim}{\rightarrow}\bB_{\pi_1(X,\ox)}
\end{equation}
définie dans \eqref{higgs2-Kp2c}; d'où la proposition.

\begin{prop}\label{higgs2-Kpp1}
Soit $X$ un schéma cohérent (resp. un schéma ayant un nombre fini de composantes con\-nexes). 
Alors tout objet de $\Et_{\rf/X}$ est cohérent dans $X_\fet$; en particulier, 
le topos $X_\fet$ est cohérent.  
\end{prop}

Pour la première proposition, comme les produits fibrés sont représentables dans $\Et_{\rf/X}$, 
il suffit de montrer que tout objet de $\Et_{\rf/X}$ est quasi-compact (\cite{sga4} VI 2.1). 
Si $X$ est cohérent, tout objet de $\Et_{\rf/X}$ est un schéma cohérent, et est donc
quasi-compact en tant qu'objet du site $\Et_{\rf/X}$. 
Supposons ensuite que $X$ ait un nombre fini de composantes connexes $X_1,\dots,X_n$. 
Il suffit de montrer que pour tout objet $Y$ de $\Et_{\rf/X}$ 
et tout $1\leq i\leq n$, $Y\times_XX_i$ est quasi-compact (\cite{sga4} VI 1.3). 
On peut donc supposer que l'image de $Y$ dans $X$ est égale à $X_1$. 
Quitte à remplacer $X$ par $X_1$ \eqref{higgs2-Kpp32}, on peut se borner au cas où $X$ est connexe.   
Il résulte alors de l'équivalence de catégories \eqref{higgs2-Kp2a} que $Y$ est quasi-compact.

Comme les limites projectives finies sont représentables dans $\Et_{\rf/X}$, 
la seconde proposition résulte de la première et de (\cite{sga4} VI 2.4.5).

\begin{cor}\label{higgs2-Kpp102}
Pour tout schéma $X$ dont l'ensemble des composantes connexes est localement fini, 
le topos $X_\fet$ est algébrique. 
\end{cor}

En effet, les composantes connexes $(X_i)_{i\in I}$ de $X$ sont ouvertes et fermées dans $X$ \eqref{higgs2-elc0}. 
Par suite, $(X_i\rightarrow X)_{i\in I}$ est un recouvrement de $\Et_{\rf/X}$. D'autre part, pour tout $i\in I$, 
le topos $(X_\fet)_{/X_i}=(X_i)_{\fet}$ est cohérent en vertu de \ref{higgs2-Kpp1}, et le morphisme
$X_i\rightarrow X$ est clairement quasi-séparé dans $X_\fet$. 
Il résulte alors de la définition (\cite{sga4} VI 2.3) que $X_\fet$ est algébrique. 

\begin{cor}\label{higgs2-Kpp111}
Pour tout schéma cohérent $X$, le morphisme 
$\rho_X\colon X_\et\rightarrow X_\fet$ est cohérent.
\end{cor}
En effet, pour tout objet $Y$ de $\Et_{\rf/X}$, $\rho_X^*(Y)=Y$ est un objet cohérent de $X_\et$. 
Par suite, $\rho_X$ est cohérent en vertu de \ref{higgs2-Kpp1} et (\cite{sga4} VI 3.2).

\begin{lem}\label{higgs2-Kpp33}
Soient $X$ un schéma, $F$ un faisceau abélien localement constant et constructible de $X_\et$. Alors~:
\begin{itemize}
\item[{\rm (i)}] Il existe un revêtement étale $Y\rightarrow X$ et un homomorphisme surjectif 
$u\colon \mZ_{Y_\et}\rightarrow F$, où $\mZ_{Y_\et}$ est le $\mZ$-module libre de $X_\et$ engendré par $Y$
{\rm (\cite{sga4} IV 11.3.3)}. 
\item[{\rm (ii)}] Si $\rho_{X*}(F)=0$, alors $F=0$. 
\end{itemize}
\end{lem}

(i) Par descente, $F$ est représentable par un objet $Y$ de $\Et_{\rf/X}$ (\cite{sga4} IX 2.2).  
L'identité de $F$ définit alors une section $e$ de $F(Y)$ et par suite un homomorphisme 
$u\colon \mZ_{Y_\et}\rightarrow F$ qui est clairement surjectif. 

(ii) Considérons $Y$ et $u$ comme dans (i) et notons 
$\mZ_{Y_\fet}$ le $\mZ$-module libre de $X_\fet$ engendré par $Y$. 
On a $\rho_X^*(\mZ_{Y_\fet})\simeq \mZ_{Y_\et}$ en vertu de (\cite{sga4} IV 13.4(b)). 
On en déduit un homomorphisme surjectif $u\colon \rho_X^*(\mZ_{Y_\fet})\rightarrow F$. 
Soit $v\colon \mZ_{Y_\fet}\rightarrow \rho_{X*}(F)$ le morphisme adjoint de $u$. 
Si $\rho_{X*}(F)=0$, alors $v=0$ et par suite $u=0$, ce qui implique que $F=0$.

\begin{lem}\label{higgs2-Kpp2}
Soient $X$ un schéma dont l'ensemble des composantes connexes est localement fini, 
$Y$ un objet de $\Et_{\rf/X}$, $F$ un faisceau abélien de $X_\fet$, 
$e\in F(Y)$. Notons $\mZ_{Y_\fet}$ le $\mZ$-module libre de $X_\fet$ engendré par $Y$ 
{\rm (\cite{sga4} IV 11.3.3)} et $u\colon \mZ_{Y_\fet}\rightarrow F$ l'homomorphisme associé à $e$. Alors~:
\begin{itemize}
\item[{\rm (i)}] $Y$ est un faisceau localement constant et constructible de $X_\et$. 
\item[{\rm (ii)}] Si $u$ est un épimorphisme, alors $\rho_X^*(F)$ est un faisceau de $\mZ$-modules 
localement constant et constructible de $X_\et$. 
\end{itemize}
\end{lem}

(i) On peut se borner au cas où $X$ est connexe \eqref{higgs2-elc0}. D'après \eqref{higgs2-Kp2a}, il existe un revêtement étale surjectif 
$X'\rightarrow X$ tel que $Y\times_XX'$ soit $X'$-isomorphe à une somme disjointe finie de copies de $X'$, 
d'où l'assertion. 

(ii) Procédant comme dans (i), on peut se réduire au cas où $X$ est connexe et 
$Y$  est $X$-isomorphe à une somme disjointe finie de copies de $X$, en vertu
de \eqref{higgs2-Kpp31a} et (\cite{sga4} IV 13.4(b)). Donc $e$ correspond à des sections 
$e_1,\dots,e_n\in F(X)$ telles que l'homomorphisme induit $\mZ_{X_\fet}^n\rightarrow F$ soit surjectif. 
Considérons un point géométrique $\ox$ de $X$, et reprenons les notations de \ref{higgs2-Kp2}. 
On en déduit, compte tenu de l'équivalence de catégories \eqref{higgs2-Kp2c},  
que $\nu_\ox(F)$ est un groupe abélien de type fini, muni de l'action triviale de $\pi_1(X,\ox)$.  
Par suite, $F$ est un faisceau abélien constant de $X_\fet$ de valeur $\nu_\ox(F)$, 
et il en est alors de même de $\rho_X^*(F)$.

\begin{lem}\label{higgs2-Kpp3}
Soient $X$ un schéma  dont l'ensemble des composantes connexes est localement fini, 
$F$ un faisceau de $X_\fet$. Alors la catégorie $(\Et_{\rf/X})_{/F}$ est filtrante. 
\end{lem}

Il suffit de montrer que les sommes de deux objets et les conoyaux des doubles flèches 
sont représentables dans $(\Et_{\rf/X})_{/F}$ (\cite{sga4} I 2.7.1). 
Si $U$ et $V$ sont deux objets de $(\Et_{\rf/X})_{/F}$, 
la somme disjointe $U\sqcup V$ est un revêtement étale de $X$ qui représente la somme de $U$ et $V$ dans 
$(\Et_{\rf/X})_{/F}$. Soient $f,g\colon U\rightrightarrows V$ une double flèche de $(\Et_{\rf/X})_{/F}$, 
$G$ son conoyau dans $X_\fet$. 
Il suffit de montrer que $G$ est représentable par un objet de $\Et_{\rf/X}$. 
Il existe un recouvrement étale $(X_i\rightarrow X)_{i\in I}$ de $\Et_{\rf/X}$ tel que, pour tout $i\in I$, $X_i$ soit connexe
\eqref{higgs2-elc0}. Par descente, il suffit de montrer que, pour tout $i\in I$, 
$G|X_i$ est représentable par un objet de $\Et_{\rf/X_i}$ (\cite{sga1} VIII 2.1 et 5.7, \cite{giraud2} II 3.4.4). 
On peut donc se borner au cas où $X$ est connexe. 
Considérons un point géométrique $\ox$ de $X$, reprenons les notations de \ref{higgs2-Kp2}
et identifions $\Et_{\rf/X}$ avec la catégorie $\cC(\pi_1(X,\ox))$ au moyen de l'équivalence de catégories \eqref{higgs2-Kp2a}.
Comme dans toute catégorie galoisienne, les limites inductives finies sont représentables (\cite{sga1} V 4.2),  
le conoyau $W$ de la double flèche $(f,g)$ est représentable dans $\cC(\pi_1(X,\ox))$. 
Il est immédiat de voir que $W$ est aussi le conoyau de $(f,g)$ dans le topos $\bB_{\pi_1(X,\ox)}$. 
Identifiant $X_\fet$ avec $\bB_{\pi_1(X,\ox)}$ au moyen du foncteur \eqref{higgs2-Kp2c}, 
on voit alors que $W$ représente $G$.

\begin{prop}\label{higgs2-Kpp4}
Si $X$ est un schéma cohérent, ayant un nombre fini de composantes connexes, 
alors le morphisme d'adjonction $\id\rightarrow \rho_{X*}\rho_X^*$ est un isomorphisme~;
en particulier, le foncteur $\rho^*_X\colon X_\fet\rightarrow X_\et$ est pleinement fidèle.  
\end{prop}

Soit $G$ un faisceau de $X_\fet$.  D'après (\cite{sga4} II 4.1.1),
le morphisme canonique de $X_\fet$ 
\begin{equation}\label{higgs2-Kpp4a}
\underset{\underset{(\Et_{\rf/X})/G}{\longrightarrow}}{\lim}\ U\rightarrow G
\end{equation}
est un isomorphisme. Comme le foncteur $\rho_X^*$ commute aux limites inductives et qu'il prolonge 
le foncteur d'injection canonique $\Et_{\rf/X}\rightarrow \Et_{/X}$, on en déduit un isomorphisme de $X_\et$
\begin{equation}\label{higgs2-Kpp4b}
\underset{\underset{(\Et_{\rf/X})/G}{\longrightarrow}}{\lim}\ U\stackrel{\sim}{\rightarrow} \rho_X^*(G).
\end{equation}
D'autre part, les topos $X_\et$ et $X_\fet$ sont cohérents \eqref{higgs2-Kpp1} et le morphisme $\rho_X$ est cohérent \eqref{higgs2-Kpp111}.
Donc le foncteur $\rho_{X*}$ commute aux limites inductives filtrantes. En effet, 
la preuve de (\cite{sga4} VI 5.1) vaut aussi pour le foncteur image directe des faisceaux d'ensembles. 
Comme la catégorie $(\Et_{\rf/X})_{/G}$ est filtrante d'après \ref{higgs2-Kpp3}, on en déduit un isomorphisme 
\begin{equation}\label{higgs2-Kpp4c}
\underset{\underset{(\Et_{\rf/X})/G}{\longrightarrow}}{\lim}\ U\stackrel{\sim}{\rightarrow} \rho_{X*}(\rho_X^*(G)).
\end{equation}
Par fonctorialité, le morphisme d'adjonction $G\rightarrow \rho_{X*}(\rho_X^*G)$ est la limite inductive des morphismes 
identiques $\id_U$, pour les objets $U$ de $(\Et_{\rf/X})/G$. C'est donc un isomorphisme, d'où la première 
proposition. 
La seconde proposition est équivalente à la première par les propriétés générales des foncteurs adjoints. 

\begin{rema}\label{higgs2-Kpp101}
Soient $X$ un schéma, $G$ un faisceau de $X_\fet$. D'après (\cite{sga4} I 5.1 et III 1.3), 
$\rho_X^*(G)$ est le faisceau associé au préfaisceau $F$ sur $\Et_{/X}$ défini, pour tout $V\in \ob(\Et_{/X})$, par 
\begin{equation}
F(V)=\underset{\underset{(U,u)\in I_V^\circ}{\longrightarrow}}{\lim}\ G(U),
\end{equation}
où $I_V$ est la catégorie des couples $(U,u)$ formés d'un objet $U$ de $\Et_{\rf/X}$ et d'un $X$-morphisme 
$u\colon V\rightarrow U$. De plus, le morphisme d'adjonction $G\rightarrow \rho_{X*}(\rho_X^*G)$ est induit 
par le morphisme de préfaisceaux sur $\Et_{\rf/X}$ défini, pour tout $U\in \ob(\Et_{\rf/X})$, par l'isomorphisme canonique
\begin{equation}
G(U)\stackrel{\sim}{\rightarrow} F(U);
\end{equation}
en effet, $I_U$ admet un objet final, à savoir $(U,\id_U)$. On prendra garde que cela n'implique pas en général 
que le morphisme d'adjonction $\id\rightarrow \rho_{X*}\rho_X^*$ soit un isomorphisme, puisque 
$\rho_X^*(G)$ n'est pas en général égal à $F$.
\end{rema}

\begin{cor}\label{higgs2-Kpp5}
Soient $X$ un schéma cohérent, ayant un nombre fini de composantes connexes, 
$F$ un faisceau (resp. faisceau abélien de torsion) de $X_\et$.
Alors les conditions suivantes sont équivalentes~:
\begin{itemize}
\item[{\rm (i)}] Il existe un faisceau (resp. faisceau abélien de torsion) 
$G$ de $X_\fet$ et un isomorphisme $F\simeq \rho_X^*(G)$.
\item[{\rm (ii)}] $F$ est limite inductive dans $X_\et$ d'un système inductif filtrant de faisceaux 
(resp. faisceaux abéliens de torsion) localement constants et constructibles.
\end{itemize}
\end{cor}
Montrons que (i) implique (ii). Soit $G$ un faisceau de $X_\fet$. 
D'après \eqref{higgs2-Kpp4b}, on a un isomorphisme de $X_\et$
\begin{equation}\label{higgs2-Kpp5a}
\underset{\underset{(\Et_{\rf/X})/G}{\longrightarrow}}{\lim}\ U\stackrel{\sim}{\rightarrow} \rho_X^*(G).
\end{equation}
On en déduit par \ref{higgs2-Kpp2}(i) et \ref{higgs2-Kpp3} que $\rho_X^*(G)$ 
vérifie la condition non respée du (ii) .

Soit $G$ un faisceau abélien de torsion de $X_\fet$. Pour tout objet $U$ de $(\Et_{\rf/X})_{/G}$, 
notons $\mZ_{U_\fet}$ le $\mZ$-module libre de $X_\fet$ engendré par $U$  
et $H_{U}$ l'image du morphisme canonique $\mZ_{U_\fet}\rightarrow G$ (\cite{sga4} IV 11.3.3). 
L'isomorphisme \eqref{higgs2-Kpp4a} induit un isomorphisme de groupes de $X_\fet$
\begin{equation}\label{higgs2-Kpp4e}
\underset{\underset{U\in (\Et_{\rf/X})/G}{\longrightarrow}}{\lim}\ H_U\stackrel{\sim}{\rightarrow} G.
\end{equation}
On en déduit alors un isomorphisme de groupes de $X_\et$
\begin{equation}\label{higgs2-Kpp4f}
\underset{\underset{(\Et_{\rf/X})/G}{\longrightarrow}}{\lim}\ \rho_X^*(H_U)\rightarrow \rho_X^*(G).
\end{equation}
D'après \ref{higgs2-Kpp2}(ii) et (\cite{sga4} IX 1.2), $\rho_X^*(H_U)$ est un faisceau abélien de torsion, 
localement constant et constructible de $X_\et$. 
Comme la catégorie $(\Et_{\rf/X})/G$ est filtrante \eqref{higgs2-Kpp3}, $\rho_X^*(G)$ vérifie la condition respée du (ii).

Montrons que (ii) implique (i).  Considérons en premier lieu le cas où $F$ est un faisceau d'ensembles. 
Par descente, si $F$ est localement constant et constructible, 
il est représentable par un objet $Y$ de $\Et_{\rf/X}$ (\cite{sga4} IX 2.2);
on a donc un isomorphisme $F\stackrel{\sim}{\rightarrow} \rho_X^*(Y)$, ce qui démontre la propriété requise 
dans le cas envisagé. Le cas général s'en déduit puisque $\rho_X^*$ est pleinement fidèle \eqref{higgs2-Kpp4} et 
commute aux limites inductives.  

Considérons ensuite le cas où $F$ est un faisceau abélien de torsion, localement constant et constructible. 
D'après \ref{higgs2-Kpp33}(i), il existe un revêtement étale $Y\rightarrow X$ et un homomorphisme surjectif 
$u\colon \mZ_{Y_\et}\rightarrow F$, où $\mZ_{Y_\et}$ est le $\mZ$-module libre de $X_\et$ engendré par $Y$. 
Notons $\mZ_{Y_\fet}$ le $\mZ$-module libre de $X_\fet$ engendré par $Y$.
On a $\rho_X^*(\mZ_{Y_\fet})\simeq \mZ_{Y_\et}$ en vertu de (\cite{sga4} IV 13.4(b)). 
Soient $v\colon \mZ_{Y_\fet}\rightarrow \rho_{X*}(F)$ le morphisme adjoint de $u$, 
$G$ l'image de $v$, $w\colon \rho_X^*(G)\rightarrow F$ l'adjoint 
du morphisme canonique $G\rightarrow \rho_{X*}(F)$ et $H$ le noyau de $w$. Il est clair que $w$ est surjectif. 
D'autre part, $\rho_{X*}(F)$ est de torsion en vertu de \ref{higgs2-Kpp1}, \ref{higgs2-Kpp111} et (\cite{sga4} IX 1.2(v)); 
il en est alors de même de $G$ et de $\rho_X^*(G)$. 
D'après \ref{higgs2-Kpp2}(ii), $\rho_X^*(G)$ est un faisceau abélien localement constant et constructible de $X_\et$. 
Par suite, $H$ est un faisceau abélien de torsion, localement constant et constructible
en vertu de (\cite{sga4} IX 2.1(ii) et 2.6). D'une part, la suite 
\begin{equation}
0\rightarrow \rho_{X*}(H)\rightarrow \rho_{X*}(\rho_X^*(G))\rightarrow \rho_{X*}(F)
\end{equation}
est exacte. D'autre part, le morphisme d'adjonction $G\rightarrow \rho_{X*}(\rho_X^*(G))$ est un isomorphisme 
en vertu de \ref{higgs2-Kpp4}. On en déduit que $\rho_{X*}(H)=0$ et par suite que $H=0$ en vertu de \ref{higgs2-Kpp33}(ii).
Donc $w$ est un isomorphisme, ce qui démontre la propriété requise dans le cas envisagé.  

Enfin, le cas où $F$ est limite inductive dans $X_\et$ d'un système inductif filtrant de faisceaux abéliens 
de torsion, localement constants et constructibles se déduit du cas précédent puisque $\rho_X^*$ est pleinement 
fidèle \eqref{higgs2-Kpp4} et commute aux limites inductives.  

\begin{defi}\label{higgs2-Kp5}\index{Schema@Schéma!Kpiun@$K(\pi,1)$}
On dit qu'un schéma $X$ est $K(\pi,1)$ si pour tout entier $n$ inversible dans $\co_X$
et tout $(\mZ/n\mZ)$-module $F$ de $X_\fet$, 
l'homomorphisme d'adjonction $F\rightarrow \rR\rho_{X*}(\rho_X^*F)$ est un isomorphisme.
\end{defi}

Cette notion ne semble raisonnable que pour les schémas qui satisfont à la conclusion
de \ref{higgs2-Kpp4}, en particulier pour les schémas cohérents, ayant un nombre fini de composantes connexes.

\section{Site et topos de Faltings}\label{higgs2-tf}

\subsection{}\label{higgs2-tf1}\index{Site fibre de Faltings@Site fibré de Faltings}
\index{Topos!de Faltings@de Faltings}
\index{100000510@$\hE$}\index{100000512@$\tE$}\index{1000001002@$\pi\colon E\rightarrow \Et_{/X}$}\index{1000001003@$(V\rightarrow U)$}
Dans cette section, $f\colon Y\rightarrow X$ désigne un morphisme de schémas et 
\begin{equation}\label{higgs2-tf1b}
\pi\colon E\rightarrow \Et_{/X}
\end{equation}
le $\mU$-site fibré déduit du site fibré des revêtements étales $\cR/\Sch$ \eqref{higgs2-Kp6a} 
par changement de base par le foncteur 
\begin{equation}\label{higgs2-tf1c}
\Et_{/X}\rightarrow \Sch, \ \ \ U\mapsto U\times_XY.
\end{equation} 
On dit que $\pi$ est le {\em site fibré de Faltings} associé à $f$. 
On peut décrire explicitement la catégorie $E$ de la façon suivante. Les objets de $E$ 
sont les morphismes de schémas $V\rightarrow U$ au-dessus de $f\colon Y\rightarrow X$ tels que le morphisme
$U\rightarrow X$ soit étale et que le morphisme $V\rightarrow U_Y=U\times_XY$ soit étale fini. 
Soient $(V'\rightarrow U')$, $(V\rightarrow U)$ deux objets de $E$. Un morphisme 
de $(V'\rightarrow U')$ dans $(V\rightarrow U)$ est la donnée d'un $X$-morphisme $U'\rightarrow U$ et 
d'un $Y$-morphisme $V'\rightarrow V$ tels que le diagramme
\begin{equation}\label{higgs2-tf1a}
\xymatrix{
V'\ar[r]\ar[d]&U'\ar[d]\\
V\ar[r]&U}
\end{equation}
soit commutatif. Le foncteur $\pi$ est alors défini pour tout objet $(V\rightarrow U)$ de $E$, par
\begin{equation}\label{higgs2-tf1d}
\pi(V\rightarrow U)=U.
\end{equation}

On notera que le site fibré $\pi$ vérifie les conditions de \ref{higgs2-tcevg1} ainsi que la condition \ref{higgs2-tcevg18}(i'). 
On munit $E$ de la topologie co-évanescente associée à $\pi$ \eqref{higgs2-tcevg3}, autrement dit,  
la topologie engendrée par les recouvrements $\{(V_i\rightarrow U_i)\rightarrow (V\rightarrow U)\}_{i\in I}$
des deux types suivants~: 
\begin{itemize}
\item[(v)] $U_i=U$ pour tout $i\in I$, et $(V_i\rightarrow V)_{i\in I}$
est un recouvrement étale. 
\item[(c)] $(U_i\rightarrow U)_{i\in I}$ est un recouvrement étale et 
$V_i=U_i\times_UV$ pour tout $i\in I$. 
\end{itemize}
Le site co-évanescent $E$ ainsi défini est encore appelé {\em site de Faltings} associé à $f$;  c'est un $\mU$-site.  
On désigne par $\hE$ (resp. $\tE$) la catégorie des préfaisceaux (resp. le topos des faisceaux) 
de $\mU$-ensembles sur $E$. On appelle encore $\tE$ le {\em topos de Faltings} associé à $f$. 
Si $F$ est un préfaisceau sur $E$, on note $F^a$ le faisceau associé.

\begin{rema}\label{higgs2-tf7}
La catégorie $E$ a été initialement introduite par Faltings, mais 
avec une topologie qui est en général strictement plus fine que la topologie co-évanescente, à savoir, 
la topologie engendrée par les familles de morphismes $\{(V_i\rightarrow U_i)\rightarrow (V\rightarrow U)\}_{i\in I}$
telles que $(V_i\rightarrow V)_{i\in I}$ et $(U_i\rightarrow U)_{i\in I}$ soient des recouvrements étales 
(\cite{faltings2} page 214). 
En effet, si $Y$ est vide, $\tE$ est le topos vide, 
c'est-à-dire qu'il est équivalent à la catégorie ponctuelle (\cite{sga4} IV 2.2 et 4.4).
Cela résulte de \ref{higgs2-tcevg6} puisque pour tout $U\in \ob(\Et_{/X})$, $U_Y$ est vide, et donc $(U_Y)_\fet$ 
est le topos vide. Par contre, si l'on munit $E$ de la topologie considérée par Faltings, on obtient le topos $X_\et$. 
Cet exemple montre aussi que la topologie considérée par Faltings  
ne satisfait pas en général la proposition \ref{higgs2-tcevg5}, qui joue pourtant un rôle essentiel dans son approche, 
d'où la nécessité de la modifier comme nous l'avons fait dans \ref{higgs2-tf1}. Nous donnerons dans (\cite{ag3} 8.18) 
un exemple qui illustre un autre défaut de la topologie considérée par Faltings. 
\end{rema}

\subsection{}\label{higgs2-tf2}
Il résulte de \ref{higgs2-tcevg4} et du fait que $E$ admet un objet final que les limites projectives finies sont représentables 
dans $E$, que le foncteur $\pi$ est exact à gauche et que la famille des recouvrements verticaux (resp. cartésiens) de $E$ 
est stable par changement de base. En fait, la limite projective d'un diagramme  
\begin{equation}\label{higgs2-tf2a}
\xymatrix{
&{(V''\rightarrow U'')}\ar[d]\\
{(V'\rightarrow U')}\ar[r]&{(V\rightarrow U)}}
\end{equation} 
de $E$ est représentable par le morphisme $(V'\times_VV''\rightarrow U'\times_UU'')$.  
En effet, ce morphisme représente clairement la limite projective du diagramme \eqref{higgs2-tf2a}
dans la catégorie des morphismes de schémas au-dessus de $f$. 
Il suffit donc de montrer que c'est un objet de $E$, ou encore que le morphisme 
\begin{equation}\label{higgs2-tf2b}
V'\times_VV''\rightarrow  U'_Y\times_{U_Y}U''_Y
\end{equation}
est un revêtement étale. Celui-ci est clairement un morphisme de $\Et_{/Y}$, et est donc étale.  
D'autre part, le diagramme 
\begin{equation}
\xymatrix{
{V'\times_VV''}\ar[r]\ar[d]&{V'\times_{U_Y}V''}\ar[d]\\
V\ar[r]^-(0.5){\Delta_V}&{V\times_{U_Y}V}}
\end{equation}
où $\Delta_V$ est le plongement diagonal, est cartésien. Comme $\Delta_V$ est une immersion fermée, 
on en déduit aussitôt que le morphisme \eqref{higgs2-tf2b} est fini, d'où notre assertion.

\subsection{} \label{higgs2-tf3}\index{1000001005@$E_\coh$, $E_\scoh$, $\tE_\coh$, $\tE_\scoh$ ($E$ site de Faltings)}
Soit $\star$ l'un des deux symboles ``$\coh$'' pour cohérent, ou ``$\scoh$'' pour séparé et cohérent, introduits dans \ref{higgs2-not2}. 
On désigne par 
\begin{equation}\label{higgs2-tf3a}
\pi_\star\colon E_\star\rightarrow \Et_{\star/X}
\end{equation}
le site fibré déduit de $\pi$ par changement de base par le foncteur d'injection canonique \eqref{higgs2-not2}
\begin{equation}\label{higgs2-tf3b}
\varphi\colon \Et_{\star/X}\rightarrow \Et_{/X}, 
\end{equation} 
et par 
\begin{equation}\label{higgs2-tf3c}
\Phi\colon E_\star\rightarrow E
\end{equation}
la projection canonique. On munit $E_\star$ de la topologie co-évanescente définie par $\pi_\star$
et on note $\tE_\star$ le topos des faisceaux de $\mU$-ensembles sur $E_\star$. 
En vertu de \ref{higgs2-tcevg10}, si $X$ est quasi-séparé, 
le foncteur $\Phi$ induit par restriction une équivalence de catégories 
\begin{equation}\label{higgs2-tf3d}
\Phi_s\colon \tE\stackrel{\sim}{\rightarrow} \tE_\star.
\end{equation}
De plus, sous la même hypothèse, la topologie co-évanescente de $E_\star$ est induite par celle de $E$ au moyen du foncteur $\Phi$,
d'après \ref{higgs2-tcevg11}.

\begin{prop}\label{higgs2-tf4}
Supposons $X$ et $Y$ cohérents. Alors,
\begin{itemize}
\item[{\rm (i)}] Pour tout objet $(V\rightarrow U)$ de $E_{\coh}$, 
$(V\rightarrow U)^a$ est un objet cohérent de $\tE$. 
\item[{\rm (ii)}] Le topos $\tE$ est cohérent~; en particulier, il a suffisamment de points.
\end{itemize}
\end{prop}

(i) En effet, tout objet de $\Et_{\coh/X}$ est quasi-compact. D'autre part, pour tout $W\in \ob(\Et_{\coh/Y})$,
comme $W$ est un schéma cohérent, tout objet de $\Et_{\rf/W}$ est quasi-compact en vertu de \ref{higgs2-Kpp1}. 
La proposition résulte alors de \eqref{higgs2-tf3d} et \ref{higgs2-tcevg15}(iii) (appliqués au site fibré $\pi_\coh$ \eqref{higgs2-tf3a}).

(ii) Cela résulte de \ref{higgs2-tcevg15}(iv) et (\cite{sga4} VI § 9).

\subsection{}\label{higgs2-tf11}\index{100000525@$\sigma$, $\beta$}
On munit $\Et_{/X}$ de l'objet final $X$,  et $E$ de l'objet final $(Y\rightarrow X)$.
Les foncteurs $\alpha_{X!}$ \eqref{higgs2-tcevg1ab} et $\sigma^+$ \eqref{higgs2-tcevg18d} sont alors explicitement définis par 
\begin{eqnarray}
\alpha_{X!}\colon \Et_{\rf/Y}\rightarrow E,&& V\mapsto (V\rightarrow X),\label{higgs2-tf11b}\\
\sigma^+\colon \Et_{/X}\rightarrow E,&& U\mapsto (U_Y\rightarrow U).\label{higgs2-tf11a}
\end{eqnarray}
Ceux-ci sont exacts à gauche et continus \eqref{higgs2-tcevg18}. Ils définissent donc deux morphismes de topos 
\begin{eqnarray}
\beta\colon \tE\rightarrow Y_\fet,\label{higgs2-tf11d}\\
\sigma\colon \tE\rightarrow X_\et.\label{higgs2-tf11c}
\end{eqnarray} 
Pour tout faisceau $F=\{U\mapsto F_U\}$ sur $E$, on a $\beta_*(F)=F_X$.

\begin{lem}\label{higgs2-tf13}\index{1000001010@$\Psi \colon Y_\et\rightarrow \tE$}
{\rm (i)}\ Le foncteur 
\begin{equation}\label{higgs2-tf13a}
\Psi^+\colon E\rightarrow \Et_{/Y},\ \ \ (V\rightarrow U)\mapsto V
\end{equation}
est continu et exact à gauche~; il définit donc un morphisme de topos 
\begin{equation}\label{higgs2-tf13b}
\Psi\colon Y_\et\rightarrow \tE.
\end{equation}

{\rm (ii)}\ Pour tout faisceau $F$ de $Y_\et$, on a un isomorphisme canonique de $\tE$
\begin{equation}\label{higgs2-tf13c}
\Psi_*(F)\stackrel{\sim}{\rightarrow} \{U\mapsto \rho_{U_Y*}(F|U_Y)\},
\end{equation}
où pour tout objet $U$ de $\Et_{/X}$, 
$\rho_{U_Y}\colon (U_Y)_\et\rightarrow (U_Y)_\fet$ est le morphisme canonique \eqref{higgs2-Kp1a},
et pour tout morphisme $g\colon U'\rightarrow U$ de $\Et_{/X}$, le morphisme de transition 
\begin{equation}\label{higgs2-tf13d}
\rho_{U_Y*}(F|U_Y)\rightarrow (g_Y)_{\fet*}(\rho_{U'_Y*}(F|U'_Y))
\end{equation}
est le composé 
\begin{equation}\label{higgs2-tf13e}
\rho_{U_Y*}(F|U_Y)\rightarrow \rho_{U_Y*}((g_Y)_{\et*}(F|U'_Y))\stackrel{\sim}{\rightarrow}
(g_Y)_{\fet*}(\rho_{U'_Y*}(F|U'_Y)),
\end{equation}
dans lequel la première flèche est induite par le morphisme d'adjonction $\id\rightarrow (g_Y)_{\et*}(g_Y)_{\et}^*$
et la seconde flèche par \eqref{higgs2-Kpp31a}.
\end{lem}
En effet, $\Psi^+$ est clairement exact à gauche \eqref{higgs2-tf2}. D'autre part, 
pour tout faisceau $F$ de $Y_\et$, on a un isomorphisme canonique de $\hE$
\begin{equation}\label{higgs2-tf13f}
F\circ \Psi^+\stackrel{\sim}{\rightarrow} \{U\mapsto \rho_{U_Y*}(F|U_Y)\},
\end{equation}
où le membre de droite est le préfaisceau sur $E$ défini par les morphismes de transition \eqref{higgs2-tf13e}. 
Soit $(U_i\rightarrow U)_{i\in I}$ un recouvrement de $\Et_{/X}$. 
Pour tout $(i,j)\in I^2$, posons $V_i=U_i\times_XY$, $U_{ij}=U_i\times_UU_j$ et $V_{ij}=U_{ij}\times_XY$, 
et notons $h_{i}\colon V_{i}\rightarrow U_Y$ et $h_{ij}\colon V_{ij}\rightarrow U_Y$ 
les morphismes structuraux. La suite 
\begin{equation}
0\rightarrow F|U_Y\rightarrow \prod_{i\in I} (h_i)_{\et*}(F|V_i)\rightrightarrows 
\prod_{(i,j)\in I^2} (h_{ij})_{\et*}(F|V_{ij})
\end{equation}
est exacte. Comme $ \rho_{U_Y*}$ commute aux limites projectives, 
on en déduit par \eqref{higgs2-Kpp31a} que la suite 
\begin{equation}
0\rightarrow \rho_{U_V*}(F|U_Y)\rightarrow \prod_{i\in I} (h_i)_{\fet*}(\rho_{V_i*}(F|V_i))\rightrightarrows 
\prod_{(i,j)\in I^2} (h_{ij})_{\fet*}(\rho_{V_{ij}*}(F|V_{ij}))
\end{equation}
est exacte. Par suite, $F\circ \Psi^+$ est un faisceau sur $E$ en vertu de \ref{higgs2-tcevg5}, d'où la proposition. 

\subsection{}\label{higgs2-tf12}
Explicitons les constructions du \eqref{higgs2-fccp1} pour le foncteur $\Psi^+$ défini dans \eqref{higgs2-tf13a}.
Le foncteur composé 
\begin{equation}\label{higgs2-tf12a}
\Psi^+\circ \sigma^+\colon \Et_{/X}\rightarrow \Et_{/Y}
\end{equation}
n'est autre que le foncteur image inverse par $f\colon Y\rightarrow X$; on a donc $f_\et=\sigma\Psi$. 
D'autre part, pour tout objet $U$ de $\Et_{/X}$, le foncteur  \eqref{higgs2-fccp1aa}
\begin{equation}\label{higgs2-tf12b}
\Psi_U^+\colon \Et_{\rf/U_Y}\rightarrow \Et_{/U_Y}
\end{equation}
induit par $\Psi^+$, s'identifie à l'injection canonique. On peut donc identifier le morphisme $\Psi_U$ \eqref{higgs2-fccp1ac}
au morphisme canonique $\rho_{U_Y}\colon (U_Y)_\et\rightarrow (U_Y)_\fet$ \eqref{higgs2-Kp1a}; en particulier, on a
$\beta\Psi=\rho_Y$ \eqref{higgs2-fccp1ad}. 
Pour tout morphisme $g\colon U'\rightarrow U$ de $\Et_{/X}$, 
le diagramme \eqref{higgs2-fccp2a} s'identifie alors au diagramme \eqref{higgs2-Kpp31a}.  
Des isomorphismes $\Psi^*\sigma^* =f_\et^*$ et 
$\Psi^*\beta^*=\rho^*_Y$, on déduit par adjonction des morphismes 
\begin{eqnarray}
\sigma^*&\rightarrow& \Psi_* f_\et^*, \label{higgs2-tf12c}\\
\beta^*&\rightarrow&\Psi_* \rho^*_Y. \label{higgs2-tf12d}
\end{eqnarray}

\begin{prop}\label{higgs2-tf15}
Supposons $X$ quasi-séparé, et $Y$ cohérent et étale-localement connexe \eqref{higgs2-elc}. 
Pour tout $U\in \ob(\Et_{/X})$, notons $h_U\colon U_Y\rightarrow Y$ la projection canonique. 
Alors~: 
\begin{itemize}
\item[{\rm (i)}] Pour tout faisceau $F$ de $Y_\fet$, $\beta^*(F)$ est le faisceau sur $E_{\coh}$
défini par $\{U\mapsto (h_U)^*_\fet F\}$.  
\item[{\rm (ii)}] Le morphisme d'adjonction $\id\rightarrow \beta_*\beta^*$ est un isomorphisme. 
\item[{\rm (iii)}] Le morphisme d'adjonction $\beta^*\rightarrow \Psi_*\rho^*_Y$ \eqref{higgs2-tf12d} est un isomorphisme.
\end{itemize}
\end{prop}

On notera d'abord que les trois foncteurs 
\begin{eqnarray}
\alpha_{X!}^\coh\colon \Et_{\rf/Y}\rightarrow E_{\coh},&& V\mapsto (V\rightarrow X),\\
\sigma^+_\coh\colon \Et_{\coh/X}\rightarrow E_{\coh},&& U\mapsto (U_Y\rightarrow U),\\
\Psi^+_\coh\colon E_\coh\rightarrow \Et_{\coh/Y},&& (V\rightarrow U)\mapsto V,
\end{eqnarray}
sont bien définis, continus et exacts à gauche. Les morphismes de topos qu'ils définissent 
s'identifient à $\beta$, $\sigma$ et $\Psi$, respectivement, compte tenu de \eqref{higgs2-tf3d} et \ref{higgs2-not2}. 
D'autre part, pour tout objet $U$ de $\Et_{\coh/X}$, le schéma $U_Y$ est cohérent et localement connexe. 
Par suite, le morphisme d'adjonction $\id\rightarrow \rho_{U_Y*} \rho_{U_Y}^*$
est un isomorphisme en vertu de \ref{higgs2-Kpp4}. La proposition résulte donc de \ref{higgs2-fccp4}. 

\begin{prop}\label{higgs2-tf20}
Supposons $X$ et $Y$ cohérents. Alors
les morphismes $\beta$, $\sigma$ et $\Psi$ sont cohérents. 
\end{prop}

En effet, tout objet de $\Et_{\coh/X}$ est cohérent dans $X_\et$; tout objet de 
$\Et_{\rf/Y}$ est cohérent dans $Y_\fet$ \eqref{higgs2-Kpp1}; pour tout objet $(V\rightarrow U)$ de $E_{\coh}$,
$(V\rightarrow U)^a$ est un objet cohérent de $\tE$ d'après \ref{higgs2-tf4}(i).
La proposition résulte donc de (\cite{sga4} VI 3.3), compte tenu de la preuve de \ref{higgs2-tf15}.

\begin{rema}\label{higgs2-tf10}
On a un $2$-morphisme 
\begin{equation}\label{higgs2-tf10a}
\tau\colon f_\fet \beta\rightarrow \rho_X \sigma,
\end{equation}
tel que pour tout faisceau $F$ sur $E$ et tout $U\in \ob(\Et_{\rf/X})$, 
\[
f_{\fet*}(\beta_*(F))(U)\rightarrow \rho_{X*}(\sigma_*(F))(U)
\]
soit l'application canonique $F(U_Y\rightarrow X)\rightarrow F(U_Y\rightarrow U)$.
\end{rema}

\subsection{}\label{higgs2-tf19}
Considérons un diagramme commutatif de morphismes de schémas
\begin{equation}\label{higgs2-tf19a}
\xymatrix{
Y'\ar[r]^{f'}\ar[d]_{g'}&X'\ar[d]^g\\
Y\ar[r]^f&X}
\end{equation}
On désigne par 
\begin{equation}\label{higgs2-tf19b}
\pi'\colon E'\rightarrow \Et_{/X'}
\end{equation}
le site fibré de Faltings associé à $f'$ \eqref{higgs2-tf1}. On munit $E'$ de la topologie co-évanescente associée à $\pi'$
et note $\tE'$ le topos des faisceaux de $\mU$-ensembles sur $E'$. 
Pour tout objet $(V\rightarrow U)$ de $E$, le morphisme canonique $V\times_YY'\rightarrow U\times_XX'$
est un objet de $E'$. On définit ainsi un foncteur 
\begin{equation}\label{higgs2-tf19c}
\Phi^+\colon E\rightarrow E',\ \ \ (V\rightarrow U)\mapsto (V\times_YY'\rightarrow U\times_XX'), 
\end{equation}
qui est clairement exact à gauche \eqref{higgs2-tf2}. 
Pour tout faisceau $F=\{U'\mapsto F_{U'}\}$ sur $E'$, $F\circ \Phi^+$ est le préfaisceau sur $E$ défini par 
\begin{equation}\label{higgs2-tf19d}
\{U\mapsto (g'_{U})_{\fet*}(F_{U\times_XX'})\},
\end{equation} 
où pour tout $U\in \ob(\Et_{/X})$, $g'_U\colon U\times_XY'\rightarrow U\times_XY$ est le changement de base de $g'$. 
Soit $(U_i\rightarrow U)_{i\in I}$ un recouvrement de $\Et_{/X}$.  
Comme le foncteur $(g'_{U})_{\fet*}$ commute aux limites projectives, la relation
de recollement de $F$ relativement au recouvrement $(U_i\times_XX'\rightarrow U\times_XX')_{i\in I}$ \eqref{higgs2-tcevg5a} 
implique la relation analogue pour $F\circ \Phi^+$ relativement au recouvrement $(U_i\rightarrow U)_{i\in I}$. 
Par suite, $F\circ \Phi^+$ est un faisceau sur $E$ \eqref{higgs2-tcevg5}, et $\Phi^+$ est continu. 
Celui-ci définit donc un morphisme de topos 
\begin{equation}\label{higgs2-tf19e}
\Phi\colon \tE'\rightarrow \tE.
\end{equation}
Il résulte aussitôt des définitions que les carrés du diagramme
\begin{equation}\label{higgs2-tf19f}
\xymatrix{
{X'_\et}\ar[d]_{g_\et}&{\tE'}\ar[l]_-(0.5){\sigma'}\ar[d]^{\Phi}\ar[r]^-(0.5){\beta'}&
{Y'_\fet}\ar[d]^{g'_{\fet}}\\
{X_\et}&{\tE}\ar[l]_{\sigma}\ar[r]^{\beta}&{Y_\fet}}
\end{equation}
où $\beta'$ et $\sigma'$ sont les morphismes canoniques \eqref{higgs2-tf11d} et \eqref{higgs2-tf11c} relatifs à $f'$, 
sont commutatifs à isomorphismes canoniques près. D'autre part, le diagramme 
\begin{equation}\label{higgs2-tf19g}
\xymatrix{
{Y'_\et}\ar[r]^-(0.4){\Psi'}\ar[d]_{g'_\et}&{\tE'}\ar[d]^{\Phi}\\
{Y_\et}\ar[r]^-(0.4){\Psi}&{\tE}}
\end{equation}
où $\Psi'$ est le morphisme \eqref{higgs2-tf13b} relatif à $f'$, 
est commutatif à isomorphisme canonique près. 

\begin{rema}\label{higgs2-tf195}
Conservons les hypothèses de \ref{higgs2-tf19}, soient de plus $F$ un faisceau abélien de $\tE$, $i$ un entier $\geq 0$. 
Alors le diagramme
\begin{equation}\label{higgs2-tf195a}
\xymatrix{
{\rH^i(Y_\fet,\beta_*F)}\ar[r]\ar[d]_u&{\rH^i(\tE,F)}\ar[dd]^w\\
{\rH^i(Y'_\fet,g'^*_\fet(\beta_*F))}\ar[d]_v&\\
{\rH^i(Y'_\fet,\beta'_*(\Phi^*F))}\ar[r]&{\rH^i(\tE',\Phi^*F)}}
\end{equation}
où les flèches horizontales proviennent des suites spectrales de Cartan-Leray (\cite{sga4} V 5.3), $u$ et $w$ sont 
les morphismes canoniques et $v$ est induit par le morphisme de changement de base relativement 
au carré droit de \eqref{higgs2-tf19f} (\cite{egr1} (1.2.2.2)) est commutatif.
\end{rema}

\subsection{}\label{higgs2-tf17}
Soit $(f'\colon Y'\rightarrow X')$ un objet de $E$. On note
\begin{equation}\label{higgs2-tf17a}
\pi'\colon E'\rightarrow \Et_{/X'}
\end{equation}
le site fibré de Faltings associé à $f'$. Tout objet de $E'$ est naturellement un objet de $E$. 
On définit ainsi un foncteur 
\begin{equation}\label{higgs2-tf17c}
\Phi\colon E'\rightarrow E.
\end{equation}
On vérifie aussitôt que $\Phi$ se factorise à travers une équivalence de catégories 
\begin{equation}
E'\stackrel{\sim}{\rightarrow} E_{/(Y'\rightarrow X')},
\end{equation}
qui est même une équivalence de catégories sur $\Et_{/X'}$ (\cite{sga1} VI 4.3), 
où l'on considère $E_{/(Y'\rightarrow X')}$ comme une ($\Et_{/X'}$)-catégorie par changement 
de base du foncteur $\pi$. 
Il résulte alors de \ref{higgs2-tcevg71} que la topologie co-évanescente de $E'$ 
est induite par celle de $E$ au moyen du foncteur $\Phi$. 
Par suite, $\Phi$ est continu et cocontinu (\cite{sga4} III 5.2). 
Il définit donc une suite de trois foncteurs adjoints~:
\begin{equation}\label{higgs2-tf17e}
\Phi_!\colon \tE'\rightarrow \tE, \ \ \ \Phi^*\colon \tE\rightarrow \tE', \ \ \ \Phi_*\colon \tE'\rightarrow \tE,
\end{equation}
dans le sens que pour deux foncteurs consécutifs de la suite, celui de droite est
adjoint à droite de l'autre. D'après (\cite{sga4} III 5.4), le foncteur $\Phi_!$ se factorise à travers 
une équivalence de catégories 
\begin{equation}\label{higgs2-tf17f}
\tE'\stackrel{\sim}{\rightarrow} \tE_{/(Y'\rightarrow X')^a}.
\end{equation}
Le couple de foncteurs $(\Phi^*,\Phi_*)$ définit le morphisme de localisation de $\tE$ en $(Y'\rightarrow X')^a$
\begin{equation}\label{higgs2-tf17g}
\Phi\colon \tE'\rightarrow \tE.
\end{equation}
Comme $\Phi\colon E'\rightarrow E$ est un adjoint à gauche du foncteur $\Phi^+\colon E\rightarrow E'$ 
défini dans \eqref{higgs2-tf19c}, le morphisme \eqref{higgs2-tf17g}
s'identifie au morphisme défini dans \eqref{higgs2-tf19e}, en vertu de (\cite{sga4} III 2.5).

\subsection{}\label{higgs2-tf21}\index{1000001020@$\rho \colon X_\et\gtimes_{X_\et}Y_\et\rightarrow \tE$}
On désigne par 
\begin{equation}\label{higgs2-tf21a}
\varpi\colon D\rightarrow \Et_{/X}
\end{equation} 
le site fibré associé au foncteur image inverse 
$f^+\colon \Et_{/X}\rightarrow \Et_{/Y}$, défini dans \eqref{higgs2-tcevg40}.  
On munit $D$ de la topologie co-évanescente relative à $\varpi$. 
On obtient ainsi le site co-évanescent associé au foncteur $f^+$ \eqref{higgs2-co-ev1}, dont le topos des 
faisceaux de $\mU$-ensembles est $X_\et\gtimes_{X_\et}Y_\et$ \eqref{higgs2-co-ev101}. 
Tout objet de $E$ est naturellement un objet de $D$. 
On définit ainsi un foncteur pleinement fidèle et exact à gauche
\begin{equation}\label{higgs2-tf21b}
\rho^+\colon E\rightarrow D.
\end{equation} 
Pour tout $U\in \ob(\Et_{/X})$, la restriction de $\rho^+$ aux fibres au-dessus de $U$ n'est autre 
que le foncteur d'injection canonique $\Et_{\rf/U_Y}\rightarrow \Et_{/U_Y}$, autrement dit, le foncteur $\Psi_U^+$ \eqref{higgs2-tf12b}.
Le foncteur $\rho^+$ s'identifie donc au foncteur portant le même nom, associé au foncteur $\Psi^+$ \eqref{higgs2-tf13a}
et défini dans \eqref{higgs2-fccp5b}. Celui-ci est continu et exact à gauche. 
Il définit donc un morphisme de topos \eqref{higgs2-fccp5d}
\begin{equation}\label{higgs2-tf21c}
\rho\colon X_\et\gtimes_{X_\et}Y_\et\rightarrow \tE.
\end{equation}
Il résulte aussitôt des définitions que les carrés du diagramme
\begin{equation}\label{higgs2-tf21d}
\xymatrix{
{X_\et}\ar@{=}[d]&{X_\et\gtimes_{X_\et}Y_\et}\ar[l]_-(0.5){\rp_1}\ar[d]^{\rho}\ar[r]^-(0.5){\rp_2}&
{Y_\et}\ar[d]^{\rho_Y}\\
{X_\et}&{\tE}\ar[l]_{\sigma}\ar[r]^{\beta}&{Y_\fet}}
\end{equation}
où $\rp_1$ et $\rp_2$ sont les projections canoniques \eqref{higgs2-co-ev11}, 
sont commutatifs à isomorphismes canoniques près. De plus, on a un diagramme commutatif  
\begin{equation}\label{higgs2-tf21f}
\xymatrix{
{f_{\fet}\beta\rho}\ar[rr]^{\tau_E*\rho}\ar[d]&&{\rho_X\sigma\rho}\ar[d]\\
{f_{\fet}\rho_Y\rp_2}\ar[r]&{\rho_Xf_\et\rp_2}\ar[r]^{\rho_X*\tau_D}&{\rho_X\rp_1}}
\end{equation}
où $\tau_D$ est le $2$-morphisme \eqref{higgs2-co-ev11c}, $\tau_E$ est le $2$-morphisme \eqref{higgs2-tf10a}, 
les flèches verticales sont les isomorphismes sous-jacents au diagramme \eqref{higgs2-tf21d}, 
et la flèche horizontale non libellée provient de \eqref{higgs2-Kpp31a}.
D'autre part, le diagramme 
\begin{equation}\label{higgs2-tf21e}
\xymatrix{
{Y_\et}\ar[r]^-(0.5){\Psi_D}\ar[dr]_{\Psi}&{X_\et\gtimes_{X_\et}Y_\et}\ar[d]^{\rho}\\
&{\tE}}
\end{equation}
où $\Psi_D$ est le morphisme \eqref{higgs2-co-ev14a}, 
est commutatif à isomorphisme canonique près.

\begin{prop}\label{higgs2-tf212}
Supposons le schéma $X$ quasi-séparé et le schéma $Y$ cohérent. 
Alors, pour qu'une famille $((V_i\rightarrow U_i)\rightarrow (V\rightarrow U))_{i\in I}$ 
de morphismes de $E_\coh$ \eqref{higgs2-tf3} soit couvrante, il faut et il suffit qu'elle le soit dans $D$ \eqref{higgs2-tf21}.
\end{prop}

On notera d'abord que la topologie co-évanescente de $E_\coh$ étant induite par celle de $E$ \eqref{higgs2-tf3}, 
une famille de morphismes de même but de $E_\coh$ est couvrante dans $E_\coh$ si et seulement si 
elle l'est dans $E$ (\cite{sga4} III 3.3). Le morphisme $f$ induit un foncteur \eqref{higgs2-not2}
\begin{equation}
f^+_\coh\colon \Et_{\coh/X}\rightarrow \Et_{\coh/Y}.
\end{equation}
On note $D_\coh$ le site co-évanescent associé à $f^+_\coh$ \eqref{higgs2-co-ev1}. 
Le foncteur d'injection canonique $D_\coh\rightarrow D$ est continu et exact à gauche. Il induit une équivalence entre les 
topos associés, d'après \ref{higgs2-topfl11} et \ref{higgs2-co-ev101}.
Par suite, une famille de morphismes de même but de $D_\coh$ est couvrante dans $D_\coh$ si et seulement si 
elle l'est dans $D$ (\cite{sga4} II 4.4).
Par ailleurs, tout objet de $E_{\coh}$ est naturellement un objet de $D_\coh$. 
On définit ainsi un foncteur pleinement fidèle et exact à gauche
\begin{equation}
\rho^+_\coh\colon E_\coh\rightarrow D_\coh.
\end{equation} 
Celui-ci est un $(\Et_{\coh/X})$-foncteur cartésien dont la restriction 
aux fibres au-dessus de tout $U\in \ob(\Et_{\coh/X})$ n'est autre 
que le foncteur d'injection canonique $\Et_{\rf/U_Y}\rightarrow \Et_{\coh/U_Y}$  (cf. \ref{higgs2-tcevg40}).
Donc $\rho^+_\coh$ est continu en vertu de \ref{higgs2-tcevg85}. Pour montrer la proposition, 
il suffit de montrer que pour qu'une famille $\cF=((V_i\rightarrow U_i)\rightarrow (V\rightarrow U))_{i\in I}$ 
de morphismes de $E_\coh$ soit couvrante, il faut et il suffit qu'elle le soit dans $D_\coh$. 
La condition est nécessaire puisque $\rho^+_\coh$ est continu (\cite{sga4} III 1.6).
Inversement, supposons la famille $\cF$ couvrante dans $D_\coh$ et montrons qu'elle l'est dans $E_\coh$. 
Pour tout $U\in \ob(\Et_{\coh/X})$, le schéma $U_Y$ est cohérent. Par suite, tout objet de $\Et_{\coh/U_Y}$ est quasi-compact. 
D'après \ref{higgs2-tcevg44}, il existe donc un recouvrement étale $(U'_j\rightarrow U)_{j\in J}$ 
et pour tout $j\in J$, un recouvrement étale $(V'_{j,k}\rightarrow U'_j\times_UV)_{k\in K_j}$ 
tels que pour tous $j\in J$ et $k\in K_j$, il existe $i_{j,k}\in I$, un $U$-morphisme $U'_j\rightarrow U_{i_{j,k}}$ et
un $V$-morphisme $V'_{j,k}\rightarrow V_{i_{j,k}}$ tels que le diagramme
\begin{equation}
\xymatrix{
{V'_{j,k}}\ar[r]\ar[d]&{V_{i_{j,k}}}\ar[d]\\
{U'_j}\ar[r]&{U_{i_{j,k}}}}
\end{equation}
soit commutatif. 
Posons $W_{j,k}=U'_j\times_{U_{i_{j,k}}}V_{i_{j,k}}$. Pour tout $j\in J$, $(W_{j,k}\rightarrow U'_j\times_UV)_{k\in K_j}$ 
est un recouvrement de $\Et_{\rf/U'_i\times_XY}$. 
Par suite, $((V_i\rightarrow U_i)\rightarrow (V\rightarrow U))_{i\in I}$ est un recouvrement de $E_\coh$.

\begin{rema}\label{higgs2-tf211}
Sous les hypothèses de \eqref{higgs2-tf19}, le diagramme de morphismes de topos
\begin{equation}\label{higgs2-tf211a}
\xymatrix{
{X'_\et\gtimes_{X'_\et}Y'_\et}\ar[r]^{\Xi}\ar[d]_{\rho'}&{X_\et\gtimes_{X_\et}Y_\et}\ar[d]^\rho\\
{\tE'}\ar[r]^\Phi&{\tE}}
\end{equation}
où $\rho$ et $\rho'$ sont les morphismes canoniques \eqref{higgs2-tf21c}, $\Phi$ est le morphisme \eqref{higgs2-tf19e} 
et $\Xi$ est le morphisme déduit de la fonctorialité des topos co-évanescents \eqref{higgs2-topfl10}, 
est commutatif à isomorphisme canonique près. En effet, d'après \ref{higgs2-co-ev8}(i),
pour tout $(V\rightarrow U)\in \ob(D)$ \eqref{higgs2-tf21}, on a un isomorphisme canonique 
\begin{equation}\label{higgs2-tf211b}
\Xi^*((V\rightarrow U)^a)\stackrel{\sim}{\rightarrow} (V\times_YY'\rightarrow U\times_XX')^a.
\end{equation} 
\end{rema}

\subsection{}\label{higgs2-tf24}\index{1000001030@$\rho(\oy\rightsquigarrow \ox)$ point d'un topos de Faltings}
D'après \ref{higgs2-co-ev17} et (\cite{sga4} VIII 7.9), la donnée d'un point de $X_\et\gtimes_{X_\et}Y_\et$ 
est équivalente à la donnée d'une paire de points géométriques $\ox$ de $X$ et $\oy$ de $Y$
et d'une flèche de spécialisation $u$ de $f(\oy)$ vers $\ox$, c'est-à-dire, d'un $X$-morphisme 
$u\colon \oy\rightarrow X_{(\ox)}$, où $X_{(\ox)}$ désigne le localisé strict de $X$ en $\ox$. 
Un tel point sera noté $(\oy\rightsquigarrow \ox)$ ou encore $(u\colon \oy\rightsquigarrow \ox)$. 
On désigne par $\rho(\oy\rightsquigarrow \ox)$ son image par  $\rho$ \eqref{higgs2-tf21c}, qui est donc un point de $\tE$.  
Pour tous $F\in \ob(X_\et)$ et $G\in \ob(Y_\fet)$, on a des isomorphismes canoniques fonctoriels 
\begin{eqnarray}
(\sigma^*F)_{\rho(\oy\rightsquigarrow \ox)} &\stackrel{\sim}{\rightarrow}& F_{\ox}, \label{higgs2-tf24a}\\
(\beta^*G)_{\rho(\oy\rightsquigarrow \ox)} &\stackrel{\sim}{\rightarrow}& (\rho_Y^*G)_{\oy}. \label{higgs2-tf24b}
\end{eqnarray}
D'après \eqref{higgs2-tf21f}, pour tout $H\in \ob(X_\fet)$, l'application 
\begin{equation}\label{higgs2-tf24c}
(\sigma^*(\rho_X^*H))_{\rho(\oy\rightsquigarrow \ox)} \rightarrow (\beta^*(f^*_\fet H))_{\rho(\oy\rightsquigarrow \ox)} 
\end{equation}
induite par $\tau$ \eqref{higgs2-tf10a}, s'identifie canoniquement au morphisme 
de spécialisation défini par $u$
\[
(\rho_X^*H)_{\ox}\rightarrow (\rho_X^*H)_{f(\oy)}.
\] 
Pour tout objet $(V\rightarrow U)$ de $E$, on a un isomorphisme canonique fonctoriel
\begin{equation}\label{higgs2-tf24d}
(V\rightarrow U)^a_{\rho(\oy\rightsquigarrow \ox)} \stackrel{\sim}{\rightarrow} U_{\ox}\times_{U_{f(\oy)}}V_{\oy},
\end{equation}
l'application $V_{\oy}\rightarrow U_{f(\oy)}$ est induite par $V\rightarrow U$, et 
l'application $U_{\ox}\rightarrow U_{f(\oy)}$ est le morphisme de spécialisation défini par $u$. 
Cela résulte de \eqref{higgs2-co-ev17d} et du fait que $\rho^*$ prolonge $\rho^+$ (\cite{sga4} III 1.4).

\begin{rema}\label{higgs2-tf25}
Soient $U\in \ob(\Et_{/X})$, $\ox,\oz$ deux points géométriques de $X$, 
$u$ une flèche de spécialisation de $\oz$ vers $\ox$, $u^*\colon U_\ox\rightarrow U_\oz$
le morphisme de spécialisation correspondant (\cite{sga4} VIII 7.7). Alors~:
\begin{itemize}
\item[(i)] Pour tout $\ox'\in U_\ox$, $\oz'=u^*(\ox')$ est le point de $U_\oz$ correspondant au morphisme composé 
\begin{equation}
\oz\stackrel{u}{\longrightarrow} X_{(\ox)}\stackrel{i}{\longrightarrow} U,
\end{equation}
où $i$ est le $X$-morphisme défini par $\ox'$ (\cite{sga4} VIII 7.3); plus précisément, $i$ est induit 
par l'inverse de l'isomorphisme canonique $U_{(\ox')}\stackrel{\sim}{\rightarrow}X_{(\ox)}$, 
où $U_{(\ox')}$ désigne le localisé strict de $U$ en $\ox'$. 
Il existe donc une unique spécialisation $u'$ de $\oz'$ vers $\ox'$ 
qui s'insère dans un diagramme commutatif 
\begin{equation}
\xymatrix{
{\oz'}\ar[r]^-(0.5){u'}\ar[d]&{U_{(\ox')}}\ar[d]\\
\oz\ar[r]^-(0.5)u&{X_{(\ox)}}}
\end{equation}
où les flèches verticales sont les isomorphismes canoniques. 

\item[(ii)] Si $U$ est séparé sur $X$, l'application $u^*$ est injective. 
En effet, si $\ox'$ et $\ox''$ sont deux points de $U_\ox$ tels que $u^*(\ox')=u^*(\ox'')=\oz'$,
alors on a deux spécialisations de $\oz'$ vers $\ox'$ et $\ox''$, respectivement. 
En particulier, l'image du morphisme canonique $\oz'\rightarrow U\times_{X}X_{(\ox)}$ est contenue 
dans chacun des sous-schémas ouverts et fermés $U_{(\ox')}$ et $U_{(\ox'')}$  de $U\times_{X}X_{(\ox)}$  
(\cite{ega4} 18.5.11), ce qui n'est pas possible. 

\item[(iii)] Si $U$ est fini sur $X$, l'application $u^*$ est bijective puisque
$U\times_{X}X_{(\ox)}=\prod_{\ox'\in U_\ox}U_{(\ox')}$.
\end{itemize}
\end{rema}

\subsection{}\label{higgs2-tf250}
Soient $\ox$ un point géométrique de $X$, $\oy$ un point géométrique de $Y$, 
$X_{(\ox)}$ le localisé strict de $X$ en $\ox$, $u\colon \oy\rightarrow X_{(\ox)}$ un $X$-morphisme,
de sorte que $(\oy \rightsquigarrow \ox)$ est un point de $X_\et\gtimes_{X_\et}Y_\et$ \eqref{higgs2-tf24}.
On désigne par $\cP_{\rho(\oy \rightsquigarrow \ox)}$ la catégorie des objets $\rho(\oy \rightsquigarrow \ox)$-pointés de $E$,
définie comme suit. Les objets de $\cP_{\rho(\oy \rightsquigarrow \ox)}$ sont les triplets $((V\rightarrow U), \xi,\zeta)$ 
formés d'un objet $(V\rightarrow U)$ de $E$, d'un $X$-morphisme $\xi\colon \ox\rightarrow U$ et 
d'un $Y$-morphisme $\zeta\colon \oy\rightarrow V$ tels que, notant encore 
$\xi\colon X_{(\ox)}\rightarrow U$ le $X$-morphisme induit par $\xi$ (\cite{sga4} VIII 7.3), le diagramme
\begin{equation}\label{higgs2-tf250a}
\xymatrix{
\oy\ar[r]^-(0.5)u\ar[d]_{\zeta}&{X_{(\ox)}}\ar[d]^{\xi}\\
V\ar[r]&U}
\end{equation}
soit commutatif. 
Soient $((V\rightarrow U), \xi,\zeta)$, $((V'\rightarrow U'), \xi',\zeta')$ deux objets de $\cP_{\rho(\oy \rightsquigarrow \ox)}$.
Un morphisme de $((V'\rightarrow U'), \xi',\zeta')$ dans $((V\rightarrow U), \xi,\zeta)$ est la donnée d'un 
morphisme $(g\colon U'\rightarrow U, h\colon V'\rightarrow V)$ de $E$ tel que $g\circ \xi'=\xi$ et $h\circ \zeta'=\zeta$. 
Il résulte de \eqref{higgs2-tf24d} et \ref{higgs2-tf25}(i) que $\cP_{\rho(\oy \rightsquigarrow \ox)}$ est canoniquement équivalente à 
la catégorie des voisinages de 
$\rho(\oy \rightsquigarrow \ox)$ dans $E$ (\cite{sga4} IV 6.8.2). Elle est donc cofiltrante et pour tout préfaisceau 
$F=\{U\mapsto F_U\}$ sur $E$, on a un isomorphisme canonique fonctoriel (\cite{sga4} IV (6.8.4))
\begin{equation}\label{higgs2-tf250b}
(F^a)_{\rho(\oy \rightsquigarrow \ox)}\stackrel{\sim}{\rightarrow} \underset{\underset{((V\rightarrow U), \xi,\zeta)\in 
\cP^\circ_{\rho(\oy \rightsquigarrow \ox)}}{\longrightarrow}}\lim\ F_U(V).
\end{equation}
Si $X$ est quasi-séparé, on peut remplacer dans la limite ci-dessus $\cP_{\rho(\oy \rightsquigarrow \ox)}$ 
par la sous-catégorie pleine $\cP^\coh_{\rho(\oy \rightsquigarrow \ox)}$ 
formée des objets $((V\rightarrow U), \xi,\zeta)$ tels que $U$ soit de présentation 
finie sur $X$ ({\em i.e.}, soit un objet de $\Et_{\coh/X}$), qui est aussi cofiltrante (cf. \ref{higgs2-tf3}).

\begin{prop}\label{higgs2-tf251}
Supposons les schémas $X$ et $Y$ cohérents. 
Alors, lorsque $(\oy \rightsquigarrow \ox)$ décrit la famille des points de $X_\et\gtimes_{X_\et}Y_\et$, 
la famille des foncteurs fibres de $\tE$ associés aux points $\rho(\oy \rightsquigarrow \ox)$ est conservative
{\rm (\cite{sga4} IV 6.4.0)}. 
\end{prop}

Soit $u\colon F\rightarrow G$ un morphisme de $\tE$ tel que pour tout point $(\oy \rightsquigarrow \ox)$
de $X_\et\gtimes_{X_\et}Y_\et$, le morphisme correspondant 
$u_{\rho(\oy \rightsquigarrow \ox)}\colon F_{\rho(\oy \rightsquigarrow \ox)}\rightarrow 
G_{\rho(\oy \rightsquigarrow \ox)}$ soit un monomorphisme. Montrons que $u$ est un monomorphisme.
Il s'agit de montrer que pour tous $(V\rightarrow U)\in \ob(E_\coh)$ \eqref{higgs2-tf3}
et $a,b \in F_U(V)$ tels que $u(a)=u(b)$, on a $a=b$. 
Compte tenu de \ref{higgs2-tf17} et \ref{higgs2-tf211}, on peut supposer $(V\rightarrow U)=(Y\rightarrow X)$. 
Pour tout point $(\oy \rightsquigarrow \ox)$ de $X_\et\gtimes_{X_\et}Y_\et$, on a 
$a_{\rho(\oy \rightsquigarrow \ox)}=b_{\rho(\oy \rightsquigarrow \ox)}$ puisque 
$u_{\rho(\oy \rightsquigarrow \ox)}$ est un monomorphisme. D'après \eqref{higgs2-tf250b}, il existe un objet 
$((U_{(\oy \rightsquigarrow \ox)}\rightarrow  V_{(\oy \rightsquigarrow \ox)}),\xi_{(\oy \rightsquigarrow \ox)},
\zeta_{(\oy \rightsquigarrow \ox)})$ de $\cP^\coh_{\rho(\oy \rightsquigarrow \ox)}$ tel que $a$ et $b$ aient mêmes images dans
$F_{U_{(\oy \rightsquigarrow \ox)}}(V_{(\oy \rightsquigarrow \ox)})$. 
Par ailleurs, les topos $X_\et$ et $Y_\et$ étant cohérents et 
le morphisme $f\colon Y_\et\rightarrow X_\et$ étant cohérent (\cite{sga4} VI 3.3), 
la famille des foncteurs fibres de $X_\et\gtimes_{X_\et}Y_\et$ associés aux points 
$(\oy \rightsquigarrow \ox)$ est conservative, d'après \ref{higgs2-tcevg17} et \ref{higgs2-tf24}.
On en déduit que la famille des morphismes 
$((U_{(\oy \rightsquigarrow \ox)}\rightarrow  V_{(\oy \rightsquigarrow \ox)})\rightarrow 
(Y\rightarrow X))_{(\oy \rightsquigarrow \ox)}$
est couvrante dans $D$ \eqref{higgs2-tf21} (\cite{sga4} IV 6.5). Elle est donc couvrante dans $E_\coh$ en vertu de \ref{higgs2-tf212}.
Par suite, $a=b$ et $u$ est un monomorphisme. 

Supposons, de plus, que pour tout point $(\oy \rightsquigarrow \ox)$
de $X_\et\gtimes_{X_\et}Y_\et$, le morphisme 
$u_{\rho(\oy \rightsquigarrow \ox)}$ soit un épimorphisme et montrons qu'il en est de même de $u$. 
Il suffit de montrer que pour tous $(V\rightarrow U)\in \ob(E_\coh)$ et $b \in G_U(V)$, 
il existe $a\in F_U(V)$ tel que $b=u(a)$. On peut encore supposer $(V\rightarrow U)=(Y\rightarrow X)$.
D'après \eqref{higgs2-tf250b}, pour tout point $(\oy \rightsquigarrow \ox)$ de $X_\et\gtimes_{X_\et}Y_\et$, 
il existe un objet $((U_{(\oy \rightsquigarrow \ox)}\rightarrow  V_{(\oy \rightsquigarrow \ox)}),\xi_{(\oy \rightsquigarrow \ox)},
\zeta_{(\oy \rightsquigarrow \ox)})$ de $\cP^\coh_{\rho(\oy \rightsquigarrow \ox)}$ et une section 
$a_{(\oy \rightsquigarrow \ox)}\in F_{U_{(\oy \rightsquigarrow \ox)}}(V_{(\oy \rightsquigarrow \ox)})$
dont l'image par $u$ dans $G_{U_{(\oy \rightsquigarrow \ox)}}(V_{(\oy \rightsquigarrow \ox)})$ est la restriction de $b$. 
Comme $u$ est un monomorphisme, les sections $a_{(\oy \rightsquigarrow \ox)}$
coïncident sur $(V_{(\oy' \rightsquigarrow \ox')}\times_YV_{(\oy \rightsquigarrow \ox)}\rightarrow 
U_{(\oy' \rightsquigarrow \ox')}\times_XU_{(\oy \rightsquigarrow \ox)})$, pour tous les points $(\oy \rightsquigarrow \ox)$
et $(\oy' \rightsquigarrow \ox')$ de $X_\et\gtimes_{X_\et}Y_\et$. Ils proviennent donc d'une section $a\in F_X(Y)$,
et on a $u(a)=b$ puisque les restrictions 
aux $(U_{(\oy \rightsquigarrow \ox)}\rightarrow  V_{(\oy \rightsquigarrow \ox)})$ coïncident.

\subsection{}\label{higgs2-tf26}
Supposons $X$ strictement local, de point fermé $x$. Pour tout $U\in \ob(\Et_{\scoh/X})$ \eqref{higgs2-not2}, 
on désigne par $U^\rf$ la somme disjointe des localisés stricts de $U$ en les points de $U_x$;
c'est un sous-schéma ouvert et fermé de $U$, qui est fini sur $X$ (\cite{ega4} 18.5.11). 
La correspondance $U\mapsto U^\rf$ définit un foncteur 
\begin{equation}\label{higgs2-tf26a}
\iota_x^+\colon \Et_{\scoh/X}\rightarrow \Et_{\rf/X},
\end{equation}
qui est clairement exact à gauche et continu. Le morphisme de topos associé 
\begin{equation}\label{higgs2-tf26b}
\iota_x\colon X_\fet\rightarrow X_\et
\end{equation} 
s'identifie au morphisme $\Ens\rightarrow X_\et$ défini par $x$. 
En effet, le foncteur fibre en $x$ induit une équivalence 
entre les topos $X_\fet$ et $\Ens$ car le groupe $\pi_1(X,x)$ est trivial \eqref{higgs2-Kp2a}, 
et pour tout $U\in \ob(\Et_{\scoh/X})$, l'application canonique $U^\rf_x\rightarrow U_x$ est bijective.

Considérons le site fibré 
\begin{equation}\label{higgs2-tf26c}
\pi_\scoh\colon E_\scoh\rightarrow \Et_{\scoh/X}
\end{equation}
défini dans \eqref{higgs2-tf3}, et munissons $E_\scoh$ de la topologie co-évanescente associée à $\pi_\scoh$.  
On rappelle que le topos des faisceaux de $\mU$-ensembles sur $E_\scoh$ s'identifie canoniquement à $\tE$ \eqref{higgs2-tf3d}.
Pour tout objet $(V\rightarrow U)$ de $E_{\scoh}$, $V\times_UU^\rf=V\times_{U_Y}U^\rf_Y$ 
est un revêtement étale de $Y$. On obtient ainsi un foncteur 
\begin{equation}\label{higgs2-tf26d}
\theta^+\colon E_{\scoh}\rightarrow \Et_{\rf/Y}, \ \ \ (V\rightarrow U)\mapsto V\times_UU^\rf. 
\end{equation}

\begin{lem}\label{higgs2-tf261}\index{1000001040@$\theta\colon Y_\fet\rightarrow \tE$}
Conservons les hypothèses de \eqref{higgs2-tf26}. 
\begin{itemize}
\item[{\rm (i)}] Le foncteur $\theta^+$ est continu et exact à gauche~; il définit donc un morphisme de topos 
\begin{equation}\label{higgs2-tf261a}
\theta\colon Y_\fet\rightarrow \tE.
\end{equation}
\item[{\rm (ii)}] Pour tout faisceau $F$ de $Y_\fet$, on a un isomorphisme canonique de faisceaux sur $E_\scoh$
\begin{equation}\label{higgs2-tf261b}
\theta_*(F)\stackrel{\sim}{\rightarrow} \{U\mapsto (j_{U})_{\fet*}(F|U^\rf_Y)\},
\end{equation}
où pour tout $U\in \ob(\Et_{\scoh/X})$, $j_{U}\colon U^\rf_Y\rightarrow U_Y$ désigne l'injection canonique,
et pour tout morphisme $g\colon U'\rightarrow U$ de $\Et_{\scoh/X}$, le morphisme de transition 
\begin{equation}\label{higgs2-tf261c}
(j_{U})_{\fet*}(F|U^\rf_Y)\rightarrow (g_Y)_{\fet*}((j_{U'})_{\fet*}(F|U'^\rf_Y))
\end{equation}
est le composé 
\begin{equation}\label{higgs2-tf261d}
(j_{U})_{\fet*}(F|U^\rf_Y)\rightarrow (j_{U})_{\fet*}((g^\rf_Y)_{\fet*}(F|U'^\rf_Y))\stackrel{\sim}{\rightarrow}
(g_Y)_{\fet*}((j_{U'})_{\fet*}(F|U'^\rf_Y)),
\end{equation}
où $g^\rf\colon U'^\rf\rightarrow U^\rf$ est l'image de $g$ par le foncteur \eqref{higgs2-tf26a}, la première flèche 
est induite par le morphisme d'adjonction $\id\rightarrow (g^\rf_Y)_{\fet*} (g^\rf_Y)_{\fet}^*$ et la seconde flèche 
est l'isomorphisme canonique.
\end{itemize}
\end{lem}

On vérifie aisément que $\theta^+$ commute aux produits fibrés et transforme objet final en objet final. 
Il est donc exact à gauche. D'autre part, pour tout faisceau $F$ de $Y_\fet$, 
on a un isomorphisme canonique de $\hE_\scoh$
\begin{equation}\label{higgs2-tf261e}
F\circ \theta^+\stackrel{\sim}{\rightarrow}\{U\mapsto (j_{U})_{\fet*}(F|U^\rf_Y)\}, 
\end{equation}
où le membre de droite est le préfaisceau sur $E_\scoh$ défini par les morphismes de transition \eqref{higgs2-tf261d}. 
Soit $(U_i\rightarrow U)_{i\in I}$ un recouvrement de $\Et_{\scoh/X}$. 
Pour tout $(i,j)\in I^2$, posons $V_i=U_i\times_XY$, $W_i=U^\rf_i\times_XY$,
$U_{ij}=U_i\times_UU_j$, $V_{ij}=U_{ij}\times_XY$ et $W_{ij}=U^\rf_{ij}\times_XY$, 
et notons $h_{i}\colon V_{i}\rightarrow U_Y$, $g_i \colon W_{i}\rightarrow U^\rf_Y$, 
$h_{ij}\colon V_{ij}\rightarrow U_Y$ et $g_{ij}\colon W_{ij}\rightarrow U^\rf_Y$ 
les morphismes structuraux. Le foncteur \eqref{higgs2-tf26a} étant exact à gauche et continu, la suite 
\begin{equation}
0\rightarrow F|U^\rf_Y\rightarrow \prod_{i\in I} (g_i)_{\fet*}(F|W_i)\rightrightarrows 
\prod_{(i,j)\in I^2} (g_{ij})_{\fet*}(F|W_{ij})
\end{equation}
est exacte. Comme $(j_{U})_{\fet*}$ commute aux limites projectives, on en déduit que la suite 
\begin{equation}
0\rightarrow (j_{U})_{\fet*}(F|U^\rf_Y)\rightarrow \prod_{i\in I} (h_i)_{\fet*}((j_{U_i})_{\fet*}(F|W_i))\rightrightarrows 
\prod_{(i,j)\in I^2} (h_{ij})_{\fet*}((j_{U_{ij}})_{\fet*}(F|W_{ij}))
\end{equation}
est exacte. Par suite, $F\circ \theta^+$ est un faisceau en vertu de \ref{higgs2-tcevg5}. 
Donc $\theta^+$ est continu, d'où la proposition. 

\subsection{}\label{higgs2-tf262}
Conservons les hypothèses de \eqref{higgs2-tf26} et 
explicitons les constructions de \eqref{higgs2-fccp1} pour le foncteur $\theta^+$ défini dans \eqref{higgs2-tf26d}.
Le foncteur composé 
\begin{equation}\label{higgs2-tf262a}
\theta^+\circ \sigma^+\colon \Et_{\scoh/X}\rightarrow \Et_{\rf/Y}
\end{equation}
n'est autre que le foncteur $U\mapsto U^\rf_Y$; on a donc $\sigma\theta=\iota_x f_\fet$, où $\iota_x$ est 
le morphisme \eqref{higgs2-tf26b}. 
D'autre part, pour tout objet $U$ de $\Et_{\scoh/X}$, le foncteur \eqref{higgs2-fccp1aa}
\begin{equation}\label{higgs2-tf262b}
\theta^+_U\colon \Et_{\rf/U_Y}\rightarrow \Et_{\rf/U_Y^\rf}
\end{equation}
induit par $\theta^+$ n'est autre que l'image inverse par le morphisme canonique $j_U\colon U_Y^\rf\rightarrow U_Y$;
en particulier, $\theta$ est une section de $\beta$, {\em i.e.}, 
on a un isomorphisme canonique \eqref{higgs2-fccp1ad}
\begin{equation}\label{higgs2-tf262c}
\beta\theta \stackrel{\sim}{\rightarrow} \id_{Y_\fet}.
\end{equation} 
On obtient un morphisme de changement de base 
\begin{equation}\label{higgs2-tf262d}
\beta_*\rightarrow \theta^*,
\end{equation}
composé de $\beta_*\rightarrow \beta_* \theta_*\theta^*\stackrel{\sim}{\rightarrow}\theta^*$,
où la première flèche est induite par le morphisme d'adjonction $\id\rightarrow \theta_*\theta^*$,
et la seconde flèche par \eqref{higgs2-tf262c}.

\subsection{}\label{higgs2-tf264}
Conservons les hypothèses de \eqref{higgs2-tf26}. Le morphisme canonique $\rho_X\colon X_\et\rightarrow X_\fet$ \eqref{higgs2-Kp1a}
s'identifie à l'unique morphisme de topos $X_\et\rightarrow \Ens$ (\cite{sga4} IV 4.3).
Par ailleurs, le morphisme composé $\iota_x \rho_X\colon X_\et\rightarrow X_\et$ est défini par le morphisme de sites 
\begin{equation}\label{higgs2-tf264c}
\Et_{\scoh/X}\rightarrow \Et_{\scoh/X}, \ \ \ U\mapsto U^\rf.
\end{equation}
L'injection canonique $U^\rf\rightarrow U$, pour $U\in  \ob(\Et_{\scoh/X})$, définit alors un $2$-morphisme
\begin{equation}\label{higgs2-tf264a}
\id_{X_\et} \rightarrow \iota_x \rho_X. 
\end{equation}
D'après \ref{higgs2-topfl6}, les morphismes de topos 
$\iota_x \rho_X f_\et \colon Y_\et\rightarrow X_\et$ et $\id_{Y_\et}$, et le $2$-morphisme 
$f_\et \rightarrow \iota_x \rho_X f_\et $ induit par \eqref{higgs2-tf264a}, définissent un morphisme de topos 
\begin{equation}\label{higgs2-tf264b}
\gamma\colon Y_\et \rightarrow X_\et\gtimes_{X_\et}Y_\et
\end{equation}
tel que $\rp_1\gamma= \iota_x \rho_X f_\et$, $\rp_2 \gamma=\id_{Y_\et}$ et $\tau*\gamma$ est induit  
par \eqref{higgs2-tf264a}~:
\begin{equation}\label{higgs2-tf264d}
\xymatrix{
{X_\et}\ar[d]_{\iota_x \rho_X}&{Y_\et}\ar[d]_{\gamma}\ar[l]_{f_\et}\ar[rd]^\id&\\
{X_\et}\ar[rd]&{X_\et\gtimes_{X_\et}Y_\et}\ar[l]_-(0.4){\rp_1}\ar[r]^-(0.4){\rp_2}&{Y_\et}\ar[ld]^{f_\et}\\
&{X_\et}&}
\end{equation}

Reprenons les notations de \eqref{higgs2-tf21} et notons $D_\scoh$ le site fibré sur $\Et_{\scoh/X}$ déduit 
de $D$ \eqref{higgs2-tf21a} par changement de base par le foncteur d'injection canonique $\Et_{\scoh/X}\rightarrow \Et_{/X}$ 
\eqref{higgs2-tf3}. Il résulte de \ref{higgs2-co-ev8}(i) que pour tout $(V\rightarrow U)\in \ob(D_\scoh)$, 
on a un isomorphisme canonique 
\begin{equation}\label{higgs2-tf264e}
\gamma^*((V\rightarrow U)^a)\stackrel{\sim}{\rightarrow} V\times_U U^\rf.
\end{equation}
On en déduit que le diagramme 
\begin{equation}\label{higgs2-tf264f}
\xymatrix{
{Y_\et}\ar[r]^-(0.5){\gamma}\ar[d]_{\rho_Y}&{X_\et\gtimes_{X_\et}Y_\et}\ar[d]^\rho\\
{Y_\fet}\ar[r]^{\theta}&{\tE}}
\end{equation}
où $\rho$ est le morphisme canonique \eqref{higgs2-tf21c}, est commutatif à isomorphisme canonique près. 

\begin{prop}\label{higgs2-tf263}
Sous les hypothèses de \eqref{higgs2-tf26}, pour tout faisceau $F$ de $\tE$ et tout point géométrique $\oy$ de $Y$, l'application 
\begin{equation}\label{higgs2-tf263a}
(\beta_*F)_{\rho_Y(\oy)}\rightarrow (\theta^*F)_{\rho_Y(\oy)}
\end{equation}
induite par le morphisme de changement de base \eqref{higgs2-tf262d} est bijective.
\end{prop}

On observera d'abord qu'il existe une unique spécialisation de $f(\oy)$ dans $x$; 
on peut donc considérer le point $(\oy \rightsquigarrow x)$ de $X_\et\gtimes_{X_\et}Y_\et$ \eqref{higgs2-tf24}.
Il est clair que les points $\gamma(\oy)$ et $(\oy \rightsquigarrow x)$ sont canoniquement isomorphes. 
Il résulte alors de \eqref{higgs2-tf264f} que les points $\theta(\rho_Y(\oy))$ et $\rho(\oy \rightsquigarrow x)$ de $\tE$
sont canoniquement isomorphes. On désigne par $\cC_\oy$ la catégorie des $Y$-schémas étales et finis, $\oy$-pointés, 
que l'on identifie à la catégorie des voisinages de $\rho_Y(\oy)$ dans le site $\Et_{\rf/Y}$ (\cite{sga4} IV 6.8.2), 
et par $\cP^\scoh_{\rho(\oy\rightsquigarrow x)}$ la sous-catégorie pleine de la catégorie 
$\cP_{\rho(\oy\rightsquigarrow x)}$ \eqref{higgs2-tf250}
formée des objets $\rho(\oy \rightsquigarrow x)$-pointés $((V\rightarrow U),\xi,\zeta)$ de $E$ 
tels que $U$ soit séparé et de présentation finie sur $X$
({\em i.e.}, soit un objet de $\Et_{\scoh/X}$). Il résulte de \eqref{higgs2-tf24d} et \ref{higgs2-tf25}(i) 
que $\cP^\scoh_{\rho(\oy \rightsquigarrow \ox)}$ est canoniquement équivalente 
à la catégorie des voisinages du point $\rho(y \rightsquigarrow x)$ dans le site $E_{\scoh}$. Elle est donc cofiltrante. 
Notons $\xi_0 \colon x\rightarrow X$ l'injection canonique. On a un foncteur pleinement fidèle
\begin{equation}
\jmath_\oy\colon \cC_\oy\rightarrow \cP^\scoh_{\rho(\oy\rightsquigarrow x)}, \ \ \ (V,\zeta\colon \oy\rightarrow V)\mapsto 
((V\rightarrow X),\xi_0,\zeta),
\end{equation}
compatible avec  le foncteur canonique \eqref{higgs2-tf11b}
\begin{equation}
\alpha_{X!}^{\scoh}\colon \Et_{\rf/Y}\rightarrow E_{\scoh}, \ \ \ V\mapsto (V\rightarrow X).
\end{equation} 
Le morphisme d'adjonction $F\rightarrow \theta_*(\theta^*F)$ est 
défini, pour tout $(V\rightarrow U)\in \ob(E_{\scoh})$, par l'application canonique 
\[
F(V\rightarrow U)\rightarrow (\theta^*F)(V\times_UU^\rf).
\]
Par suite, \eqref{higgs2-tf263a} s'identifie à l'application
\begin{equation}\label{higgs2-tf263b}
\underset{\underset{(W,\zeta)\in \cC_\oy^\circ}{\longrightarrow}}{\lim}\ F(W\rightarrow X)
\rightarrow F_{\rho(\oy\rightsquigarrow x)}
\end{equation}
induite par le foncteur $\jmath^\circ_\oy$. Il suffit donc de montrer que  $\jmath_\oy^\circ$ est cofinal.
Soit $((V\rightarrow U),\xi,\zeta)$ un objet de $\cP^\scoh_{\rho(\oy\rightsquigarrow x)}$.
Notons encore $\xi\colon X\rightarrow U$ le $X$-morphisme induit par $\xi$, de sorte que le diagramme
\begin{equation}
\xymatrix{
\oy\ar[r]\ar[d]_{\zeta}&X\ar[d]^{\xi}\\
V\ar[r]&U}
\end{equation}
est commutatif. On en déduit un $Y$-morphisme $\zeta'\colon \oy\rightarrow V\times_{U,\xi}X$ qui s'insère dans 
un diagramme commutatif 
\begin{equation}\label{higgs2-tf263c}
\xymatrix{
\oy\ar[r]^-(0.5){\zeta'}\ar[rd]_-(0.5){\zeta}&{V\times_{U,\xi}X}\ar[r]\ar[d]&X\ar[d]^\xi\\
&V\ar[r]&U}
\end{equation}
Par suite, $(V\times_{U,\xi}X,\zeta')$ est un objet de $\cC_\oy$, et le diagramme \eqref{higgs2-tf263c} induit  un morphisme 
\[
\jmath_\oy(V\times_{U,\xi}X,\zeta')=((V\times_{U,\xi}X\rightarrow X),\xi_0,\zeta')\rightarrow ((V\rightarrow U),\xi,\zeta)
\] 
de $\cP^\scoh_{\rho(\oy\rightsquigarrow x)}$. On en déduit que  $\jmath_\oy^\circ$ est cofinal d'après (\cite{sga4} I 8.1.3(c));
d'où la proposition.

\begin{cor}\label{higgs2-tf265}
Supposons $X$ strictement local, et notons $\theta\colon Y_\fet\rightarrow \tE$ le morphisme 
de topos défini dans \eqref{higgs2-tf261a}. Alors, le morphisme de changement de base
$\beta_*\rightarrow \theta^*$ \eqref{higgs2-tf262d} est un isomorphisme~; en particulier, 
le foncteur $\beta_*$ est exact.   
\end{cor}

Cela résulte de \ref{higgs2-tf263} et \ref{higgs2-Kpp20}. 

\begin{cor}\label{higgs2-tf266}
Supposons $X$ strictement local, et notons $\theta\colon Y_\fet\rightarrow \tE$ le morphisme 
de topos défini dans \eqref{higgs2-tf261a}. Alors~:
\begin{itemize}
\item[{\rm (i)}] Pour tout faisceau $F$ de $\tE$, l'application canonique 
\begin{equation}\label{higgs2-tf266a}
\Gamma(\tE,F)\rightarrow \Gamma(Y_\fet, \theta^*F)
\end{equation}
est bijective. 
\item[{\rm (ii)}] Pour tout faisceau abélien $F$ de $\tE$, l'application canonique 
\begin{equation}\label{higgs2-tf266b}
\rH^i(\tE,F)\rightarrow \rH^i(Y_\fet, \theta^*F)
\end{equation}
est bijective pour tout $i\geq 0$.
\end{itemize}
\end{cor}

(i) En effet, le diagramme 
\begin{equation}
\xymatrix{
{\Gamma(\tE,F)}\ar[d]_v&{\Gamma(Y_\fet,\beta_*F)}\ar[d]^u\ar[l]^-(0.5)\sim_-(0.5)w\\
{\Gamma(\tE,\theta_*(\theta^*F))}\ar[rd]_{v'}&{\Gamma(Y_\fet,\beta_*(\theta_*(\theta^*F)))}\ar[d]^{u'}
\ar[l]_-(0.4){w'}^-(0.4)\sim\\
&{\Gamma(Y_\fet,\theta^*F)}}
\end{equation}
où $w$, $w'$ et $v'$ sont les bijections canoniques, $u$ et $v$ sont 
induits par le morphisme d'adjonction $\id\rightarrow \theta_*\theta^*$ et
$u'$ est induit par l'isomorphisme \eqref{higgs2-tf262c}, est  commutatif. Comme $u'\circ u$ est bijectif en vertu de \ref{higgs2-tf265}, 
il en est de même de $v'\circ v$, d'où la proposition.  

(ii) En effet, le diagramme 
\begin{equation}
\xymatrix{
{\rH^i(\tE,F)}\ar[d]_v&{\rH^i(Y_\fet,\beta_*F)}\ar[d]^u\ar[l]_-(0.5)w\\
{\rH^i(\tE,\theta_*(\theta^*F))}\ar[rd]_{v'}&{\rH^i(Y_\fet,\beta_*(\theta_*(\theta^*F)))}\ar[d]^{u'}\ar[l]_-(0.4){w'}\\
&{\rH^i(Y_\fet,\theta^*F)}}
\end{equation}
où $w$, $w'$ et $v'$ sont induits par la suite spectrale de Cartan-Leray (\cite{sga4} V 5.3), $u$ et $v$ sont 
induits par le morphisme d'adjonction $\id\rightarrow \theta_*\theta^*$ et
$u'$ est induit par l'isomorphisme \eqref{higgs2-tf262c}, est  commutatif. D'autre part, $u'\circ u$ est bijectif,
et le foncteur $\beta_*$ est exact en vertu de \ref{higgs2-tf265}; donc $w$ est bijectif, d'où la proposition. 

\subsection{}\label{higgs2-tf267}\index{1000001040@$\varphi_\ox\colon \tE\rightarrow \uY_\fet$}
Soient $\ox$ un point géométrique de $X$, $\uX$ le localisé strict de $X$ en $\ox$, 
$\uY=Y\times_X\uX$, $\uf\colon \uY\rightarrow \uX$ la projection canonique. 
On désigne par $\tuE$ le topos de Faltings associé à $\uf$ \eqref{higgs2-tf1} et par
\begin{equation}\label{higgs2-tf267a}
\theta\colon \uY_\fet\rightarrow \tuE
\end{equation}
le morphisme défini dans \eqref{higgs2-tf261a}. Le morphisme canonique $\uX\rightarrow X$ induit par fonctorialité un morphisme
\eqref{higgs2-tf19e}
\begin{equation}\label{higgs2-tf267c}
\Phi\colon \tuE\rightarrow \tE.
\end{equation}
On note 
\begin{equation}\label{higgs2-tf267d}
\varphi_\ox\colon \tE\rightarrow \uY_\fet
\end{equation}
le foncteur composé $\theta^*\circ \Phi^*$.

\begin{prop}\label{higgs2-tf268}
Conservons les hypothèses de \eqref{higgs2-tf267}, supposons de plus que $f$ soit cohérent. Alors~:
\begin{itemize}
\item[{\rm (i)}] Pour tout faisceau $F$ de $\tE$, on a un isomorphisme canonique fonctoriel 
\begin{equation}\label{higgs2-tf268a}
\sigma_*(F)_\ox\stackrel{\sim}{\rightarrow}\Gamma(\uY_\fet,\varphi_\ox(F)).
\end{equation}
\item[{\rm (ii)}] Pour tout faisceau abélien $F$ de $\tE$ et tout entier $i\geq 0$, on a un isomorphisme canonique fonctoriel 
\begin{equation}\label{higgs2-tf268b}
\rR^i\sigma_*(F)_\ox\stackrel{\sim}{\rightarrow}\rH^i(\uY_\fet,\varphi_\ox(F)).
\end{equation}
\item[{\rm (iii)}] Pour toute suite exacte de faisceaux abéliens 
$0\rightarrow F'\rightarrow F\rightarrow F''\rightarrow 0$ de $\tE$ et tout entier $i\geq 0$, le diagramme 
\begin{equation}\label{higgs2-tf268c}
\xymatrix{
{\rR^i\sigma_*(F'')_\ox} \ar[r]\ar[d]&{\rR^{i+1}\sigma_*(F')_\ox}\ar[d]\\
{\rH^i(\uY_\fet,\varphi_\ox(F''))}\ar[r]&{\rH^{i+1}(\uY_\fet,\varphi_\ox(F'))}}
\end{equation}
où les flèches verticales sont les isomorphismes canoniques \eqref{higgs2-tf268b} et les flèches horizontales  
sont les bords des suites exactes longues de cohomologie, est commutatif. 
\end{itemize}
\end{prop}

Cette proposition sera démontrée dans \ref{higgs2-lptf6}. 

\begin{rema}\label{higgs2-tf278}
Conservons les hypothèses et notations de \eqref{higgs2-tf267}, soient de plus $\oy$ un point géométrique de $Y$, 
$u\colon \oy\rightarrow \uX$ un $X$-morphisme, de sorte que   
$(\oy \rightsquigarrow \ox)$ est un point de $X_\et\gtimes_{X_\et}Y_\et$ \eqref{higgs2-tf24}.
Notons $\tx$ le point fermé de $\uX$, $v\colon \oy\rightarrow \uY$ le morphisme induit par $u$,
$\ty$ le point géométrique de $\uY$ défini par $v$, $\rho_\uY\colon \uY_\et\rightarrow \uY_\fet$ 
le morphisme canonique \eqref{higgs2-Kp1a} et $\psi_\ty\colon \uY_\fet\rightarrow \Ens$ 
le foncteur fibre associé au point $\rho_{\uY}(\ty)$ de $\uY_\fet$. Le foncteur composé 
\begin{equation}\label{higgs2-tf278a}
\psi_\ty\circ \varphi_\ox \colon \tE\rightarrow \Ens
\end{equation}
est alors canoniquement isomorphe au foncteur fibre associé au point $\rho(\oy \rightsquigarrow \ox)$ de $\tE$ \eqref{higgs2-tf21c}. 
En effet, les carrés du diagramme de morphismes de topos 
\begin{equation}\label{higgs2-tf278b}
\xymatrix{
{\uY_\et}\ar[r]^-(0.5){\gamma}\ar[d]_{\rho_{\uY}}&{\uX_\et\gtimes_{\uX_\et}\uY_\et}
\ar[r]^-(0.5){\Xi}\ar[d]_{\urho}&{X_\et\gtimes_{X_\et}Y_\et}\ar[d]^{\rho}\\
{\uY_\fet}\ar[r]^-(0.4){\theta}&{\tuE}\ar[r]^{\Phi}&{\tE}}
\end{equation}
où $\rho$ et $\urho$ sont les morphismes canoniques \eqref{higgs2-tf21c}, $\gamma$ est le morphisme \eqref{higgs2-tf264b}
et $\Xi$ est le morphisme déduit de la fonctorialité des topos co-évanescents \eqref{higgs2-topfl10},
sont commutatifs à isomorphismes canoniques près~:
le carré de gauche correspond au diagramme \eqref{higgs2-tf264f} et le carré de droite
correspond au diagramme \eqref{higgs2-tf211a}. 
Par ailleurs, il existe une unique flèche de spécialisation de $\uf(\ty)$ vers $\tx$; on peut donc considérer le 
point $(\ty\rightsquigarrow \tx)$ de $\uX_\et\gtimes_{\uX_\et}\uY_\et$. Il est clair que 
$\gamma(\ty)$ est canoniquement isomorphe à $(\ty\rightsquigarrow \tx)$ 
et que $\Xi(\ty\rightsquigarrow \tx)=(\oy \rightsquigarrow \ox)$. 
Par suite, $\rho(\oy \rightsquigarrow \ox)$ est canoniquement isomorphe à $\Phi(\theta(\rho_{\uY}(\ty)))$~; d'où l'assertion. 
\end{rema}

\begin{prop}\label{higgs2-tf29}
Supposons les schémas $X$ et $Y$ cohérents. Pour tout point géométrique $\ox$ de $X$, 
soient $X_{(\ox)}$ le localisé strict de $X$ en $\ox$, $Y_{(\ox)}=Y\times_XX_{(\ox)}$,
\begin{equation}\label{higgs2-tf29a}
\varphi_\ox\colon \tE\rightarrow (Y_{(\ox)})_\fet
\end{equation}
le foncteur défini dans \eqref{higgs2-tf267d}. Alors, la famille des foncteurs $(\varphi_\ox)$,
lorsque $\ox$ décrit l'ensemble des points géométriques de $X$, est conservative. 
\end{prop}

En effet, la famille des foncteurs fibres de $\tE$ associés 
aux points de la forme $\rho(\oy \rightsquigarrow \ox)$, lorsque $(\oy \rightsquigarrow \ox)$
décrit les points de $X_\et\gtimes_{X_\et}Y_\et$, est conservative en vertu de \ref{higgs2-tf251}. 
La proposition s'ensuit compte tenu de \ref{higgs2-tf278}.

\begin{cor}\label{higgs2-tf30}
Sous les hypothèses de \eqref{higgs2-tf29}, pour qu'un morphisme $u$ de $\tE$ soit un monomorphisme (resp. épimorphisme), 
il faut et il suffit que pour tout point géométrique $\ox$ de $X$, il en soit ainsi pour $\varphi_\ox(u)$ \eqref{higgs2-tf29a}.
\end{cor}

Cela résulte de \ref{higgs2-tf29} et (\cite{sga4} I 6.2(ii)).

\subsection{}\label{higgs2-tf275}
Reprenons les hypothèses et notations de \eqref{higgs2-tf267}; de plus, on note $\cC_\ox$ la catégorie des $X$-schémas étales 
$\ox$-pointés (\cite{sga4} VIII 3.9), que l'on identifie à la catégorie des voisinages de $\ox$ dans le site 
$\Et_{/X}$ (\cite{sga4} IV 6.8.2). C'est une catégorie cofiltrante. 
Pour tout objet $(U,\xi\colon \ox\rightarrow U)$ de $\cC_\ox$, 
on désigne encore par $\xi\colon \uX\rightarrow U$ le morphisme déduit de
$\xi$ (\cite{sga4} VIII 7.3) et par
\begin{equation}\label{higgs2-tf275b}
\xi_Y\colon \uY\rightarrow U_Y
\end{equation}
son changement de base par $f$. Le topos $\tE_{/(U_Y\rightarrow U)^a}$ est canoniquement équivalent 
au topos de Faltings associé au morphisme $f_U\colon U_Y\rightarrow U$ d'après \eqref{higgs2-tf17f}. Notons 
\begin{equation}\label{higgs2-tf275d}
\jmath_{U}\colon \tE_{/(U_Y\rightarrow U)^a}\rightarrow \tE
\end{equation}
le morphisme de localisation de $\tE$ en $(U_Y\rightarrow U)^a$. 
Le morphisme $\xi\colon \uX\rightarrow U$ induit alors par fonctorialité un morphisme
de topos \eqref{higgs2-tf19e}
\begin{equation}\label{higgs2-tf275f}
\Phi_{\xi}\colon \tuE\rightarrow \tE_{/(U_Y\rightarrow U)^a}.
\end{equation}
D'après \eqref{higgs2-tf19f}, le diagramme 
\begin{equation}\label{higgs2-tf275aa}
\xymatrix{
{\tuE}\ar[r]^-(0.5){\Phi_\xi}\ar[d]_{\ubeta}&{\tE_{/(U_Y\rightarrow U)^a}}\ar[d]^{\beta_U}\\
{\uY_\fet}\ar[r]^-(0.5){(\xi_Y)_\fet}&{(U_Y)_\fet}}
\end{equation} 
où $\ubeta$ et $\beta_U$ sont les morphismes canoniques \eqref{higgs2-tf11d},
est commutatif à isomorphisme canonique près. On en déduit un morphisme de changement de base 
\begin{equation}\label{higgs2-tf275ab}
(\xi_Y)^*_\fet\beta_{U*}\rightarrow \ubeta_*\Phi_\xi^*.
\end{equation}
D'après (\cite{egr1} 1.2.4(i)), le composé 
\begin{equation}\label{higgs2-tf275bc}
(\xi_Y)^*_\fet\beta_{U*}\rightarrow \ubeta_*\Phi_\xi^* \rightarrow \theta^*\Phi_\xi^*,
\end{equation}
où la seconde flèche est induite par \eqref{higgs2-tf262d}, est le morphisme de changement de base déduit de 
l'isomorphisme canonique \eqref{higgs2-tf262c}
\begin{equation}\label{higgs2-tf275bd}
\beta_U\circ \Phi_\xi\circ \theta\stackrel{\sim}{\rightarrow}(\xi_Y)_\fet.
\end{equation}

Soient $F=\{U\mapsto F_U\}$ un objet de $\hE$, 
$F^a=\{U\mapsto G_U\}$ le faisceau de $\tE$ associé à $F$, et pour tout $U\in \ob(\Et_{/X})$, 
$F_U^a$ le faisceau de $(U_Y)_\fet$ associé à $F_U$. D'après \ref{higgs2-tcevg8}, 
$\{U\mapsto F^a_U\}$ est un préfaisceau sur $E$ et on a un morphisme canonique 
$\{U\mapsto F_U\}\rightarrow \{U\mapsto F_U^a\}$ de $\hE$, induisant un isomorphisme entre les faisceaux associés. 
On associe au préfaisceau $\{U\mapsto F_U^a\}$ le foncteur
\begin{equation}\label{higgs2-tf275c}
\cC_\ox^\circ \rightarrow \uY_\fet, \ \ \ (U,\xi)\mapsto (\xi_Y)^*_\fet(F^a_U),
\end{equation} 
qui à tout morphisme $t\colon (U',\xi')\rightarrow (U,\xi)$ de $\cC_\ox$ fait correspondre le morphisme composé
\[
(\xi_Y)^*_\fet(F^a_U)\stackrel{\sim}{\rightarrow} (\xi'_Y)^*_\fet((t_Y)_\fet^*F^a_U)\rightarrow (\xi'_Y)^*_\fet(F^a_{U'}),
\]
où la première flèche est induite par la relation $\xi=t\circ \xi'$ et la seconde flèche provient du
morphisme de transition $F^a_U\rightarrow (t_Y)_{\fet*}(F^a_{U'})$ de $\{U\mapsto F_U^a\}$. 
De même, on associe au faisceau $F^a=\{U\mapsto G_U\}$ le foncteur
\begin{equation}\label{higgs2-tf275cc}
\cC_\ox^\circ \rightarrow \uY_\fet, \ \ \ (U,\xi)\mapsto (\xi_Y)^*_\fet(G_U).
\end{equation} 
Le morphisme canonique $\{U\mapsto F_U^a\}\rightarrow \{U\mapsto G_U\}$ induit alors un morphisme de foncteurs 
de $\cC_\ox^\circ$ dans $\uY_\fet$~:
\begin{equation}\label{higgs2-tf275cd}
(\xi_Y)^*_\fet(F^a_U)\rightarrow (\xi_Y)^*_\fet(G_U), \ \ \ (U,\xi)\in \ob(\cC_\ox).
\end{equation}

D'après (\cite{sga4} III 5.3), on a un isomorphisme canonique 
\begin{equation}\label{higgs2-tf275db}
\beta_{U*}(\jmath_U^*(F^a))\stackrel{\sim}{\rightarrow} G_U. 
\end{equation}
Les morphismes \eqref{higgs2-tf275ab} et \eqref{higgs2-tf275bc} induisent donc deux morphismes fonctoriels 
\begin{eqnarray}
(\xi_Y)^*_\fet(G_U)&\rightarrow& \ubeta_*(\Phi^*F^a),\label{higgs2-tf275ac}\\
(\xi_Y)^*_\fet(G_U)&\rightarrow& \varphi_\ox(F^a).\label{higgs2-tf275ca}
\end{eqnarray}
D'après (\cite{egr1} 1.2.4(i)), ce sont des morphismes de systèmes projectifs sur la catégorie $\cC_\ox^\circ$ \eqref{higgs2-tf275cc}. 
Compte tenu de \eqref{higgs2-tf275cd}, on en déduit deux morphismes fonctoriels en $F$
\begin{eqnarray}
\underset{\underset{(U,\xi)\in \cC_\ox^\circ}{\longrightarrow}}{\lim}\ (\xi_Y)_\fet^*(F^a_U) &\rightarrow&
\ubeta_*(\Phi^*F^a),\label{higgs2-tf275ad}\\
\underset{\underset{(U,\xi)\in \cC_\ox^\circ}{\longrightarrow}}{\lim}\ (\xi_Y)_\fet^*(F^a_U) &\rightarrow&
\varphi_\ox(F^a).\label{higgs2-tf275da}
\end{eqnarray}

\begin{remas}\label{higgs2-tf274}
Conservons les hypothèses de \eqref{higgs2-tf275}. 
\begin{itemize} 
\item[(i)] D'après (\cite{sga4} XVII 2.1.3), le morphisme \eqref{higgs2-tf275bc} est égal au composé 
\begin{equation}\label{higgs2-tf274a}
(\xi_Y)^*_\fet\beta_{U*}\stackrel{\sim}{\rightarrow}
\theta^*\Phi_\xi^*\beta_U^*\beta_{U*}\rightarrow \theta^*\Phi_\xi^*,
\end{equation}
où la première flèche est induite par \eqref{higgs2-tf275bd} et la seconde flèche est induite 
par le morphisme d'adjonction $\beta^*_U\beta_{U*}\rightarrow\id$.
\item[(ii)] Choisissons un objet affine $(X_0,\xi_0)$ de $\cC_\ox$ et notons $I$ la catégorie
des $X_0$-schémas étales $\xi_0$-pointés qui sont affines au-dessus de $X_0$. 
Le foncteur canonique $I\rightarrow \cC_\ox$ est alors cofinal (\cite{sga4} VIII 4.5). 
On peut donc remplacer dans les limites inductives dans \eqref{higgs2-tf275ad} et \eqref{higgs2-tf275da} la catégorie $\cC_\ox$ par $I$. 
\end{itemize}
\end{remas}

\begin{prop}\label{higgs2-tf279}
Les hypothèses étant celles de \eqref{higgs2-tf275}, soient de plus $\oy$ un point géométrique de $Y$, 
$u\colon \oy\rightarrow \uX$ un $X$-morphisme, de sorte que   
$(\oy \rightsquigarrow \ox)$ est un point de $X_\et\gtimes_{X_\et}Y_\et$ \eqref{higgs2-tf24}. 
On note $v\colon \oy\rightarrow \uY$ le $Y$-morphisme induit par $u$, $\ty$ le point géométrique de $\uY$
défini par $v$, $\rho_\uY\colon \uY_\et\rightarrow \uY_\fet$ le morphisme canonique \eqref{higgs2-Kp1a} et 
$\psi_\ty\colon \uY_\fet\rightarrow \Ens$ le foncteur fibre associé au point $\rho_\uY(\ty)$. 
Soient $F=\{U\mapsto F_U\}$ un objet de $\hE$, 
$F^a$ le faisceau de $\tE$ associé à $F$, et pour tout $U\in \ob(\Et_{/X})$, 
$F_U^a$ le faisceau de $(U_Y)_\fet$ associé à $F_U$. Alors, on a un isomorphisme 
canonique et fonctoriel 
\begin{equation}\label{higgs2-tf279a}
(F^a)_{\rho(\oy\rightsquigarrow \ox)} \stackrel{\sim}{\rightarrow} 
\underset{\underset{(U,\xi)\in \cC_\ox^\circ}{\longrightarrow}}{\lim}\ 
\psi_\ty((\xi_Y)_\fet^*(F^a_U)),
\end{equation}
dont l'inverse s'identifie à l'image du morphisme canonique \eqref{higgs2-tf275da} par le foncteur $\psi_\ty$. 
\end{prop}

On notera d'abord que $\{U\mapsto F_U^a\}$ est naturellement un objet de $\hE$ et que le morphisme 
canonique $\{U\mapsto F_U\}\rightarrow \{U\mapsto F_U^a\}$ induit un isomorphisme entre les faisceaux
associés, d'après \ref{higgs2-tcevg8}. Soit $\cP_{\rho(\oy \rightsquigarrow \ox)}$ la catégorie des objets 
$\rho(\oy \rightsquigarrow \ox)$-pointés de $E$ \eqref{higgs2-tf250}. On a un foncteur 
\begin{equation}
\phi\colon \cP_{\rho(\oy \rightsquigarrow \ox)}\rightarrow \cC_\ox, \ \ \ ((V\rightarrow U),\xi,\zeta)\mapsto (U,\xi).
\end{equation}
Pour tout $(U,\xi)\in \ob(\cC_\ox)$, la fibre de $\phi$ au-dessus de $(U,\xi)$ est canoniquement équivalente 
à la catégorie $\cD^\ty_{(U,\xi)}$ des $U_Y$-schémas finis étales, $\xi_Y(\ty)$-pointés \eqref{higgs2-tf275b}. 
L'isomorphisme \eqref{higgs2-tf250b} induit donc un isomorphisme canonique et fonctoriel 
\begin{equation}\label{higgs2-tf279b}
(F^a)_{\rho(\oy\rightsquigarrow \ox)} \stackrel{\sim}{\rightarrow} 
\underset{\underset{(U,\xi)\in \cC_\ox^\circ}{\longrightarrow}}{\lim}\ 
\underset{\underset{(V,\zeta)\in (\cD^\ty_{(U,\xi)})^\circ}{\longrightarrow}}{\lim}\ F_U(V).
\end{equation}
Compte tenu de \eqref{higgs2-Kpp31a} et (\cite{sga4} IV (6.8.4)), pour tout $(U,\xi)\in \ob(\cC_\ox)$,  
on a un isomorphisme canonique et fonctoriel
\begin{equation}\label{higgs2-tf279c}
\psi_\ty((\xi_Y)_\fet^*(F^a_U))\stackrel{\sim}{\rightarrow} 
\underset{\underset{(V,\zeta)\in (\cD^\ty_{(U,\xi)})^\circ}{\longrightarrow}}{\lim}\ F_U(V);
\end{equation}
d'où l'isomorphisme \eqref{higgs2-tf279a}.
Par ailleurs, $\psi_\ty\circ \varphi_\ox$ est le foncteur fibre de $\tE$ associé au point 
$\rho(\oy \rightsquigarrow \ox)$, en vertu de \ref{higgs2-tf278}. Pour établir la seconde assertion, il suffit donc de montrer 
que pour tout objet $((V\rightarrow U),\xi,\zeta)$ de $\cP_{\rho(\oy\rightsquigarrow \ox)}$, l'application canonique \eqref{higgs2-tf279b}
\begin{equation}
F_U(V)\rightarrow (F^a)_{\rho(\oy\rightsquigarrow \ox)}
\end{equation}
est le composé
\begin{equation}
F_U(V)\rightarrow \psi_\ty((\xi_Y)^*_\fet(F^a_U))\rightarrow \psi_\ty(\varphi_\ox(F^a)),
\end{equation}
où la première flèche est définie par l'objet $(V,\zeta)$ de $\cD^\ty_{(U,\xi)}$ et la seconde flèche est l'image par $\psi_\ty$
du composé de \eqref{higgs2-tf275cd} et \eqref{higgs2-tf275ca}. 
On peut se borner au cas où $F$ est un faisceau de $\tE$, de sorte que $F_U^a=F_U$. 
Par localisation \eqref{higgs2-tf17}, on peut supposer $U=X$. Notons $\pr_Y\colon \uY\rightarrow Y$
la projection canonique. D'après \ref{higgs2-tf274}(i), 
l'image par $\psi_\ty$ du morphisme $(\pr_Y)^*_\fet(F_X)\rightarrow \varphi_\ox(F)$ \eqref{higgs2-tf275ca} 
coïncide avec l'application 
\begin{equation}
(F_X)_{\rho_{Y}(\oy)}\rightarrow F_{\rho(\oy \rightsquigarrow \ox)}
\end{equation}
induite par le morphisme d'adjonction $\beta^*(\beta_*(F))\rightarrow F$ et l'isomorphisme \eqref{higgs2-tf24b}, d'où l'assertion recherchée.

\begin{cor}\label{higgs2-tf290}
Conservons les hypothèses de \eqref{higgs2-tf275}, soient, de plus, $F=\{U\mapsto F_U\}$ un objet de $\hE$, 
$F^a$ le faisceau de $\tE$ associé à $F$, et pour tout $U\in \ob(\Et_{/X})$, 
$F_U^a$ le faisceau de $(U_Y)_\fet$ associé à $F_U$. Alors, le morphisme canonique \eqref{higgs2-tf275da}
\begin{equation}\label{higgs2-tf290a}
\underset{\underset{(U,\xi)\in \cC_\ox^\circ}{\longrightarrow}}{\lim}\ (\xi_Y)_\fet^*(F^a_U) 
\rightarrow \varphi_\ox(F^a)
\end{equation}
est un isomorphisme.
\end{cor}

Cela résulte de \ref{higgs2-tf279} et \ref{higgs2-Kpp20}.

\begin{prop}\label{higgs2-tf276}
Sous les hypothèses de \eqref{higgs2-tf275}, pour tout faisceau $F=\{U\mapsto F_U\}$ de $\tE$, le morphisme canonique \eqref{higgs2-tf275ad}
\begin{equation}\label{higgs2-tf276a}
\underset{\underset{(U,\xi)\in \cC_\ox^\circ}{\longrightarrow}}{\lim}\ (\xi_Y)_\fet^*(F_U) \rightarrow \ubeta_*(\Phi^*F)
\end{equation}
est un isomorphisme. 
\end{prop}

Cela résulte  de \ref{higgs2-tf265} et \ref{higgs2-tf290}.

\subsection{}\label{higgs2-tf27}
Soit $F=\{U\mapsto F_U\}$ un préfaisceau en groupes abéliens sur $E$ tel que pour tout $U\in \ob(\Et_{/X})$, 
$F_U$ soit un faisceau (en groupes abéliens) de $(U_Y)_\fet$ (par exemple, $F$ est un faisceau abélien de $\tE$), 
et soit $i$ un entier $\geq 0$. On désigne par $\cH^i(F)$
le faisceau de $X_\et$ associé au préfaisceau sur $\Et_{/X}$ défini pour tout $U\in \ob(\Et_{/X})$ par  le groupe 
\begin{equation}\label{higgs2-tf27a}
\rH^i((U_Y)_\fet,F_U), 
\end{equation}
et pour tout morphisme $g\colon U'\rightarrow U$ de $\Et_{/X}$, par l'application composée  
\begin{equation}\label{higgs2-tf27b}
\rH^i((U_Y)_\fet,F_U)\rightarrow \rH^i((U_Y)_\fet,(g_{Y})_{\fet*}(F_{U'}))\rightarrow \rH^i((U'_Y)_\fet,F_{U'}),
\end{equation}
où la première flèche est induite par le morphisme de transition de $F$
et la seconde flèche est induite par la suite spectrale de Cartan-Leray (\cite{sga4} V 5.3). 

On note $F^a=\{U\mapsto G_U\}$ le faisceau (en groupes abéliens) de $\tE$ associé à $F$ (\cite{sga4} III 6.4). 
Pour tout $U\in \ob(\Et_{/X})$, le topos $\tE_{/(U_Y\rightarrow U)^a}$ est canoniquement équivalent au 
topos de Faltings associé au morphisme $U_Y\rightarrow U$ d'après \eqref{higgs2-tf17f}. 
On a donc un morphisme canonique de topos \eqref{higgs2-tf11d}
\begin{equation}\label{higgs2-tf27c}
\beta_{U}\colon \tE_{/(U_Y\rightarrow U)^a}\rightarrow (U_Y)_\fet.
\end{equation}
Par définition du foncteur de restriction (\cite{sga4} III 5.3),  
on a un isomorphisme canonique 
\begin{equation}\label{higgs2-tf27d}
\beta_{U*}(F^a|(U_Y\rightarrow U)^a)\stackrel{\sim}{\rightarrow} G_U. 
\end{equation}
Par suite, la suite spectrale de Cartan-Leray induit une application fonctorielle en $F$
\begin{equation}\label{higgs2-tf27e}
\rH^i((U_Y)_\fet,G_U)\rightarrow \rH^i((U_Y\rightarrow U)^a,F^a).
\end{equation}
Composant avec l'application $\rH^i((U_Y)_\fet,F_U)\rightarrow \rH^i((U_Y)_\fet,G_U)$ induite par le morphisme 
canonique $F\rightarrow F^a$, on obtient une application fonctorielle en $F$
\begin{equation}\label{higgs2-tf27ef}
\rH^i((U_Y)_\fet,F_U)\rightarrow \rH^i((U_Y\rightarrow U)^a,F^a).
\end{equation}

Soit $g\colon U'\rightarrow U$ un morphisme de $\Et_{/X}$. Le diagramme de morphismes de topos
\begin{equation}\label{higgs2-tf27g}
\xymatrix{
{\tE_{/(U'_Y\rightarrow U')^a}}\ar[r]^{\beta_{U'}}\ar[d]_{j}&{(U'_Y)_\fet}\ar[d]^{(g_Y)_\fet}\\
{\tE_{/(U_Y\rightarrow U)^a}}\ar[r]^{\beta_U}&{(U_Y)_\fet}}
\end{equation}
où $j$ est le morphisme de localisation de $\tE_{/(U_Y\rightarrow U)^a}$ en $(U'_Y\rightarrow U')^a$, 
est commutatif, à isomorphisme canonique près. Cela se vérifie aussitôt sur le diagramme de morphismes de sites correspondant. Par ailleurs, le diagramme
\begin{equation}\label{higgs2-tf27h}
\xymatrix{
{\beta_{U*}(F^a|(U_Y\rightarrow U)^a)}\ar[r]^-(0.5)u\ar@{=}[d]&{\beta_{U*}(j_*(F^a|(U'_Y\rightarrow U')^a))}
\ar[r]_-(0.5)\sim^-(0.5)v&{(g_Y)_{\fet*}(\beta_{U'*}(F^a|(U'_Y\rightarrow U')^a))}\ar@{=}[d]\\
{G_U}\ar[rr]^w&&{(g_Y)_{\fet*}(G_{U'})}}
\end{equation}
où $u$ est induit par le morphisme d'adjonction $\id\rightarrow j_*j^*$, $v$ est induit 
par le diagramme \eqref{higgs2-tf27g}, $w$ est le morphisme de transition de $F^a$, et les identifications verticales 
proviennent de l'isomorphisme \eqref{higgs2-tf27d}, est commutatif. On en déduit que le diagramme 
\begin{equation}
\xymatrix{
{\rH^i((U_Y)_\fet,G_U)}\ar[r]\ar[d]&{\rH^i((U_Y\rightarrow U)^a,F^a)}\ar[dd]\\
{\rH^i((U_Y)_\fet,(g_Y)_{\fet*}(G_{U'}))}\ar[d]&\\
{\rH^i((U'_Y)_\fet,G_{U'})}\ar[r]&{\rH^i((U'_Y\rightarrow U')^a,F^a)}}
\end{equation}
où les flèches horizontales sont les applications \eqref{higgs2-tf27e}, les flèches verticales de gauche sont 
définies dans \eqref{higgs2-tf27b} (relativement à $F^a$) et la flèche verticale de droite est l'application canonique, est commutatif. 
Par suite, l'application \eqref{higgs2-tf27ef} définit un morphisme de préfaisceaux sur $\Et_{/X}$. Prenant les faisceaux associés,
on obtient un morphisme de groupes abéliens de $X_\et$
\begin{equation}\label{higgs2-tf27f}
\cH^i(F)\rightarrow \rR^i\sigma_*(F^a).
\end{equation}

\begin{teo}\label{higgs2-tf28}
Soient $F=\{U\mapsto F_U\}$ un préfaisceau en groupes abéliens sur $E$, 
$\ox$ un point géométrique de $X$, $X_{(\ox)}$ le localisé strict de $X$ en $\ox$. 
Supposons que $f$ soit cohérent, et que pour tout $U\in \ob(\Et_{/X})$, $F_U$ soit un faisceau de $(U_Y)_\fet$. 
Alors, pour tout entier $i\geq 0$, la fibre 
\begin{equation}
\cH^i(F)_\ox\rightarrow \rR^i\sigma_*(F^a)_\ox
\end{equation} 
du morphisme \eqref{higgs2-tf27f} en $\ox$ est un isomorphisme.
\end{teo}

Ce théorème sera démontré dans \ref{higgs2-lptf8}.

\section{Limite projective de topos de Faltings}

\subsection{}\label{higgs2-lptf1} 
On désigne par $\fM$ la catégorie des morphismes de schémas et par $\fD$ la catégorie des morphismes 
de $\fM$. Les objets de $\fD$ sont donc des diagrammes commutatifs de morphismes de schémas 
\begin{equation}\label{higgs2-lptf1a} 
\xymatrix{
V\ar[r]\ar[d]&U\ar[d]\\
Y\ar[r]&X}
\end{equation} 
où l'on considère les flèches horizontales comme des objets de $\fM$
et les flèches verticales comme des morphismes de $\fM$; un tel objet sera noté $(V,U,Y,X)$.   
On désigne par $\fE$ la sous-catégorie pleine de $\fD$ formée des objets $(V,U,Y,X)$ tels que le morphisme
$U\rightarrow X$ soit étale de présentation finie et que le morphisme $V\rightarrow U_Y=U\times_XY$ soit étale fini. 
Le ``foncteur but''
\begin{equation}\label{higgs2-lptf1b} 
\fE\rightarrow \fM, \ \ \ (V,U,Y,X)\mapsto (Y\rightarrow X),
\end{equation}
fait de $\fE$ une catégorie fibrée, clivée et normalisée.  
La catégorie fibre au-dessus d'un objet $f\colon Y\rightarrow X$ de $\fM$ est la catégorie notée $E_\coh$ dans \ref{higgs2-tf3}. 
Pour tout diagramme commutatif de morphismes de schémas
\begin{equation}\label{higgs2-lptf1c}
\xymatrix{
Y'\ar[r]^{f'}\ar[d]_{g'}&X'\ar[d]^g\\
Y\ar[r]^f&X}
\end{equation}
le foncteur image inverse de \eqref{higgs2-lptf1b} associé au morphisme $(g',g)$ de $\fM$ 
est le foncteur \eqref{higgs2-tf19c}
\begin{equation}\label{higgs2-lptf1d}
\Phi^+\colon \fE_f\rightarrow \fE_{f'}, \ \ \ (V\rightarrow U)\mapsto (V\times_YY'\rightarrow U\times_XX').
\end{equation}
Munissant chaque fibre de $\fE/\fM$ de la topologie co-évanescente \eqref{higgs2-tf3}, 
$\fE$ devient un $\mU$-site fibré (cf. \ref{higgs2-tf19} et \cite{sga4} VI 7.2.4). On désigne par 
\begin{equation}\label{higgs2-lptf1e}
\fF\rightarrow \fM
\end{equation}
le $\mU$-topos fibré associé à $\fE/\fM$ (\cite{sga4} VI 7.2.6)~: 
la catégorie fibre de $\fF$ au-dessus d'un objet $f\colon Y\rightarrow X$
de $\fM$ est le topos $\tfE_f$ des faisceaux de $\mU$-ensembles sur le site co-évanescent $\fE_f$, 
et le foncteur image inverse relatif au morphisme défini par 
le diagramme \eqref{higgs2-lptf1c} est le foncteur $\Phi^*\colon \tfE_{f}\rightarrow \tfE_{f'}$ image inverse par le morphisme 
de topos $\Phi\colon \tfE_{f'}\rightarrow \tfE_{f}$ associé au morphisme de sites $\Phi^+$ \eqref{higgs2-lptf1d}. 
On note 
\begin{equation}\label{higgs2-lptf1f}
\fF^\vee\rightarrow \fM^\circ
\end{equation}
la catégorie fibrée obtenue en associant à tout objet $f\colon Y\rightarrow X$
de $\fM$ la catégorie $\fF_{f}=\tfE_{f}$, et à tout morphisme défini par 
un diagramme \eqref{higgs2-lptf1c} le foncteur $\Phi_*\colon \tfE_{f'}\rightarrow \tfE_{f}$ image directe par le morphisme 
de topos $\Phi\colon \tfE_{f'}\rightarrow \tfE_{f}$. 
 
\subsection{}\label{higgs2-lptf2}
Soient $I$ une catégorie cofiltrante essentiellement petite (\cite{sga4} I 2.7 et 8.1.8), 
\begin{equation}\label{higgs2-lptf2a} 
\varphi\colon I\rightarrow \fM, \ \ \ i\mapsto (f_i\colon Y_i\rightarrow X_i)
\end{equation}
un foncteur tel que pour tout morphisme $j\rightarrow i$ de $I$, 
les morphismes $Y_j\rightarrow Y_i$ et $X_j\rightarrow X_i$ soient affines. 
On suppose qu'il existe $i_0\in \ob(I)$ tel que $X_{i_0}$ et $Y_{i_0}$ soient cohérents.
On désigne par 
\begin{eqnarray}
\fE_\varphi&\rightarrow& I\label{higgs2-lptf2aa}\\
\fF_\varphi&\rightarrow& I\\
\fF^\vee_\varphi&\rightarrow& I^\circ
\end{eqnarray}
les site, topos et catégorie fibrés déduits de $\fE$ \eqref{higgs2-lptf1b}, $\fF$ \eqref{higgs2-lptf1e} et $\fF^\vee$ \eqref{higgs2-lptf1f}, respectivement,
par changement de base par le foncteur $\varphi$. On notera que $\fF_\varphi$ est le topos fibré associé à $\fE_\varphi$
(\cite{sga4} VI 7.2.6.8). D'après (\cite{ega4} 8.2.3), les limites projectives 
\begin{equation}\label{higgs2-lptf2c}
X=\underset{\underset{i\in \ob(I)}{\longleftarrow}}{\lim}\ X_i \ \ {\rm et}\ \ 
Y=\underset{\underset{i\in \ob(I)}{\longleftarrow}}{\lim}\ Y_i
\end{equation}
sont représentables dans la catégorie des schémas. Les morphismes $(f_i)_{i\in I}$ induisent 
un morphisme $f\colon Y\rightarrow X$, qui représente la limite projective du foncteur \eqref{higgs2-lptf2a}. 

Pour tout $i\in \ob(I)$, on a un diagramme commutatif canonique
\begin{equation}\label{higgs2-lptf2e}
\xymatrix{
Y\ar[r]^f\ar[d]&X\ar[d]\\
Y_i\ar[r]^{f_i}&X_i}
\end{equation}
Il lui correspond un foncteur image inverse \eqref{higgs2-lptf1d}
\begin{equation}\label{higgs2-lptf2ea}
\Phi_i^+\colon \fE_{f_i}\rightarrow \fE_f,
\end{equation}
qui est continu et exact à gauche, et par suite un morphisme de topos 
\begin{equation}\label{higgs2-lptf2eb}
\Phi_i\colon \tfE_{f}\rightarrow \tfE_{f_i}.
\end{equation}

On a un foncteur naturel 
\begin{equation}\label{higgs2-lptf2d}
\fE_\varphi \rightarrow \fE_f,
\end{equation}
dont la restriction à la fibre au-dessus de tout $i\in \ob(I)$ est le foncteur $\Phi_i^+$ \eqref{higgs2-lptf2ea}. 
Ce foncteur transforme morphisme cartésien en isomorphisme. Il se factorise donc de façon unique à travers 
un foncteur (\cite{sga4} VI 6.3)
\begin{equation}\label{higgs2-lptf2f}
\underset{\underset{I^\circ}{\longrightarrow}}{\lim}\ \fE_\varphi \rightarrow \fE_{f}.
\end{equation}
Le $I$-foncteur $\fE_\varphi\rightarrow \fE_f\times I$ déduit de \eqref{higgs2-lptf2d} est un morphisme cartésien
de sites fibrés (\cite{sga4} VI 7.2.2). Il induit donc un morphisme cartésien de topos fibrés (\cite{sga4} VI 7.2.7)
\begin{equation}\label{higgs2-lptf2g}
\tfE_f\times I\rightarrow \fF_\varphi.
\end{equation}

\begin{prop}\label{higgs2-lptf3}
Le couple formé du topos $\tfE_f$ et du morphisme \eqref{higgs2-lptf2g} 
est une limite projective du topos fibré $\fF_\varphi/I$ {\rm (\cite{sga4} VI 8.1.1)}.
\end{prop} 

On notera d'abord que le foncteur \eqref{higgs2-lptf2f} est une équivalence de catégories en vertu de 
(\cite{ega4} 8.8.2,  8.10.5 et 17.7.8). Soient $T$ un $\mU$-topos, 
\begin{equation}
h\colon T\times I\rightarrow \fF_\varphi
\end{equation}
un morphisme cartésien de topos fibrés au-dessus de $I$. Notons $\varepsilon_I\colon \fE_\varphi\rightarrow \fF_\varphi$
le foncteur cartésien canonique (\cite{sga4} VI (7.2.6.7)), et posons 
\begin{equation}
h^+= h^*\circ \varepsilon_I\colon \fE_\varphi\rightarrow T\times I.
\end{equation}
Pour tout $i\in \ob(I)$, on désigne par
\begin{equation}
h_i^+\colon \fE_{f_i}\rightarrow T
\end{equation}
la restriction de $h^+$ aux fibres au-dessus de $i$. 
Compte tenu de l'équivalence de catégories \eqref{higgs2-lptf2f} et de  (\cite{sga4} VI 6.2), il existe un et essentiellement 
un unique foncteur 
\begin{equation}
g^+\colon \fE_{f}\rightarrow T
\end{equation}
tel que $h^+$ soit isomorphe au composé 
\begin{equation}
\xymatrix{
{\fE_\varphi}\ar[r]&{\fE_{f}\times I}\ar[rr]^{g^+\times \id_I}&&{T\times I}},
\end{equation}
où la première flèche est le foncteur déduit de \eqref{higgs2-lptf2d}. Montrons que $g^+$ est un morphisme de sites. 
Pour tout objet $(V\rightarrow U)$ de $\fE_f$, il existe $i\in \ob(I)$, un objet $(V_i\rightarrow U_i)$ de $\fE_{f_i}$ et un isomorphisme de $\fE_f$
\begin{equation}
(V\rightarrow U)\stackrel{\sim}{\rightarrow} \Phi_i^+(V_i\rightarrow U_i).
\end{equation}
Comme les foncteurs $h_i^+$ et $\Phi_i^+$ sont exacts à gauche, 
on en déduit que $g^+$ est exact à gauche. 
D'autre part, tout recouvrement cartésien (resp. vertical) {\em fini} de $\fE_f$ \eqref{higgs2-tcevg3}  
est l'image inverse d'un recouvrement cartésien (resp. vertical) de $\fE_{f_i}$ pour un objet $i\in I$,  
en vertu de (\cite{ega4} 8.10.5(vi)). Comme les schémas $X$ et $Y$ sont cohérents, on en déduit 
que $g^+$ transforme les recouvrements cartésiens (resp. verticaux) de $\fE_f$ en familles épimorphiques de $T$. 
Par suite, $g^+$ est continu en vertu de \ref{higgs2-tcevg5}.  Il définit donc un morphisme de topos 
\begin{equation}
g\colon T\rightarrow \tfE_f
\end{equation}
tel que $h$ soit isomorphe au composé 
\begin{equation}
\xymatrix{
{T\times I}\ar[rr]^{g\times \id_I}&&{\tfE_f\times I}\ar[r]&{\fF_\varphi}},
\end{equation}
où la seconde flèche est le morphisme \eqref{higgs2-lptf2g}. Un tel morphisme $g$ est essentiellement unique car la ``restriction'' 
$g^+\colon \fE_f\rightarrow T$ du foncteur $g^*$ est essentiellement unique d'après ce qui précède, d'où la proposition. 

\subsection{}\label{higgs2-lptf4}
Munissons $\fE_\varphi$ de la topologie totale (\cite{sga4} VI 7.4.1) 
et notons $\Top(\fE_\varphi)$ le topos des faisceaux de $\mU$-ensembles sur $\fE_\varphi$. 
D'après (\cite{sga4} VI 7.4.7), on a une équivalence canonique de catégories \eqref{higgs2-not3}
\begin{equation}\label{higgs2-lptf4b}
\Top(\fE_\varphi)\stackrel{\sim}{\rightarrow}\bHom_{I^\circ}(I^\circ, \fF^\vee_\varphi). 
\end{equation}
D'autre part, le foncteur naturel $\fE_\varphi\rightarrow \fE_f$ \eqref{higgs2-lptf2d} est un morphisme de sites (\cite{sga4} VI 7.4.4)
et définit donc un morphisme de topos 
\begin{equation}\label{higgs2-lptf4c}
\varpi\colon \tfE_f\rightarrow \Top(\fE_\varphi). 
\end{equation}
En vertu de \ref{higgs2-lptf3} et (\cite{sga4} VI 8.2.9), il existe une équivalence de catégories $\Theta$
qui s'insère dans un diagramme commutatif 
\begin{equation}\label{higgs2-lptf4d}
\xymatrix{
{\tfE_f}\ar[r]^-(0.5){\Theta}_-(0.5)\sim\ar[d]_{\varpi_*}&{\bHom_{\cart/I^\circ}(I^\circ, \fF^\vee_\varphi)}\ar@{^(->}[d]\\
{\Top(\fE_\varphi)}\ar[r]^-(0.5)\sim&{\bHom_{I^\circ}(I^\circ, \fF^\vee_\varphi)}}
\end{equation}
où la flèche horizontale inférieure est l'équivalence de catégories \eqref{higgs2-lptf4b} et 
la flèche verticale de droite est l'injection canonique. 

Pour tout objet $F$ de $\Top(\fE_\varphi)$, si 
$\{i\mapsto F_i\}$ est la section correspondante de $\bHom_{I^\circ}(I^\circ, \fF^\vee_\varphi)$, on a un 
isomorphisme canonique fonctoriel (\cite{sga4} VI 8.5.2)
\begin{equation}\label{higgs2-lptf4e}
\varpi^*(F)\stackrel{\sim}{\rightarrow} \underset{\underset{i\in I^\circ}{\longrightarrow}}{\lim}\ \Phi_i^*(F_i).
\end{equation}

\begin{cor}\label{higgs2-lptf55}
Soient $F$ un faisceau de $\Top(\fE_\varphi)$, 
$\{i\mapsto F_i\}$ la section de $\bHom_{I^\circ}(I^\circ, \fF^\vee_\varphi)$ qui lui est 
associée par l'équivalence de catégories \eqref{higgs2-lptf4b}.
Alors on a un isomorphisme canonique fonctoriel
\begin{equation}\label{higgs2-lptf55a}
\underset{\underset{i\in I^\circ}{\longrightarrow}}{\lim}\ \Gamma(\tfE_{f_i},F_i)\stackrel{\sim}{\rightarrow}
\Gamma(\tfE_{f},\underset{\underset{i\in I^\circ}{\longrightarrow}}{\lim}\ \Phi_i^*(F_i)).
\end{equation}
\end{cor}

\begin{cor}\label{higgs2-lptf5}
Soient $F$ un faisceau abélien de $\Top(\fE_\varphi)$, 
$\{i\mapsto F_i\}$ la section de $\bHom_{I^\circ}(I^\circ, \fF^\vee_\varphi)$ qui lui est 
associée par l'équivalence de catégories \eqref{higgs2-lptf4b}.
Alors pour tout entier $q\geq 0$, on a un isomorphisme canonique fonctoriel
\begin{equation}\label{higgs2-lptf5a}
\underset{\underset{i\in I^\circ}{\longrightarrow}}{\lim}\ \rH^q(\tfE_{f_i},F_i)\stackrel{\sim}{\rightarrow}
\rH^q(\tfE_{f},\underset{\underset{i\in I^\circ}{\longrightarrow}}{\lim}\ \Phi_i^*(F_i)).
\end{equation}
\end{cor}

Les corollaires \ref{higgs2-lptf55} et \ref{higgs2-lptf5} résultent de \ref{higgs2-lptf3} et (\cite{sga4} VI 8.7.7). 
On notera que les conditions requises dans 
(\cite{sga4} VI 8.7.1 et 8.7.7) sont satisfaites en vertu de \ref{higgs2-tf4} et (\cite{sga4} VI 3.3, 5.1 et 5.2).

\subsection{}\label{higgs2-lptf51}
On désigne par 
\begin{eqnarray}
\cR_\varphi&\rightarrow& I,\label{higgs2-lptf51a}\\
\cG_\varphi&\rightarrow& I,\label{higgs2-lptf51b}\\
\cG^\vee_\varphi&\rightarrow& I^\circ\label{higgs2-lptf51c}
\end{eqnarray} 
les site, topos et catégorie fibrés déduits, respectivement, du site fibré des revêtements étales $\cR/\Sch$ \eqref{higgs2-Kp6a}, 
du topos fibré $\cG/\Sch$ \eqref{higgs2-Kp6b} et de la catégorie fibrée $\cG^\vee/\Sch^\circ$ \eqref{higgs2-Kp6c}, 
par changement de base par le foncteur 
\begin{equation}\label{higgs2-lptf51d}
I\rightarrow \Sch, \ \ \ i\mapsto Y_i
\end{equation}
induit par $\varphi$ \eqref{higgs2-lptf2a}. Pour tout $i\in \ob(I)$, on note 
\begin{equation}\label{higgs2-lptf51de}
t_i\colon Y\rightarrow Y_i
\end{equation}
le morphisme canonique \eqref{higgs2-lptf2c}. On a un foncteur naturel 
\begin{equation}\label{higgs2-lptf51e}
\cR_\varphi \rightarrow \cR_Y
\end{equation}
dont la restriction à la fibre au-dessus de $i\in \ob(I)$ est donnée par le foncteur de changement 
de base par le morphisme $t_i$
\[
\cR_{Y_i}\rightarrow \cR_Y,\ \ \ Y'_i\mapsto Y'_i\times_{Y_i}Y.
\] 
Ce foncteur transforme morphisme cartésien en isomorphisme, et se factorise donc de façon unique à travers 
un foncteur  
\begin{equation}\label{higgs2-lptf51f}
\underset{\underset{I^\circ}{\longrightarrow}}{\lim}\ \cR_\varphi \rightarrow \cR_Y.
\end{equation}
Le $I$-foncteur $\cR_\varphi \rightarrow \cR_Y\times I$ déduit de \eqref{higgs2-lptf51e} est un morphisme cartésien de sites fibrés 
(\cite{sga4} VI 7.2.2). Il induit donc un morphisme de topos fibrés 
\begin{equation}\label{higgs2-lptf51g}
Y_\fet\times I\rightarrow \cG_\varphi.
\end{equation}

\begin{lem}\label{higgs2-lptf52}
Le foncteur \eqref{higgs2-lptf51f} est une équivalence de sites lorsque l'on munit 
la source de la topologie de la limite inductive du site fibré $\cR_\varphi$
{\rm (\cite{sga4} VI 8.2.5)} et le but de la topologie étale.
\end{lem} 

En effet, le foncteur \eqref{higgs2-lptf51f} est une équivalence de catégories en vertu de (\cite{ega4} 8.8.2, 8.10.5 et 17.7.8).  
Soient $i\in \ob(I)$, $g_i\colon Y'_i\rightarrow Y_i$ un revêtement étale, $g\colon Y'\rightarrow Y$ le revêtement 
étale déduit de $g_i$ par changement de base par le morphisme $Y\rightarrow Y_i$. D'après (\cite{ega4} 8.10.5), 
pour que $g$ soit surjectif, il faut et il suffit qu'il existe un morphisme $j\rightarrow i$ de $I$ 
tel que le revêtement étale $g_j\colon Y'_j\rightarrow Y_j$ déduit de $g_i$ par changement de base, soit surjectif. 
L'assertion relative aux topologies s'ensuit compte tenu de (\cite{sga4} VI 8.2.2 et III 1.6). 

\begin{prop}\label{higgs2-lptf53}
Le couple formé du topos $Y_\fet$ et du morphisme \eqref{higgs2-lptf51g} 
est une limite projective du topos fibré $\cG_\varphi/I$. 
\end{prop}

Cela résulte de \ref{higgs2-lptf52} et (\cite{sga4} VI 8.2.3).

\begin{cor}\label{higgs2-lptf535}
Soient $F$ un faisceau abélien du topos total de $\cR_\varphi$, 
$\{i\mapsto F_i\}$ la section de $\bHom_{I^\circ}(I^\circ, \cG^\vee_\varphi)$ qui lui est associée.
Alors pour tout entier $q\geq 0$, on a un isomorphisme canonique fonctoriel
\begin{equation}\label{higgs2-lptf535a}
\underset{\underset{i\in I^\circ}{\longrightarrow}}{\lim}\ \rH^q((Y_i)_\fet,F_i)\stackrel{\sim}{\rightarrow}
\rH^q(Y_\fet,\underset{\underset{i\in I^\circ}{\longrightarrow}}{\lim}\ (t_i)^*_\fet(F_i)).
\end{equation}
\end{cor}

Cela résulte de \ref{higgs2-lptf53} et (\cite{sga4} VI 8.7.7). 
On notera que les conditions requises dans 
(\cite{sga4} VI 8.7.1 et 8.7.7) sont satisfaites en vertu de \ref{higgs2-Kpp1} et (\cite{sga4} VI 3.3, 5.1 et 5.2).

\subsection{}\label{higgs2-lptf54}
On désigne par 
\begin{equation}\label{higgs2-lptf54a}
\beta \colon \tfE_f\rightarrow Y_\fet
\end{equation}
et, pour tout $i\in \ob(I)$, par
\begin{equation}\label{higgs2-lptf54b}
\beta_i \colon \tfE_{f_i}\rightarrow (Y_i)_\fet
\end{equation}
les morphismes canoniques \eqref{higgs2-tf11d}. Il résulte de \eqref{higgs2-tf19f} et (\cite{sga1} VI 12; cf. aussi \cite{egr1} 1.1.2) 
qu'il existe essentiellement un unique morphisme cartésien de topos fibrés
\begin{equation}
\beta_\varphi \colon \fF_{\varphi}\rightarrow \cG_\varphi
\end{equation}
dont la fibre au-dessus de tout $i\in \ob(I)$ est le morphisme $\beta_i$. 
De plus, le diagramme de morphismes de topos fibrés
\begin{equation}
\xymatrix{
{\tfE_f\times I}\ar[r]^{\beta\times \id}\ar[d]&{Y_\fet\times I}\ar[d]\\
{\fF_{\varphi}}\ar[r]^{\beta_\varphi}&{\cG_\varphi}}
\end{equation}
où les flèches verticales sont les morphismes \eqref{higgs2-lptf2g} et \eqref{higgs2-lptf51g}, 
est commutatif à isomorphisme canonique près. On peut donc identifier $\beta$ au morphisme 
déduit de $\beta_\varphi$ par passage à la limite projective dans le sens de (\cite{sga4} VI 8.1.4).

\begin{prop}\label{higgs2-lptf56}
Soient $F$ un faisceau de $\Top(\fE_\varphi)$, 
$\{i\mapsto F_i\}$ la section de $\bHom_{I^\circ}(I^\circ,\fF_\varphi^\vee)$
qui lui est associée par l'équivalence de catégories \eqref{higgs2-lptf4b}.
Alors, on a un isomorphisme canonique fonctoriel
\begin{equation}\label{higgs2-lptf56a} 
\underset{\underset{i\in I^\circ}{\longrightarrow}}{\lim}\
(t_i)^*_\fet(\beta_{i*}(F_i))\stackrel{\sim}{\rightarrow} \beta_*(\varpi^*F),
\end{equation}
où $\varpi$ est le morphisme \eqref{higgs2-lptf4c} et $t_i$ est le morphisme \eqref{higgs2-lptf51de}.
\end{prop}

Cela résulte de \ref{higgs2-lptf54} et (\cite{sga4} VI 8.5.9). On notera que les conditions requises dans {\em loc. cit.}
sont satisfaites en vertu de \ref{higgs2-tf4}, \ref{higgs2-tf20} et (\cite{sga4} VI 3.3 et 5.1). 

\begin{rema}\label{higgs2-lptf565}
Soient $F$ un faisceau abélien de $\Top(\fE_\varphi)$, 
$\{i\mapsto F_i\}$ l'objet de $\bHom_{I^\circ}(I^\circ,\fF_\varphi^\vee)$
qui lui est associé par l'équivalence de catégories \eqref{higgs2-lptf4b}, $q$ un entier $\geq 0$. 
Il résulte alors de \ref{higgs2-tf195} que le diagramme 
\begin{equation}\label{higgs2-lptf565a}
\xymatrix{
{\underset{\underset{i\in I^\circ}{\longrightarrow}}{\lim}\ \rH^q((Y_i)_\fet,\beta_{i*}(F_i))}\ar[r]\ar[d]_u&
{\underset{\underset{i\in I^\circ}{\longrightarrow}}{\lim}\ \rH^q(\tfE_{f_i},F_i)}\ar[dd]^w\\
{\underset{\underset{i\in I^\circ}{\longrightarrow}}{\lim}\ \rH^q(Y_\fet,(t_i)_\fet^*(\beta_{i*}(F_i)))}\ar[d]_v&\\
{\rH^q(Y_\fet,\beta_*(\varpi^*F))}\ar[r]&{\rH^q(\tfE_{f},\varpi^*F)}}
\end{equation}
où les flèches horizontales proviennent des suites spectrales de Cartan-Leray, 
$u$ est le morphisme canonique, $v$ est induit par \eqref{higgs2-lptf56a}, et $w$ est induit par \eqref{higgs2-lptf4e}, est commutatif. 
\end{rema}

\subsection{}\label{higgs2-lptf6}
Nous pouvons maintenant démontrer la proposition \ref{higgs2-tf268}. Choisissons un voisinage étale affine 
$X_0$ de $\ox$ dans $X$. Notons $I$ la catégorie des $X_0$-schémas étales $\ox$-pointés
qui sont affines au-dessus de $X_0$ (cf. \cite{sga4} VIII 3.9 et 4.5), et $\varphi \colon I\rightarrow \fM$ 
le foncteur qui à un objet $U$ de $I$, associe la projection canonique $f_{U}\colon U_Y\rightarrow U$. 
Alors, $\uf$ s'identifie canoniquement à la limite projective du foncteur $\varphi$. 
Pour tout $U\in \ob(I)$, le topos $\tE_{/(U_Y\rightarrow U)^a}$ est canoniquement équivalent au 
topos de Faltings associé au morphisme $f_{U}$ d'après \eqref{higgs2-tf17f}. Par suite, avec les notations 
de cette section, pour tout faisceau $F$ de $\tE$, 
$\{U\mapsto F|(U_Y\rightarrow U)^a\}$ est naturellement une section de $\bHom_{I^\circ}(I^\circ,\fF_\varphi^\vee)$.
Elle définit donc un faisceau de $\Top(\fE_\varphi)$ \eqref{higgs2-lptf4b}. On a un isomorphisme canonique fonctoriel \eqref{higgs2-lptf4e}
\begin{equation}
\Phi^*(F)\stackrel{\sim}{\rightarrow} \varpi^*(\{U\mapsto F|(U_Y\rightarrow U)^a\}).
\end{equation}

(i) D'après (\cite{sga4} VIII 3.9 et 4.5), on a un isomorphisme canonique 
\begin{equation}\label{higgs2-lptf6a}
\sigma_*(F)_\ox\stackrel{\sim}{\rightarrow} \underset{\underset{U\in I^\circ}{\longrightarrow}}{\lim}\ 
\Gamma((U_Y\rightarrow U)^a,F).
\end{equation}
Celui-ci induit un isomorphisme fonctoriel
\begin{equation}\label{higgs2-lptf6b}
\sigma_*(F)_\ox\stackrel{\sim}{\rightarrow} \Gamma(\tuE,\Phi^*(F)),
\end{equation}
en vertu de \ref{higgs2-lptf55}.
La proposition s'en déduit compte tenu de \ref{higgs2-tf266}(i).

(ii) Cela résulte, comme pour (i), de \ref{higgs2-lptf5}, \ref{higgs2-tf266}(ii) et de l'isomorphisme canonique (\cite{sga4} V 5.1(1))
\begin{equation}
\rR^i\sigma_*(F)_\ox\stackrel{\sim}{\rightarrow} \underset{\underset{U\in I^\circ}{\longrightarrow}}{\lim}\ 
\rH^i((U_Y\rightarrow U)^a,F).
\end{equation}

(iii) Cela résulte aussitôt de la preuve de (ii).

\subsection{}\label{higgs2-lptf8}
Nous pouvons enfin démontrer le théorème \ref{higgs2-tf28}.
Choisissons un voisinage étale affine $X_0$ de $\ox$ dans $X$. Notons $I$ la catégorie des $X_0$-schémas étales $\ox$-pointés
qui sont affines au-dessus de $X_0$, et $\varphi \colon I\rightarrow \fM$ 
le foncteur qui à un objet $U$ de $I$, associe la projection canonique $f_{U}\colon U_Y\rightarrow U$. 
D'après (\cite{sga4} IV (6.8.4)) et \ref{higgs2-lptf535}, on a des isomorphismes canoniques 
\begin{equation}\label{higgs2-lptf8a}
\cH^i(F)_\ox\stackrel{\sim}{\rightarrow}\underset{\underset{U\in I^\circ}{\longrightarrow}}{\lim}\ \rH^i((U_Y)_\fet,F_U)
\stackrel{\sim}{\rightarrow}\rH^i(\uY_\fet,\underset{\underset{U\in I^\circ}{\longrightarrow}}{\lim}\ (\xi_Y)_\fet^*F_U).
\end{equation}
En vertu de de \ref{higgs2-tf268}(ii), on a un isomorphisme 
\begin{equation}\label{higgs2-lptf8c}
\rR^i\sigma_*(F^a)_\ox\stackrel{\sim}{\rightarrow} \rH^i(\uY_\fet,\varphi_\ox(F^a)).
\end{equation}
D'autre part, il résulte de \ref{higgs2-lptf565} et des définitions que le diagramme 
\begin{equation}\label{higgs2-lptf8d}
\xymatrix{
{\cH^i(F)_\ox}\ar[r]^-(0.4)\sim\ar[d]_u&{\rH^i(\uY_\fet,\underset{\underset{U\in I^\circ}{\longrightarrow}}{\lim}\ 
(\xi_Y)_\fet^*F_U)}\ar[d]^v\\
{\rR^i\sigma_*(F^a)_\ox}\ar[r]^-(0.4)\sim&{\rH^i(\uY_\fet,\varphi_\ox(F^a))}}
\end{equation}
où les flèches horizontales sont les isomorphismes \eqref{higgs2-lptf8a} et \eqref{higgs2-lptf8c}, 
$u$ est la fibre du morphisme \eqref{higgs2-tf27f} en $\ox$, et $v$ est induit par l'isomorphisme \eqref{higgs2-tf290a}, est commutatif.  
Par suite, $u$ est un isomorphisme, d'où le théorème.

\printindex

\end{document}